\pdfoutput=1
\RequirePackage{ifpdf}
\ifpdf 
\documentclass[pdftex]{sigma}
\else
\documentclass{sigma}
\fi

\makeindex

\usepackage{tabularx,mdwtab}
\usepackage{bbm}

\usepackage{lscape}

\def\cal{\relax}
\usepackage{DLdef1,makeidx}

\def\mmat #1,#2,#3,#4,{\text{\small\arraycolsep=3pt $
\begin{pmatrix}#1&#2\\#3&#4\end{pmatrix}$}}

\newcommand {\fH}{{\mathfrak{H}}}
\newcommand {\fme}{{\mathfrak{me}}}

\newcommand{\un}{\protect\underline{N}}
\newcommand{\uM}{\protect\underline{M}}

\newcommand{\del}{\partial}
\newcommand{\Size}{\text{Size}}

\newcommand {\lra} {\longleftrightarrow}
\newcommand {\htr}{{\text{htr}}}
\newcommand {\fba}{{\mathfrak{ba}}}

\numberwithin{equation}{section}

\newtheorem{Theorem}{Theorem}[section]

\newtheorem{Lemma}[Theorem]{Lemma}
\newtheorem{Proposition}[Theorem]{Proposition}
\newtheorem{Conjecture}[Theorem]{Conjecture}
\newtheorem{Problem}[Theorem]{Problem}
\newtheorem{Claim}[Theorem]{Claim}

 { \theoremstyle{definition}

\newtheorem{Remark}[Theorem]{Remark}
\newtheorem{Remarks}[Theorem]{Remarks} }

\begin{document}

\allowdisplaybreaks

\newcommand{\arXivNumber}{1510.07255}

\renewcommand{\thefootnote}{}

\renewcommand{\PaperNumber}{089}

\FirstPageHeading

\ShortArticleName{Simple Vectorial Lie Algebras in Characteristic~2 and their Superizations}

\ArticleName{Simple Vectorial Lie Algebras in Characteristic~2\\ and their Superizations\footnote{This paper is a~contribution to the Special Issue on Algebra, Topology, and Dynamics in Interaction in honor of Dmitry Fuchs. The full collection is available at \href{https://www.emis.de/journals/SIGMA/Fuchs.html}{https://www.emis.de/journals/SIGMA/Fuchs.html}}}

\Author{Sofiane BOUARROUDJ~$^{\dag^1}$, Pavel GROZMAN~$^{\dag^2}$, Alexei LEBEDEV~$^{\dag^2}$,\newline
Dimitry LEITES~$^{\dag^1\dag^3}$ and Irina SHCHEPOCHKINA~$^{\dag^4}$}

\AuthorNameForHeading{S.~Bouarroudj, P.~Grozman, A.~Lebedev, D.~Leites and I.~Shchepochkina}

\Address{$^{\dag^1}$~New York University Abu Dhabi, Division of Science and Mathematics,\\
\hphantom{$^{\dag^1}$}~P.O.~Box 129188, United Arab Emirates}
\EmailDD{\href{mailto:sofiane.bouarroudj@nyu.edu}{sofiane.bouarroudj@nyu.edu}, \href{mailto:dl146@nyu.edu}{dl146@nyu.edu}}

\Address{$^{\dag^2}$~Equa Simulation AB, R{\aa}sundav\"{a}gen 100, Solna, Sweden}
\EmailDD{\href{mailto:pavel.grozman@bredband.net}{pavel.grozman@bredband.net}, \href{mailto:alexeylalexeyl@mail.ru}{alexeylalexeyl@mail.ru}}

\Address{$^{\dag^3}$~Department of Mathematics, Stockholm University, SE-106 91 Stockholm, Sweden}
\EmailDD{\href{mailto:mleites@math.su.se}{mleites@math.su.se}}

\Address{$^{\dag^4}$~Independent University of Moscow, Bolshoj Vlasievsky per.~11, 119002 Moscow, Russia}
\EmailDD{\href{mailto:irina@mccme.ru}{irina@mccme.ru}}

\ArticleDates{Received September 25, 2019, in final form August 25, 2020; Published online September 24, 2020}

\Abstract{We overview the classifications of simple finite-dimensional modular Lie algebras. In characteristic~2, their list is wider than that in other characteristics; e.g., it contains desuperizations of modular analogs of complex simple vectorial Lie superalgebras.
We consider odd parameters of deformations. For all 15 Weisfeiler gradings of the 5 exceptional families, and one Weisfeiler grading for each of 2 serial simple complex Lie superalgebras (with 2 exceptional subseries), we describe their characteristic-2 analogs~-- new simple Lie algebras. Descriptions of several of these analogs, and of their desuperizations, are far from obvious. One of the exceptional simple vectorial Lie algebras is a~previously unknown deform (the result of a~deformation) of the characteristic-2 version of the Lie algebra of divergence-free vector fields; this is a~new simple Lie algebra with no analogs in characteristics distinct from~2. In characteristic~2, every simple Lie superalgebra can be obtained from a~simple Lie algebra by one of the two methods described in arXiv:1407.1695. Most of the simple Lie superalgebras thus obtained from simple Lie algebras we describe here are new.}

\Keywords{modular vectorial Lie algebra; modular vectorial Lie superalgebra}

\Classification{17B50; 17B20; 70F25}

\begin{flushright}
\it To Dmitry Borisovich Fuchs with admiration
\end{flushright}

\setcounter{tocdepth}{3}
\tableofcontents

\printindex

\renewcommand{\thefootnote}{\arabic{footnote}}
\setcounter{footnote}{0}

\section{Notation and background}\label{Sbackgr}
To make the text understandable for the uninitiated, we place the most basic facts before the Introduction which we have divided into two parts to make it more readable; for the same reason we divided the background to accommodate both students and experts. Statements proved directly, or by means of \textit{Mathematica}-based \textit{SuperLie} package \cite{Gr}, are called Claims.

We recall the basics and show how to modify familiar
formulas in order to pass from $\Cee$ to fields of positive characteristic, especially
characteristic~2. In some formulas given for $p=2$, we retain notation
convenient for comparison with the cases where $p\neq 2$.

\subsection{Main points}
\begin{enumerate}\itemsep=0pt
\item[\textbf{1.}] We give an overview of the classification of simple finite-dimensional modular Lie algebras and Lie superalgebras over an algebraically closed field $\Kee$ of characteristic~$p>0$. We update the conjectures for various values of $p>0$.
\item[\textbf{2.}] We use our results on classification of Lie superalgebras of vector fields with polynomial coefficients over $\Cee$ to describe their characteristic-$p$ versions, especially, their desuperizations, for all 15~Weisfeiler gradings of all 5~exceptional simple vectorial Lie superalgebras, and for several serial ones, also exceptional in a~sense.
\item[\textbf{3.}] One of the deforms\footnote{\textit{Deform} is the result of a~deformation, like transform is the result of a~transformation.} of the divergence-free Lie algebra $\fsvect(5;\un)$ which exists only if $p=2$. It is one of the exceptional simple vectorial Lie algebras~-- a~desuperization of an exceptional simple vectorial Lie superalgebra. \textit{This is the most unexpected result of this paper}.
\end{enumerate}

\subsection{Generalities} As is now customary, we denote the elements of $\Zee/2$ by $\ev$ and
$\od$, to distinguish them from integers. For us, $\Nee:=\{1,
2,\dots\}$,\index{$Nmathbb$@$\Nee$} as it used to be in the past, and still is in some countries; we set $\Zee_+:=\Nee\cup\{0\}$.\index{$Zmathbb_+:=\Nee\cup\{0\}$@$\Zee_+:=\Nee\cup\{0\}$} The parity $p$ of a~non-zero element $v$ of a~$\Zee/2$-graded space $V$, called a~\textit{superspace},\index{(Superspace@Superspace} is equal to $i$ if and only if $v\in V_i$. Any $\Zee/2$-graded algebra is called a~\textit{superalgebra}.\index{(Superalgebra@Superalgebra}

Hereafter $\Kee$ is an algebraically closed field of characteristic~
$p>0$; usually, $p=2$ and all Lie (super)algebras are finite-dimensional
unless otherwise stated.
Mostly (exceptions indicated), $\Pi$\index{((ZZPi, change of parity@$\Pi$, change of parity}~denotes the change of parity functor, i.e., tensoring by
$\Pi(\Zee)$.

The superization of most formulas of algebra is achieved via the following \textit{sign rule}\index{(Sign Rule@Sign Rule}
\begin{equation*}\label{SignR}
\begin{minipage}[c]{14cm}
``If something of parity $a$ is moved past something of parity $b$, the sign $(-1)^{ab}$ accrues. Formulas defined only on homogeneous elements are extended to arbitrary elements via linearity."
\end{minipage}
\end{equation*}

Note that only even homomorphisms are considered as
\textit{morphisms of superalgebras}.\index{(Morphism of superalgebras@Morphism of superalgebras}

Observe that sometimes applying the Sign Rule requires some dexterity. For example,
we have to distinguish between two versions both of which turn in the
non-super case into one, called skew- or anti-commutativity, which are synonyms only in the \textit{non-super} case; for two elements $a$ and $b$ of a~superalgebra we call the following conditions
\begin{alignat*}{3}
&ba=(-1)^{p(b)p(a)}ab \quad &&\text{super commutativity},& \nonumber \\
&ba=-(-1)^{p(b)p(a)}ab \quad &&\text{super anti-commutativity},& \nonumber \\
&ba=(-1)^{(p(b)+1)(p(a)+1)}ab \quad &&\text{super skew-commutativity},& \nonumber \\
&ba=-(-1)^{(p(b)+1)(p(a)+1)}ab
\quad &&\text{super antiskew-commutativity}.& 
\end{alignat*}
Examples: the bracket in any Lie superalgebra is \textit{super anti-commutative}; the anti-bracket $\{-,-\}_{\rm B.b.}$, see~\eqref{BB}, being anti-commutative relative the parity in the Lie superalgebra is, however,
\textit{super antiskew-commutative} relative to the \textit{natural parity of generating functions}.

\subsubsection{Conventions and notation often used}\label{Conv1}\index{(Convention@Convention}
 In what follows, we assume that every \textit{supercommutative} superalgebra is associative with 1; their morphisms send~1 to~1.\index{(Morphisms of supercommutative superalgebras@Morphisms of supercommutative\newline superalgebras}

We denote by $\fc$ the center of a~given Lie (super)algebra; both
$\fc(\fg)$ and $\fc\fg:=\fc\oplus \fg$ denote a~trivial central extension of~$\fg$. \index{$cmathfrak$, $\fc\fg:=\fc(\fg)$@$\fc$, $\fc\fg:=\fc(\fg)$}

Let $\fa\ltimes \fb$\index{$amathfrak\ltimes \fb$@$\fa\ltimes \fb$}
or $\fb\rtimes \fa$ denote the semi-direct sum of modules (algebras)
in which $\fa$ is a~submodule (ideal).

Let $\fd(\fg):=\fg\ltimes\Kee D$, \index{$dmathfrak(\fg)$@$\fd(\fg)$} where $D$ is an outer derivation of~$\fg$; unless specified otherwise, $D$ is the grading operator of~$\fg$. For example, for $\fd(\fo_\Pi(2k))$, we take $D=\diag(I_k, 0_k)$.

Let $\fg^{\bcdot k}:=[\fg, [\fg,
\dots [\fg, \fg]\dots]]$,\index{$gmathfrak^{\bcdot k}:=[\fg, [\fg,
\dots [\fg, \fg]\dots]]$@$\fg^{\bcdot k}:=[\fg, [\fg,
\dots [\fg, \fg]\dots]]$} the $k$-fold bracket.

We denote the
functor of raising to the $n$th symmetric (resp.\ exterior) power by\index{$S^n$, the $n$th symmetric power}\index{$E^n=\Lambda^n$, the $n$th exterior power}
$S^n$ (resp.~$E^n$, often denoted also by $\Lambda^n$); sometimes we denote the exterior (Grassmann)
algebra by~$\Lambda[\theta]$ or~$\Lambda(r)$\index{((ZZLambda[\theta]=\Lambda(r)$, the exterior algebra in $\theta=(\theta_1,\dots, \theta_r)$@$\Lambda[\theta]=\Lambda(r)$, the exterior algebra\newline in $\theta=(\theta_1,\dots, \theta_r)$}
in generators $\theta=(\theta_1,\dots, \theta_r)$ satisfying anti-commutativity
relations (and, additionally, $\theta_i^2=0$ for all $i$ if $p=2$).

The symbol $\id$ (sometimes $\id_n$, $\id_{a|b}$)\index{$Zmi$, $\id_n$, $\id_{a\vert b}$, the identity operator or the tautological module @$\id$, $\id_n$, $\id_{a\vert b}$, the identity operator\newline or the tautological module} denotes
not only the identity operator (in the space of dimension $n$, resp.\ superspace of superdimension $a|b$), but also the tautological module $V$ over the linear Lie superalgebra $\fg\subset\fgl(V)$; sometimes we write
$\id_\fg:=V$ for clarity.

$L_D$ is the Lie derivative along the vector field~$D$.\index{$L_D$, the Lie derivative}

\subsubsection[Definition of Lie superalgebras for $p\neq 2, 3$]{Definition of Lie superalgebras for $\boldsymbol{p\neq 2, 3}$}

The ``naive'' definition of Lie superalgebras for $p\neq 2, 3$ is obtained by applying the Sign Rule to anti-commutativity and Jacobi identities. To understand deformations with odd parameter, we need a~more sophisticated approach using the functor of points. The multiplication in the Lie superalgebra will be called \textit{super-bracket} or just \textit{bracket}.\index{(Bracket@Bracket}

\subsubsection[Definition of Lie superalgebras for $p=2$]{Definition of Lie superalgebras for $\boldsymbol{p=2}$}\label{defLieS2}\index{(Lie superalgebra for $p=2$@Lie superalgebra for $p=2$}

\underline{If $p=2$}, the \textit{antisymmetry}\index{(Antisymmetry@Antisymmetry} condition for Lie algebra $\fg_\ev$ should be replaced by an equivalent for $p\neq 2$, but otherwise
stronger, \textit{alternating} or \textit{antisymmetry} condition
\[
[x,x]=0 \  \text{for any $x\in\fg_\ev$}.
\]
\underline{If $p=2$}, a~\textit{Lie superalgebra}\index{(Lie superalgebra@Lie superalgebra}
is a~superspace $\fg=\fg_\ev\oplus\fg_\od$ such that $\fg_\ev$ is a
Lie algebra, $\fg_\od$ is a~$\fg_\ev$-module (made two-sided by
symmetry), with a~\textit{squaring}\index{(Squaring@Squaring} $x\mapsto x^2$ and a~\textit{bracket}\index{(Bracket@Bracket} of odd elements, which are defined via a~linear map $s\colon S^2(\fg_\od)\tto\fg_\ev$, where $S^2$ denotes the operator of raising to symmetric square, as follows:
\begin{gather*}
x^2 := s(x\otimes x),\\ 
{}[x,y] := s(x\otimes y + y\otimes x)\text{~for any $x,y\in\fg_\od$.}
\end{gather*}
The linearity of the $\fg_\ev$-valued function $s$ implies that
\begin{gather*}
\text{$(ax)^2=a^2x^2$ for any $x\in \fg_\od$ and $a\in \Kee$, and}\\ 
{}[x,y] \ \text{is a~bilinear form on $\fg_\od$ with values in $\fg_\ev$.} 
\end{gather*}

The \textit{Jacobi identity} involving odd elements takes the form of the following two conditions:
\begin{gather}
 \big[x^2,y\big]=[x,[x,y]] \ \text{for any~} x\in\fg_\od, y\in\fg_\ev,\quad 
 \big[x^2,x\big]=0 \ \text{for any $x\in\fg_\od$.}\label{6}
\end{gather}

The (super)algebra satisfying only Jacobi identity, without any symmetry conditions, is called a~\textit{Leibniz $($super$)$algebra}.\index{(Leibniz (super)algebra@Leibniz (super)algebra}

Over $\Zee/2$, the condition \eqref{6} must (for a~reason, see \cite{KLLS}) be replaced with a~more general one:{\samepage
\begin{gather}\label{JIgen}
\big[x^2,y\big]=[x,[x,y]] \quad \text{for any $x \in \fg_\od$ and $y\in \fg$}.
\end{gather}
For any other ground field, this condition is equivalent to condition~\eqref{6}.}

More generally, for any Lie superalgebra $\fg$, since we want the Lie superalgebra $\fder\ \fg$ of all derivations of $\fg$ to be a~Lie superalgebra, we have to add (to the Leibniz rule) the following condition on derivations (it becomes \eqref{JIgen} for $D=\ad_y$)
\[
D\big(x^{2}\big) = [D(x),x] \ \text{for odd elements $x\in\fg_\od$ and any $D\in\fder\ \fg$.}
\]

By an \textit{ideal} $\fii$ of a~Lie superalgebra $\fg$, one always means $\fii$ \textit{homogeneous} with respect to parity, i.e., equal to ${\fii\cap\fg_\ev\bigoplus \fii\cap\fg_\od}$; for $p=2$, the ideal should be closed with respect to squaring.\index{(Ideal, of the Lie superalgebra@Ideal, of the Lie superalgebra}

Recall that a~given Lie (super)algebra $\fg$ is said to be
\textit{simple} if $\dim\fg>1$ and $\fg$ has no proper ideals;\index{(Lie superalgebra, simple@Lie superalgebra, simple} $\fg$
is \textit{semisimple}\index{(Lie superalgebra, semisimple@Lie superalgebra, semisimple} if its radical is zero; $\fg$ is
\textit{almost simple}\index{(Lie superalgebra, almost simple@Lie superalgebra, almost simple} if it can be sandwiched (non-strictly) between
a~simple Lie superalgebra $\fs$ and the Lie superalgebra
$\fder(\fs)$ of derivations of~$\fs$, i.e.,
$\fs\subset\fg\subset\fder(\fs)$.

The definition of the \textit{derived of the Lie superalgebra}\index{(Lie superalgebra, derived@Lie superalgebra, derived} $\fg$ changes when $p=2$: let $\fg^{(0)}:=\fg$ and for any $i\geq 0$, set
\begin{equation*}\label{deralg}
\fg^{(i+1)}=(\fg^{(i)})^{(1)}:=\begin{cases}\big[\fg^{(i)},\fg^{(i)}\big]&\text{if $p\neq 2$},\\
\big[\fg^{(i)},\fg^{(i)}\big]+\Span\big\{g^2\,|\,
g\in\fg^{(i)}_\od\big\}&\text{if $p=2$}.
\end{cases}
\end{equation*}

An even linear mapping $r\colon \fg\tto\fgl(V)$ is said to be a
\textit{representation of the Lie superalgebra}\index{(Representation of the Lie superalgebra@Representation of the Lie superalgebra}~$\fg$, and $V$ is said to be a~\textit{$\fg$-module} if
\begin{gather*}
r([x, y])=[r(x), r(y)]\quad \text{for any $x, y\in \fg$,} \\
r(x^2)=(r(x))^2\quad \text{for any $x\in\fg_\od$.}
\end{gather*}

\subsubsection[Definition of Lie superalgebras for $p=3$]{Definition of Lie superalgebras for $\boldsymbol{p=3}$}\label{defLieP3}\index{(Lie superalgebra for $p=3$@Lie superalgebra for $p=3$} Since we are giving a~review of the context for any $p>0$, we have to note a~peculiarity of
$p=3$, where the \textit{Jacobi identity}\index{(Jacobi identity@Jacobi identity} for \textit{Lie superalgebras}\index{(Lie superalgebra@Lie superalgebra} entails, additionally, that
\begin{equation}\label{p3super}
[x, [x,x]]=0\quad \text{for any $x\in\fg_\od$}.
\end{equation}
The super anti-commutative algebra satisfying the Jacobi identity, but not \eqref{p3super} is called a~\textit{pre-Lie superalgebra}.\index{(Pre-Lie superalgebra@Pre-Lie superalgebra} For interesting examples of pre-Lie superalgebras, see~\cite{BBH}.

\subsection{Basics on the functor of points}\label{FoP} In this subsection, we follow \cite{Ma}; we advise the reader interested in subtleties that we, like most authors, do not dwell on, to read the Appendices to~\cite{Mo}.

For a~fixed object $M$ and any object $X$ of a~category $\sfC$, the association $X\longmapsto \Hom_\sfC(X,M)$ defines a~functor $F\colon \sfC\leadsto\Sets$. The idea is

 (A) To consider $\Hom_\sfC(X,M)$ as the set of points of $M$, which is indeed the case for any set~$M$ if $X$ is a~point of $M$, and $\sfC=\Sets$;

 (B) Considering objects of the category $\sfC$ of sets endowed with a~structure (of a~group, algebra, module over a~fixed algebra, topological space, etc.), and the morphisms in $\sfC$ being the maps of these sets preserving (exactly, or up to an equivalence, see \cite[Section~1.16.3]{Ma}) a~certain structure (that of a~group, or of an associative algebra, or of Lie algebra, etc.), as a~model, we'd like to imitate these sets-with-a-structure by objects of another category.

For example, a~\textit{Lie supergroup} is any group in the category of supermanifolds, see~\cite{Ma}.

Likewise, a~\textit{Lie superalgebra}\index{(Lie superalgebra in the category of linear supervarieties@Lie superalgebra in the category \newline of linear supervarieties} is any Lie algebra in the category of linear supervarieties, see~\cite{KLLS}. There we reformulate the naive definition of Lie superalgebras, which are $\Zee/2$-graded linear spaces with multiplication satisfying certain identities, in terms of supervarieties.

\subsubsection{PBW-theorem for Lie superalgebras}\label{PBW}
In \cite{Et}, an interesting description of conditions when the Poincar\'e--Birkhoff--Witt theorem for Lie superalgebras holds (or not) is offered for $p>0$. Note, although we will not use this in this paper, that for Lie superalgebras understood ``naively'', the PBW theorem holds.

\subsubsection{Deformations of the brackets}\label{sssDeform}
Let $C$ be a~supercommutative superalgebra, let  $\Spec C$ be the affine super scheme.

Recall, see \cite{Ru}, where the non-super case is considered, that a~\textit{deformation} of
a Lie superalgebra $\fg$ over $\Spec C$, is a~Lie algebra $\fG$ such that $\fG\simeq\fg\otimes C$, as superspaces. The deformation is \textit{trivial} if $\fG\simeq\fg\otimes C$, as Lie superalgebras, not just as superspaces, and \textit{non-trivial} otherwise.

Generally, the \textit{deforms}\index{(Deform, the result of deformation@Deform, the result of deformation} of
a Lie superalgebra $\fg$ over $\Kee$ are Lie superalgebras $\fG\otimes_I\Kee$, where~$I$ is any closed point in $\Spec C$.

In particular, consider a~deformation with an odd parameter $\tau$. This is a~Lie superalgebra $\fG$ isomorphic to $\fg\otimes\Kee[\tau]$ as a~\textit{super space}; if, moreover, $\fG\simeq\fg\otimes\Kee[\tau]$ as a~\textit{Lie superalgebra}, i.e.,
\[
[a\otimes f, b\otimes g]=(-1)^{p(f)p(b)}[a,b]\otimes fg \quad \text{for all $a,b\in \fg$ and $f,g\in\Kee[\tau]$},
\]
then the deformation is considered \textit{trivial} (and \textit{non-trivial} otherwise).\index{(Deformation, (non)trivial@Deformation, (non)trivial} Observe that $\fg\otimes \tau$ is not an ideal of $\fG$: the ideal should be a~free $\Kee[\tau]$-module.

\textbf{Comment.} Consider formal deformations over $\Kee[[\tau]]$. If
the formal series in $\tau$ converges in a~domain $D$, we can evaluate $\tau$ for any $\tau\in D$ and~-- if $\dim \fg<\infty$~-- consider copies $\fg_\tau$, where $\tau\in D$, of the same dimension as~$\fg$. If the parameter is formal or odd, such an evaluation is possible only trivially: $\tau\mapsto 0$.

\subsection{Linear (matrix) Lie superalgebras}\label{SS:2.3} Certain basics of linear superalgebra are not well-known, or given wrongly in the literature; no harm in recalling about them.

The \textit{general linear} Lie
superalgebra of all supermatrices of size $\Size$\index{$Size$@$\Size$} corresponding to linear operators in the superspace $V=V_{\bar 0}\oplus V_{\bar 1}$ over the ground field $\Kee$ is denoted by
$\fgl(\Size)$, where $\Size=(p_1, \dots, p_{|\Size|})$ is an ordered
collection of parities of the basis vectors of $V$ for which we take only vectors \textit{homogeneous with respect to parity} and $|\Size|:=\dim V$;
usually, for the \textit{standard}\index{(Format, standard (of basis, of supermatrix)@Format, standard (of basis, of supermatrix)} (simplest from a~certain point of view) format, $\fgl(\ev,
\dots, \ev, \od, \dots, \od)$ is abbreviated to $\fgl(\dim V_{\bar
0}|\dim V_{\bar 1})$. Any supermatrix from $\fgl (\Size)$\index{$glmathfrak(\Size)$@$\fgl(\Size)$} can
be uniquely expressed as the sum of its even and odd parts; in the
standard format this is the following block expression; on non-zero summands the parity is defined:
\[
\mmat A,B,C,D,=\mmat A,0,0,D,+\mmat 0,B,C,0,,\qquad
 p\left(\mmat A,0,0,D,\right)=\ev, \qquad p\left(\mmat 0,B,C,0,\right)=\od.
\]
The \textit{supertrace}\index{(Supertrace@Supertrace}\index{$Zstr$, supertrace@$\str$, supertrace} is the map $\fgl (\Size)\tto \Kee$,
$(X_{ij})\longmapsto \sum (-1)^{p_i(p(X)+1)} X_{ii}$.
Thus, in the standard format, $\str \left(\begin{smallmatrix}A&B\\ C&D\end{smallmatrix}\right)=\tr A- \tr D$.
Observe that for the Lie superalgebra $\fgl_\cC(p|q)$ over a~supercommutative superalgebra $\cC$, i.e., for supermatrices with elements in $\cC$, we have
\begin{gather*}
\str X=\tr A- (-1)^{p(X)}\tr D\text{~~for any $X=\begin{pmatrix}A&B\\ C&D\end{pmatrix}$,}\nonumber\\
\text{where $p(X)=p(A_{ij})= p(D_{kl})=p(B_{il})+\od=p(C_{kj})+\od$}.
\end{gather*}
So if $\cC_\od\neq 0$, then on odd supermatrices the supertrace coincides with the trace.

Since $\str\ [x, y]=0$, the subsuperspace of supertraceless
matrices constitutes a~Lie subsuperalgebra of $\fgl(\Size)$ called \textit{special linear} and denoted
$\fsl(\Size)$.\index{$slmathfrak(\Size)$@$\fsl(\Size)$}

\subsubsection[The queer version of $\fgl(n)$]{The queer version of $\boldsymbol{\fgl(n)}$}\label{sssQue} There are at least two super versions of $\fgl(n)$, not
one; for reasons, see \cite[Chapters~1 and~7]{Lsos}. The other version~-- $\fq(n)$~-- is called the \textit{queer}\index{$qmathfrak(n)$, the general queer
Lie superalgebra @$\fq(n)$, the general queer
Lie superalgebra}
Lie superalgebra and is defined as the one that preserves~-- if $p\neq 2$~-- the
complex structure given by an \textit{odd} operator $J$, i.e.,
$\fq(n)$ is the centralizer $C(J)$ of $J$:
\[
\fq(n)=C(J)=\{X\in\fgl(n|n)\,|\,[X, J]=0 \}, \quad \text{where } J^2=-\id.
\]
It is clear that by a~change of basis we can reduce $J$ to the form (shape)
$J_{2n}$\index{$J_{2n}$}
in the standard
format, and then $\fq(n)$ takes the form
\begin{equation*}\label{q}
\fq(n)=\left \{(A,B):=\begin{pmatrix} A&B\\B&A\end{pmatrix}, \ \text{where $A, B\in\fgl(n)$ and $J_{2n}:=\begin{pmatrix}0&1_n\\-1_n&0\end{pmatrix}$}\right\}.
\end{equation*}
(Over any algebraically closed field $\Kee$, instead of $J$ we can take any odd operator $K$ such that $K^2=a\id_{n|n}$, where $a\in \Kee^\times$; and the Lie superalgebras $C(K)$ are isomorphic for distinct~$K$; if $p=2$, it is natural to select $K^2=\id$, and hence $\Pi_{n|n}:=\Pi_{2n}=\left(\begin{smallmatrix}0&1_n\\1_n&0\end{smallmatrix}\right)$ can serve as the normal shape of~$K$.)\index{((ZZPi or $\Pi_{k}$ or $\Pi_{n\vert n}$, matrix@$\Pi$ or $\Pi_{k}$ or $\Pi_{n\vert n}$, matrix}

On $\fq(n)$, the supertrace vanishes, but the \textit{queertrace} is defined: $\qtr\colon (A,B)\longmapsto \tr B$.\index{(Queertrace@Queertrace}\index{$Zqtr$, queertrace@$\qtr$, queertrace} Denote by $\fsq(n)$ the Lie superalgebra
of \textit{queertraceless} matrices; set $\fp\fsq(n):=\fsq(n)/\Kee 1_{2n}$.

If $p=2$, on $\mathfrak{q}(n)$ there is another (even) trace $\htr\colon (A,B)\mapsto \tr(A)$, see~\cite{KLLS}.\index{$\htr$}

\subsubsection{The supermatrix of the dual operator (after \cite{Lsos})}\label{sssStr} Let $F\in\End_{\cC}(V)$. The passage from the matrix of $F$ in a~basis of $V$ to the matrix of $F^*$ in the dual basis of $V^*$ is performed by means of the \textit{supertransposition}\index{$({}^{\rm st}$, supertransposition@${}^{\rm st}$, supertransposition}\index{(Supertransposition@Supertransposition} which in the standard format is of the shape
\[
X=\begin{pmatrix} R&S\\
T&U\end{pmatrix}\mapsto X^{\rm st}:=\begin{cases}
\begin{pmatrix} R^{\rm t}&T^{\rm t}\\
-S^{\rm t}&U^{\rm t}\end{pmatrix}&\text{if $p(X)=\ev$},\\
\begin{pmatrix} R^{\rm t}&-T^{\rm t}\\
S^{\rm t}&U^{\rm t}\end{pmatrix}&\text{if $p(X)=\od$},
\end{cases}
\]
where $M^{\rm t}$\index{(Transposition, of matrices@Transposition, of matrices} is the transposed of the matrix $M$.

\subsubsection{Lie superalgebras preserving bilinear forms}\label{sssBilMatr}
The supermatrices $X\in\fgl(\Size)$ such that
\[
X^{\rm st}B+(-1)^{p(X)p(B)}BX=0\quad \text{for an homogeneous matrix
$B\in\fgl(\Size)$}
\]
constitute the Lie superalgebra $\faut (\cB)$\index{$autmathfrac(Bcal)@$\faut (\cB)$} that preserves the
bilinear form $\cB$ on $V$ whose Gram matrix $B=(B_{ij})$ is given by the formula
\begin{gather}\label{martBil}
B_{ij}=(-1)^{p(B)p(v_i)}\cB(v_{i}, v_{j}) \quad \text{for the basis vectors $v_{i}\in V$.}
\end{gather}

In order to identify a~bilinear form $B(V, W)$ with an operator, an element of $\Hom(V, W^*)$, the matrix~$B$ of the bilinear form~$\cB$ is defined in \cite[Chapter~1]{Lsos} by equation~\eqref{martBil}, not by \textit{seemingly} natural~-- but inappropriate for such an identification~-- formula
\begin{gather}\label{wrong}
B_{ij}=\cB(v_i,v_j)\quad \text{for the basis vectors $v_i\in V$.}
\end{gather}
Moreover, the would-be definition \eqref{wrong} contradicts the manifest symmetry of the odd bilinear form $\qtr$ on $\fq(n)$. To correctly define symmetry of bilinear forms, consider the \textit{upsetting}\index{(Upsetting, of bilinear forms@Upsetting, of bilinear forms} of bilinear forms
$u\colon\Bil (V, W)\tto\Bil(W, V)$, see \cite[Chapter~1]{Lsos}, given by the formula
\begin{gather*}
u(\cB)(w,v)=(-1)^{p(v)p(w)}\cB(v,w)\quad \text{for any $v \in V$ and $w\in W$.}
\end{gather*}

If $V=W$, we say that
$\cB$ is \textit{symmetric} if \index{(Form, bilinear (anti)symmetric@Form, bilinear (anti)symmetric}
\begin{gather*}
u(B)=B,\quad \text{where $u(B)=
\mmat R^{\rm t},(-1)^{p(B)}T^{\rm t},(-1)^{p(B)}S^{\rm t},-U^{\rm t},$ for $B=\mmat R,S,T,U,$.}
\end{gather*}
Similarly, $\cB$ is \textit{anti-symmetric} if $u(B)=-B$.

\subsubsection{Notational convention}\label{ConvB}\index{(Convention@Convention} By abuse of notation we will often denote the bilinear form~$\cB$ by its Gram matrix $B$ in a~normal shape.

\subsection[Analogs of polynomials for $p>0$]{Analogs of polynomials for $\boldsymbol{p>0}$}\label{AlgDivP} Let $\Cee[x]:=\Cee[x_1,\dots,x_a]$ denote the supercommutative superalgebra of polynomials
in indeterminates $x$ in their \textit{standard order}; i.e.,\index{(Order, standard of indeterminates@Order, standard of indeterminates} let the
first~$m$ indeterminates be even and the other $n$ be odd
($m+n=a$). Among the bases of $\Cee[x]$ in which the structure
constants are integers, the two bases are usually considered: the monomial
basis and the basis of \textit{divided powers}\index{(Divided powers@Divided powers} constructed as follows.

For any multi-index $\underline{r}=(r_1, \dots , r_a)$, where
$r_1,\dots,r_m\in\Zee_+$ and $r_{m+1},\dots, r_a \in\{0,1\}$, we set\index{$Zu_i^{(r_{i})}$, $u^{(r)}$@$u_i^{(r_{i})}$, $u^{(r)}$}
\begin{gather}\label{udiv}
u_i^{(r_{i})} := \frac{x_i^{r_{i}}}{r_i!}\quad \text{and}\quad
u^{(\underline{r})} := \prod\limits_{1\leq i\leq a} u_i^{(r_{i})}.
\end{gather}

Clearly, we have
\begin{gather}
\label{divp}
u^{(\underline{r})} \cdot u^{(\underline{s})} =
\left(\mathop{\prod}\limits_{m+1\leq i\leq a}
\min(1,2-r_i-s_i)\cdot(-1)^{\mathop{\sum}\limits_{m<i<j\leq a}
r_js_i}\right)\cdot \binom {\underline{r} + \underline{s}}
{\underline{r}}
u^{(\underline{r} + \underline{s})}, \\
\text{where}\quad \binom {\underline{r} + \underline{s}}
{\underline{r}}:=\mathop{\prod}\limits_{1\leq i\leq m}\binom {r_{i}
+ s_{i}} {r_{i}}.\nonumber
\end{gather}
These $u_i^{(r_{i})}$ form an ``integer basis'' (i.e., a~basis in which all structure constants with respect to the product~\eqref{divp} are integers) of $\Cee[x]$.

\subsubsection{Notational convention}\label{divP}\index{(Convention@Convention} In what follows, for clarity, we will write exponents of divided powers in parentheses.

Over any field $\Kee$ of characteristic~$p>0$, we consider the
supercommutative superalgebra (now we do not have any elements $x$, only the $u_i^{(r_{i})}$)
\begin{gather*}
\label{u;N} \cO(m; \un| n):=\Kee[u;
\un]:=\Span_{\Kee}\left(u^{(\underline{r})}\,|\,r_i
\begin{cases}< p^{N_{i}}&\text{for $i\leq m$}\\
=0\text{ or 1}&\text{for $i>m$}\end{cases}\right),
\end{gather*}
with multiplication given by formula \eqref{divp}
where $\un = (N_1,\dots, N_m)$ is the \textit{shearing
vector}\index{(Vector, shearing@Vector, shearing} with $N_i\in\Zee_+\cup\infty$ (we assume that
$p^\infty=\infty$). Important particular cases of shearing vectors:\index{$(k$, linear form@$\One:=(1,\dots,1)$}\index{$N_\infty:=(\infty, \dots, \infty)$@$\un_\infty:=(\infty, \dots, \infty)$}
\begin{gather}\label{un_s}
\One:=(1,\dots,1) \quad \text{and} \quad \un_\infty:=(\infty, \dots,
\infty); \quad \text{we set~} \widehat\cO(m):=\cO(m; \un_\infty).
\end{gather}
The algebra \index{$Fcal$, the algebra of functions@$\cF$, the algebra of functions}\index{$Ocal(m):=\cO(m; \un_\infty\vert n)$, the completion of $\cO(m; \un\vert n)$@$\widehat\cO(m):=\cO(m; \un_\infty\vert n)$, the completion\newline of $\cO(m; \un\vert n)$}
\be\label{functions}
\cF:=\cO(m; \un| n)=\Kee[u; \un] \quad \text{and its completion $\widehat\cO(m; \un_\infty|n)$}
\ee
are called the \textit{algebras of
divided powers}.\index{$Ocal(m; \un\vert n)$, the algebra of
divided powers@$\cO(m; \un\vert n)$, the algebra of
divided powers} We will sometimes need completion $\widehat\cO\big(\widehat \un_i\big)$\index{$Ocal\big(\widehat \un_i\big)$@$\widehat\cO\big(\widehat \un_i\big)$} \index{$N_i$@$\widehat \un_i$} with respect to only one indeterminate, where $\big(\widehat \un_i\big)_j=\un_j$ except for $\big(\widehat \un_i\big)_i=\infty$.

Clearly, $\cO(a; \One)=\Kee[u; \One]$ is the algebra of truncated
polynomials. Only $\Kee[u; \One]$ is indeed generated by the
declared indeterminates whereas the list of generators of $\Kee[u; \un]$
consists of~$u_i^{(p^{k_i})}$ for all $i$ and all $k_i$ such that
$0\leq k_i<N_i$ if $u_i$ is even.

\subsection{The (generalized) Cartan prolongation}\label{CartProl} \index{(Cartan prolongation@Cartan prolongation}
Let
$\fg_-=\mathop{\oplus}\limits_{-d\leq i\leq -1}\fg_i$ be a~nilpotent
$\Zee$-graded Lie (super)algebra and $\fg_0$ a~Lie sub(super)algebra of the Lie
(super)algebra $\mathfrak{der}_0(\fg_-)$ of degree 0 derivations of $\fg_-$.
Recall that the graded Lie superalgebra
$\fb=\oplus_{p \geqslant -d}\ \fb_p$
is said to be \textit{transitive}\index{(Lie (super)algebra, transitive@Lie (super)algebra, transitive} if for all $p\geq 0$ we have
\begin{gather*}
\{x\in \fb_p \,|\,[x,\fb_{-}]=0\}=0, \quad \text{where $\mathfrak{b}_-:=\oplus_{i<0}\mathfrak{b}_i$}.
\end{gather*}

The maximal transitive $\Zee$-graded Lie (super)algebra whose non-positive part is $\fg_-\oplus \fg_0$ is called the (generalized) \textit{Cartan prolong}\index{(Cartan prolongation@Cartan prolongation} of the pair $(\fg_-,\fg_0)$ and is denoted by $(\fg_-,\fg_0)_*$.\index{$(gmathfrak_-,\fg_0)_*$, Cartan prolongation@$(\fg_-,\fg_0)_*$, Cartan prolongation}

If $p=0$, we can realize $\fg_-$ by elements of negative degree of $\fvect(n|m; \vec r)$ and $\fg_0$ by elements of 0th degree of $\fvect(n|m; \vec r)$ in a non-standard (see Section~\ref{Zgrad}) grading of $\fvect(n|m)$, where $n|m={\rm sdim} \fg_-$. Then the Cartan prolong $(\fg_-, \fg_0)_{*}:=\mathop{\oplus}\limits_{k\geq -d}\fg_{k}$
 of the pair $(\fg_-, \fg_0)$ is obtained for any $k>0$ by
\[
\fg_k:=\{D\in \fvect(n|m; \vec r)_k\,|\,[D, \fg_i]\subset \fg_{k+i}\text{~~for any $i<0$}\}.
\]
The above-described procedure is called \textit{generalized} prolongation because the initial Cartan prolongation was defined for $d=1$ only.

\subsubsection{Partial Cartan prolongation involving positive components}\label{CartProlPart}\index{(Cartan prolongation, partial@Cartan prolongation, partial} Let $\fh_1\subset \fg_1$ be
a~proper $\fg_0$-sub\-module such that $[\fg_{-1}, \fh_1]=\fg_0$. If such
$\fh_1$ exists (usually, $[\fg_{-1}, \fh_1]\subset\fg_0$), define
the 2nd prolongation of $(\mathop{\oplus}\limits_{i\leq 0}\fg_i)\oplus
\fh_1$ to be
\begin{equation*}
\label{partprol} \fh_{2}:=\{D\in\fg_{2}\,|\,[D, \fg_{-1}]\subset
\fh_1\}.
\end{equation*}
The terms $\fh_{i}$, where $i>2$, are similarly defined. Set
$\fh_i:=\fg_i$ for $i\leq 0$ and call $\fh_*:=\oplus\fh_i$ the \textit{partial Cartan prolong involving positive components}.

\textbf{Examples.} The Lie superalgebra $\fvect(1|n; n)$ is a~subalgebra of $\fk(1|2n;
n)$. The former is obtained as the Cartan prolong of the same
nonpositive part as $\fk(1|2n; n)$ and a~submodule of $\fk(1|2n;
n)_1$. The simple exceptional superalgebra $\fk\fas$ discovered in
\cite{Sh5, Sh14} is another example.

\subsection{Vectorial Lie algebras and algebras of divided powers}
The~Cartan prolong of $(\fg_-, \fg_0)$, where $\fg_0$ acts faithfully on $\fg_-$ and $\sdim\fg_- = m|n$, can be embedded into the superalgebra of polynomial vector fields of $m$ even and $n$ odd indeterminates, i.e., into $\fder\ \Cee[x_1, \dots, x_a]$ (where $a = m+n$, the first $m$ indeterminates are even, and the rest are odd), see~\cite{Shch}.

Over a~field $\Kee$ of characteristic~$p>0$, if one tries to follow the recipe of Section~\ref{CartProl} naively and use derivations of usual polynomials, instead of divided powers, it would not work. For example, let us consider the prolong $(\fg_-, \fg_0)_{*}$, where $\fg_- = \fg_{-1}$, $\sdim \fg_{-1} = \sdim \fg_0 = 1|0$, and the action of $\fg_0$ on $\fg_{-1}$ is non-trivial. It has the form $\bigoplus\limits_{i=-1}^\infty \fg_i$ such that $\sdim \fg_i = 1|0$ and $\fg_i=[\fg_{-1}, \fg_{i+1}]$ for all $i\geq -1$.

The corresponding prolong over $\Cee$ would be embedded isomorphically into $\fder\ \Cee[x]$ so that~$\fg_i$ would be mapped into $\Span\big(x^{i+1}\del_x\big)$. Over $\Kee$, the construction of embedding would fail for $\fg_{p-1}$, because $[\del_x, x^p\del_x] = 0$, and there is no element $X$ such that $[\del_x, X] = x^{p-1}\del_x$.

However, over $\Kee$, Cartan prolongs can be embedded into the superalgebra of derivations of the algebra of divided powers. Let us first say a~bit about these derivations.

Over $\Cee$, consider the action of derivation $\del_{x_i}$ of $\Cee[x_1, \dots, x_a]$ in the basis of divided powers. It is given by (recall the definition \eqref{udiv} of $u^{(r)}$)
\begin{gather}\label{distPD}
\del_{x_i} u^{(\underline{r})} =
\begin{cases}
0 &\text{if~~} r_i = 0,\\
(-1)^{\max(0, i-m-1)} u^{((r_1,\dots,r_{i-1}, r_i-1, r_{i+1},\dots, r_a))} &\text{otherwise.}
\end{cases}
\end{gather}
Since all the coefficients are integer, the map given by this formula is a~derivation of $\Kee[m;\un|n]$. We will denote this map $\del_i:=\del_{x_i} $ and call the maps $\del_1, \dots, \del_a$ \textit{distinguished} partial derivatives.\index{(Distinguished partial derivative@Distinguished partial derivative}

The \textit{general Lie algebra of vector fields}
consists of the following derivations expressed in terms of \textit{distinguished} partial derivatives\index{$vectmathfrak(m;\un\vert n)=W(m;\un\vert n)$@$\fvect(m;\un\vert n)=W(m;\un\vert n)$}
\[
\fvect(m;\un|n)=\left\{\sum_{1\leq i\leq a} f_i\partial_i\,|\,f_i\in \cO(m;
\un|n)\right\}.
\]

Note that if $\un\neq \One$, then $\fvect(m;\un|n)$ is not the whole $\fder\, \Kee[m;\un|n]$. Maps $\del_i^{p^k}$, where $1\leq i\leq m$ and $1\leq k<N_i$, are also derivations of $\Kee[m;\un|n]$, and a~general derivation of $\Kee[m;\un|n]$ has the form\footnote{It is easy to see that $\Kee[m;\un|n]$ is isomorphic to $\Kee[\sum N_i;\One|n]$, so their algebras of derivations are isomorphic as well, so it is not surprising that a~general element of $\fder\, \Kee[m;\un|n]$ has $\sum N_i +n$ functional parameters.}
\[
\sum\limits_{1\leq i\leq m}\ \sum\limits_{0\leq k\leq N_i-1} f_{i,k}\del_i^{p^k} + \sum\limits_{m+1 \leq i \leq a} f_{i}\del_i, \quad \text{where $f_{i,k}, f_i\in \Kee[m;\un|n]$.}
\]

The Lie superalgebra
$\fvect(m;\un|n)$ and its subalgebras, are called \textit{vectorial} Lie
superalgebras\index{(Lie superalgebra, vectorial@Lie superalgebra, vectorial} (cf.\ with \textit{matrix} or \textit{linear} Lie superalgebras). Cartan prolongs can be embedded into $\fvect(m;\un|n)$; in particular, the above Cartan prolong would be isomorphic to $\fvect(1;\un_\infty|0)$, with $\fg_i$ corresponding to $\Span\big(u_1^{(i+1)}\del_1\big)$.

\subsubsection{Notation, again} Hereafter, the symbol $\fg(a|b)$ or $\fg(a;\un |b)$\index{$gmathfrak(a\vert b)$ or $\fg(a;\un \vert b)$, a~vectorial Lie superalgebra of type $\fg$@$\fg(a\vert b)$ or $\fg(a;\un \vert b)$, a~vectorial Lie\newline superalgebra of type $\fg$} will designate the vectorial Lie superalgebra
with given name $\fg$ realized by vector fields on the linear supermanifold $\cK^{a|b}$ (the one corresponding to the superspace $\Kee^{a|b}$, see~\cite{KLLS}), and endowed with a~W-grading, see Table~\eqref{nonstandgr}. The \textit{standard}\index{(Grading, standard@Grading, standard} grading is taken as a~point of reference for regradings governed
by the vector $\vec r$\index{$Zr$ or $\vec r$, the vector of degrees @$r$ or $\vec r$, the vector of degrees} of degrees, which often can be described by one number $r$ that usually (for details, see \cite{LS1, Sh14})
is equal to the number of odd indeterminates of degree 0. The regraded
Lie superalgebra is denoted by~$\fg(a|b;r)$. In the standard grading
the parameter $r$ is usually omitted, see Table~\eqref{nonstandgr} and tables in Section~\ref{ssExcAsCTS}.

The module $\cF$ of ``functions''\index{$Fcal$, the algebra of functions@$\cF$, the algebra of functions} over $\fvect(m;\un|n)$\index{$vectmathfrak(m;\un\vert n)$@$\fvect(m;\un\vert n)$} and its
subalgebras (usually with the same negative part) is an analog of the
\textit{tautological} module $V$ over $\fgl(V)$ and its subalgebras.

\subsubsection{Names} The Lie algebra $\fvect(1;\un)$ is called a~\textit{Zassenhaus algebra}. For $p=2$ it is not simple. Observe that
$\fvect(1;\un)\simeq\fk(1;\un)$ (indeed, $f\partial_x\longleftrightarrow K_f$, see definition~\eqref{K_f}; clearly, $\partial_x$~is the distinguished derivative with respect to the only indeterminate). The simple derived algebra $\fvect^{(1)}(1;\un)\simeq\fk^{(1)}(1;\un)$ is also called a~\textit{Zassenhaus algebra}\index{(Zassenhaus algebra@Zassenhaus algebra} causing confusion, while $\fvect(1;\One)$ is lately called (even for $p=0$) the
\textit{Witt algebra} in honor of Witt\index{$W(m;\un)$, Witt algebra}\index{(Witt algebra@Witt algebra} who was the first to study one of its modular incarnations, see Introduction to the first volume of~\cite{S}.

In the old literature, $\fvect(m;\un)$, like its version for $p=0$, was called \textit{the general Lie algebra of Cartan type}; lately, it is called the \textit{Jacobson--Witt algebra},\index{(Jacobson--Witt algebra@Jacobson--Witt algebra} whereas the name \textit{Witt algebra}\index{(Witt algebra@Witt algebra} is reserved for the particular case $\fvect(1;\underline{1})$ for $p>2$.

\subsection{Traces and divergencies on vectorial Lie superalgebras}\label{traceDiv}
On any Lie (super)algebra~$\fg$ over a~supercommutative superalgebra $\cC$, e.g., over a~field
$\cC=\Kee$, a~\textit{trace}\index{(Trace@Trace}\index{$Ztr$, trace@$\tr$, trace}\index{$Zstr$, supertrace@$\str$, supertrace} is any linear mapping $\tr\colon \fg\tto \cC$ such that
\begin{equation}\label{deftr}
\tr \big(\fg^{(1)}\big)=0.
\end{equation}

Now let $\fg$ be a~$\Zee$-graded vectorial Lie (super)algebra with
$\fg_{-}:=\mathop{\oplus}\limits_{i<0}\fg_i$ generated by
$\fg_{-1}$, and let $\tr$ be a~trace on $\fg_0$. Recall that any
$\Zee$-grading of a~given vectorial Lie (super)algebra is given by degrees
of the indeterminates, so the space of functions $\cF$, see equation~\eqref{functions}, is also
$\Zee$-graded.

The \textit{divergence}\index{(Divergence@Divergence}\index{$Div$, divergence@$\Div$, divergence} $\Div\colon \fg\tto\cF$ is a~degree-preserving
$\ad_{\fg_{-1}}$-invariant extension of the trace to the Cartan prolong; this extension should satisfy the following conditions,
so $\Div\in Z^1(\fg; \cF)$, i.e., is a~cocycle:
\begin{gather*}
X_i(\Div D)=\Div([X_i,D])\text{~~for all elements $X_i$ that span $\fg_{-1}$},\\
\Div|_{\fg_0}=\tr,\\
\Div|_{\fg_{-}}=0.
\end{gather*}

We denote by $\Vol(u;\un)$ or simply
$\Vol_u:=\cF^*$ the $\fvect(m;\un|n)$-module of \textit{volume forms}\index{$Vol(u;\un)$, module of volume forms@$\Vol(u;\un)$, module of volume forms} dual to $\cF$ over $\cF$. As an $\cF$-module, $\Vol_u$ is generated by the \textit{volume element} $\vvol_u=1^*$ with
fixed indeterminates (``coordinates") $u$ which we often do not
indicate. On the rank-1 $\cF$-module of \textit{weighted}\index{(Density, weighted@Density, weighted}
$\lambda$-densities $\Vol^{\lambda}(m;\un|n)$ with generator
$\vvol_u^{\lambda}$\index{$Vol^{\lambda}(m;\un\vert n)$, module of weighted densities@$\Vol^{\lambda}(m;\un\vert n)$, module of weighted densities} over $\cF$, the $\fvect(m;\un|n)$-action is given for any $f\in\cF$ and $D\in\fvect(m;\un|n)$ by the
\textit{Lie derivative}
\begin{equation}\label{LieDer}
L_D\big(f\vvol_u^{\lambda}\big)=\big(D(f)+(-1)^{p(f)p(D)}\lambda f\Div(D)\big)\vvol_u^{\lambda}.
\end{equation}

The \textit{special} Lie algebra
$\fs\fg:=\Ker\Div$ of \textit{divergence-free} elements of
$\fg$ is the Cartan prolong of $(\fg_{-}, \Ker\tr|_{\fg_0})$. For example,
$\fsvect(m;\un|n)=(\id_{\fsl(m|n)}, \fsl(m|n))_{*,\un}$.

The $\fvect(0|n)$-module $\Vol(0|n)$ contains a~submodule $\Vol_0(0|n)$\index{$Vol_0(0\vert n)$@$\Vol_0(0\vert n)$} of
codimension~1:
\begin{equation}
\label{Vol_0} \Vol_{0}(0|n):=\left\{f\vvol\,|\,\int f\vvol=0,\  \text{where $f\in\cF=\Lambda(n)$}\right\},
\end{equation}
where the \textit{Berezin integral}\index{(Berezin integral@Berezin integral} $\int f\vvol_\theta$ is equal to the coefficient of the monomial of the highest in~$\theta$s degree.

Over $\fsvect(0|n)$, we often identify $\Vol_0(0|n)$ with a~submodule of $\cF$ and omit $(0|n)$; set \index{$T_0^0(0\vert n)$}
\begin{equation*}
\label{T_0^0}
T_0^0(0|n):=\Vol_0(0|n)/\Kee\cdot \vvol.
\end{equation*}

\subsubsection{Examples of several divergences}\label{sssDiv} On $\fvect(m;\un|n)$, the explicit expression of the \textit{standard} divergence\index{(Divergence, standard@Divergence, standard} is as follows
\begin{equation}\label{DivVect}
\Div\colon \ \sum f_i \del_{x_i}\tto \sum (-1)^{p(x_i)p(f_i)} \del_{x_i} f_i.
\end{equation}

The supertrace restricted from $\fgl$ vanishes on $\fq$, but there is an ``indigenous'' queer trace on $\fq$; analogously, the standard divergence~\eqref{DivVect} vanishes on certain Lie subsuperalgebras of $\fvect(m;\un|n)$ on which there might be defined an ``indigenous'' divergence.
This happens, e.g., with $\fk(2n+1|2n+2)$ and $\fm(n)$ as will be shown later on.

If there are several traces on $\fg_0$, and
hence divergences on $\fg=(\fg_-, \fg_0)_{\ast, \un}$, there are several types of special
subalgebras, and we need an individual name for each.

If $\fg$ is a~Lie \textit{super}algebra, then the linear functional $\tr$
satisfying condition \eqref{deftr} is often called, for emphasis,
\textit{super}trace and denoted by $\str$. If we were consistent, we
should, accordingly, use the term \textit{super}divergence but
instead we drop the preface ``super'' in both cases.

\subsection{Critical coordinates and unconstrained shearing vectors}\label{CritCoor} \index{(Coordinate, critical@Coordinate, critical}The coordinate of the shearing vector $\un$
corresponding to an even indeterminate of the $\Zee$-graded vectorial Lie (super)algebra~$\fg$ is said to
be \textit{critical} if it cannot take an arbitrarily big value.

The shearing vector without any imposed restrictions on its
coordinates is said to be \textit{unconstrained}; we denote
it by $\un^u$.\index{$N^u$@$\un^u$} Let \index{$dim\un$@$\dim\un$}\index{$Par\un^u$@$\Par\un^u$}
\begin{gather*}
\text{$\dim\un$} \ \text{be the number of coordinates of $\un$},\nonumber\\
\text{$\Par\un^u$} :=\dim\un^u-\card(\{\text{critical coordinates})\}\nonumber\\
\hphantom{\text{$\Par\un^u$} :=}{} \text{be the number of parameters $\un^u$ depends on}.
\end{gather*}

We established the (non)critical
coordinates of the shearing vectors of the $\Zee$-graded vectorial
Lie (super)algebra $\fg$ with a~computer's aid by
explicitly computing the bases of the first several
terms~$\fg_i$ for $i\geq 0$ without imposing any constraints on $\un$.

\begin{Conjecture}\label{CritCoorConj}\index{(Conjecture@Conjecture} If the value of the coordinate $\un_i$ \textup{(of the $\Zee$-graded vectorial Lie
(super)algebra $\fg$)} can be $>1$, then it can take any value.
\end{Conjecture}

\section{Background continued. Subtleties}\label{SbackSub}

\subsection[The serial simple vectorial Lie superalgebras over $\Cee$ as prolongs]{The serial simple vectorial Lie superalgebras over $\boldsymbol{\Cee}$ as prolongs}\label{sssSerialProl} When we only need the vectorial Lie
superalgebras considered as abstract, not realized by vector fields,
we may consider their simplest filtrations with the smallest
codimension of their maximal subalgebras, and gradings associated with
such filtrations.

\subsubsection[Convention: on central element $z\in\fg_0$]{Convention: on central element $\boldsymbol{z\in\fg_0}$}\label{Conv4.2.1}\index{(Convention@Convention} We chose the central element $z\in\fg_0$
so that it acts on $\fg_i$ as $i\cdot\id$. The irreducible 1-dimensional module over the commutative Lie algebra spanned by $z$ which acts as $i\cdot\id$ is denoted by $\Kee[i]$.\index{$Kmathbb[i]$@$\Kee[i]$}

\subsection[The two types of superizations of the contact series over $\Cee$]{The two types of superizations of the contact series over $\boldsymbol{\Cee}$}\label{ssTwoTypesKM}{}~{}

\underline{The type $\fk$}: Define the Lie superalgebra
$\fhei(2n|m)$\index{$heimathfrak(2n\vert m)$@$\fhei(2n\vert m)$} on the direct sum of a~$(2n, m)$-dimensional
superspace $W$, endowed with a~non-degenerate antisymmetric bilinear
form $B$, and a~$(1, 0)$-dimensional space spanned by $z$. Clearly,
we have
\begin{equation*}
\label{hei} \fk(2n+1|m)=(\fhei(2n|m), \fc\fosp(m|2n))_*
\end{equation*}
and, given $\fhei(2n|m)$ and a~subalgebra $\fg$ of
$\fc\fosp(m|2n)$, we call $(\fhei(2n|m), \fg)_*$ the
\textit{$k$-prolong} of $(W, \fg)$, where $W$ is the tautological
$\fosp(m|2n)$-module.

\underline{The type $\fm$}: The ``odd'' analog of $\fk$ is
associated with the following ``odd'' analog of $\fhei(2n|m)$.
Denote by $\fba(n)$\index{$bamathfrak(n)$@$\fba(n)$} the \textit{antibracket} Lie superalgebra
($\fba$ is Anti-Bracket read backwards). Its space is $W\oplus
\Cee\cdot z$, where $W$ is an $n|n$-dimensional superspace endowed
with a~non-degenerate antisymmetric odd bilinear form $B$; the
bracket in $\fba(n)$ is given by the following relations:
\begin{equation*}
\label{defBAneHei} \text{$z$ is odd and lies in the center; $[v,
w]=B(v, w)\cdot z$ for any $v, w\in W$.}
\end{equation*}

Given $\fba(n)$ and a~subalgebra $\fg$ of
$\fc\fpe(n)$, we call $(\fba(n), \fg)_*$ the
\textit{$m$-prolong} of $(W, \fg)$, where~$W$ is the tautological
$\fpe(n)$-module.

\subsection[Generating functions over $\Cee$]{Generating functions over $\boldsymbol{\Cee}$}\label{GenFun} A~laconic way
of describing $\fk$, $\fm$ and their subalgebras is via generating functions.

On the $2n+1|m$-dimensional superspace
with even coordinates $t$, and $p=(p_1,\dots, p_n)$, $q=(q_1,\dots, q_n)$, and odd indeterminates (``coordinates'') $(\xi, \eta, \theta)$, where $\xi=(\xi_1,\dots, \xi_k)$, $\eta=(\eta_1,\dots, \eta_k)$ and $\theta=(\theta_1,\dots, \theta_s)$, the \textit{odd contact form $\alpha_1$}\index{((ZZZalpha_1, odd contact form@$\alpha_1$, odd contact form} is defined to be
\begin{gather}
\label{alpha1} \alpha_1=dt+\mathop{\sum}\limits_{1\leq i\leq
n}(p_idq_i-q_idp_i)+ \mathop{\sum}\limits_{1\leq j\leq
k}(\xi_jd\eta_j+\eta_jd\xi_j)+
\begin{cases}0&
\text{if}\ m=2k,\\
\mathop{\sum}\limits_{1\leq\ell\leq s}\theta_\ell d\theta_\ell&\text{if}\ m=2k+s.\end{cases}\!\!\!\!\!
\end{gather} For any $f\in\Cee [t, p, q, \xi, \eta,
\theta]$, set\index{$K_f$}
\begin{equation}
\label{K_f} K_f=(2-E)(f)\pder{t} + H_f + \pderf{f}{t} E,
\end{equation}
where $E=\mathop{\sum}\limits_i y_i \pder{y_{i}}$ (here the $y_{i}$
are all the coordinates except $t$) is the \textit{Euler
operator},\index{$E$, Euler operator} and
\begin{equation*}
\label{H_f} H_f=\mathop{\sum}\limits_{i\leq n}\left(\pderf{f}{p_i}
\pder{q_i}-\pderf{f}{q_i} \pder{p_i}\right )
-(-1)^{p(f)}\left(\mathop{\sum}\limits_{\ell\leq s}\pderf{ f}{\theta_\ell}
\pder{\theta_\ell}+ \mathop{\sum}\limits_{j\leq k}\left(\pderf{f}{\xi_j}
\pder{\eta_j}+ \pderf{f}{\eta_j} \pder{\xi_j}\right)\right ).
\end{equation*}
The \textit{Hamiltonian vector
field} $H_f$\index{$H_f$, Hamiltonian vector
field} with Hamiltonian $f$ preserves the symplectic form\index{((ZZZomega_0, symplectic form@$\omega_0$, symplectic form} $\omega_0:=d \alpha_1$.

On the $(n|n+1)$-dimensional\index{(Vector field, (peri)contact@Vector field, (peri)contact}
superspace with even coordinates $q=(q_1,\dots, q_n)$, and odd indeterminates (``coordinates'')
$\xi=(\xi_1,\dots, \xi_n)$, and $\tau$, the \textit{even contact form $\alpha_0$},\index{((ZZZalpha_0, even contact form@$\alpha_0$, even contact form} is defined to be
\begin{equation}\label{alp0}
\alpha_0=d\tau+\mathop{\sum}\limits_{i\leq n}(\xi_idq_i+q_id\xi_i).
\end{equation} For any $f\in\Cee [q,
\xi, \tau]$, set\index{$M_f$}\index{$Le_f$@$\Le_f$}
\begin{gather*}
 M_f=(2-E)(f)\del_{\tau}- \Le_f -(-1)^{p(f)} \del_{\tau}(f) E,\quad
\text{where $E=\mathop{\sum}\limits_iy_i \del_{y_i}$ and
where $y=(q,\xi)$},\\
  \Le_f=\mathop{\sum}\limits_{i\leq n}\left( \pderf{f}{q_i}\
\pder{\xi_i}+(-1)^{p(f)} \pderf{f}{\xi_i}\ \pder{q_i}\right ).
\end{gather*}
Since
\begin{gather}
 L_{K_f}(\alpha_1)=2 \pderf{f}{t}\alpha_1=K_1(f)\alpha_1, \nonumber\\
L_{M_f}(\alpha_0)=-(-1)^{p(f)}2 \pderf{
f}{\tau}\alpha_0=-(-1)^{p(f)}M_1(f)\alpha_0.\label{LieK_f}
\end{gather}
Let $\fk (2n+1|m)$\index{$kmathfrak (2n+1\vert m)$@$\fk (2n+1\vert m)$} be the (contact) Lie superalgebra preserving the distribution given by the Pfaff equation with the 1-form $\alpha_1$; let
$\fm(n):=\fm(n+1|n)$\index{$mmathfrak(n):=\fm(n+1\vert n)$, pericontact Lie superalgebra@$\fm(n):=\fm(n+1\vert n)$, pericontact Lie\newline superalgebra} be the pericontact (``odd'' contact) Lie superalgebra preserving the distribution given by the Pfaff equation with the 1-form $\alpha_0$.

Equation~\eqref{LieK_f} implies that $K_f\in \fk (2n+1|m)$ and $M_f\in \fm (n)$. Observe that
\[
p(\Le_f)=p(M_f)=p(f)+\od.
\]

\textbf{Contact brackets, Poisson bracket, antibracket a.k.a.\ Buttin (Schouten) bracket.} To the
(super)com\-mu\-tators $[K_f, K_g]$ and $[M_f, M_g]$ there correspond
\textit{contact brackets}\index{$(\{- , -\}_{\rm k.b.}$, the contact bracket@$\{- , -\}_{\rm k.b.}$, the contact bracket}\index{$(\{- , -\}_{\rm m.b.}$, the (peri)contact bracket@$\{- , -\}_{\rm m.b.}$, the (peri)contact bracket} of the generating functions:
\begin{gather*}
 [K_f, K_g]=K_{\{f, \; g\}_{\rm k.b.}},\\
 [M_f, M_g]=M_{\{f, \; g\}_{\rm m.b.}}.
\end{gather*}
The explicit expressions for the contact brackets are as follows.
Let us first define the brackets on functions that do not depend on~$t$ (resp.~$\tau$).

The \textit{Poisson bracket} $\{- , -\}_{\rm P.b.}$\index{$(\{- , -\}_{\rm P.b.}$, the Poisson bracket@$\{- , -\}_{\rm P.b.}$, the Poisson bracket} (in the
realization with the form
 $\omega_0:=d\alpha_1$ for
$m=2k+1$ it is given by the formula
\begin{gather*}
\{f, g\}_{\rm P.b.}=\mathop{\sum}\limits_{i\leq n}
\left( \pderf{f}{p_i} \pderf{g}{q_i}- \pderf{f}{q_i}
\pderf{g}{p_i}\right)\\
{} -
(-1)^{p(f)}\left( \mathop{\sum}\limits_{\ell\leq s}\pderf{ f}{\theta_\ell}
\pder{\theta_\ell}+\mathop{\sum}\limits_{j\leq k}\left(
\displaystyle\pderf{f}{\xi_j} \pderf{ g}{\eta_j}+\pderf{f}{\eta_j}
\pderf{ g}{\xi_j}\right)\right)\quad \text{for any }f, g\in \Cee [p, q, \xi, \eta,
\theta].
\end{gather*}

The \textit{Buttin\footnote{C.~Buttin was the first to publish that
the Schouten bracket satisfies super Jacobi identity.} bracket}\index{(Buttin bracket@Buttin bracket}\index{$(\{- , -\}_{\rm B.b.}$, the Buttin bracket, a.k.a. antibracket, Schouten bracket@$\{- , -\}_{\rm B.b.}$, the Buttin bracket, a.k.a.\newline antibracket, Schouten bracket}\index{(Antibracket@Antibracket}\index{(Schouten bracket@Schouten bracket}
$\{- , -\}_{\rm B.b.}$, discovered by Schouten and initially
known as the \textit{Schouten bracket}, is very popular in physics
under the name \textit{antibracket}, see \cite{GPS}. It is given by
the formula
\begin{gather}\label{BB}
\{f, g\}_{\rm B.b.}=\mathop{\sum}\limits_{i\leq n}
\bigg(\pderf{f}{q_i} \pderf{g}{\xi_i}+(-1)^{p(f)} \pderf{f}{\xi_i}
\pderf{g}{q_i}\bigg) \quad \text{for any }f, g\in \Cee [q, \xi].
\end{gather}

In terms of the Poisson and Buttin brackets, respectively, the
contact brackets are as follows:
\begin{gather}
\{f, g\}_{\rm k.b.}= (2-E) (f)\pderf{g}{t}-\pderf{f}
{t}(2-E) (g)-\{ f, g\}_{\rm P.b.}, \nonumber\\
\{ f, g\}_{\rm m.b.}= (2-E) (f)\pderf{g}{\tau}+(-1)^{p(f)}
\pderf{f}{\tau}(2-E) (g)-\{ f, g\}_{\rm B.b.}.\label{KB}
\end{gather}

The Lie superalgebras of \textit{Hamiltonian vector
fields}\index{(Hamiltonian vector fields@Hamiltonian vector fields} (or \textit{Hamiltonian
superalgebras}) and their special subalgebras (defined only if $n=0$)
are\index{$hmathfrak (2n\vert m)$@$\fh (2n\vert m)$}\index{$hmathfrak^{(1)} (0\vert m)$@$\fh^{(1)} (0\vert m)$}
\begin{gather*}
\begin{split}
&\fh (2n|m)=\{D\in \fvect (2n|m)\,|\,L_D\omega_0=0, \text{~~where
$\omega_0:=d\alpha_1$}\},\\
&\fh^{(1)} (0\vert m) =\left\{H_f\in \fh (0|m)\,|\,\int f\vvol=0\right\}.
\end{split}
\end{gather*}
The ``odd'' analogues of the Lie superalgebra of Hamiltonian fields
are the Lie superalgebra of vector fields $\Le_{f}$ introduced in
\cite{Le2}, and its special subalgebra: \index{$lemathfrak(n)$@$\fle(n)$}\index{$slemathfrak(n)$@$\fsle(n)$}
\begin{gather}
\fle (n)=\{D\in \fvect (n|n)\,|\,L_D\omega_1=0, \
\text{where $\omega_1:=d\alpha_0$}\},\nonumber\\
\fsle (n)=\{D\in \fle (n)\,|\,\Div D=0\}.\label{le_n}
\end{gather}

It is not difficult to prove the following isomorphisms as Lie
superalgebras with the brackets on the right-hand sides given by the above-described brackets k.b.\ (where $K_f$ and $H_f$ are involved) and m.b.\ (where $M_f$ and $Le_f$ are involved) \index{$bmathfrak(n)$@$\fb(n)$}
\begin{gather*}
\fk (2n+1|m)=\{K_f\,|\,f\in \Cee[t, p, q, \xi]\} \cong \Cee[t, p, q, \xi],\\
\fpo (2n|m):=\left\{K_f\,|\,f\in \Cee[t, p, q, \xi]
\  \text{such that} \ \pderf f t=0\right\}\cong \Cee[p, q, \xi],\\
\fh (2n|m)=\{H_f\,|\,f\in \Cee [p, q, \xi]\}\simeq\Cee[p,
q,\xi]/\Cee\cdot 1,\\
\fm (n)=\{M_f\,|\,f\in \Cee [\tau, q, \xi]\}\cong \Cee [\tau, q, \xi],\\
\fb (n):=\left\{M_f\,|\,f\in \Cee [\tau, q, \xi),
\  \text{such that} \ \pderf f \tau=0\right\}\cong\Pi(\Cee[q, \xi]),\\
\fle(n) = \{\Le_f\,|\,f\in \Cee [q, \xi]\}\cong\Pi(\Cee[q,
\xi]/\Cee\cdot 1).
\end{gather*}
We have\index{$pomathfrak^{(1)} (0\vert m)$@$\fpo^{(1)} (0\vert m)$}\index{$hmathfrak^{(1)} (0\vert m)$@$\fh^{(1)} (0\vert m)$}
\begin{gather*}
\fpo^{(1)} (0|m)=\left\{ K_f\in \fpo (0|m)\,|\,\int f\vvol =0\right\},\\
\fh^{(1)} (0|m)=\fpo^{(1)} (0|m)/\Cee\cdot K_1.
\end{gather*}

\subsubsection[Generating functions over $\Kee$]{$\!\!\!$Generating functions over $\boldsymbol{\Kee}$}\label{sssGenF} Recall that the contact Lie
superalgebra $\fk(n_\ev+1; \un |n_\od)$ consists of the vector fields $D$ that preserve the
contact structure (a~non-integrable distribution given by a~contact
form $\alpha_1$, cf.~\eqref{alpha1}) on the supervariety $\cM$
associated with the superspace~$\Kee^{n_\ev+1|n_\od}$:
\begin{gather}
\label{L_D}
L_D(\alpha_1)=F_D\alpha_1\ \text{for some $F_D\in\cF$,
where $\cF$ is the space of functions on $\cM$.}
\end{gather}
Consider the form (here $n_\ev=2k_\ev$; if $n_\ev$ is odd, no contact
form exists)
\begin{gather}
\alpha_1=dx_0+\mathop{\sum}\limits_{1\leq i\leq k}
x_idx_{k+i}\nonumber\\
\hphantom{\alpha_1=}{} + \begin{cases}0&\text{if
$n=(n_\ev+1)+n_\od=(2k_\ev+1)+2k_\od=2k+1$}, \\
x_{2k+1}dx_{2k+1}&\text{if $n=(n_\ev+1)+n_\od=(2k_\ev+1)+(2k_\od+1)=2k+2$},
\end{cases}\label{alpha_1}
\end{gather}
where in order to make expressions for brackets simpler, we consider the following nonstandard order of indeterminates
the constituents of dual pairs one above/under the other:
\[
\begin{array}{@{}llll}
\text{even}\colon \ x_0,& \text{even}\colon \ x_1, \phantom{k_\ev+k}\dots, x_{k_\ev},&\text{odd}\colon \ x_{k_\ev+1},\dots,x_{k_\ev+k_\od},&\\
& \text{even}\colon \ x_{k_\ev+k_\od+1},
\dots, x_{2k_\ev+k_\od},& \text{odd}\colon \
 x_{2k_\ev+k_\od+1},\dots,x_{2k}; &\text{odd}\colon \ x_{2k+1}.
\end{array}
\]
The vector fields $D$ satisfying \eqref{L_D} for some function
$F_D$ look differently for different characteristics:

\textit{For $p\neq 2$, and also if $p=2$ and $n=2k+1$}, the fields $D$
satisfying \eqref{L_D} have, for any $f\in\cF$, the following form
(compare with \eqref{K_f}):\index{$K_f$}
\begin{gather}
K_f=(1-E')(f)\displaystyle\pder{x_0} + \pderf{f}{x_0}
E'+\mathop{\sum}\limits_{1\leq i\leq k_\ev} \left(\frac{\del f}{\del
x_{k+i}}\pder{x_{i}}-\frac{\del f}{\del x_i}\pder{x_{k+i}}\right)\label{K_f1}\\
-(-1)^{\Pi(f)}\displaystyle\left(\mathop{\sum}\limits_{k_\ev+1\leq
i\leq k}\left(\frac{\del
f}{\del x_{k+i}}\pder{x_{i}}+\frac{\del f}{\del x_i}\pder{x_{k+i}}\right)+\begin{cases}0&\text{if
$n=2k+1$}\\
\frac12\frac{\del f}{\del x_{2k+1}}\pder{x_{2k+1}}&\text{if
$n=2k+2$}\end{cases}\right),\nonumber
\end{gather}
where \index{$E'$, Euler operator}
\[
E':= \sum _{1\leq i\leq k} x_i
\del_{x_{i}}+\begin{cases}0&\text{if $n=2k+1$},\\
\frac12 x_{2k+1} \del_{x_{2k+1}}&\text{if $n=2k+2$}.\end{cases}
\]

\textit{For $p=2$ and $n=2k+2$}, we cannot use formula~\eqref{K_f1}
anymore (at least, not for arbitrary~$f$) since it contains
$\frac12$. In this case, the elements of the contact algebra are of
the following three types, and their linear combinations, where $k=k_\ev+k_\od$:

a) For any $f\in\cF$ such that $\pderf{f}{x_0}=\pderf{f}{x_{2k+1}}=0$,
we have \index{$K_f$ for $p=2$}
\begin{equation}\label{K_f2}
 K_f=(1+E')(f)\pder{x_0} + \sum_{1\leq i\leq k} \left(\frac{\del f}{\del x_i}\pder{x_{k+i}}+\frac{\del f}{\del
x_{k+i}}\pder{x_{i}}\right),
\end{equation}
where
$E':=\mathop{\sum}\limits_{1\leq i\leq k} x_i \pder{x_{i}}$.

\textit{Observe} that in \eqref{K_f1} and \eqref{K_f2} we can also take $E':=\mathop{\sum}\limits_{k+1\leq i\leq 2k} x_i \pder{x_{i}}$.

b) Set
\index{$Fcal$@$\widehat\cF$}
\begin{equation}\label{hatF}
\widehat\cF:=\mathcal{O}(x_0, x_1,\dots,x_{k_\ev}, x_{k+1},\dots,x_{k+k_\ev}; \un|x_{k_\ev+1}, \dots, x_{k}, x_{k+k_\ev+1}, \dots, x_{2k}).
\end{equation} For any $g\in\widehat\cF$, or equivalently, for any $g\in\cF$ such that $\pderf{g}{x_{2k+1}}=0$, we set\index{$A_g$}\index{$B_g$}
\begin{align*}
\text{b1)}\quad &A_g:=g\left(x_{2k+1}\del_{x_0}+ \del_{x_{2k+1}}\right),\\
\text{b2)}\quad &B_g:=gx_{2k+1}\del_{x_{2k+1}}.
\end{align*}

For the pericontact Lie superalgebra $\fm(n;\un|n)$ the analog of the formula \eqref{alpha_1} takes the form
\begin{equation*}
\label{alpha_0} \alpha_0=dx_0+\mathop{\sum}\limits_{1\leq i\leq k} x_idx_{k+i},
\end{equation*}
where the parities of indeterminates are such that $p(x_i)=p(x_{i+k})+\od$; e.g., they are as follows
\[
 \text{even}\colon \ x_1, \dots, x_{k}; \qquad \text{odd}\colon \ x_{k+1},\dots,x_{2k}, \text{~and~~} x_0.\\
\]

\begin{Theorem}[on explicit squaring and contact brackets for $p=2$]\label{brackets} Brackets and squares of contact vector fields, and the corresponding contact brackets of generating functions, are given by formulas~\eqref{K_fbr}. 
Both the contact brackets $\{-, -\}_{\rm k.b.}$ and $\{-, -\}_{\rm m.b.}$ are of the shape
\begin{gather}\{f,f_1\}=\pderf{f}{x_0}(1+E')(f_1)+(1+E')(f)\pderf{f_1}{x_0}+
\mathop{\sum}\limits_{1\leq i\leq k}
\left(\pderf{f}{x_i}\pderf{f_1}{x_{k+i}}+
\pderf{f}{x_{k+i}}\pderf{f_1}{x_{i}}\right),\nonumber\\
\text{where~~$E':=\mathop{\sum}\limits_{1\leq
i\leq k} x_i \pder{x_{i}}$ or $E':=\mathop{\sum}\limits_{k+1\leq
i\leq 2k} x_i \pder{x_{i}}$ for any $f, f_1\in \widehat\cF$}.\label{Konbr}
\end{gather}
Then, for any $f, f_1, g, g_1\in \widehat\cF$, we deduce
\begin{alignat}{3}
& [K_f,K_{f_1}]=K_{\{f,f_1\}_{\rm k.b.}},\qquad && [M_f,M_{f_1}]=M_{\{f,f_1\}_{\rm m.b.}}, &\nonumber\\
& (K_f)^2=K_{\mathop{\sum}\limits_{1\leq i\leq
k}\pderf{f}{x_i}\pderf{f}{x_{k+i}}},\qquad &&
  (M_f)^2=M_{\mathop{\sum}\limits_{1\leq i\leq
k}\pderf{f}{x_i}\pderf{f}{x_{k+i}}}, &\nonumber\\
& [K_f,A_{g}]=A_{\{f,g\}_{\rm k.b.}},\qquad &&(A_g)^2=
B_{g\pderf{g}{x_0}}+
g^2K_1, &\nonumber\\
& [K_f,B_{g}]=B_{\{f,g\}_{\rm k.b.}},\qquad &&[A_g,B_{g_1}]=A_{gg_1}, &\nonumber\\
& [A_g,A_{g_1}]=B_{\pderf{(gg_1)}{x_0}},\qquad &&[B_g,B_{g_1}]=(B_g)^2=0. & \label{K_fbr}
\end{alignat}
\end{Theorem}

\begin{proof}[Proof {\normalfont is based on direct computations, see formulas}~\eqref{K_fbr}]\end{proof}

In particular, for $\fk(1;\un|0)$, we have (unless $a=b=0$)
\begin{equation*}\label{K(1vert 0}
\{x^{(a)}, x^{(b)}\}=\left(\binom{a+b-1}{a-1}+\binom{a+b-1}{b-1}\right)x^{(a+b-1)}=\binom{a+b}{a}x^{(a+b-1)}.
\end{equation*}

\begin{Lemma}[A helpful lemma]\label{long} For any $g\in\widehat\cF$, see~\eqref{hatF}, we have $g^2\in\Kee$, see the expressions for $(A_g)^2$ and $[A_g, B_g]$ in equation~\eqref{K_fbr}.
\end{Lemma}

\begin{proof} Indeed, $g = \sum_{r} g_r x^{(r)}$, where $r$ is a~$(2k+1)$-tuple of non-negative numbers
and the sum runs over a~set of such tuples, and $g_r\in\Kee$ for all $r$. Then,
\[
g\cdot g = \sum_{r,s} g_r g_s x^{(r)}\cdot x^{(s)},
\]
i.e., the terms with $r\neq s$, are encountered 2 times, so what remains is
\be\label{prodBinom}
\sum_{r} g_r^2 x^{(r)}\cdot x^{(r)} = \sum_{r} g_r^2 \binom{2r}{r}x^{(2r)} = \sum_{r} g_r^2 \left(\prod_{i=0}^{2k} \binom{2r_i}{r_i}\right) x^{(2r)}.
\ee
If $a>0$, then $\binom{2a}{a} = 0 \pmod 2$.
By definition, $\binom{0}{1} := 0$; hence the product in equation~\eqref{prodBinom} vanishes even if there are only even indeterminates involved. Therefore, only the summand with $r=(0,\dots,0)$, i.e., $(g_{(0,\dots,0)})^2$, survives. \end{proof}

\begin{Claim}[Grading operators in $\fk(2n_\ev+1;\un|m)$ and $\fm(n;\un|n+1)$]\label{sssGrOp} \underline{For $p\neq 2$} and the bracket~\eqref{Konbr}, the $x_0$-action gives
a grading of $\fg$ by the formula ${\ad_{x_0}|_{\fg_i}=i\id}$; this action also
defines the $1$-dimensional $\fg_0$-module we denote by $\Kee[i]$.\index{$Kmathbb[i]$@$\Kee[i]$}

\underline{For $p=2$}, let\index{$Kmathbb[*]$@$\Kee[*]$}\index{((ZZPhi@$\Phi$}
\begin{equation}\label{Kee*}
\Kee[*] \text{~~denote the $\Kee x_0$-module analogous to
$\Kee[i]$}; \ \ \Phi:=\mathop{\sum}\limits_{1\leq i\leq k}
x_ix_{k+i}.
\end{equation}

\underline{For $p=2$}, the element $x_0$ annihilates a~subspace
$\textrm{ann}(x_0)|_{\fg_{-1}}$ of $\fg_{-1}$ and acts as
multiplication by $1$ on both
$\fg_{-1}/(\textrm{ann}(x_0)|_{\fg_{-1}})$, and $\fg_{-2}$.

\underline{For $p=2$}, the operators $\ad_{x_0}|_{\fg_0}$ and
$\ad_{\Phi}|_{\fg_0}$ interchange their roles as compared with $p\neq 2$: now $\ad_{\Phi}$ commutes
with $\fg_0$. \end{Claim}

\begin{Claim}[brackets in $\fk(2k_\ev+1;\un|2k_\od+1)$]\label{sssGrOp1}
 In items $1)$--$3)$ the isomophisms are between the LHS, defined as the Lie superalgebras of vector fields that multiply the contact form $\alpha$ by a~function, and the RHS, the brackets and squarings in which is given by formulas~\eqref{K_fbr}.
\begin{enumerate}\itemsep=0pt
\item[$1)$] $\fk(1;\un|0)\simeq\fvect(1;\un|0)$.

\item[$2)$] $\fk(1;\un|1)\simeq
(\Kee K_1\oplus \cO(x_0;\un))\bigoplus \Pi(\cO(x_0;\un))$. The even part of the simple Lie superalgebra $\fk^{(1)}(1;\un|1)$ is solvable. \textup{(For other examples of the same phenomenon indigenous to $p=2$, see \cite[Section~$16.2$]{BGL1}.)} This Lie superalgebra is the superization of $\fvect^{(1)}(1;\underline{N+1})$ by ``method $2$'', see $\cite{BLLS2}$ and $\cite{KrLe}$, and the
desuperization
$\textbf{F}(\fk^{(1)}(1;\un))$ is simple Lie algebra $\fvect^{(1)}(1;\underline{N+1})$.

\item[$3)$] $\fk(2k_\ev+1;\un|2k_\od+1)\simeq
\big(\fpo(2k_\ev;\hat\un|2k_\od)\oplus \widehat\cF\big)\bigoplus \Pi\big(\widehat\cF\big)$, where $\un=(N_0, \hat\un)$ and at least one of $k_\ev$ and $k_\od$ is non-zero. These Lie superalgebras and their desuperizations are not simple, the ideal $\fii$ is generated by the $A_g$ and $B_g$; recall \eqref{hatF}. We have
\[
\fk(2k_\ev+1;\un|2k_\od+1)/\fii\simeq\fh(2k_\ev;\un|2k_\od).
\]
\end{enumerate}
\end{Claim}

\subsection{Weisfeiler gradings}\label{WeisF} For vectorial Lie
superalgebras, the invariant notion is filtration, not grading. In characteristic~
$0$, the \textit{Weisfeiler} filtrations\index{$W$-filtration, Weisfeiler filtration}
\index{(W-grading, Weisfeiler grading@W-grading, Weisfeiler grading} were used in the description
of the infinite-dimensional Lie (super)algebras $\cL$ by selecting a
maximal subalgebra $\cL_{0}$ of finite codimension; for the simple
vectorial Lie algebra, there is only one such $\cL_{0}$. (Dealing
with finite-dimensional algebras for $p>0$, we can confine ourselves
to maximal subalgebras of \textit{least} codimension, or ``almost
least''.)

Let $\cL_{-1}$ be a~minimal $\cL_{0}$-invariant subspace strictly
containing $\cL_{0}$; for $i\geq 1$, set:
\begin{gather*}
\cL_{-i-1}=\begin{cases}[\cL_{-1}, \cL_{-i}]+\cL_{-i}&\text{unless $p=2$ and $i=1$},\\
[\cL_{-1}, \cL_{-i}] +\cL_{-i}+\Span\big(X^2\,|\,X\in \cL_{-1}\big)&\text{for $p=2$ and $i=1$},\\
\end{cases}\nonumber\\
\cL_i =\{D\in \cL_{i-1}\,|\,[D, \cL_{-1}]\subset\cL_{i-1}\}.\label{WeiGr}
\end{gather*}
We thus get a~filtration:
\begin{equation}
\label{Wfilt} \cL= \cL_{-d}\supset \cL_{-d+1}\supset \dots \supset
\cL_{0}\supset \cL_{1}\supset \cdots.
\end{equation}
The $d$ in \eqref{Wfilt} is called the \textit{depth} of $\cL$, and
of the associated \textit{Weisfeiler-graded} Lie superalgebra
$\fg=\mathop{\oplus}\limits_{-d\leq i}\fg_i$, where
$\fg_{i}=\cL_{i}/\cL_{i+1}$. We will for brevity say
\textit{W-graded} and \textit{W-filtered}.

For the list of simple \textit{W-graded} vectorial Lie superalgebras
$\fg=\mathop{\oplus}\limits_{-d\leq i}\fg_i$ over $\Cee$, see
\cite{LS} reproduced in Tables \ref{forgetseries} and \ref{excepC}.

\subsubsection{The $\Zee$-gradings of vectorial Lie superalgebras}\label{Zgrad} These gradings are defined by the vector $\vec r$ of degrees of the indeterminates, but this
vector can be shortened \textit{for W-gradings} to a~number $r$,
or a~symbol, which we do not indicate for $r=0$. Let the indeterminates $t$, $p_i$, $q_j$, and $u_\ell$ be even, while $\tau$, $\xi_i$, $\eta_j$, and $\theta_\ell$ be
odd. Let the contact Lie superalgebra $\fk (2n+1|m)$ preserve the distribution given by
the Pfaff equation
\[
\text{$\alpha_1(X) =0$ \ \ for $X\in
\fvect(2n+1|m)$,}
\]
where\index{((ZZZalpha_1$, odd contact form@$\widetilde\alpha_1$, odd contact form}
the form $\alpha_1$ is given by \eqref{alpha1}. For the $\fk$ series, let $u=(t; p,q)$ be even indeterminates, the odd indeterminates being the $\theta$ (resp.~$\theta$, $\xi$, $\eta$), see \eqref{alpha1}.

 For the $\fm$ series, the
indeterminates in Table \eqref{nonstandgr} are denoted as in formula \eqref{alp0}, i.e., the~$q_i$ even, the $\xi_i$, and $\tau$ odd.

In Table \eqref{nonstandgr}, the ``standard" gradings\index{(Grading, standard@Grading, standard} correspond to
$r=0$, they are marked by an asterisk~$(*)$. For $r=0$, the
codimension of ${\cal L}_0$ is the smallest.
\begin{equation}\label{nonstandgr}\footnotesize
\renewcommand{\arraystretch}{1.3}
\begin{tabular}{|c|c|}
\hline Lie superalgebra & its $\Zee$-grading \\ \hline

$\fvect (n|m; r)$, & $\deg u_i=\deg \theta_j=1$ for any $i, j$
\hfill $(*)$\\

\cline{2-2} where $0\leq r\leq m$ & $\deg \theta_j=0$ for $1\leq j\leq r$; \\
&$\deg u_i=\deg \theta_{r+s}=1$ for any $i, s$
\\

\hline $\fm(n; r),$ & $\deg \tau=2$, $\deg q_i=\deg \theta_i=1$ for
any $i$ \hfill $(*)$\\

\cline{2-2} where $0\leq r< n-1$& $\deg \tau=\deg
q_i=2$, $\deg \xi_i=0$ for $1\leq i\leq r$;\\

and one more grading (next line):& $\deg q_{r+j}=\deg \xi_{r+j}=1$ for any $j$\\
\hline

$\fm(n; n)$ & $\deg \tau=\deg q_i=1$, $\deg \xi_i=0$ for $1\leq
i\leq n$ \\ \hline

\cline{2-2} $\fk (2n+1|m; r)$, & $\deg t=2$,  $\deg
p_i=\deg q_i= \deg \xi_j=\deg \eta_j=\deg \theta_\ell=1$\\
where $0\leq r\leq [\frac{m}{2}]$,& for any $i$, $j$, $\ell$ \hfill $(*)$ \\

\cline{2-2} $r\neq k-1$ if $(n,m)=(0,2k)$& $\deg t=\deg
\xi_i=2$,
$\deg \eta_{i}=0$ for $1\leq i\leq r$; \\
and one more grading (next line):&$\deg p_i=\deg q_i=\deg
\theta_{\ell}=1$
for $\ell\geq 1$ and all $i$\\

\hline $\fk(1|2m; m)$ & $\deg t =\deg \xi_i=1$, $\deg \eta_{i}=0$
for $1\leq i\leq m$ \\ \hline
\end{tabular}
\end{equation}

\subsection{Divergence-free and traceless subalgebras}\label{SSdivFree} In this subsection, the ground field is any $\Kee$ for $p\neq 2$. The peculiarities of $p=2$ are considered in Sections~\ref{mandbab} and \ref{kandpoab}. Here we will not mention $\un$ if $p>2$.

\subsubsection{$\fk$ series}\label{SSdivFreeK12} Since the restriction of the standard divergence \eqref{DivVect} to the subalgebra of degree 0 is (super)trace, and since the space $\fg_0/[\fg_0, \fg_0]$, where $\fg:=\fk(2n+1|m)$, is spanned by~$K_t$ for $(n,m)\neq (0,2)$, it is easy to calculate that
\begin{equation}
\label{div2}
 \Div K_f =(2n+2-m)\del_t (f) \quad \text{if $2n+2-m\neq 0$},
\end{equation}
it follows that the divergence-free (relative the restriction of the divergence \eqref{DivVect} to $\fk(2n+1|m)$) subalgebra of the contact Lie
superalgebra either coincides with it for $m=2n+2$ or is the Poisson superalgebra singled our by the condition $\del_t(f)=0$.

On $\fk(2n+1|2n+2)$ there is its own, ``indigenous'' divergence $K_f\mapsto \del_t(f)$; it also singles out the
Poisson superalgebra. This, however, is not the whole story: the case $\fk(1|2)$ is exceptional.

{\bf The case of $\boldsymbol{\fk(1|2)}$.} 
\underline{Let $\alpha_1=dt+\xi d\eta+\eta d\xi$}. Since $\fk(1|2)_0$ is commutative and 2-di\-men\-sio\-nal, there are 2 linearly independent traces on it: one~-- $\tr$~-- is equal to 1 at $t$ and vanishes at $\xi\eta$, the other one~-- call it $\tr_{(2)}$~-- is equal to 1 at $\xi\eta$ and vanishes at $t$.

Clearly, the condition $K_1(f)=0$ singles out the subalgebra $\fk_-\oplus\Kee\xi\eta$
of $\fk(1|2)$. In other words, the operator $\del_t=\frac12K_1$ in the adjoint representation is an analog of the divergence~-- the prolong of the trace on $\fk_0$; this analog is equal to 1 at $t$ and vanishes at $\xi\eta$.

The divergence-free condition $\Div D=0$, where $D\in\fg$ for a~$\Zee$-graded vectorial Lie super\-algebra $\fg$, should single out the complete prolong of $(\fg_-, \fs)$, where $\fs=\{g\in \fg_0\,|\,\tr g=0\}$. Therefore, the condition that determines the divergence is
\begin{gather}\label{(1)}
X(\Div D)=\Div([X,D]) \quad \text{for any} \ X\in \fg_-.
\end{gather}
Since $\fg_{-1}$ generates the negative part, it suffices to require fulfillment of the condition~\eqref{(1)} for any $X\in \fg_{-1}$.

Therefore, we have to express the divergence not in terms of partial derivatives, but in terms of the operators commuting (not supercommuting) with $\fg_-$ (recall that in~\cite{Shch}, the operators that span $\fg_-$ are denoted by $X_i$, and the operators commuting with $\fg_-$ are denoted by $Y_i$).

To write the second divergence $\Div_{(2)}$, which is the prolong of $\tr_{(2)}$, we need two operators commuting (not supercommuting!) with $\fk(1|2)_-$. In our case, the $X$-operators are
\[
K_1=2\del_t, \quad  K_\xi=\del_\eta+\xi\del_t \quad \text{and} \quad K_\eta=\del_\xi+\eta\del_t,
\]
then the needed $Y$-operators are
\begin{gather}\label{Yop}
\tilde K_1=\partial_t,\quad \tilde K_\xi(f)=(-1)^{p(f)}(\del_\eta-\xi\del_t)(f) \quad \text{and} \quad \tilde K_\eta(f)=(-1)^{p(f)}(\del_\xi-\eta\del_t)(f).
\end{gather}

\underline{Let $\alpha_1=dt+\xi d\eta$} which works for any characteristic.
Then, the $X$-operators are
\[
K_1=\del_t, \quad K_\xi=\del_\eta \quad \text{and} \quad K_\eta=\del_\xi+\eta\del_t,
\]
and the $Y$-operators are
\begin{gather}\label{YopTil}
\tilde K_1=\del_t, \quad \tilde K_\xi(f)=(-1)^{p(f)}(\del_\eta-\xi\del_t)(f) \quad \text{and}\quad \tilde K_\eta(f)=(-1)^{p(f)}\del_\xi(f).
\end{gather}

\begin{Claim}[the second divergence on $\fk(1|2)$]\index{(Divergence@Divergence} The prolong of $\tr_{(2)}$ composed of the $Y$-opera\-tors~\eqref{Yop} or \eqref{YopTil} is the same \index{$Div_{(2)}$@$\Div_{(2)}$}
\begin{gather*}
\Div_{(2)}:=\tilde K_\eta\tilde K_\xi-\tilde K_1,
\end{gather*}
and $\Div_{(2)}(\xi\eta)=1$ while $\Div_{(2)}(t)=0$.
\end{Claim}

{\bf The $\boldsymbol{\fk(1|2)}$-module of weighted densities.} 
 Over contact Lie superalgebras $\fk(2n+1|m)$ it is natural to express the spaces of weighted densities\index{(Density, weighted@Density, weighted} in terms of the conformally preserved form~$\alpha_1$. This recalculation is well-known for $m=0$, where $\vvol=\alpha_1\wedge (d\alpha_1)^n$. The general case follows from equation~\eqref{div2}: from the point of view of the $\fk(2n+1|m)$-action
\[
\vvol=\begin{cases}(\alpha_1)^{(2n+2-m)/2}&\text{if $(m, n)\neq (0,0)$},\\
\alpha_1&\text{if $(m, n)= (0,0)$}.
\end{cases}
\]

Since the center of $\fk(1|2)_0$ is of dimension 2, the weights of the spaces of weighted densities have 2 parameters, not one: $\cF_{a,b}:=\cF\alpha_1^a\beta^b$, where $a,b\in\Kee$.
\index{$Fcal_{a,b}:=\cF\alpha_1^a\beta^b$, the spaces of weighted densities over $\fk(1\vert 2)$@$\cF_{a,b}:=\cF\alpha_1^a\beta^b$, the spaces of weighted densities over $\fk(1\vert 2)$}

Let $\beta$ be the symbol of the \textit{class} of the differential form $d\xi$ (or, equivalently, $(d\eta)^{-1}$ )in the quotient space $\Omega^1/\cF\alpha_1$ of 1-forms. The Lie derivative acts as follows
\begin{align}
L_{K_f}\big(\alpha_1^a\beta^b\big)&=\big(a\del_t(f)+(-1)^{p(f)}b\Div_2(K_f)\big)\big(\alpha_1^a\beta^b\big)\nonumber\\
&=\big(\big(a-(-1)^{p(f)}b\big)\del_t(f)+(-1)^{p(f)} b\tilde K_\eta\tilde K_\xi(f)\big)\big(\alpha_1^a\beta^b\big).\label{L_Dk12}
\end{align}
The space $\cF_{a,b}$ of weighted densities over $\fk(1|2)$ is a~rank-1 module generated by $\alpha_1^a\beta^b$ over the algebra of functions $\cF=\cF_{0,0}$.

\subsubsection[$\fm$ series, its simple subalgebras, and weighted densities]{$\boldsymbol{\fm}$ series, its simple subalgebras, and weighted densities}\label{SSdivFreeM} For the pericontact series, the situation is
more interesting than that for contact series: the divergence-free subalgebra is simple and new (only as compared with the above-described algebras; it is known since ca~1978, see \cite{ALSh}).

Let $p\neq 2$. Since
\begin{equation*}
\label{Div} \Div M_f =(-1)^{p(f)}2\left ((1-E)\pderf{f}{\tau} -
\mathop{\sum}\limits_{i\leq n}\frac{\del^2 f}{\del q_i
\del\xi_i}\right ),
\end{equation*}
it follows that the divergence-free subalgebra of the pericontact
superalgebra is\index{$smmathfrak (n)$@$\fsm (n)$}
\begin{equation*}
\label{sm2} \fsm (n) = \Span\left (M_f \in \fm (n)\,|\,
(1-E)\pderf{f}{\tau} =\mathop{\sum}\limits_{i\leq n}\frac{\del^2
f}{\del q_i \del\xi_i}\right ).
\end{equation*}
In particular,\index{((ZZDelta@$\Delta$}
\begin{equation*}
 \Div \Le_f = (-1)^{p(f)}2\mathop{\sum}\limits_{i\leq
n}\frac{\del^2 f}{\del q_i \del\xi_i}=(-1)^{p(f)}2\Delta(f), \quad \text{where $\Delta:=\mathop{\sum}\limits_{i\leq
n}\frac{\del^2}{\del q_i \del\xi_i}$}.
\end{equation*}
The divergence-free vector fields from $\fsle (n)$ are generated by
\textit{harmonic} functions, i.e., such that $\Delta(f)=0$.

Rank 1 over the algebra $\cF$ modules $\cF_{a,b}^\fm:=\cF\alpha_0^a\gamma^b$,\index{$Fcal_{a,b}^\fm:=\cF\alpha_0^a\gamma^b$, the space of weighted densities over $\fm(n)$@$\cF_{a,b}^\fm:=\cF\alpha_0^a\gamma^b$, the space of weighted densities over $\fm(n)$} where $a,b\in\Kee$,
 are generated by $\alpha_0^a\gamma^b$, where $\gamma$ is a~symbol of the class of differential forms (whose explicit expression is irrelevant, same as that of $\beta$, see equation~\eqref{L_Dk12}). The Lie derivative acts as follows:
\begin{equation*}\label{DeltaM_f2}
L_{M_f}\big(\alpha_0^a\gamma^b\big) =\big((-1)^{p(f)}b\partial_\tau(f)+ a~\Delta^\fm(f)\big)\big(\alpha_0^a\gamma^b\big).
\end{equation*}

The divergence-free relative the standard divergence Lie superalgebras $\fsle (n)$, $\fsb (n)$ and $\fsvect (1|n)$
have traceless ideals $\fsle^{(1)}(n)$, $\fsb^{(1)}(n)$ and
$\fsvect^{(1)}(n)$ of codimension 1; they are defined from the exact
sequences\index{$slemathfrak^{(1)}(n)$@$\fsle^{(1)}(n)$}\index{$bmathfrak^{(1)}(n)$@$\fb^{(1)}(n)$}\index{$svectmathfrak^{(1)}(1\vert n)$@$\fsvect^{(1)}(1\vert n)$}
\begin{gather*}
0\tto \fsle^{(1)}(n)\tto \fsle(n)\tto \Kee\cdot \Le_{\xi_1\dots\xi_n} \tto 0, \\
0\tto \fsb^{(1)}(n)\tto \fsb(n)\tto \Kee\cdot M_{\xi_1\dots\xi_n} \tto 0, \\
 0\tto \fsvect^{(1)}(1|n)\tto \fsvect (1|n)\tto \Kee
\cdot\xi_1\cdots\xi_n\partial_{t}\tto 0.
\end{gather*}

\subsubsection[A~deform of the series $\fb$]{A~deform of the series $\boldsymbol{\fb}$}\label{deforBab} Let $p\neq 2$. For an explicit
form of $M_f$, see Section~\ref{GenFun}. Set\index{$bmathfrak_{a, b}(n)$@$\fb_{a, b}(n)$}
\begin{equation*}
\fb_{a, b}(n) =\left\{M_f\in \fm (n)\,|\,a\; \Div
M_f=(-1)^{p(f)}2(aE-bn)\pderf{f}{\tau}\right\}.
\end{equation*}
We denote the operator that singles out $\fb_{\lambda}(n)$ in $\fm
(n)$ as follows, cf.~\eqref{LieDer}:\index{$Div_{\lambda}$@$\Div_{\lambda}$}\index{((ZZDelta@$\Delta$}
\begin{equation*}
\Div_{\lambda}=(bn-aE)\pder{\tau}-a\Delta,\quad \text{for $\lambda
=\frac{2a}{n(a-b)}$ \quad and \quad $\Delta=\mathop{\sum}\limits_{i\leq
n}\frac{\del^2 }{\del q_i \del\xi_i}$}.
\end{equation*}
Taking the explicit form of the divergence of $M_{f}$ into account, we get
\begin{align}
\fb_{a, b}(n) &=\left\{M_f\in \fm (n)\,|\,(bn-aE)\pderf{f}{\tau} =
a\Delta f\right\}\nonumber \\
&=\big\{D\in\fvect(n|n+1) \,|\,L_{D}\big(\vvol_{q, \xi, \tau}^a\alpha_{0}^{a-bn}\big)=0\big\}. \label{b_a,b'}
\end{align}
It is subject to a~direct verification that $\fb_{a, b}(n)\simeq
\fb_\lambda(n)$ for $\lambda =\frac{2a}{n(a-b)}\in\Kee P^1$.
\index{$bmathfrak_{\lambda}(n)$@$\fb_{\lambda}(n)$}
Obviously, if $\lambda=0, 1, \infty$ (where $\fb_{0}:=\fb$ and
$\fb_\infty:=\fb_{a,a}$) the structure of
$\fb_\lambda(n)$ differs from the other members of the parametric
family: the following exact sequences single out simple Lie
superalgebras (the quotient $\fle(n)$ and ideals, the first derived
subalgebras):\index{$bmathfrak_1^{(1)}(n)$@$\fb_1^{(1)}(n)$}\index{$bmathfrak_\infty^{(1)}(n)$@$\fb_\infty^{(1)}(n)$}
\begin{gather}
0\tto \Kee M_{1} \tto \fb(n)\tto \fle(n)\tto 0,\nonumber\\
0\tto \fb_{1}^{(1)}(n)\tto \fb_{1}(n)\tto \Kee\cdot M_{\xi_1\dots\xi_n} \tto 0,\nonumber\\
0\tto \fb_{\infty}^{(1)}(n)\tto \fb_{\infty}(n)\tto \Kee\cdot M_{\tau\xi_1\dots\xi_n} \tto 0.\label{specValuesOfLamb}
\end{gather}

\begin{Problem}\label{OPfb}\index{(Problem, open@Problem, open}
The Lie superalgebras $\fb_{\lambda}(n)$ can be further deformed at
certain points $\lambda$, see~{\rm \cite{LSh}}, where $\Kee=\Cee$; the Lie superalgebras of series $\fh$ and $\fle$ also have extra deformations. Describe the deformations of $\fb_{\lambda}(n;\un)$, as well as $\fh$ and $\fle$ for all $p>0$.
\end{Problem}

\subsection[Passage from $p=0$ to $p>0$]{Passage from $\boldsymbol{p=0}$ to $\boldsymbol{p>0}$}\label{CtoK} Here we have
collected answers to several questions that stunned us while we were
writing this paper. We hope that even the simplest of these answers
will help the reader familiar with representations of Lie algebra
over $\Cee$, but with no experience of working with characteristic~$p>0$.

\textit{For $p=2$, several of our definitions are new}, see Sections~\ref{mandbab} and~\ref{kandpoab}.

\subsubsection[The Lie (super)algebras preserving symmetric non-degenerate\\ bilinear forms~$\cB$]{The Lie (super)algebras preserving symmetric non-degenerate\\ bilinear forms~$\boldsymbol{\cB}$}

We often denote the Gram matrix of the bilinear form
$\cB$ also by $\cB$, let $\faut(\cB)$\index{$autmathfrak(\cB)$@$\faut(\cB)$} be the Lie (super)algebra
preserving $\cB$. If $\cB$ is odd and the superspace, on which it is
defined, is of superdimension $n|n$, we write $\fpe_\cB(n)$\index{$pemathfrak_\cB(n)$@$\fpe_\cB(n)$} instead of
$\faut(\cB)$.

{\bf Let $\boldsymbol{p\neq 2}$ and $\boldsymbol{\fg=\fpe_\cB(n)}$.} The Lie superalgebra
$\fg$ consists of the supermatrices of the form
\begin{gather*}
X=\begin{pmatrix} A & B\\C & -A^{\rm t}
 \end{pmatrix}, \text{~~where $B$ is symmetric and $C$ is
antisymmetric}\\
\text{if the form $\cB$ is in its
normal shape $\Pi_{n|n}:=\Pi_{2n}=\begin{pmatrix}0&1_n\\1_n&0\end{pmatrix}$}.\label{matrgen4}
\end{gather*}
Clearly, $\str X= 2\tr A$. We also have $\fg^{(1)}=\fspe(n)$, i.e., $\fspe(n)$
is of codimension 1; it is singled out by the condition $\str X=0$,
which is equivalent to $\tr A=0$.

The Lie superalgebra $\fle(n;\un|n)$ is, by definition, the Cartan
prolong $(\id, \fpe(n))_{*,\un}$.

Over $\Cee$, there is no shearing vector, and $\fle(n):=\fle(n|n)$
is spanned by the elements $\Le_f$, where $f\in\Cee[q,\xi]$.

If $p>2$, the elements of $\cO(q;\un|\xi)\oplus
\Span\big(q_i^{(p^{\un_i})}\,|\,N_i<\infty\big)$, or $\cO(q;\un_\infty|\xi)$
for $\un=\un_\infty$, see~\eqref{un_s}, generate $\fle(n;\un)$. If
$N_i<\infty$ for at least one $i$, the additional part
$\text{Irreg}$ does not change while the regular part
looks the same for any $p>2$:
\begin{gather}
\text{Reg}=\Span (\Le_f\,|\,f\in \cO(q;\un|\xi) )\oplus
\Span\big(q_i^{(p^{\un_i})}\,|\,N_i<\infty\big), \nonumber\\ \text{Irreg}
=\Span(\xi_i\del_{q_i})_{i=1}^n.\label{irreg}
\end{gather}
In other words: \textit{there are vector fields corresponding to
non-existing generating functions, like $q_i^{(p^{\un_i})}$ and~$\xi_j^2$}.
The prolong $\fsle(n;\un):=(\id, \fspe(n))_{*, \un}$ is singled out by the
condition\index{((ZZDelta@$\Delta$}
\[
\Div \Le_f=0\Longleftrightarrow \Delta f=0, \quad \text{where} \ \Delta
=\mathop{\sum}\limits_{i\leq n}\frac{\del^2 }{\del q_i \del\xi_i}.
\]
The operator $\Delta$ is, therefore, the ``Cartan prolong of the
supertrace on $\fg_0$'' expressed as an operator acting on the space of
generating functions.

{\bf Modifications in the above description for $\boldsymbol{p=2}$.} 
If $p=2$, the analogs of symplectic (resp.\ periplectic) Lie
(super)algebras accrue additional elements: if the matrix of the
bilinear form $\cB$ is $\Pi_{2n}$
(resp.~$\Pi_{n|n}$), then
$\faut(\cB)$ consists of the (super)matrices of the form
\begin{equation}
\label{matrgen} X=\begin{pmatrix} A & B\\C & A^{\rm t}
\end{pmatrix}, \quad \text{where $B$ and $C$ are symmetric, $A\in\fgl(n)$}.
\end{equation}
Denote the \textit{general} Lie (super)algebra preserving the form $\cB$ as follows:\index{$omathfrak_{\rm gen}(2n)$@$\fo_{\rm gen}(2n)$}
\index{$pemathfrak_{\rm gen}(n)$@$\fpe_{\rm gen}(n)$}
\begin{equation*}\label{AutGen}
\faut(\cB)=\begin{cases}\fo_{\rm gen}(2n)&\text{for
$\cB=\Pi_{2n}$},\\
\fpe_{\rm gen}(n)&\text{for $\cB=\Pi_{n|n}$}.
\end{cases}
\end{equation*}

Let\index{(ZD@ZD}
\begin{equation*}\label{ZD}
ZD\text{~~denote the space of symmetric matrices with zeros on their
main diagonals.}
\end{equation*}
The derived Lie (super)algebra $\faut^{(1)}(\cB)$
consists of the (super)matrices of the form \eqref{matrgen}, where
$B,C\in ZD$. In other words, these Lie (super)algebras resemble the
orthogonal Lie algebras.

On these Lie (super)algebras $\faut^{(1)}(\cB)$ the following
(super)trace (\textit{half-trace}) is defined:\index{$Zmh$, trace or half-trace@$\htr$, trace or half-trace}
\begin{equation*}
\label{matrgen3} \htr\colon \ \begin{pmatrix} A~& B\\C & A^{\rm t}
\end{pmatrix}\tto \tr A.\end{equation*}
The half-traceless Lie sub(super)algebra of $\faut^{(1)}(\cB)$\index{$autmathfrak^{(1)}(\cB)$@$\faut^{(1)}(\cB)$} is
isomorphic to $\faut^{(2)}(\cB)$.

There is, however, an algebra $\widetilde\faut(\cB)$,\index{$autmathfrak(\cB)$@$\widetilde\faut(\cB)$} such that
$\faut^{(1)}(\cB)\subset\widetilde\faut(\cB)\subset\faut(\cB)$,
consisting of (super)matrices of the form \eqref{matrgen}, where
$B\in ZD$, and any symmetric~$C$ (or isomorphic to it version
of the Lie superalgebra with
$C\in ZD$, and any symmetric $B$). We suggest that it be denoted as follows:
\begin{gather*}
\widetilde \faut(\cB)=\begin{cases}
\textbf{F}(\fpe)(2n)&\text{for $\cB$ even},\\
\fpe(n)&\text{for $\cB$ odd},
\end{cases}\\ \faut^{(1)}(\cB)=\{X\in \widetilde \faut(\cB)
\,|\,\htr X=0\}=\begin{cases}
\textbf{F}(\fspe)(2n)&\text{for $\cB$ even},\\
\fspe(n)&\text{for $\cB$ odd}.
\end{cases}
\end{gather*}

\subsection{Central extensions}\label{ssAS} There is only one non-trivial
central extension of $\fspe(n)$ for ${p\neq 2, 3}$ existing only for
$n=4$. We denote it $\fas$\index{$asmathfrak$@$\fas$} because
it was discovered by A.~Sergeev (1970s, unpublished). For numerous non-trivial central extensions of versions of $\fspe(n)$ and its simple subquotients for ${p=2, 3}$, see~\cite{BGLL3}.

Let us represent an arbitrary
element $A\in\fas$ as a~pair $A=x+d\cdot z$, where $x\in\fspe(4)$,
$d\in{\Cee}$, and $z$ is the central element. The bracket in $\fas$ is
\begin{equation}\label{2.1.4}
\left[\mat {a &b\\ c&-a^{\rm t}} +d\cdot z,
\mat{a' & b' \cr c' & -(a')^{\rm t}} +d'\cdot
z\right]= \left[\mat{ a& b \cr c & -a^{\rm t}},
\mat{ a' & b' \cr c' & -(a')^{\rm t}
}\right]+\tr~c\widetilde c'\cdot z,
\end{equation}
where $\ \widetilde {}\ $ is extended via linearity from matrices
$c_{ij}=E_{ij}-E_{ji}$ on which $\widetilde c_{ij}=c_{kl}$ for any
even permutation $(1234)\longmapsto(ijkl)$. Recall that $b=b^{\rm t}$ and $b'=(b')^{\rm t}$ in \eqref{2.1.4}, whereas
$c=-c^{\rm t}$ and $c'=-(c')^{\rm t}$.

The Lie superalgebra $\fas$ can also be described in terms of
the spinor representation. For this, we need several vectorial
superalgebras. Consider $\fpo(0|6)$, the Lie superalgebra whose
superspace is the Grassmann superalgebra $\Lambda(\xi, \eta)$
generated by $\xi=(\xi_1, \xi_2, \xi_3)$ and $\eta=(\eta _1, \eta _2, \eta_3)$ with the
Poisson bracket.

Recall that $\fh(0|6)=\Span (H_f\,|\,f\in\Lambda (\xi, \eta))$. Now,
observe that $\fspe(4)$ can be embedded into $\fh(0|6)$. Indeed,
setting $\deg \xi_i=\deg \eta _i=1$ for all $i$ we introduce
a~$\Zee$-grading on $\Lambda(\xi, \eta)$ which, in turn, induces
a~$\Zee$-grading on $\fh(0|6)$ of the form
$\fh(0|6)=\mathop{\oplus}\limits_{i\geq -1}\fh(0|6)_i$. Since
$\fsl(4)\cong\fo(6)$, we can identify $\fspe(4)_0$ with
$\fh(0|6)_0$.

It is not difficult to see that the elements of degree $-1$ in the
standard gradings of $\fspe(4)$ and $\fh(0|6)$ constitute isomorphic
$\fsl(4)\cong\fo(6)$-modules. It is subject to a~direct verification
that it is possible to embed $\fspe(4)_1$ into $\fh(0|6)_1$.

Sergeev's extension $\fas$ is the result of the restriction to
$\fspe(4)\subset\fh(0|6)$ of the cocycle that turns $\fh(0|6)$ into
$\fpo(0|6)$. The quantization deforms $\fpo(0|6)$ into
$\fgl(\Lambda(\xi))$; the through maps $T_\lambda\colon
\fas\tto\fpo(0|6)\tto\fgl(\Lambda (\xi))$ are representations of
$\fas$ in the $4|4$-dimensional modules $\spin_\lambda$ isomorphic
to each other for all $\lambda\neq 0$. The explicit form of
$T_\lambda$ is as follows:
\begin{equation*}
T_\lambda\colon \ \mat{ a & b \cr c & -a^{\rm t} }+d\cdot
z\longmapsto \mat{ a & b-\lambda \widetilde c \cr c & -a^{\rm t}
}+\lambda d\cdot 1_{4|4},
\end{equation*}
where $1_{4|4}$ is the unit matrix and $\widetilde c$ is defined in
the line under equation~\eqref{2.1.4}. Clearly, $T_\lambda$ is an
irreducible representation for any $\lambda$.

\subsection{Prolongs} The Lie superalgebra $\fpe_{\rm gen}(n)$ is larger than
$\fpe(n)$: both $B$ and $C$ in $\fpe_{\rm gen}(n)$ are symmetric, see
\eqref{matrgen}. Observe that $\fpe_{\rm gen}(n)\subset\fsl(n|n)$.
Denote\index{$lemathfrak_{\rm gen}(n;\un\vert n)$@$\fle_{\rm gen}(n;\un\vert n)$}
\begin{equation*}
\fle_{\rm gen}(n;\un|n):=(\id,
\fpe_{\rm gen}(n))_{*,\un}.
\end{equation*}
Clearly, if $\un=\un_\infty$,
see \eqref{un_s}, then $\fle_{\rm gen}(n;\un|n)$ consists of the following two parts, cf.\ equation~\eqref{irreg}:
\begin{equation}\label{reg_gen}
\text{Reg}_{\rm gen}=\Span(\Le_f \,|\,f\in \cO(q;\un| \xi)), \quad
\text{Irreg}_{\rm gen} =\Span(B_i:=\xi_i\del_{q_i})_{i=1}^n.
\end{equation}
The part $\text{Irreg}_{\rm gen}$ corresponds to the \textit{nonexisting}
generating functions $\xi_i^2$. Clearly, $\fle_{\rm gen}(n;\un|n)$ is
contained in $\fsvect(n;\un|n)$, and therefore
\[
\fle_{\rm gen}(n;\un|n)=\fsle_{\rm gen}(n;\un|n).
\]

The difference between $\fle_{\rm gen}(n;\un|n)$ and
$\fle(n;\un|n)$ is constituted by the space $\text{Irreg}_{\rm gen}$. The
nonexisting generating functions $\xi_i^{(2)}$ generate linear
vector fields corresponding to the diagonal elements of the matrices
$B$ in \eqref{matrgen}, like the $q_i^{(2)}$ generate linear
vector fields corresponding to the diagonal elements of the matrices
$C$ in \eqref{matrgen}, but these two sets of elements are different
in their nature: there are no elements of degree $>0$ in
$(\id,\fpe_{\rm gen})_{*,\un}$ whose brackets with $\fg_{-1}$ give the
$B_i$ in \eqref{reg_gen}.

The \underline{\textit{correct} $p=2$ analogs} of the complex Lie superalgebras
$\fsle(n)$ and $\fspe(n)$ are, respectively, $\big(\id,(\fpe(n))^{(1)}\big)_{*,\un}$ and
$\fpe(n)^{(1)}$.

In \cite{Leb}, Lebedev considered $\fg=\fpe(n)$, the derived
algebras $\fg^{(1)}$ and $\fg^{(2)}$, and the Cartan prolongs of
these derived algebras playing the role of $\fg_0$, whereas for $\fg_{-1}$ he considered the tautological $\fg_0$-module
$\id$. Clearly, $\fg^{(1)}$
consists of supermatrices of the form \eqref{matrgen} with
zero-diagonal matrices $B$ and $C$, whereas $\fg^{(2)}$ is singled
out of $\fg^{(1)}$ by the condition $\htr =0$. The corresponding Cartan prolongs
only have the regular parts:
\begin{gather*}
\big(\id,\fg^{(1)}\big)_{*,\un}=\Span\big(\Le_f\,|\,f\in
\cO(q;\One|\xi)\big);\\
\big(\id,\fg^{(2)}\big)_{*,\un}=\Span\big(\Le_f\,|\,f\in \cO(q;\One|\xi)
\  \text{and} \ \Delta f =0\big).
\end{gather*}
Let a~non-degenerate (anti)symmetric bilinear form $\cB$ be defined on a~superspace $V$; let
$\textbf{F}(\cB)$ be the same form considered on $\textbf{F}(V)$, the same space with superstructure forgotten. Let
$\fh_\cB(a;\un|b)$ denote the Hamiltonian Lie superalgebra~-- the Cartan prolong of the ortho-orthogonal Lie superalgebra
$\fo\fo_\cB(a|b)$ preserving the non-degenerate form $B$; its
desuperization is $\fh_{\textbf{F}(\cB)}(a+b;\widetilde\un)$, where
$\widetilde\un$ has no critical coordinates.

\begin{Remark} For $\un$ with $N_i<\infty$ for all $i$
and $p=2$, the Lie superalgebra $\fle^{(1)}(n;\un|n)$ is spanned by
the elements $f\in \cO(q;\un|\xi)$, whereas each of the ``virtual''
generating functions $q_i^{(2^{\un_i})}\not\in \cO(q;\un|\xi)$ determines an outer
derivation of $\fle^{(1)}(n;\un|n)$.
\end{Remark}

\subsubsection[Divergence-free subalgebras $\fg$ of series $\fh$ and $\fle$ in the standard W-grading]{Divergence-free subalgebras $\boldsymbol{\fg}$ of series $\boldsymbol{\fh}$ and $\boldsymbol{\fle}$ in the standard W-grading}

These
subalgebras are prolongations of subalgebras of 0th components of
$\fh$ and $\fle$ consisting of traceless subalgebras; that is how these (super)algebras were described in~\cite{Leb}.

It is possible, however, to describe various subalgebras of $\fh_0$ or $\fle_0$, generated by (linear combinations of)
quadratic monomials, by eliminating squares of indeterminates from the set of functions generating $\fg_0$. In other words,
\begin{gather*}
\begin{minipage}[l]{13cm}
\textit{constraints imposed on the shearing vector $\un$
corresponding to the space of generating functions determine various
divergence-free subalgebras of $\fh(n;\un)$ and $\fle(n;\un)$}.
\end{minipage}
\end{gather*}

\subsubsection[$\fsvect_{a, b}(0|n)$]{$\boldsymbol{\fsvect_{a, b}(0|n)}$}\label{AandB} For $p>0$, let\index{$svectmathfrak_{a, b}(0\vert n)$@$\fsvect_{a, b}(0\vert n)$}
$\fsvect_{a, b}(0|n)$ denote $\fsvect(0|n)\ltimes
\Kee(az+bd)$, where the element ${d:=\sum\xi_i\partial_{\xi_i}}$ determines the standard
$\Zee$-grading of $\fsvect(0|n)$, while
$z$ is an element generating the trivial center commuting with
$\fsvect(0|n)\ltimes \Kee\cdot d$.

\subsubsection[$\fspe_{a, b}(n)$]{$\boldsymbol{\fspe_{a, b}(n)}$}

\label{AandBspe} For $p=0$, the meaning of\index{$spemathfrak_{a, b}(n)$@$\fspe_{a, b}(n)$}
$\fspe_{a, b}(n)$ is similar to that of $\fsvect_{a, b}(0|n)$, but
with $d:=\diag(1_n,-1_n)$. To define the analog of $\fspe_{a, b}(n)$
for $p=2$, see line $N=7$ in Table~\ref{cts-prol}, observe that the
codimension of $\fspe(n)$ in $\fm_0$, where $\fm:=\fm(n)$ is considered
in its standard $\Zee$-grading, is equal to~2.

So, to pass from
$\fspe(n)$ to $\fm_0$, we have to add \textit{two linearly independent} elements, whereas
to pass to $\fspe_{a,b}$ we have to add a~linear combination of
these elements with coefficients~$a$ and~$b$. The question is: ``can
we single out these elements in a~canonical way?"

\underline{For $p=0$}. The identity operator (in matrix realization)
is one of these elements. How to select the other element? There is
no distinguished element in $\fpe(n)\setminus \fspe(n)$. But, if $p=0$, there is an element $\diag(1_n,-1_n)$
corresponding to a~``most symmetric" generating function $\sum
q_i\xi_i$.

\underline{For $p>2$} this ``most symmetric" element
lies in $\fspe(n)$ if $p$ divides $n$ and the choice of the linearly independent second element
from $\fpe(n)\setminus \fspe(n)$ becomes a~matter of taste.

\underline{For $p=2$}, the situation becomes completely miserable. Now, the
restriction of $M_{\sum q_i\xi_i}$ to~$\fm_{-1}$ not only lies in
$\fspe(n)$ for $n$ even, it coincides with the identity operator.
So, in this case, there is no distinguished operator not lying in
$\fspe(n)$. What to do?

We suggest considering the elements of $\fm_0$ as operators acting
not just on $\fm_{-1}$, but on the whole $\fm_-$. If $\fm_0$ is thus
understood, there are two well-defined linear forms \index{$(l$, linear form@$\ell$, linear form}\index{((ZZZmu, linear form@$\mu$, linear form}
$\ell $ and $\mu$ that single out
$\fspe(n)$ in $\fm_0$:
\begin{gather}
\text{for any operator $A\in\fm_0$, let
$A_i=\ad_A|_{\fm_i}$; then $A_{-2}=\ell (A)\cdot\id$ and $A_{-1}$
is as in \eqref{matrgen}},\nonumber\\
\mu(A)=\begin{cases}\htr (A_{-1})&\text{for $p=2$},\\
\str (A_{-1})&\text{for $p>2$}.
\end{cases}\label{2traces}
\end{gather}

Now, $\fspe(n)$ is singled out by conditions $\ell (A)=\mu(A)=0$,
while\index{$spemathfrak_{a,b}(n)$@$\fspe_{a,b}(n)$}
\begin{equation*}\label{spe_a,b}
\fspe_{a,b}(n):=\{X\in\fpe(n)=\fm_0\,|\,(a\mu+b\ell )(X)=0\}.
\end{equation*}

\subsection[On $\fm$ and $\fb$]{On $\boldsymbol{\fm}$ and $\boldsymbol{\fb}$}\label{mandbab} To pass from $\fb(n;\un|n)$
to $\fm(n;\un|n+1)$, we have to add to $\fb(n)_0=\fpe(n)$ the
central element; it will serve as a~grading operator of the prolong.
We see that $\fm$ is the generalized Cartan prolong of $(\fb(n)_-,\fc\fb(n)_0)$.

The commutant of $\fm(n;\un|n+1)_0$ is the like that of
$\fb(n)_0=\fpe(n)$, so is of codimension~2. Hence there are two
traces on $\fm(n;\un|n+1)_0$, namely $\htr$ and $\ell$, see~\eqref{2traces}, and therefore there are two
divergences on $\fm$.\index{(Divergence@Divergence} One of them is given by the operator\index{$D_\tau$}
\begin{gather*}
\del_\tau, \quad \text{more precisely}  \ D_\tau:=\del_\tau \circ \sign,\\ \text{i.e.,} \quad D_\tau(f)=(-1)^{p(f)}\del_\tau(f)
\quad \text{for any} \
f\in\cO(q;\un|\tau,\xi)
\end{gather*}
since this should be the mapping \textit{commuting} (not
\textit{super}commuting) with $\fm_-$, see \cite{Shch}. The
condition $D_\tau(f)=0$, i.e., just $\del_\tau(f)=0$ singles out
precisely $\fb(n)$.

The other divergence is given by the operator \eqref{DeltaM_f1}.

\subsubsection[$\fsb$]{$\boldsymbol{\fsb}$}\label{sb} The definition of $\fsb(n;\un)$\index{$sbmathfrak(n;\un)$@$\fsb(n;\un)$} is the same
for any characteristic~$p$ (in terms of generating ``functions" from
an appropriate space $\cF$, see \eqref{functions}):
\begin{equation*}\label{sbf}
\fsb(n;\un)=\Span(f\in\cF\,|\,\Delta(f)=0).
\end{equation*}

\subsubsection[$\fb_{a,b}$ for $p=2$]{$\boldsymbol{\fb_{a,b}}$ for $\boldsymbol{p=2}$}\label{bab} The direct analog of trace on $\fm_0$
is $\htr$. On $\fle$, the prolong of $\htr$ is the operator $\Delta$. But
$\Delta$ does not commute with the whole of $\fm_-$. To obtain the
$\fm_-$-invariant prolong of this trace on $\fm_0$, we have to
express $\htr$ in terms of the operators commuting with~$\fm_-$
(\textit{$Y$-type vectors} in terms of~\cite{Shch}). Taking $\fm_-$
spanned by the elements
\[
\fm_{-2}=\Kee\cdot \del_\tau, \quad
 \fm_{-1}=\Span(\del_{q_i}+\xi_i\del_\tau,\;
 \del_{\xi_i})_{i=1}^n,
\]
we see that the operators commuting with $\fm_-$ are spanned by
\[
\del_\tau, \quad \del_{q_i}, \quad \del_{\xi_i}+q_i\del_\tau.
\]
In terms of these operators, the vector field $M_f$ takes the form:\index{$M_f$}
\begin{equation}\label{M_f1}
M_f= f\del_\tau+\mathop{\sum}\limits_i\big(
\del_{q_i}(f)(\del_{\xi_i}+q_i\del_\tau)+ (\del_{\xi_i}+
q_i\del_\tau)(f)\del_{q_i}\big)
\end{equation}
and the invariant prolong of $\htr$~-- the direct analog of \textit{divergence}~-- takes the form:\index{((ZZDelta^\fm$@$\Delta^\fm$}\index{$E_q=\mathop{\sum}\limits q_i\del_{q_i}$}

\begin{equation}\label{DeltaM_f1}
\Delta^\fm(f)=\mathop{\sum}\limits_i((\del_{\xi_i}+
q_i\del_\tau)\del_{q_i}(f)=\Delta(f)+E_q\del_\tau(f), \quad \text{where} \ E_q:=\mathop{\sum}\limits_i q_i\del_{q_i}.
\end{equation}
The condition $\Delta^\fm(f)=0$ singles out the $p=2$ analog of
$\fsm$, whereas the condition
\begin{equation}\label{DeltaM_f22}
b\del_\tau(f)+ a \Delta^\fm(f)=(b+aE_q)\del_\tau(f)+ a \Delta(f)=0
\end{equation}
singles out the $p=2$ analog of $\fb_{a,b}$, cf.~\eqref{b_a,b'}.

Setting \index{$pomathfrak_{a,b}(2n;\widetilde\un)$@$\fpo_{a,b}(2n;\widetilde\un)$}
\[
\fpo_{a,b}\big(2n;\widetilde\un\big):=\textbf{F}(\fb_{a,b}(n;\un))
\]
we single out a~subalgebra in the Lie algebra of contact vector
fields \textit{which has no analogs for $p\neq 2$}.

Let us figure out how the parameter $\lambda$ of the regrading
$\fpo_{\lambda}(2n;\un):=\textbf{F}(\fb_{\lambda}(n;\un))$ depends
on parameters $a$, $b$ above; for a~summary, see $N=6, 7$ in Table~\ref{cts-prol}.
The space of $\fb_{a,b}(n;\un)$ consists of
vector fields~\eqref{M_f1} whose generating functions satisfy equation~\eqref{DeltaM_f22}; the
regrading
\[
\deg \tau=\deg q_i=1, \quad \deg \xi_i=0\text{~~for all $i$}
\]
turns $\fb_{a,b}(n;\un)$ into the Lie superalgebra $\fb_{\lambda}(n;\un; n)$\index{$bmathfrak_{\lambda}(n;\un; n)$@$\fb_{\lambda}(n;\un; n)$}\index{$pomathfrak_{\lambda}(2n+1;\un)$@$\fpo_{\lambda}(2n+1;\un)$} whose 0th
component is isomorphic to $\fvect(0|n)$ and the $(-1)$st component
is isomorphic to the $\fvect(0|n)$-module $\Vol^\lambda$ of weighted
$\lambda$-densities.
Set
\[
\fpo_{\lambda}\big(2n+1;\widetilde\un\big):=\textbf{F}(\fb_{\lambda}(n;\un; n)).
\]

To express $\lambda$ in terms of the parameters
$a,b$, we take an element in the 0th component of
$\fb_{a,b}(n;\un)$ not lying in $\fb(n;\un)$ and see
how it acts on $M_1$.

Let $f=\alpha \tau+\beta q_1\xi_1$. Equation~\eqref{DeltaM_f22} implies that $\alpha b+\beta a=0$, so we can take $f=a
\tau+b q_1\xi_1$. Observe that $M_{b q_1\xi_1}$ acts on the $(-1)$st
component as the vector field $D=b\xi_1\partial_{\xi_1}$ and
\begin{equation*}\label{lambdaIab}
{}[M_{a\tau}, M_1] = [a\tau \partial_{\tau},
\partial_{\tau}] = a\partial_{\tau} =\nfrac ab (\Div D) M_1,
\quad \text{hence $\lambda=\nfrac ab$}.
\end{equation*}

The $p=2$ version of equation~\eqref{specValuesOfLamb} are the following
exact sequences that single out the simple Lie superalgebras (recall
that $\fb(n;\un)=\fb_{\lambda}(n;\un)$ for $\lambda=0$):
\begin{gather}
0\tto \Kee M_{1} \tto \fb(n;\un)\tto \fle(n;\un)\tto 0,\nonumber\\
0\tto \fb_{1}^{(1)}(n;\un)\tto \fb_{1}(n;\un)\tto \Kee\cdot M_{\xi_1\dots\xi_n} \tto 0,\nonumber\\
0\tto \fb_{\infty}^{(1)}(n;\un)\tto \fb_{\infty}(n;\un)\tto \Kee\cdot M_{\tau\xi_1\dots\xi_n} \tto 0.\label{specValOfLambForP2}
\end{gather}

The following statement follows directly.

\begin{Proposition}[two exceptional deforms of the Poisson algebra]\label{Desup} Desuperization of the simple Lie superalgebras introduced in~\eqref{specValOfLambForP2} yields $2$ exceptional serial Lie algebras that have no analogs for $p\neq 2$:
 \index{$pomathfrak_{\infty}^{(1)}(n;\un):=\textbf{F}\big(\fb_{\infty}^{(1)}(n;\un)\big)$@$\fpo_{\infty}^{(1)}(n;\un):=\textbf{F}\big(\fb_{\infty}^{(1)}(n;\un)\big)$}\index{$pomathfrak_{1}^{(1)}(n;\un):=\textbf{F}\big(\fb_{1}^{(1)}(n;\un)\big)$@$\fpo_{1}^{(1)}(n;\un):=\textbf{F}\big(\fb_{1}^{(1)}(n;\un)\big)$}
\[
\fpo_{1}^{(1)}(n;\un):=\textbf{F}\big(\fb_{1}^{(1)}(n;\un)\big)\quad \text{and} \quad \fpo_{\infty}^{(1)}(n;\un):=\textbf{F}\big(\fb_{\infty}^{(1)}(n;\un)\big).
\]
\end{Proposition}

\begin{Problem}\label{OP}\index{(Problem, open@Problem, open} In~{\rm \cite{LSh}}, we described the deformations of the Buttin algebra $\fb(n)$ over~$\Cee$, having corrected a~result due to Kochetkov~{\rm \cite{Koch}} who described exceptional deformations of~$\fb_{\lambda}(n)$ at certain values of $\lambda$; some of these deformations having an odd parameter. It is an open problem to obtain a~version of~{\rm \cite{LSh}} in the modular case.
\end{Problem}

\subsection[On $\fk$ and $\fpo$ for $p=2$]{On $\boldsymbol{\fk}$ and $\boldsymbol{\fpo}$ for $\boldsymbol{p=2}$}\label{kandpoab} Observe that $\fpo_I(n;\un)$, a~central extension of the Lie algebra $\fh_I(n;\un)$, \textit{is not a~Lie algebra}.\index{$hmathfrak_I(n;\un)$@$\fh_I(n;\un)$} Indeed, the bracket should be anti-symmetric, i.e.,
alternate, while $\{x_i,x_i\}_{I}=1$, not 0, so $\fpo_I(n;\un)$\index{$pomathfrak_I(n;\un)$@$\fpo_I(n;\un)$} is a~\textit{Leibniz} algebra,\index{(Leibniz (super)algebra@Leibniz (super)algebra} not Lie algebra. Only $\fh_\Pi(2n;\un)$
has an analog of the familiar central extension; \textit{this} nontrivial central extension is a~correct
direct analog of the complex Poisson Lie (super)algebra
$\fpo(2n|0)$.

To pass from $\fpo(0|n)$ to $\fk(1;\un|2n)$, we have to add, as a
direct summand, a~central element to $\fpo(0|n)_0=\fo_\Pi^{(1)}(n)$;
it will act
on the prolong of $\big(\fpo(0|n)_-, \fc\fo_\Pi^{(1)}(n)\big)$ as a~grading operator. We see that the generalized
Cartan prolong of $\big(\fpo(0|n)_-, \fc\fo_\Pi^{(1)}(n)\big)$ is
$\fk(1;\un|2n)$.

$\bullet$ The commutant of $\fk_0(1;\un|2n)$ is isomorphic to that of $\fpo_0(0|n)=\fo_\Pi^{(1)}(n)$,
so it is of codimension 2 in $\fk_0(1;\un|2n)$. Thus, there are \textit{two} traces on
$\fk_0(1;\un|2n)$, and hence there are \textit{two} divergences on $\fk(1;\un|2n)$,
like on $\fm(n;\underline{L}|n)$. These divergences are given by almost the same formulas as for $\fm(n;\underline{L}|n)$, where $\underline{L}=(N,\One)$ and $2^N$ is the hight of $t$, ``almost'' because $\partial_t$ should replace $\partial_\tau$.\index{(Divergence@Divergence}

For $p\neq 2$, we have two divergences on $\fk(1;\un|2n)$ only if $n=1$, see Section~\ref{SSdivFreeK12}.

\subsection[Exceptional simple vectorial Lie superalgebras for $p=2$ analogous\\
to their namesakes over $\Cee$]{Exceptional simple vectorial Lie superalgebras for $\boldsymbol{p=2}$\\ analogous
to their namesakes over $\boldsymbol{\Cee}$}\label{ssExDef} We give detailed description of all exceptional simple vectorial Lie superalgebras over $\Cee$ and fields of characteristic~2 in the main text; for a~summary, see Section~\ref{STables}. These Lie superalgebras constitute two non-intersecting sets as follows.

\underline{The
complete Cartan prolong of its \textit{negative} part}: such is every Lie superalgebra of series $\fvect$, $\fk$ and $\fm$ in the standard
grading, see \eqref{nonstandgr} and each simple exceptional Lie
superalgebra~$\fg$ of depth $>1$, whose negative part in its
W-grading is different from the negative part of the Lie superalgebras
of series $\fk$ or $\fm$ in their respective standard gradings.

\underline{The
complete Cartan prolong of its \textit{nonpositive} part}: such are the
exceptional vectorial Lie superalgebras, and their desuperizations, see
Tables~\ref{excepC} and~\ref{forgetexcept1}, other than in the above paragraph;
the corresponding gradings
are explicitly given in Table~\eqref{regrad}.

The desuperizations of two nonisomorphic Lie superalgebras realized by
vector fields on supervarieties of different superdimension might
turn out to be vectorial Lie algebras realized on varieties of the same
dimension. We distinguish these cases by indicating their depths as an index at the name
($\fm\fb_2(11;\un)$ and $\fm\fb_3(11;\un)$, and also
$\fk\fle_2(20;\un)$ and $\fk\fle_3(20;\un)$); for the case of equal
depths, we distinguish non-isomorphic algebras by a~tilde: $\fv\fle(9;\un)$ and $\widetilde{\fv\fle}(9;\underline{L})$, as well as
$\fk\fas(7;\un)$ and $\widetilde{\fk\fas}(7;\underline{L})$; for details, see respective sections.

\subsection{A~technical remark: natural generators of vectorial Lie superalgebras}\label{natGener} This subsection is needed for calculations
only. Let $\fg=\oplus\fg_i$ be a~Weisfeiler
grading of a~given simple vectorial Lie superalgebra. We see that
$\fg_{-1}$ is an irreducible $\fg_0$-module with highest-weight
vector $H$, and $\fg_{1}$ is the direct sum of indecomposable
(sometimes, irreducible) $\fg_0$-modules with lowest-weight vectors~$v_i$.

Over $\Cee$, and over $\Kee$ for $\un=\One$, the simple Lie
superalgebra $\fg$ is generated (bar a~few exceptions) by the
generators of $\fg_0$, the vector $H$, and the $v_i$. (For other
values of $\un$, we have to add the $\fg_0$-lowest-weight vectors
$v_j^k\in\fg_j$ for some $j>1$ to the above generators;
these cases are not considered.) So we have to describe the generators of
$\fg_0$, or rather of its quotient modulo its center.

If $\fg_0$ is of the form $\fg(A)$ or its simple subquotient, we select its
Chevalley generators, see~\cite{BGL1}.

If $\fg_0$ is an almost simple ``lopsided", see Section~\ref{sssChunk} (in particular,
of type $\fpe$, $\fspe$), but $\Zee$-graded Lie superalgebra, we
apply the above-described procedure to $\fg_0$: first, take
\textit{its} 0th components and its generators, then the highest and
lowest-weight vectors in \textit{its} components of degree $\pm1$,
etc.

If $\fg_0$ is semi-simple of the form
$\fs\otimes\Lambda(r)\ltimes\fvect(0|k)$, where $\fs$ is almost
simple, then we take the already described generators of
$\fvect(0|k)$ and apply the above procedures to~$\fs$.

For a~list of defining relations for many simple Lie superalgebras
over~$\Cee$, and their relatives, see~\cite{GL1,GLP}. For defining relations for Lie algebras with Cartan
matrix over~$\Kee$, see~\cite{BGLL1}.

\section{Introduction: overview of the scenery}\label{Sintro}

In the Introduction (divided into two parts to ease digesting it) we give a
brief sketch of the main constructions and ideas; for basic
background, see Section~\ref{Sbackgr}. For further details, see \cite{Leb, LS, LS1}. All voluminous computations are performed with the help of the \textit{SuperLie} package, see \cite{Gr}.

\subsection[Goal: classification of simple finite-dimensional Lie algebras
over $\Kee$ a.k.a.\ \textit{modular}]{Goal: classification of simple finite-dimensional Lie algebras\\
over $\boldsymbol{\Kee}$ a.k.a.\ \textit{modular}}
In 1960s, Kostrikin and Shafarevich suggested a~method
for producing simple\ finite-dimensional Lie algebras over $\Kee$ for any $p>0$,
together with the final list for $p>5$. This list is explicit for simple $\Zee$-graded algebras;
for the rest, it is somewhat implicit (``and deforms of $\Zee$-graded algebras"), see \cite{Kos1}. The above-mentioned deforms are often deforms of non-simple algebras the stock of which was not clearly described; this made this part of the KSh-method rather vague.

\subsubsection{The original KSh-method}\index{(KSh-method@KSh-method} The initial ingredients are simple Lie algebras over
$\Cee$ of two types:
\begin{gather}
 \text{finite-dimensional, i.e., of the form $\fg(A)$,\index{$gmathfrak(A)$ Lie (super)algebra with Cartan matrix $A$@$\protect\mathfrak{g}(A)$ Lie (super)algebra with Cartan matrix $A$} where $A$ is a~Cartan matrix},\label{1} \\
 \text{infinite-dimensional vectorial types ($\fvect$, $\fsvect$, $\fh$, and $\fk$)}\nonumber\\
  \text{with polynomial coefficients.}\label{2.}
\end{gather}
Next, one, respectively,
\begin{gather}
\text{takes a~$\Zee$-form $\fg(A)_{\Zee}$ of $\fg(A)$ corresponding to the Chevalley basis,}\nonumber\\
\text{and tensors with $\Kee$ to get $\fg(A)_{\Zee}\otimes_{\Zee} \Kee$},\label{3}\\
\text{considers modular, over $\Kee$, analogs of simple vectorial Lie algebras over $\Cee$ with divided}\nonumber\\
\text{powers as coefficients of \emph{distinguished} partial derivatives, see equation~\eqref{distPD}.}\label{4.}
\end{gather}
The ingredient \eqref{1} yields via \eqref{3}
one finite-dimensional Lie algebra; the ingredient \eqref{2.} yields via \eqref{4.} an infinite family of finite-dimensional Lie algebras over $\Kee$ depending on the shearing vector $\un$. Each of the finite-dimensional Lie algebras thus obtained is either
simple, or a~``re\-la\-tive" of the~simple Lie algebra over $\Kee$ (a~central extension or a~subalgebra in the algebra of derivations).\index{(Relative of a~simple Lie (super)algebra@Relative of a~simple Lie (super)algebra} Some of these simple Lie algebras can be deformed.
\be\label{vague}
\text{To describe the \textit{deforms} is a~\textit{rather complicated} part of the KSh-method}.
\ee

Let us clarify claim \eqref{vague}. Tables \eqref{dimSvect} and \eqref{dimHam} show that some of simple Lie algebras are filtered deforms not of the simple $\Zee$-graded algebras, but of certain non-simple subalgebras of their Cartan prolongs (since their dimensions differ from those of simple algebras). The list of deforms was obtained in a~roundabout way, avoiding computing the cohomology that describes a~filtered deformation:
\begin{enumerate}\itemsep=0pt
\item[1)] In \cite{W}, Wilson sharpened a~result due to Kac \cite[Proposition~7.2]{Kfil} and classified equivalence classes of volume forms for $p>5$; later, it turned out that the description works for $p>2$. (Earlier Tyurin published a~solution of the same problem~\cite{Tyu}, but got more classes than Wilson: Tyurin missed some equivalences.)
\item[2)] Skryabin \cite{Sk0,Sk1} classified (for $p>2$) all equivalence classes of symplectic forms (Skryabin called them Hamiltonian forms); some of Skryabin's difficult-to-obtain results hold for $p=2$ as well.
\end{enumerate}

{\bf Types of Lie algebras $\fsvect$ described by Tyurin and Wilson \cite{Tyu,W}.}
In the mid-1970s, Kac observed in \cite{Kfil} that the Lie algebra that preserves the volume element of the form $h\vvol$, where $h\in\widehat \cO\big(\widehat \un_i\big)$ is invertible, can be a~subalgebra of $\fvect(m;\un)$ with finite coordinates of $\un$. Let $p>2$ and suppose that
\begin{gather}\label{Ni}
N_1=\dots=N_{t_1}<N_{t_1+1}=\dots=N_{t_2}<\dots<N_{t_{s-1}+1}=\dots=N_{t_s}=N_m.
\end{gather}
The results of Tyurin and Wilson, correct for $p>3$, state that there are only the following three types of non-equivalent classes of volume forms, and hence filtered deforms with parameter $\eps\in\Kee^\times$ of divergence-free algebras preserving them:\index{$svectmathfrak_{1+\eps\bar u}(m; \un)$@$\fsvect_{1+\eps\bar u}(m; \un)$}
\begin{gather}
\fsvect_{h}(m; \un):=\{D\in\fvect(m; \un)\,|\,
L_D(h\vvol)=0\}, \quad \text{where $h$ is one of the following:}\nonumber\\
h=\begin{cases}1,\\
1+\eps\bar u, \text{~~where $\bar u:=\prod \bar u_i$\index{$Zu$@$\bar u$}\index{$Zu_i:=u_i^{(p^{N_i}-1)}$@$\bar
u_i:=u_i^{(p^{N_i}-1)}$} and $\bar
u_i:=u_i^{(p^{N_i}-1)}$,}\\
\exp\Big(\eps u_{t_i}^{\big(p^{\un_{t_i} }\big)}\Big):=\mathop{\sum}\limits_{j\geq
0}\Big(\eps u_{t_i}^{\big(p^{\un_{t_i}}\big) }\Big)^{(j)}\in\widehat \cO\big(\widehat \un_{t_i}\big).
\end{cases}\label{svectH}
\end{gather}
For brevity, set $\fsvect_{\exp_i}(m; \un):=\fsvect_{h}(m; \un)$,\index{$svectmathfrak_{\exp_i}(m; \un)$@$\fsvect_{\exp_i}(m; \un)$} where
$h=\exp\Big(\eps u_{t_i}^{\big(p^{\un_{t_i}}\big)}\Big)$.

\begin{Remarks}\quad
\begin{enumerate}\itemsep=0pt
\item For $p=3$ and 2, these deformations of $\fsvect$ are also possible. For $p=3$, nobody knows if there are other deforms, whereas for $p=2$, there definitely is at least one more deform: its existence is the most spectacular result of this paper, see Section~\ref{S_F(mb)}.
\item S.~Tyurin described the Lie algebras of divergence-free type and got an extra type of volume form, as compared with Wilson's list~\eqref{svectH}, cf.~\cite{Tyu}.
\item S.~Kirillov~\cite{Kir} verified Skryabin's remark in passing ~\cite{Sk2} for which $i$ the $i$th derived algebra from Wilson's list~\eqref{svectH} is simple, and found the dimensions of these simple Lie algebras:
\begin{gather}\label{dimSvect}
\begin{tabular}{|l|l|}
\hline
$\dim\fsvect_{\exp_j}(m; \un)=(m-1)p^{|\un|}$&$i=0$\tsep{3pt}\\
$\dim\fsvect_{1+\bar u}^{(1)}(m; \un)=\dim\fsvect^{(1)}(m; \un)=(m-1)p^{|\un|}-m+1$&$i=1$\\
\hline
\end{tabular}
\end{gather}
\end{enumerate}
\end{Remarks}

{\bf Hamiltonian Lie algebras $\fh$ described and classified by S.~Skryabin \cite{Sk0,Sk1}.} Let $\fh(2k;\un)$ or $\fh_{\omega_{0}}(2k;\un)$\index{$hmathfrak_{\omega_{0}}(2k;\un)$@$\fh_{\omega_{0}}(2k;\un)$} be the $\Zee$-graded Lie algebra preserving the symplectic form\index{((ZZZomega_0@$\omega_0$}
\[
{\omega_0 = \mathop{\sum}\limits_{1\leq i\leq k}d u_i\wedge du_{k+i}}.
\]
The only non-isomorphic filtered deforms of $\fh_{\omega_{0}}(2k;\un)$ with parameter $\eps\in\Kee^\times$ are $\fh_{\omega_{i}}(2k;\un)$, where $i=1,2$, which preserve the following respective forms (of \textit{type $1$ and $2$} in Skryabin's terminology):\index{((ZZZomega_{2,j}@$\omega_{2,j}$}\index{((ZZZomega_{1, A}@$\omega_{1, A}$}\index{$hmathfrak_{\omega_{1, A}}(2k;\un)$@$\fh_{\omega_{1, A}}(2k;\un)$}\index{$hmathfrak_{\omega_{2,j}}(2k;\un)$@$\fh_{\omega_{2,j}}(2k;\un)$}
\begin{gather}
\omega_{2,j} = d\left(\exp (\eps u_j )\mathop{\sum}\limits_{1\leq i\leq k}u_idu_{k+i}\right), \text{where $j = t_1,\dots , t_s$, see \eqref{Ni}},\label{ham}\\
\omega_{1, A} = \omega_0 + \eps\mathop{\sum}\limits_{1\leq i, j\leq 2k} A_{i,j}d(\bar u_i )\wedge d(\bar u_{j}),\text{~~where $\bar
u_i:=u_i^{(p^{N_i})}$ for the shearing vector $\un$},\nonumber
\end{gather}
and where the non-equivalent normal shapes of the indecomposable matrices $A=(A_{i,j})$ can only be equal for $p>2$ to one of the following:\index{$J_k(0)$}\index{$J_{k,r}(\lambda)$}\index{$C_{k}$}
\begin{equation*}
\renewcommand{\arraystretch}{1.4}
\begin{tabular}{|l|l|l|}\hline
type of $A$&form of $A$&detailed notation of $\omega_1$\\
\hline
$J_k(0)$&$\antidiag\big(J_k(0),
-J_k(0)^{\rm T}\big)$&$\omega_{1,0}$ for $k>1$\\
\hline
$J_{k,r}(\lambda)$, where $\lambda\neq 0$&$\antidiag\big(J_{k,r}(\lambda), -(J_{k,r}(\lambda))^{\rm T}\big)$&$\omega_{1,r,\lambda}$ for $k=r n$ for $r,n\geq 1$\\
\hline
$C_{k}$&$\antidiag\big(C_{k}, -C_{k}^{\rm T}\big)$&$\omega_{1,C}$ for $k>1$ \\
\hline
\end{tabular}
\end{equation*}
where $J_k(\lambda)$ is a~Jordan $k\times k$ block with eigenvalue $\lambda$, and $J_{k,r}(\lambda)$ is a~$k\times k$ block matrix with blocks of size $r\times r$, so $k=r\times n$ for some $r,n\geq 1$:
\begin{gather*}
J_{k,r}(\lambda)=\begin{pmatrix}0_r&1_r&&0\\
&&\ddots&\\
0&0& &1_r\\
J_r(\lambda)&0&&0_r\\
\end{pmatrix},
\end{gather*}
and
\begin{gather*} C_{k}=\begin{pmatrix}0&1&&0\\
&&\ddots&\\
0&0& &1\\
1&0&&0\\
\end{pmatrix} \ \text{is of size $k\times k$ for $k>1$.}
\end{gather*}

\textbf{The two conditions on $\boldsymbol{J_{k,r}(\lambda)}$ and $\boldsymbol{C_k}$}. 1) The case with $J_{k,r}(\lambda)$ occurs only when
\begin{gather}\label{=}
N_1 + \dots + N_{nr} = N_{nr+1} + \dots + N_{2nr} \quad \text{(recall that $k=r n$)}
\end{gather}
and, furthermore, $N_{ir - j } = N_{ir}$ for all $i = 1,\dots, 2n$ and
all $j = 1,\dots, r-1$; i.e., $r$ indeterminates in
each of the $2n$ successive groups have equal heights.

The case with $C_k$ occurs only when condition \eqref{=} is violated.

2) Let $G$ be the group generated by the cyclic
permutations of the row vectors of length~$k$. Then, the
identity element is the only permutation in $G$ that simultaneously fixes the two vectors
\[
a=(N_1,\dots,N_k) \quad \text{and} \quad  b=(N_{k+1},\dots,N_{2k}) .
\]
 It suffices to consider representatives of equivalence classes of pairs $(a,b)$ under the $G$-action.

\begin{Remarks} \quad
\begin{enumerate}\itemsep=0pt
\item Over $\Cee$ the supervarieties of parameters of deformations of Poisson and Hamiltonian Lie superalgebras can differ, see \cite{LSh}. For $p=2$, there is at least one new type of deform: a~1-parametric family of non-isomorphic deforms different from the above~-- desuperisations of $\fb_{a,b}(n;\un)$.\index{$bmathfrak_{a,b}(n;\un)$@$\fb_{a,b}(n;\un)$}

\item S.~Kirillov~\cite{Kir} determined the $i$ for which the $i$th derived algebra of the Hamiltonian Lie algebra from Skryabin's list \cite{Sk0} is simple and what its dimension is equal to:
\begin{gather}\label{dimHam}
\dim\fh_\omega^{(i)}(2k;\un)=\left\{\begin{array}{@{}lll@{}}p^{|\un|}-2&\text{if~}\omega=\omega_0&i=2,\\
 p^{|\un|}&\text{if~}\omega=\omega_2, \text{where $k+1\not\equiv 0\pmod p$}&i=0,\\
 p^{|\un|}-1&\text{if~}\omega=\omega_2, \text{where $k+1\equiv 0\pmod p$}&i=1,\\
 p^{|\un|}-1&\text{if~}\omega=\omega_{1, A}, \text{where $\det A\neq0$}&i=1,\\
p^{|\un|}-2&\text{if~}\omega=\omega_{1, A}, \text{where $\det A=0$} & \\
& \text{(type $J_s(0)$)}&i=2.\end{array}\right.
\end{gather}
\end{enumerate}
\end{Remarks}

\subsubsection{True and semi-trivial deforms}\label{ssTrue}
In particular, the amount of
infinitesimal deformations is overwhelming and even frightening as~$p$
becomes small ($p=3$ or~-- a~horrible case~-- ${p=2}$). We recall reasons not to be too frightened;
besides, the KSh-method had been considerably improved over the past years.

The abundance of deforms of simple
Lie (super)algebras for $p>0$, especially overwhelming for $p=2$, is
somewhat misleading. It is occasioned by \textit{semi-trivial
deforms}\index{(Deform, semi-trivial@Deform, semi-trivial} each of which is given by a~cocycle representing a
nontrivial cohomology class but, though inte\-grable, yields a~deform isomorphic to the initial algebra. For a~description of many
semi-trivial deforms, see~\cite{BLLS1}. We say that a~nontrivial and
nonsemi-trivial deform is a~\textit{true deform}.\index{(Deform, true@Deform, true}

The Lie (super)algebra $\fg$ is said to be \textit{rigid} if it has no true
deforms; until recently we thought that semi-trivial deforms existed
only if $p>0$, but a~more careful study of the literature shows they
are a~universal phenomenon~\cite{Ri}.

\textit{If $p>3$, the classification has been completed}, mainly due to Premet and
Strade~\cite{BGP, S}, based on explicit description of
deforms \cite{BW, Sk1, W}.

\textit{If $p=3$, we conjecture the classification}: the examples
obtained by Cartan prolongation (see Section~\ref{CartProl}) of appropriate
parts of Lie algebras with Cartan matrix~\cite{BGLLS, GL4}, \index{(Lie algebra, standard@Lie algebra, standard}
exhaust the list of ``standard'' examples some of which were discovered by Frank, Ermolaev and, mainly, Skryabin. For an (incomplete at the moment) list of true deforms of several
``standard'' algebras, see~\cite{BGL4, LaF, LaY, LaZ, Sk1}, and
\cite{BLW} in which an earlier claim concerning deforms is corrected.

\textit{If $p=2$, we are still completing the stock of ``standard"
examples.}

\subsubsection{Amendments to the formulation of the goal} On several
occasions P.~Deligne told us what we understood as follows (for Deligne's
own words, and several open problems, see \cite{LL}):
\begin{equation*}\label{Deligne}
\begin{minipage}[c]{14cm}
``In positive characteristic, the problem ``classify ALL simple Lie (super)algebras,
and their representations'' is, perhaps, not very reasonable, and
definitely very tough; investigate first the \textit{restricted}
case related to geometry, and hence meaningful.''
\end{minipage}
\end{equation*}
Following Deligne's advice, we investigated several plausible notions of restrictedness for $p=2$ in~\cite{BLLS2} and gave explicit expressions of the restriction maps for several types of simple Lie algebras and superalgebras in \cite{BKLLS}. Nevertheless, even to describe restricted Lie (super)algebras one often
needs nonrestricted ones; for more serious examples of their usage, see \cite{Kos1}.

In this paper, we concentrate on \textit{simple} Lie (super)algebras,
keeping in mind that algebras of the following types are no less
important than simple ones:
\begin{itemize}\itemsep=0pt
\item Lie (super)algebras of the form $\fg(A)$ where $A$ is
indecomposable, see \cite{BGLLS, BGL1}.

\item Central extensions and algebras of derivations of the known
simple Lie (super)algebras (see Section~\ref{ssNonSemiSvect} and
\cite{BGLL3}). The algebras of these two types will be called
\textit{relatives}\index{(Relative of a~simple Lie (super)algebra@Relative of a~simple Lie (super)algebra} of the corresponding simple Lie (super)algebras and each other.

\item The generalized Cartan
prolongs $(\fg_-,\fg_0)_{*,\un}$ with $\fg_0$ close to simple, see Sections~\ref{CartProl} and~\ref{ssNonSemiSvect}.

\item True deforms (for definition, see Section~\ref{ssTrue}) of Lie (super)al\-gebras,
see \cite{KD,BGL4}.

\item Restricted closures of nonrestricted simple Lie algebras.
\end{itemize}

\subsection{Improvements of the KSh-method}\index{(KSh-method@KSh-method}
\subsubsection{``Standard'' modular Lie algebras}\label{ssStand} Dzhumadildaev and Kostrikin \cite{KD} suggested simplifying the KSh-method by
skipping the step over $\Cee$ and considering certain ``standard''\index{(Lie algebra, standard@Lie algebra, standard}
modular Lie algebras from the very beginning, further deforming them and
their ``relatives". On the other hand, the stock of ``standard''
examples should include, if $p<7$, certain non-simple Lie algebras,
see \cite{GZ, KD, SkT1}. The snag is: we have no idea how to select them.

Until the year 2000 or so, it was believed that the
initial KSh-method produces all simple Lie algebras only if $p>5$.
This belief was based on insufficient study of deformations and too
narrow a~choice of ``standard'' examples: as shown in \cite{KD}, the
Melikyan algebra,\index{(Melikyan algebra@Melikyan algebra} indigenous for $p=5$, are \textit{deforms} of
Poisson Lie algebra which should be considered ``standard" and processed via the KSh-scheme \eqref{1}--\eqref{vague}.

{\bf What examples should qualify as ``standard''?} 
In \cite{Lei}, the improvement of the
KSh-method suggested in \cite{KD} was developed further by
eliminating the vectorial simple Lie algebras from the input of the
KSh-method thus diminishing the stock of ``standard'' simple Lie
algebras. In the new procedure, the role of \textit{generalized Cartan prolongation}
(complete or partial), see Section~\ref{CartProl} and especially Section~\ref{CartProlPart}, becomes even more important than in the
KSh-procedure. This approach definitely works for $p>3$, and conjecturally works for $p=3$.

The stock of ``standard'' (not necessarily simple) Lie algebras must be enlarged with
examples found after~\cite{KD} was published; for $p=3$, see \cite{GL4}; for $p=2$, see \cite{BGLLS, BLLS1, Ei, GZ,Leb, SkT1}, and this paper.

\subsubsection{Splitting the problem into smaller chunks}\label{sssChunk} All simple Lie
algebras are of the following two types: the root system of a
``\textit{symmetric}'' algebra contains the root $-\sigma$ of the
same multiplicity as that of $\sigma$ for any root $\sigma$; the
algebras with root systems without this property are said to be
``\textit{lopsided}''.\index{(Lie (super)algebra, lopsided@Lie (super)algebra, lopsided} \index{(Lie (super)algebra, symmetric@Lie (super)algebra, symmetric}

This paper is devoted to the study
of lopsided algebras, but ``symmetric" Lie (super)algebra will be needed in the process

\textit{Symmetric algebras}. A~significant quantity of symmetric simple Lie
algebras consists of algebras \textit{$\fg(A)$ with indecomposable
Cartan matrix $A$ or simple ``relatives'' of such algebras} of the form
$\fg^{(i)}(A)/\fc$, where\footnote{This is shorter and more
graphic than the correct notation $(\fg(A))^{(i)}$; usually, we will
similarly place subscripts designating the degree (closer to the ``family name" $\fg$).}
$\fg^{(i)}(A)$ is the $i$th derived algebra of $\fg(A)$ and $\fc$ is the
center of $\fg^{(i)}(A)$.

\textit{For any $p$}, finite-dimensional Lie algebras $\fg(A)$ with
indecomposable Cartan matrix $A$, and their simple relatives, were
classified in \cite{WK} with an omission; for corrections, see
\cite{KWK, Sk2}, where no claim was made that these were the only
corrections needed; for this claim with a~proof, a~classification of Lie
superalgebras of the form $\fg(A)$ with indecomposable Cartan matrix
$A$, and their simple relatives, and precise definitions of related
notions, see \cite{BGL1}.

\textit{Lopsided algebras}: the set they constitute is a~virtually virgin territory a~part of
which~-- vectorial Lie (super)algebras~-- we investigate through the whole of this text.

\subsubsection{Cartan prolongations of Lie (super)algebras with Cartan
matrix}\label{ssCartProlG(A)} It turns out that every known $\Zee$-graded simple Lie algebra for $p>2$ is obtained as a~(generalized, perhaps partial) Cartan prolong of the non-positive part of a~Lie algebra $\fg(A)$ with an indecomposable Cartan matrix $A$. For $p>3$ this follows from the classification.

\begin{Conjecture} The above $\Zee$-graded simple Lie algebras and simple relatives of their deforms constitute the list of simple finite-dimensional Lie algebras for $p=3$ as well.
\end{Conjecture}

\subsection{Super goal} Although Lie \textit{super}algebras appeared in
topology in the 1940s (over finite fields, often over $\Zee/2$), the
understanding of their importance dawned only in the 1970s, thanks to
their applications in physics. This understanding put the problem ``classify simple Lie
superalgebras" on the agenda of researchers. Over $\Cee$,
the finite-dimensional simple Lie superalgebras were classified by
several teams of researchers, see reviews \cite{K2, Kapp}. The
classification of certain types of simple
vectorial Lie superalgebras was explicitly announced in \cite{K2},
together with a~conjecture listing all
\textit{primitive} vectorial Lie superalgebras; for the first
counterexamples, see~\cite{ALSh, Le2}.

A~classification of the simple vectorial Lie superalgebras over
$\Cee$ was implicitly announced when the first \textit{exceptional}
examples were given \cite{Sh3,Sh5,Sh14} and explicitly at a
conference in honor of Buchsbaum~\cite{LS0}. The claim of~\cite{K2} was corrected in
\cite{LS} (the correction contained both the complete list of simple vectorial
algebras, bar one exception later described in~\cite{Sh5}, and the
method of classification of simple $\Zee$-graded Lie superalgebras
of depth 1) and in a~series of papers \cite{CCK, CK1, CK2, CK1a,K, K10}, where the proof in the case
of $\Zee$-grading compatible with parity was given; for further
corrections and proofs, see Section~\ref{Smb}. The classification is not completed till today: there is no
classification of deformations with odd parameters. 

Although to complete the classification of the simple
finite-dimensional Lie superalgebras over $\Kee$ for $p$
``sufficiently big" (say $>7$) will be a~more cumbersome and
excruciating task than that for Lie algebras, the answer
(conjectural, but doubtless) is obvious: to get \textit{restricted}
superalgebras, take the obvious modular analogs of the complex
simple Lie superalgebras (of both finite-dimensional and
of infinite-dimensional vectorial considered for the shearing vector $\un=\One$,
see definition~\eqref{un_s}) passing to the derived algebra and quotients modulo center
if needed; to get \textit{nonrestricted} superalgebras,
consider \textit{true deforms}, see Section~\ref{ssTrue}, of the
above-mentioned analogs (for $\un$ unconstrained, speaking about
vectorial algebras). For $p$ ``small'',
the classification problem becomes more and more involved, see, e.g., \cite{BGL3, ILL}.
Nevertheless, in the two cases the classification is obtained:

\begin{itemize}\itemsep=0pt
\item \textit{For any $p$, the super goal is reached for Lie
superalgebras of the form $\fg(A)$ with in\-de\-composable Cartan
matrices $A$ or its ``relative"}, see \cite{BGL1}. Either
$\fg=\fg(A)$ or its ``relative" of the form $\fg^{(i)}/\fc$, where
$\fg^{(i)}$ is the $i$th derived algebra of $\fg$ and $\fc$ its center, is
simple. For deforms, see~\cite{BGL4}.

\item \textit{Amazingly, the super goal is reached if $p=2$},
see \cite{BLLS2}, with a~catch: \textit{modulo} the classification
of simple Lie algebras, i.e., without an explicit list of all examples.\footnote{Thus, \cite{BLLS2} resembles the classification of restricted Lie algebras for $p>5$ in the paper \cite{BW}: there are no explicit formulas for $p$-structures of simple Lie algebras of Hamiltonian series to this day, see Strade's lamentations in~\cite[Vol.~1, p.~357]{S}:
 ``The problem of restrictedness is approached. \dots\ [But] the family of Hamiltonian algebras~\dots\  is not yet handable''. This is no wonder: although Skryabin classified symplectic forms in 1985, the answer was published only in 1991, see~\cite{Sk1}, three years after~\cite{BW} appeared (and the details of~\cite{Sk1}, obtained in 1985, became available only recently, see~\cite{Sk0}). The explicit formula for the bracket in any of the deformed Lie algebra of Hamiltonian vector fields is not published to this day and can be found only in Kirillov's Ph.D.~Thesis only (in Russian), not in its published summary \cite{Kir}.}
Here, we contribute to a~conjectural list of ``standard'' simple Lie
algebras (conjecturally a~tame problem); in particular, we
explicitly describe simple vectorial Lie algebras analogous to those
over $\Cee$.

\item $p>5$. The conjecture~\cite{Lei} is easy to formulate (``take direct characteristic-$p$ versions of the simple complex Lie superalgebras of the form $\fg(A)$, queer, and vectorial with polynomial coefficients, and their deformations"), but to describe deformations of even the symmetric ones (of the form $\fg(A)$ and queer) is not easy (for partial results, see \cite{BGL4}).

\item $p=5$. Most plausible conjecture is like for $p>5$. Observe that there are indigenous $p=5$ examples of the form $\fg(A)$,

\item $p=3$. We have discovered several new vectorial Lie superalgebras, see \cite{BGL4, BL}, and of the form $\fg(A)$, see~\cite{BGL1}.
\end{itemize}

\subsection[Getting simple Lie algebras from simple Lie superalgebras if $p=2$]{Getting simple Lie algebras from simple Lie superalgebras if $\boldsymbol{p=2}$}

If $p=2$, there are two methods for constructing
a simple Lie superalgebra from a~simple Lie algebra, and every simple Lie superalgebra is obtained by
one of these two methods; for an amazingly short proof, see \cite{BLLS2}. Reversing the process we recover a~simple Lie algebra given any simple Lie superalgebra.

Even before these two methods were known, it was clear that one can get a~Lie algebra from any Lie superalgebra as follows. Observe that for any odd element $x\in \fg$ in any Lie superalgebra~$\fg$ over any field $\Kee$, we have $[x,x]:=2x^2\in U(\fg)$. That is why if $p=2$, then one needs a~squaring $x\longmapsto x^2$ for any odd $x\in\fg$;
together with the brackets of even elements with all other
elements, it is the squaring that defines the multiplication in any Lie superalgebra (for
details, see Section~\ref{defLieS2}), while the bracket of odd
elements is the polarization of the squaring. Hence,
\begin{equation}\label{KLobs}
\begin{minipage}[c]{14cm}
\textit{For $p=2$, every Lie superalgebra with the bracket as
multiplication~--\\
 we forget the squaring~-- is a~$\Zee/2$-graded Lie algebra}.
\end{minipage}
\end{equation}

To classify simple Lie superalgebras is a~much more difficult task
than to classify simple Lie algebras of the same type: the former is
based on the latter as well as on careful study of the representation theory
of Lie algebras. In \cite{KL}, it was, nevertheless, suggested~-- for $p=2$~-- to reverse
the process:
\begin{equation}\label{KLconj}
\begin{minipage}[l]{14cm}
\textit{Let $\textbf{F}$ be the \textit{desuperization}\index{$F$, the desuperization functor@$\textbf{F}$, the desuperization functor} functor forgetting squaring, see \eqref{KLobs}.}\\
\textit{To obtain simple Lie algebras for $p=2$,\\ $(A)$ apply the functor
$\textbf{F}$ to every simple Lie superalgebra $\fg$;\\ $(B)$ single out the simple Lie subalgebra $\fs(\textbf{F}(\fg))$ of $\textbf{F}(\fg)$}.
\end{minipage}
\end{equation}
Clearly, $\fs(\textbf{F}(\fg))$ is uniquely recoverable by inverting one of the two superization processes (either queerification or ``method 2'') that had lead to $\fg$, see \cite{BLLS2}.

Observe immediately that the idea of~\cite{KL} just to apply $\textbf{F}$ to the simple Lie superalgebra $\fg$ to get a~simple Lie algebra, see~\eqref{KLobs}, was naive and partly wrong: the example of $\fpsq(n)$ should have hinted at importance of item (B) in the process~\eqref{KLconj}. Understanding of this subtlety came together with the description of the two methods of superization of any simple Lie algebra as the only means to obtain any simple Lie superalgebra, see \cite{BLLS2}. For the simple vectorial Lie superalgebras, in particular, exceptional ones, just by forgetting squaring we get a~simple Lie algebra.
Let
\begin{gather*}
\text{$\textbf{F}^{-1}=\fs(\textbf{F}(-))$\index{$F^{-1}$, complete desuperization@$\textbf{F}^{-1}$, complete desuperization}\index{(Desuperization, complete@Desuperization, complete} denote the \textit{complete desuperization}} -
\end{gather*}
the composition of $\textbf{F}$ and application of item (B) of \eqref{KLconj}.

\textit{Two reasons} to take the direction of study opposite to a~seemingly reasonable one:
\begin{enumerate}\itemsep=0pt
\item[(a)] Although the classification of the
simple vectorial Lie superalgebras over $\Cee$ was only conjectured at the time \cite{KL} was written,
the list of known examples was already wider than that of known simple vectorial
Lie algebras for $p=2$, and was
(and \textit{is}, as we demonstrate in this paper) able to provide new simple examples.

\item[(b)] The results of \cite{GZ, SkT1} show that a~``frontal attack'' on the classification for $p=2$
is likely to be much more excruciating than that performed for $p>3$ by Premet and Strade. Even
to classify \textit{restricted} Lie algebras for $p=2$ will be much more
difficult problem than that Block and Wilson solved for $p>5$, see~\cite{BW}. (Even if we confine ourselves to the classical definition of restrictedness, while certain examples, which should be considered as ``classical'', have another version of restrictedness, see~\cite{BLLS2}.) So, a~plausibly complete
inventory of simple examples will be helpful. Our interpretations of the Lie (super)algebras are of independent
interest.
\end{enumerate}

Here, after a~long break, we continue exploring method
\eqref{KLconj}. It provides us with new examples of simple vectorial
Lie algebras of the form $\textbf{F}(\fg)$, where $\fg$ is a
modular, indigenous for $p=2$, version of a~simple vectorial Lie
superalgebra over $\Cee$. The two methods of superization (see~\cite{BLLS2}) applied to $\textbf{F}(\fg)$ bring many more simple Lie superalgebras than $\fg$, most of them new.

\subsection[Forgetting the superstructure if $p=2$]{Forgetting the superstructure if $\boldsymbol{p=2}$} Applying $\textbf{F}$
to the serial vectorial Lie superalgebra $\fg(m;\un|n)$ we get the Lie
algebra $\textbf{F}(\fg)\big(m+n;\widetilde\un\big)$, see Table~\ref{forgetseries}; these Lie algebras are not necessarily simple,
but their simple derived algebras are; here $\widetilde\un=(\un, 1,\dots, 1)$
with the last $n$ coordinates equal to 1.

\subsubsection{Parameters of the Lie superalgebra that change under desuperization} Here are several examples:
\begin{itemize}\itemsep=0pt
\item The unconstrained shearing vector $\widetilde \un^u$, see
Section~\ref{CritCoor}, of the vectorial Lie algebra
$\textbf{F}(\fg)$ may depend on more parameters than the shearing
vector $\un^u$ of $\fg$.

For example, $\dim\un=\Par\un^u=a$ for $\fvect(a;\un|b)$, whereas
\begin{gather*}
\dim\widetilde \un^u:=\Par\widetilde \un^u=a+b,\text{~~where $\widetilde \un:=(\un,1\dots,1)$}
\end{gather*}
for
$\textbf{F}(\fvect(a;\un|b))=\fvect\big(a+b;\widetilde \un\big)$. In all cases, except for $\fvle\big(4; \widetilde\uM|3\big)$, the tilde over any shearing vector $\underline{L}$ is understood in the above sense: it enlarges the set of coordinates of $\underline{L}$ acquiring the coordinates of the desuperized odd indeterminates.

The same applies to the desuperizations of the Lie
superalgebras of the series $\fk$, $\fh$, $\fm$, $\fle$ and their divergence-free subalgebras.

The abstract Lie superalgebra $\fg$ realized as vectorial Lie superalgebra,
$\fg(a;\un|b)$, depending on $a$ even and $b$ odd indeterminates,
can be realized in several ways as $\fg(a;\un|b;r)$ by means of
Weisfeiler filtrations or associated regradings $r$, see
Section~\ref{WeisF}. This $\fg(a;\un|b;r)$ can be interpreted as
the (generalized) Cartan prolong of the nonpositive
part of $\fg$ in the corresponding grading, see Section~\ref{CartProlPart}.

\item The Lie algebra obtained by desuperization might acquire
new properties which its namesakes for $p\neq2$ do not have. For
example, the Lie algebra $\fpo(2n;\un)=\textbf{F}(\fb_{\lambda}(n;\allowbreak n;\un))$ has a~deformation
depending on a~parameter
$\lambda\in\Kee P^1$; the corresponding non-isomorphic (for different values of $\lambda$) deforms are additional to the
well-known one, which for $p=0$ is called the result of the \textit{quantization}.\index{(Quantization@Quantization}

\item Desuperization \textbf{F} might turn distinct (types of)
Lie superalgebras into one (type of) Lie algebras:
\begin{itemize}\itemsep=0pt
\item[-- ]the Lie superalgebras of types $\fk$ and $\fm$ in the standard
grading \eqref{nonstandgr} turn into $\fk$;

\item[--] the Lie superalgebras of types $\fh_\Pi$\index{$hmathfrak_\Pi$@$\fh_\Pi$} and $\fle$ in the standard
grading \eqref{nonstandgr} turn into $\fh_\Pi$.
\end{itemize}
\end{itemize}

\subsection{Comment on Volichenko algebras in characteristic~2} The
notion of\index{(Volichenko algebra@Volichenko algebra}
\textit{Volichenko algebra}\footnote{In memory of I.~Volichenko who
was the first to study inhomogeneous subalgebras in Lie
superalgebras.}, which is an inhomogeneous (relative to parity) subspace
$\fh\subset\fg$ of a~Lie superalgebra $\fg$ closed with respect to
the superbracket in $\fg$, was introduced in~\cite{KL}. For
a~classification of simple (without any nontrivial ideals, both
one-sided and two-sided) finite-dimensional Volichenko algebras over
$\Cee$ under a~certain (hopefully, inessential) technical
assumption, and examples of certain infinite-dimensional algebras,
see~\cite{Lsos2, LSerg}.

The results of \cite{ILMS} suggested that we look at the definition of
the Lie superalgebra for $p=2$, see Section~\ref{defLieS2}, more
carefully. If one does this, it is not difficult to deduce that
\begin{equation}\label{volichISlie} \textit{if $p=2$,
the Volichenko algebras are, actually, Lie algebras}.
\end{equation}

In \cite{ILMS, KL}, the fact \eqref{volichISlie} had not
been understood, and therefore there is no need to consider these papers or Volichenko algebras in characteristic~2 while searching for simple Lie algebras. Unlike desuperizations of Lie superalgebras, which are worth investigating.

\section[Introduction, continued. Our strategy, main results and open problems]{Introduction, continued. Our strategy, main results\\ and open problems}\label{SResults}

\subsection{Generalized Cartan prolongation}\label{SSmainProced} \textit{This is a~principal procedure for
getting vectorial Lie $($super$)$algebras} over $\Cee$, the following fact is well-known \cite{Y}:
\begin{equation}\label{fact}
\begin{minipage}[l]{14cm}
\textit{given a~simple Lie algebra of the form $\fg(A)$, and its
$\Zee$-grading, the generalized prolong of the nonpositive $($with
respect to that grading$)$ part of $\fg(A)$ is isomorphic to~$\fg(A)$,
bar two series of exceptions corresponding to certain simplest
gradings of the embedded algebras~-- $\fsl(n+1)\subset\fvect(n)$
and $\fsp(2n+2)\subset\fk(2n+1)$~-- and the ambients are the
exceptional prolongs.}
\end{minipage}
\end{equation}

In \cite{Lei}, it is shown how to obtain simple Lie algebras over
$\Cee$ of the two types of prime importance for the classification
procedure over~$\Kee$: finite-dimensional and vectorial. Namely, by induction
and using (generalized, in particular, partial) Cartan prolongation,
see Section~\ref{CartProl}. First, one thus obtains all
finite-dimensional simple Lie algebras (each of them has a~Cartan matrix);
during the next step one obtains all four series of simple vectorial Lie algebras, by considering not
only complete
Cartan prolongs as in equation~\eqref{fact}, but also partial ones.

The same method works to obtain $\Zee$-graded simple Lie algebras
for $p$ sufficiently big and with new
standard examples added. After that, there still remain
considerable technical problems: namely, to describe the deforms and to classify
non-isomorphic deforms.

For any characteristic, the super version of classification of simple Lie algebras is much
more complicated than its non-super counterpart: we have to supply
the input with several more types of algebras, but the main
procedures are still the same: generalized, especially partial,
prolongations and deformations.

In several papers (e.g., \cite{BGL4, BGLLS, BL}), we have considered simple Lie (super)algebras, and have
investigated the prolongs of the nonpositive parts relative to their
$\Zee$-gradings with one (or two if $p=2$) pair(s) of Chevalley
generators being of degree $\pm1$, and the other generators being of degree~0.

Here we consider serial and exceptional simple vectorial Lie
superalgebras over $\Cee$ and desuperizations of their analogs for $p=2$.
Realization of a~given Lie (super)algebra $\fg$ in terms of vector fields
implies that $\fg$ is endowed with a~filtration; one of these filtrations, called the
\textit{Weisfeiler filtration}, is the most important, see
equation~\eqref{Wfilt}. Associated with the filtration is the grading; for brevity,
the \textit{Weisfeiler} filtrations and gradings are referred to, respectively, as
\textit{$W$-filtration} and \textit{W-grading}. For several
vectorial Lie algebras over fields of characteristic~$2$, we investigate the following problem answered by the fact
\eqref{fact} over $\Cee$:
\begin{equation*}
\begin{minipage}[l]{15cm}
\textit{when $(\fg_-,\fg_0)_{*,\un^u}\simeq\fg$ and when the prolong strictly contains $\fg$?}
\end{minipage}
\end{equation*}
We consider only the $\Zee$-gradings of the finite-dimensional vectorial algebras corresponding to the
W-gradings of their infinite-dimensional versions corresponding to $\un^u$. For examples of Lie (super)algebras $\fg$ that differ from
the prolong of the nonpositive part of a~regrading of $\fg$, see~\cite{BGLLS}.

\subsection{``Hidden supersymmetries'' of Lie algebras} It is sometimes
possible to endow the space of a~given simple Lie algebra $\fg$ with
(several) Lie superalgebra structures. For example, for
$\fg=\fsl(n)$ over any ground field, consider any distribution of
parities (of the pairs corresponding to positive and negative simple
roots) of the Chevalley generators; we thus get Lie superalgebras
$\fsl(a|b)$ for $a+b=n$ in various supermatrix formats. The sets of
defining relations between the Chevalley generators corresponding to
different formats are different.

It is, clearly, possible to perform such changes of parities of
(pairs of) Chevalley generators for any simple Lie algebra but,
except for $\fsl(n)$, the simple Lie superalgebras obtained by
factorization modulo the ideal of relations~\cite{BGL1} are
infinite-dimensional unless $p=2$.

If $p=2$, the following fact is obvious (here $x\longmapsto x^{[2]}$
is the restriction mapping and $x\longmapsto x^{2}$ is the squaring):\footnote{In \cite{BLLS2}, in addition to the
classical restrictedness, we introduced other, indigenous to
$p=2$, versions of restrictedness; their meaning is yet unclear, but
since they pertain to classical and often used algebras, e.g.,
$\fo(2n+1)$ and $\fh(2n+1;\One)$, we believe that these new
``restrictednesses'' are important.}
\begin{gather}
\textit{any classically restricted and $\Zee/2$-graded
Lie algebra $\fg=\fg_+\oplus\fg_-$ can be turned}\nonumber\\
\textit{into a~Lie
superalgebra by setting $x^2:=x^{[2]}$ for any $x\in\fg_-$.}\label{superiza}
\end{gather}
In~\cite{BLLS2}, a~(rather unexpected) generalization of the
possibility~\eqref{superiza} is described: \textit{every} simple Lie
algebra $\fg$ can be turned into a~simple Lie superalgebra by
slightly enlarging its space if $\fg$ is not restricted. This
generalization, and a~``queerification", are the two methods which,
from every simple Lie algebra, produce simple Lie superalgebras, and
every simple Lie superalgebra can be obtained in one of these two
ways, as proved in~\cite{BLLS2}. These two methods applied to the
simple Lie algebras we describe in this paper yield a~huge quantity
of new simple Lie superalgebras, both serial and exceptional. We do
not list them; this is routine, modulo the far from routine, as shown in~\cite{KrLe}, job of describing all $\Zee/2$-gradings of our newly
found simple Lie algebras.

\subsection[From $\Cee$ to $\Kee$]{From $\boldsymbol{\Cee}$ to $\boldsymbol{\Kee}$} We consider the W-grading (of the desuperized
$p=2$ version of the simple vectorial Lie superalgebra over $\Cee$)
for which the nonpositive part, especially the 0th component, is
most clear. Then, we consider the regradings described in Table~\eqref{regrad}, and perform generalized Cartan prolongation of the
nonpositive (or negative) part of the regraded algebra in the hope
of getting a~new simple Lie (super)algebra, as in~\eqref{fact}.

When this approach is inapplicable since there is no visible analog,
suitable for $p=2$, of the Lie (super)algebras we worked with over
$\Cee$, we consider the description of the Lie superalgebras as the
sum of its even and odd parts, desuperize this description and only
after this consider the W-grading of the desuperization.

\subsection{Our main results} Observe that every modular analog of one vectorial Lie (super)algebra over $\Cee$
is usually a~\textit{family} of algebras depending on $\un$; by abuse of speech we often skip the word
``family" (of algebras) and talk about \textit{one} algebra having in mind the extra parameter~$\un$.

Of the simple Lie superalgebras that can be obtained
by the two methods described in~\cite{BLLS2} from the simple Lie algebras we describe here all but one (the initial one, the one we desuperized) are new. In more details the Lie algebras obtained by desuperization are described in the following Theorems~\ref{Bjthm1}--\ref{Bjthm4} summarizing respective sections with proofs.

\begin{Theorem}[an exceptional deform of $\fsvect(5;\un)$]\label{Bjthm1} All W-regradings of
$\fm\fb_3(11;\un)$ are W-regradings of a~previously unknown true deform of
$\fsvect(5;\un)$.
\end{Theorem}

\begin{proof}
See Section~\ref{S_F(mb)}.
\end{proof}

\begin{Theorem}[desuperizations of $\fb_{a,b}(n)$ and $\widetilde{\fsb}_{\nu}\big(2^{n-1}-1|2^{n-1}\big)$ for $n$ even, as well as\linebreak  $\widetilde{\fsb}_{\nu}\big(2^{n-1}|2^{n-1}-1\big)$ for $n$ odd]\label{Bjthm2} As $($generalized$)$ Cartan
prolongs, the desuperizations of the characteristic-$2$ analogs
of complex Lie superalgebras $\fb_{a,b}(n)$ and $\widetilde{\fsb}_{\nu}\big(2^{n-1}-1|2^{n-1}\big)$ for $n$ even, as well as $\widetilde{\fsb}_{\nu}\big(2^{n-1}|2^{n-1}-1\big)$ for $n$ odd, yield $\widetilde{\fsb}_{\nu}\big(2^{n}-1;\un\big)$ and
\begin{gather*}
\text{$\fpo_{\lambda}(2n+1;\un)$}, \ \text{the serial simple $($for $\frac ab$ generic$)$ Lie algebras $\fpo_{a,b}(2n;\un)$,}\\
\text{their simple relatives for $\frac a
b=0$, $1$, and $\infty$ $($for $a\neq 0, b=0)$, see~\eqref{specValOfLambForP2}}.
\end{gather*}
\end{Theorem}

\begin{proof}
See Sections \ref{bab} and~\ref{StildeSB}.
\end{proof}

\begin{Theorem}[desuperizations of the 15 exceptional simple vectorial Lie superalgebras] \label{Bjthm21} The generalized Cartan prolongs of the nonpositive $($or negative$)$
parts of all $15$ W-gradings, see Table~\eqref{excepC}, of the~$5$
exceptional simple vectorial Lie algebras~-- desuperizations of the characteristic~$2$ analogs
of complex exceptional simple
vectorial Lie superalgebras, see Section~{\rm \ref{ssExDef}}~-- yield three simple Lie algebras, see Table~\eqref{forgetexcept1} while the other two, $\fv\fle(9;\un)$ and
$\fm\fb_3(11;\un)$, are described from another point of view in~{\rm \cite{BGLLS}}.
\end{Theorem}

\begin{Theorem}[isomorphisms between desuperizations of the 15 exceptional simple vectorial Lie superalgebras]\label{Bjthm3}\quad
\begin{enumerate}\itemsep=0pt
\item[$(i)$] For $p\neq 2$, the characteristic-$p$ analogs of the $15$ W-graded analogs of the complex exceptional simple
vectorial Lie superalgebras constitute, as \textit{abstract}
Lie superalgebras, $5$ Lie superalgebras.

\item[$(ii)$]  For $p=2$, all their finite-dimensional analogs, and their desuperizations, remain regradings of each other,
\textit{with one exception}: the $3$ indeterminates of
$\widetilde{\fk\fas}\big(7;\widetilde{\underline{K}}\big)$ are
\textit{defined} to be constrained,
$\Par\widetilde{\underline{K}}^u=3$.
\end{enumerate}

\textit{One -- for $p\neq 2$ -- exceptional family $(\fk\fas)$ yields -- for $p=2$ -- two families}:
\begin{equation*}\label{IsoAlg}\footnotesize
\renewcommand{\arraystretch}{1.4}
\begin{tabular}{|l|l|}
\hline
$\fv\fle\big(7;\widetilde{\underline{L}}\big)\simeq\fv\fle\big(9;\widetilde{\underline{M}}\big)
\simeq\widetilde{\fv\fle}\big(9;\widetilde\un\big)$&$
\Par\widetilde{\underline{L}}^u=\Par\widetilde{\underline{M}}^u=\Par\widetilde{\un}^u=3$\\
\hline
$\fm\fb_3\big(9;\widetilde{\underline{L}}\big)\simeq\fm\fb_3\big(11;
\widetilde{\underline{M}}\big)\simeq\fm\fb_2\big(11;\widetilde\un\big)$&$
\Par\widetilde{\underline{L}}^u=\Par\widetilde{\underline{M}}^u=\Par\widetilde{\un}^u=5$\\
\hline

$\fk\fle\big(15;\widetilde{\underline{K}}\big)\simeq\widetilde{\fk\fle}\big(15;\widetilde{\underline{L}}\big)
\simeq\fk\fle_3\big(20;\widetilde{\underline{M}}\big)\simeq \fk\fle_2\big(20;\widetilde\un\big)$&$\Par\widetilde{\underline{K}}^u=
\Par\widetilde{\underline{L}}^u=\Par\widetilde{\underline{M}}^u=\Par\widetilde{\un}^u=5$\\
\hline
$\fk\fas^{(1)}\big(7;\widetilde{\underline{L}}\big)
\simeq\fk\fas^{(1)}\big(8;\widetilde{\underline{M}}\big)\simeq
\fk\fas^{(1)}\big(10;\widetilde\un\big)$&$
\Par\widetilde{\underline{L}}^u=\Par\widetilde{\underline{M}}^u=\Par\widetilde{\un}^u=7$\\
\hline
$\widetilde{\fk\fas}^{(1)}\big(7;\widetilde{\underline{K}}\big)$&$\Par\widetilde{\underline{K}}^u=3$\\
\hline
$\fvas^{(1)}(8;\un)$&$\Par\un^u=4$\\
\hline
\end{tabular}
\end{equation*}
\end{Theorem}

\begin{proof}[Proof of Theorems \ref{Bjthm21} and~\ref{Bjthm3}]
For $\fv\fle\big(9; \widetilde{\underline{N}}\big)$, see Section~\ref{Svle7}. For $\fv\fle(5; \underline{N}|4)$ and $\widetilde{\fv\fle}\big(9; \widetilde{\underline{N}}\big)$, see Section~\ref{SwidetildeVLE9}. For $\widetilde{\fv\fle}\big(15; \widetilde{\underline{N}}\big)$, see Section~\ref{SKleGradK}. For $\fk\fle_3\big(20;\widetilde\un\big)$, see Section~\ref{Skle_3}. For $\fm\fb\big(9;\widetilde\un\big)$, see Section~\ref{Smb4/5}. For $\widetilde{\fk\fas}\big(5;\widetilde{\underline{N}}|5\big)$, see Section~\ref{Skas55}.

For $\fv\fle(7;\widetilde{\underline{L}})$, see Section~\ref{Svle43Kee2}. For $\fm\fb_2\big(11;\widetilde\un\big)$, see Section~\ref{Smb11_2}. For $\fk\fle_2\big(20;\widetilde\un\big)$, see Section~\ref{Skle_2}. For $\fk\fas\big(8;\widetilde{\underline{M}}\big)$, see Section~\ref{Skas8}. For $\widetilde{\fk\fas}\big(7;\widetilde{\underline{M}}\big)$, see Section~\ref{Skas7}. For $\fvas (4;\un|4)$ and $\fvas (8;\un)$, see Section~\ref{Svas}.
\end{proof}

\begin{Theorem}[desuperizations of $\fkas$]\label{Bjthm4}
For $p=2$, the analogs of the W-graded simple vectorial Lie
superalgebra $\fkas$ over $\Cee$ are not simple; they contain a~simple ideal of codimension~$1$, the derived algebra.
\end{Theorem}

\begin{proof}
See Section~\ref{ssSimpleIdeal}.
\end{proof}

\textbf{Other results.} We explicitly describe characteristic-2 Shen's version\index{(Shen algebra@Shen algebra}
(see \cite{Shen} and Section~\ref{SShen} here) of both $\fg(2)$ and the
Melikyan algebra\index{(Melikyan algebra@Melikyan algebra} as vectorial Lie algebras.

We single out the divergence-free subalgebra
$\fs\fh_\Pi\big(2n;\widetilde\un\big)$\index{$smathfrak\fh_\Pi\big(2n;\widetilde\un\big)$, divergence-free subalgebra of
 $\fh_\Pi(2n;\uM)=\textbf{F}(\fle(n;\un\vert n))$@$\fs\fh_\Pi\big(2n;\widetilde\un\big)$, divergence-free subalgebra\newline of
 $\fh_\Pi(2n;\uM)=\textbf{F}(\fle(n;\un\vert n))$}
of the Lie algebra of Hamiltonian vector fields $\fh_\Pi(2n;\uM)=\textbf{F}(\fle(n;\un|n))$, discovered in \cite{Leb}, by imposing
constraints on $\uM$.

\subsubsection{Open problems}\label{sssOpPr} \index{(Problem, open@Problem, open}

\begin{enumerate}\itemsep=0pt

\item In this paper, for serial simple vectorial Lie
superalgebras, we considered 2 W-gradings of $\textbf{F}\big(\widetilde\fsb(2^n-1)\big)$ (the standard and of depth 1). We did not consider~32 W-gradings of the remaining simple vectorial Lie superalgebras, see Section~1.16 in~\cite{LS1}.
The fact~\eqref{fact} suggests that we should investigate if these W-gradings yield new simple Lie algebra.

\item  Describe all the $\Zee/2$-grading of the newly found simple Lie
algebras in characteristic~2 (see \cite{BGLLS, BLLS1, Ei, GZ, Leb, SkT1}, and this paper) to obtain new simple Lie superalgebras.
For first results, see~\cite{KrLe}.

\item Are the partial prolongs described in
Section~\ref{ss59} and in equation~\eqref{dimVin=31} isomorphic to known simple Lie algebras?

\item See Problems \ref{OPfb},~\ref{OP}, \ref{OPas}.
\end{enumerate}

\section[The Lie superalgebra $\fvle(4\vert 3)$ over $\Cee$ and its $p>2$ versions]{The Lie superalgebra $\boldsymbol{\fvle(4\vert 3)}$ over $\boldsymbol{\Cee}$ and its $\boldsymbol{p>2}$ versions}\label{Svle43Cee}

\subsection{Recapitulations, see \cite{ShP}}\label{ssRecapShP} In the realization of
$\fle(3)$ by means of generating functions, we identify the space of
$\fle(3)$ with $\Pi(\Cee[\theta, q]/\Cee\cdot 1)$, where before the functor $\Pi$ is applied
$\theta=(\theta_1,\theta_2,\theta_3)$ are odd and $q=(q_1,q_2,q_3)$ are even, see~\eqref{le_n}. In the \textit{standard grading} $\deg_{\rm Lie}$\index{$deg_{\rm Lie}$@$\deg_{\rm Lie}$} of
$\fle(3)$, we assume that $\deg q_i=\deg \theta_i=1$ for all $i$,
and that the grading is given by the formula
\begin{equation*}
\deg_{\rm Lie}(f):=\deg \Le_f=\deg f-2\quad \text{for any monomial
$f\in\Cee[\theta, q]$.}
\end{equation*}
The \textit{nonstandard} grading $\deg_{{\rm Lie};3}$\index{$deg_{{\rm Lie}; 3}$@$\deg_{{\rm Lie}; 3}$} of $\fg=\fle(3;3)$
is determined by the formulas
\begin{gather*}
\deg_3 \theta_i=0\quad \text{and} \quad \deg_3 q_i=1\quad \text{for} \ i=1,2,3,
\\ \deg_{{\rm Lie};3} (f)=\deg_3 f-1\quad \text{for any monomial $f\in\Cee[\theta, q]$}.
\end{gather*}
This grading of $\fg=\fle(3;3)$ is of depth 1, and its homogeneous components are of the form:
\[
\fg_{-1}=\Pi(\Cee[\theta_1,\theta_2,\theta_3]/\Cee\cdot 1),\quad
\fg_k=\Pi(\Cee[\theta_1,\theta_2,\theta_3])\otimes
S^{k+1}(q_1,q_2,q_3) \quad \text{for $k\geq 0$}.
\]
In particular, $\fg_0\simeq\fvect(0|3)$, and $\fg_1$ is an
irreducible $\fg_0$-module with the lowest-weight vector~$q_1^{2}$.
The whole Lie superalgebra $\fle(3;3)$ is the Cartan prolong of its
nonpositive part and the component $\fg_1$ generates the whole
positive part.

To obtain $\fvle(4|3)$, we add the central element $d$ to the 0th
component of $\fle(3;3)$; so that $\ad_d$ is the grading operator on
the Cartan prolong of its nonpositive part; this prolong is strictly
bigger than $\fle(3;3)\ltimes \Cee\cdot d$.

This Cartan prolong is the exceptional simple Lie superalgebra
$\fvle(4|3)$.\index{$vlemathfrak(4\vert 3)$@$\fvle(4\vert 3)$}

Its component $\fvle_1$ is reducible, but indecomposable $\fvle_0$-module such that
\[
\fvle_1/\fle(3;3)_1\simeq (\fle(3;3)_{-1})^*.
\]
The other, not lying in $\fle(3;3)_1$, lowest-weight vector in $\fvle_1$ is the element dual to the highest-weight vector in $\fg_{-1}$, i.e., to $\Pi(\theta_1\theta_2\theta_3)$.

\subsubsection[Introducing indeterminate $y$ (as well as $x_i$ and $\xi_i$)]{Introducing indeterminate $\boldsymbol{y}$ (as well as $\boldsymbol{x_i}$ and $\boldsymbol{\xi_i}$)}

Under the identification
\begin{equation*}
\Pi(\theta_1\theta_2\theta_3)\longmapsto -\del_y, \quad
\Pi(\theta_i)\longmapsto -\del_{x_i},\quad
\Pi\left(\frac{\del(\theta_1\theta_2\theta_3)}{\del
\theta_i}\right)\longmapsto -\del_{\xi_i}
\end{equation*}
each vector field $D\in\fvle(4|3)$ is of the form\index{$D_{f,\;g}$}
\begin{gather}
D_{f,\;g}=\Le_f+yB_f-(-1)^{p(f)}\left(y\Delta(f)+y^2\frac{\del^3f}{\del
\xi_1\del \xi_2\del \xi_3 }\right)\del_y \nonumber\\
\hphantom{D_{f,\;g}=}{} +B_g-(-1)^{p(g)}\left(\Delta(g)+2y\frac{\del^3g}{\del \xi_1\del
\xi_2\del \xi_3}\right)\del_y,\label{D_{f,g}}
\end{gather}
where
$f,g\in\Cee[x,\xi]$; the operators $B_g$ 
\index{$B_g$} and $\Delta$\index{((ZZDelta@$\Delta$} are given by the formulas
\begin{equation}\label{BandDelta}
B_g=\frac{\del^2g}{\del \xi_2\del \xi_3}\frac{\del}{\del
\xi_1}+\frac{\del^2g}{\del \xi_3\del \xi_1}\frac{\del}{\del
\xi_2}+\frac{\del^2g}{\del \xi_1\del \xi_2}\frac{\del}{\del \xi_3},
\quad \Delta = \mathop{\sum}\limits_{1\leq i\leq 3}
\frac{\del^2}{\del x_i\del \xi_i}.
\end{equation}

There are two embeddings of $\fle(3)$ into $\fvle(4|3)$. The
embedding $i_1\colon \fle(3)\tto\fvle(4|3)$ corresponds to the grading
$\fle(3;3)$. Let us reproduce the explicit formulas from~\cite{ShP}:

Let us clarify notation of indeterminates. At the beginning, the indeterminates describing $\fle$ in any grading are denoted by $q$ and $\theta$, while the indeterminates describing $\fvle$ are denoted by~$x$ and~$\xi$. Hence, introduce the passage $i$ from one set of indeterminates to the other, and $\hat f$:
\begin{gather*}
i(q_1,q_2,q_3,\theta_1,\theta_2,\theta_2):=(x_1,x_2,x_3,\xi_1,\xi_2,\xi_3),\\
\hat f(x_1,x_2,x_3,\xi_1,\xi_2,\xi_3):=f(i(q_1,q_2,q_3,\theta_1,\theta_2,\theta_2)).
\end{gather*}
\begin{itemize}\itemsep=0pt
\item[a)] If $f=f(q)$, then
\[
 i_1(\Le_{f})=\Le_{\sum \big(\pderf{\hat f}{x_i}\big)\xi_j\xi_k-y\hat f},
\]
\textit{where $y$ is treated as a~parameter} and $(i,j, k)\in A_3$ (even
permutations of $\{1,2,3\})$.
\item[b)] If $f = \sum f_i(q)\theta_i$, then
\[
 i_1(\Le_f)=\Le_{\hat f} - \varphi(x)\sum \xi_i\partial_{\xi_i}
 +\left(-\varphi(x)y + \Delta (\varphi(x)\xi_1\xi_2\xi_3)\right)\partial_y,
\]
where $\varphi(x)=\Delta\big(\hat f\big)$.
\item[c)] If $f=\psi_1(q)\theta_2\theta_3 + \psi_2(q)\theta_3\theta_1 +
 \psi_3(q)\theta_1\theta_2$, then
\[
 i_1(\Le_f) = -\Delta\big(\hat f\big)\partial_y -\mathop{\sum}\limits_{1\leq i\leq 3}\psi_i(x)\pder{\xi_{i}}.
\]
\item[d)] If $f=\psi(q)\theta_1\theta_2\theta_3$, then
\[
 i_1(\Le_f)=-\psi(x)\partial_y.
\]
\end{itemize}

The embedding $i_2\colon \fle(3)\tto\fvle(4|3)$ corresponds to the
standard grading of $\fle(3)$. In terms of generating functions this
embedding is of the form
\begin{equation}\label{i2f}
 i_2(f(q,\theta))\tto D_{f(i(q,\theta)),\;0}.
\end{equation}

As vector spaces, we have
\begin{equation}\label{i1plusi2}
\fvle(4|3)=i_1(\fle(3;3))+i_2(\fle(3)) \quad \text{while} \quad i_1(\fle(3;3))\cap i_2(\fle(3))\simeq\fsle^{(1)}(3).
\end{equation}

By abuse of notation, denote the operator $\mathop{\sum}\limits_{1\leq i\leq 3}
\frac{\del^2}{\del q_i\del \theta_i}$ acting on the space of functions (divided powers if $p>0$) in $q_i, \theta_i$ also by $\Delta$. In this
notation, we have
\begin{equation}\label{sleaut}
i_1(f(q))=i_2(\Delta(f(q)\theta_1\theta_2\theta_3)), \quad
i_1(\Delta(f(q)\theta_1\theta_2\theta_3))=-i_2(f(q)),
\end{equation}
and
\begin{equation} \label{sleaut2}
i_1(f)=i_2(f) \quad \text{if} \quad f=\mathop{\sum}\limits_{1\leq i\leq 3}
f_i(q)\theta_i \quad \text{and}\quad \Delta f=0.
\end{equation}

The formulas \eqref{i2f}, \eqref{sleaut}, and \eqref{sleaut2} are valid for any $p$, in particular, for $p=2$.

The lowest-weight vectors in $\fvle_1$ are $i_1\big(q_1^{2}\big)$ and $i_2(\theta_1\theta_2\theta_3)$. We have
\begin{equation} \label{sdimVle1}
 \sdim \fvle_1=28|27.
\end{equation}

\begin{Proposition}[passage from $\Cee$ to $\Kee$ for $p>2$]\label{vlePgeq2} The situation described in the previous subsection does not change under passage from the ground field $\Cee$ to any field of characteristic~$0$ and also to $\Kee$ if $\Char\Kee=p>2$ provided the coordinates of the shearing vector~$\uM$ of the algebra of coefficients of the vector fields are such that $\uM_i=\infty$ for
each even indeterminate $x_i$.
\end{Proposition}

In all these cases, the Lie superalgebra $\fg$~-- the Cartan
prolong of the nonpositive part of~$\fvle$~-- remains simple and of
infinite dimension. The component $\fg_1$ also retains its structure
as a~$\fg_0$-module, but it does not generate the whole of $\fg$ (since
the $x_i$ do not generate $\cO(x;\uM)$ if~$\uM\neq\One$,
see~\eqref{un_s}).

\begin{proof} Proof follows directly from the formulas of Section~\ref{ssRecapShP}.\end{proof}

In what follows we will investigate the case of shearing vectors with finite coordinates.

\begin{Theorem}[$\fvle$ exists for any $p>0$]\label{ThVle} A~characteristic-$p$ version $\fvle\big(4;\widetilde\uM|3\big)$ of the Lie superalgebra $\fvle(4|3)$ exists for any $p>0$ and any shearing vector $\widetilde\uM=(\uM,M_y)$ provided $M_y=1$.
\end{Theorem}

For the case $p=2$, see Sections \ref{Svle43Kee2}, \ref{Svle7}, \ref{SwidetildeVLE9}. For the case $p>2$, see the next proposition.

\begin{Proposition}[description of $\fvle\big(4;\widetilde\uM|3\big)$]\label{prVle} If $p>2$ and $\uM_i<\infty$ for $i=1,2,3$, then $\fvle\big(4;\widetilde\uM|3\big)= \Span(D_{f,\;g})$, where $\widetilde\uM=(\uM,M_y)$ and $M_y=1$, and
\[
 f\in \cO(x;\uM| \xi)\oplus \Span \big(x_i^{(s_i)}\big)_{i=1}^3
\quad \text{and}\quad g\in \cO(x;\uM| \xi)\oplus \Kee\cdot
x_1^{(s_1)}x_2^{(s_2-1)}x_3^{(s_3-1)}\xi_1.
\]
\end{Proposition}

\begin{proof} For $p>2$, the component $\fg_1$ also retains its structure as
$\fg_0$-module even if any (or all) of the coordinates of the
shearing vector $\uM=(\uM_1,\uM_2,\uM_3)$ become finite. If
$\uM_i<\infty$ for all~$i$, the Cartan prolong is finite-dimensional. It can be described by means of equation~\eqref{D_{f,g}}, but we
should investigate when $D_{f,\;g}\in\fvect\big(4;\widetilde\uM|3\big)$,
where\index{$M=(\uM, M_y)$, speaking only of $\fvle\big(4; \widetilde\uM\vert 3\big)$@$\widetilde\uM=(\uM, M_y)$, speaking only of $\fvle\big(4; \widetilde\uM\vert 3\big)$}
\begin{gather*}
\widetilde\uM=(\uM, M_y)\quad \text{and} \quad M_y=1.
\end{gather*}
Direct observation gives the answer:
\begin{equation*}
f\in\cO(x;\uM|\xi)\oplus \Span
\big(x_i^{(s_i)}\big)_{i=1}^3,\quad \text{where $s_i=p^{M_i}$;}
\end{equation*}
i.e., we should add ``virtual'' (non-existing for the given $\uM$)
elements $f=x_i^{(s_i)}$ for $i=1,2,3$. Since due to \eqref{sleaut}
\[
D_{x_i^{(s_i)},0}=i_2\big(q_i^{(s_i)}\big)=i_1\big(\Delta\big(q_i^{(s_i)}
\theta_1\theta_2\theta_3\big)\big),
\]
we see that $i_2\big(q_i^{(s_i)}\big)\in \fvle^{(1)}\big(4;\widetilde\uM|3\big)$
though $q_i^{(s_i)}\notin \fle^{(1)}(3;\uM)$.

Now, let us investigate the generating functions $g$. We have\index{$D_{0,-}: g\longmapsto D_{0,\ g}$}
\begin{gather}
\Ker(D_{0,-}\colon g\longmapsto D_{0, g})=\Span\bigg(D_{0,g}=0\,|\,\text{$g=g(x)$ or $g=\mathop{\sum}\limits_i
g_i(x)\xi_i$ with $\Delta g=0$}\bigg).\!\!\!\label{ker}
\end{gather}
If $g=\mathop{\sum}\limits_i g_i(x)\xi_i$, but $\Delta g=h(x)\ne 0$,
then $D_{0,g}=h(x)\del_y$ depends on $h$ only. It is clear that
any function $h\in \cO(x;\uM)$ can be expressed as $h=\Delta g$ for
some $g\in \cO(x;\uM|\xi)$, except for
\[
h_{s}=x_1^{(s_1-1)}x_2^{(s_2-1)}x_3^{(s_3-1)},\quad \text{where $s=(s_1,
s_2, s_3)$.}
\]
To obtain $D={h_{s}}\del_y$, we should add to the space of
generating functions any of the ``virtual" functions
\[g_{s,i}=
\xi_i\del_j\del_k\big(x_1^{(s_1)}x_2^{(s_2)}x_3^{(s_3)}\big)
\quad \text{for $i=1, 2, 3$, and $j,k\ne i$, $j\ne k$.}
\]
Modulo the kernel \eqref{ker} of the map $D_{0,-}$ only one ``extra''
generator suffices; for definiteness, we select
$g_s:=x_1^{(s_1)}x_2^{(s_2-1)}x_3^{(s_3-1)}\xi_1$. Formula~\eqref{D_{f,g}} shows that $D_{0,g_s}$ lies in the homogeneous
component of degree
\[
s_1+(s_2-1)+(s_3-1)+1-3\equiv -4 \mod p.
\]
Since
the Lie superalgebra $\fvle(4;\widetilde\uM|3)$ has a~grading
operator, it follows that
\[
D_{0,g_s}\in\fvle^{(1)}\big(4;\widetilde\uM|3\big)
\]
 for $p>2$. Moreover, as was shown
in \cite{ShP} for $p=0$ (but the formulas remain true for any
$p>0$),
\[D_{0, g_s}=i_1(-h_s\xi_1\xi_2\xi_3), \quad \text{and hence $
D_{0,  g_s}\in i_1\big(\fle^{(1)}(3;\uM)\big)$ for $p>2$}.\tag*{\qed}
\]\renewcommand{\qed}{}
\end{proof}

If $p>2$, the fact \eqref{i1plusi2} does not hold. We have
\[
D_{0,\;h_s\xi_1\xi_2\xi_3}\in
\fvle\big(4;\widetilde\uM|3\big),\quad \text{but~~}D_{0,\;h_s\xi_1\xi_2\xi_3}=
i_1(f)-i_2(f),
\] where
\[
f=q_1^{(s_1)}q_2^{(s_2-1)}q_3^{(s_3-1)}\theta_1+
q_1^{(s_1-1)}q_2^{(s_2)}q_3^{(s_3-1)}\theta_2+
q_1^{(s_1-1)}q_2^{(s_2-1)}q_3^{(s_3)}\theta_3\notin
\cO(q;\uM|\theta).
\]

If $\uM_i>1$ for $i=1, 2,3$, then $\sdim\fg_1=28|27$ remains the same as over $\Cee$ and any other $p>2$.

\section[A~description of $\fvle\big(7;\widetilde{\uM}\big):=\textbf{F}(\fvle(4;\uM\vert 3))$ for $p=2$]{A~description of $\boldsymbol{\fvle\big(7;\widetilde{\uM}\big):=\textbf{F}(\fvle(4;\uM\vert 3))}$ for $\boldsymbol{p=2}$}\label{Svle43Kee2}

In this and two next sections we prove Theorem~\ref{ThVle} for the case $p=2$.

There are three W-gradings of $\fvle$ with unconstrained shearing vector. In this and two next sections, we consider each of these gradings for $p=2$, describe the corresponding $0$th and $1$st components of the regraded Lie superalgebra $\fvle$, and calculate partial prolongs. Next, we consider desuperizations of each of them. As a~result, we get a~new simple Lie superalgebra and a~new simple Lie algebra in $p=2$. Unfortunately (we'd like to get new simple algebras!), there are no partial prolongs.

The Lie superalgebra $\fvle\big(4; \widetilde\uM|3\big)$\index{$vlemathfrak\big(7;\widetilde{\uM}\big):=
\textbf{F}(\fvle(4;\uM\vert 3))$@$\fvle\big(7;\widetilde{\uM}\big):=
\textbf{F}(\fvle(4;\uM\vert 3))$}
for $p=2$ is a~direct reduction modulo~2 of the integer
form, with divided powers as coefficients, of the complex vectorial
Lie superalgebra~$\fvle(4|3)$.

First of all, let us define squares of odd elements for the Lie
superalgebra $\fle(n; \uM)$, cf.\ equation~\eqref{K_fbr}:
\begin{equation*}
 \label{(sq)}
 f^{2}:=
\mathop{\sum}\limits_{1\leq i\leq n} \frac{\del f}{\del
q_i}\frac{\del f}{\del \theta_i}
\end{equation*}
and have in mind that if $p=2$ and $\uM_i<\infty$ for all $i$, the
Lie superalgebra $\fg=\fle(n;\uM)$ is not simple: the generating
function of the maximal degree $q_1^{(s_1-1)}q_2^{(s_2-1)}\cdots
q_n^{(s_n-1)}\theta_1\theta_2\cdots\theta_n$
does not belong to $\fg^{(1)}$, the latter being simple.

\underline{For $p>2$}, we just reduce the expression~\eqref{D_{f,g}} modulo~$p$.

\underline{For $p=2$}, we cannot just reduce the expression \eqref{D_{f,g}} modulo
2; we should modify it. Indeed, the system of equations on the
coefficients of the field $D\in\fvle(4|3)$ whose solution is given
by the formula \eqref{D_{f,g}} contain coefficients $\frac12$, see
\cite{ShP}. The vector field $D\in\fvle\big(4; \widetilde\uM|3\big)$ is of
the form (recall formula \eqref{BandDelta} for $B_g$):\index{$D_{f,\;g}$}
\begin{equation}\label{D_{f,g}p=2}
D_{f,\;g}=\Le_f+yB_f+y\Delta(f)\del_y + B_g+\Delta(g)\del_y,
\quad \text{where $\widetilde\uM=(\uM,M_y)$ and $M_y=1$}.
\end{equation}

Let us explain how we got this formula: we just rewrote equations from~\cite{ShP} without $\frac12$ (multiplied the equations by 2) and solved them in the same way as it was done in~\cite{ShP}.

\underline{For $p=2$}, unlike the case $p\neq 2$, this Lie superalgebra $\fg=\fvle\big(4;\widetilde\uM|3\big)$ is not simple,
but $\fg^{(1)}$ is simple, its codimension in $\fg$ is equal to~2:
for
$f=q_1^{(s_1-1)}q_2^{(s_2-1)}q_3^{(s_3-1)}\theta_1\theta_2\theta_3$,
we have
\[
D_{f,0}=i_2(f)\notin\fg^{(1)}, \quad
D_{0,x_1^{(s_1)}x_2^{(s_2-1)}x_3^{(s_3-1)}\xi_1}=i_1(f)\notin\fg^{(1)}.
\]

\underline{For $p=2$}, the structure of the $\fg_0$-module $\fg_1$ differs
drastically from that for $p\neq 2$. Instead of two lowest-weight
vectors, we have three of them. Besides, these three lowest-weight
vectors do not describe the whole complexity of the module.

The submodule $i_1(\fle(3;3)_1)$ has a~complicated structure. To
describe it, observe that for any vectorial Lie (super)algebra
expressed in terms of generating functions, the shearing vector can
be considered either
\begin{enumerate}\itemsep=0pt
\item[(1)] on the level of generating functions (let us denote it $\uM$ in this case) or
\item[(2)] on the level of coefficients of vector fields they generate (let us denote it $\widetilde\uM$ in this case).
\end{enumerate}

In case(1), we obtain the ``underdeveloped'' Lie superalgebra
\begin{equation*}\label{leSupM}
\fle(n; \uM):=\Span(\Le_f\,|\,f\in\cO(n;\uM|n));
\end{equation*}
in case (2), we get the correct
$\fle\big(n;\widetilde\uM\big)=\Span\big(\Le_f\in \fvect\big(n;\widetilde\uM|n\big)\big)$.
We have
\begin{equation*}\label{leBezM}
\fle\big(n;\widetilde\uM\big)=\fle(n;\uM)\ltimes\Span\big(q_i^{(s_i)}\,|\,
i=1,\dots,n\big).
\end{equation*}

Accordingly, for $\widetilde\uM$ unconstrained, the component
$\fle\big(3;\widetilde\uM;3\big)_1$ is of the form:
\[
\fle\big(3;\widetilde\uM;3\big)_1=\Pi\big(\Kee[\theta_1\theta_2\theta_3]\big)\otimes S^2(q_1,q_2,q_3).
\]
This component contains submodules corresponding to the minimal
values $\uM_i=1$ for some $i$. To describe $\fg_1$ as $\fg_0$-module, consider the submodule
\[
W_0:=\fle^{(1)}(3;3)_1=\Span(q_iq_j\varphi(\theta)\text{~for any
$i$, $j$ and any function $\varphi$}).
\]
It is irreducible. It is glued to the submodules
\[
W_i:=W_0\ltimes
\Kee\cdot q_i^{(2)}
\]
in each of which $W_0$ is a~submodule, but not a~direct summand. Each $W_i$ can be further enlarged to the module
\[
W_{i,\theta}:=W_i\ltimes\Span\big(q_i^{(2)}\vf(\theta)\big)\quad \text{corresponding to $\uM_i>1$, where $\uM_j=1$ for $j\ne i$}
\]
with shearing performed on the level of generating functions.

Let us describe the subalgebras of $\fvle\big(4;\widetilde\One|3\big)$, the partial prolongs. In what follows we will often use the following

\subsection{Notational convention: on partial prolongs}\index{(Convention@Convention}
\begin{gather}
\text{Let $v_i$ be a~lowest-weight vector of the $\fg_0$-module $\fg_1$}\nonumber\\
\text{and $V_i$ the submodule generated by $v_i$.}\nonumber\\
\text{Let $\fg_{k,{(i)}}$ be the
$k$th prolong ``in the direction of $V_i\subset\fg_1$",\index{$gmathfrak_{k,{(i)}}$, the
$k$th prolong ``in the direction of $V_i\subset\fg_1$''@$\fg_{k,{(i)}}$, the
$k$th prolong ``in the direction\newline of $V_i\subset\fg_1$''} i.e.,}\nonumber\\
\text{$k$th prolong of $(\fg_-, \fg_0, V_i)$.}\label{convent}
\end{gather}

Consider the $\fg_0$-submodules $W\subset\fg_1$ not contained in
$i_1(\fle(3;3))$. There is only one such submodule $V_3\subset W$
generated by $v_3$, see \eqref{lwv1}. The $\fg_0$-module $\fg_{1,
\uM}$ has the following three lowest-weight vectors expressed in the
form $D_{f,\;g}$, and also as $i_1(-)$ or $i_2(-)$:
\begin{equation}\label{lwv1}
\footnotesize
\renewcommand{\arraystretch}{1.4}
\begin{tabular}{|l|l|l|l|} \hline
$v_1$ & $ x_1 y\del_{\xi_2} + x_1 \xi_3\del_{x_1}+ x_1
\xi_1\del_{x_3} +x_2 y\del_{\xi_1}$&
$i_1(q_1q_2)$ & $D_{x_1\xi_1\xi_3+x_2\xi_2\xi_3,\;0}$\\
 & $+x_2\xi_3\del_{x_2}
+x_2 \xi_2\del_{x_3} +\xi_3 \xi_2\del_{\xi_2} +
\xi_3 \xi_1\del_{\xi_1} $ & &\\
\hline $v_2$ & $ x_1 y\del_{\xi_1} +x_1 \xi_3\del_{x_2} +x_1
\xi_2\del_{x_3} +\xi_3 \xi_2\del_{\xi_1}$ &$i_1\big(q_1^{(2)}\big)$
 & $D_{x_1\xi_2\xi_3,\;0} $\\
 \hline
 $v_3$ & $ y \xi_3\del_{\xi_3} +y \xi_2\del_{\xi_2}
 +y \xi_1\del_{\xi_1} +\xi_3 \xi_2\del_{x_1}$
& $ i_2(\theta_1\theta_2\theta_3)$ & $D_{\xi_1\xi_2\xi_3,\;0} $\\
 &${}+\xi_3 \xi_1\del_{x_2} +\xi_2\xi_1\del_{x_3}$ & &\\
\hline
\end{tabular}
\end{equation}

By increasing the value of some of the coordinates $\uM_i$ we
enlarge $\fg_{1,(3)}=\fvle\big(4;\widetilde\One|3\big)_1$. As
$\fg_0$-module, $\fg_{1, (3)}$ is of the following form:
\[
W_0\subset (W_1+W_2+W_3)\subset \fg_{1, (3)},
\]
where $\fg_{1, (3)}/(W_1+W_2+W_3)\simeq (\fg_{-1})^*$, and
$(W_1+W_2+W_3)/W_0$ is the trivial $0|3$-dimensional module, and
where $\sdim W_0=12|12$.

The superdimensions of the positive components of $\fvle(4;\One|3)$ (and of its derived algebra in parentheses) are given in the following table:
\begin{equation*}\label{sdims1}
\footnotesize
\renewcommand{\arraystretch}{1.4}
\begin{tabular}{|l|c|c|c|c|}
\hline &$\fg_1$ &
$\fg_2$& $\fg_3$&$\fg_4$ \\
\hline $\sdim$ &$16|18$ & $10(9) |12$& $4|3$& $1(0)|0$\\ \hline
\end{tabular}
\end{equation*}

\subsection[Partial prolongs as subalgebras of $\fvle\big(4;\widetilde\uM|3\big)$]{Partial prolongs as subalgebras of $\boldsymbol{\fvle\big(4;\widetilde\uM|3\big)}$}

We have $[\fg_{-1},\fg_{1, (i)}]\simeq \fvect (0|3)$ for $i=1, 2$.
(Actually, $\fg_{1, (2)}=\fg_{1, (1)}\ltimes \Kee\cdot i_1\big(q_1^{(2)}\big)$.)

\subsubsection{Convention: on partial prolongs ``of no interest''}\label{noprol}\index{(Convention@Convention}

 \textit{In what follows, we do not investigate partial prolongs with
$[\fg_{-1},\fg_{1, (i)}]$, see \eqref{convent}, smaller than~$\fg_0$
if the $[\fg_{-1},\fg_{1, (1)}]$-module $\fg_{-1}$ is not
irreducible}: no such prolong can be a~simple Lie (super)algebra
with the given nonpositive part.

\subsection{Desuperization} \label{vle7} We have $\fg_{0}\simeq
\fc(\fvect(3;\One))$, and $\fg_{-1}\simeq \cF/\Kee$.

For $\un$ unconstrained, we have $\dim\fg_1=55$. (Compare with \eqref{sdimVle1} and \eqref{sdims12}: for $\underline{N}=\One$, the dimension drops.)

\textit{Critical coordinates of the unconstrained shearing vector}:
$\widetilde M_1$, $\widetilde M_2$, $\widetilde M_3$, $\widetilde M_7$.

The dimensions of the positive components of $\fvle(7;\widetilde
\One)$ and its simple derived algebra (in parentheses) are given in the
following table; so $\dim \fvle^{(1)}\big(7;\widetilde \One\big)=94$:
\begin{equation}\label{sdims12}
\footnotesize
\renewcommand{\arraystretch}{1.4}
\begin{tabular}{|l|c|c|c|c|} \hline
&$\fg_1$ & $\fg_2\text{ (or }\fg_1^{\bcdot 2})$&
$\fg_3=\fg_1^{\bcdot 3}$&$\fg_4\text{ or }-$\\
\hline $\dim$ &$34$ & $22\;(21)$& $7$& $1\;(-)$\\ \hline
\end{tabular}
\end{equation}

\subsection[Partial prolongs as subalgebras of $\fvle\big(7;\widetilde \uM\big)$]{Partial prolongs as subalgebras of $\boldsymbol{\fvle\big(7;\widetilde \uM\big)}$}

\begin{enumerate}\itemsep=0pt
\item[(i)] We have $\dim(\fg_1')=34$. (If $\un=\One$, then
$\dim(\fg_{1, (1)})=\dim\fg_1$; otherwise,
$\fg_{1, (1)}\subsetneq\fg_1$.) The partial Cartan prolong
\begin{equation*}\label{vle'}
\fvle'\big(7;\widetilde \uM\big):=\big(\fg_{-},\fg_0,\fg_{1, (1)}\big)_{*;\widetilde
\uM}
\end{equation*}
is such that $[\fg_{-1},\fg_1]\simeq \fcvect(3;\One)$; this prolong is
$\fvle\big(7;\widetilde \One\big)$.

\item[(ii)] The partial Cartan prolong $(\fg_{-},\fg_0,\fg_{1, (i)})_{*;\uM}$ is
such that $[\fg_{-1},\fg_{1, (i)}]\simeq \fvect(3;\One)$ for $i=2,
3$.
{\it By Convention $\ref{noprol}$, we do not investigate this partial
prolong.}
\end{enumerate}

\textit{Conclusion.} We have found a~new simple Lie superalgebra $\fvle^{(1)}\big(4;\widetilde\uM|3\big)$ and a~new simple Lie algebra $\fvle^{(1)}\big(7;\widetilde \uM\big)$. Partial prolongs do not yield new simple Lie (super)algebras.

\section[$\fvle\big(9;\widetilde{\un}\big):= \textbf{F}(\fvle(3;\un\vert 6))$, where $\fvle(3;\un\vert 6):=\fvle(4;\uM\vert 3;K)$ for $p=2$]{$\boldsymbol{\fvle\big(9;\widetilde{\un}\big):= \textbf{F}(\fvle(3;\un\vert 6))}$,\\ where $\boldsymbol{\fvle(3;\un\vert 6):=\fvle(4;\uM\vert 3;K)}$ for $\boldsymbol{p=2}$}\label{Svle7}

The Lie superalgebra $\fvle(3;\un\vert 6):=\fvle(4;\uM\vert 3;K)$\index{$vlemathfrak(9;\widetilde{\un})$@$\fvle(9;\widetilde{\un})$} is the complete prolong of its negative part,
see Section~\ref{ssExDef}. A~realization of the weight basis of the
nonpositive components by vector fields is as follows, where the
$w_i$ is a~shorthand notation for convenience:
\begin{equation*}
\footnotesize
\renewcommand{\arraystretch}{1.4}
\begin{tabular}{|c|l|}
\hline
$\fg_{i}$&the basis elements \\ \hline
&\\[-1.5em]
\hline
$\fg_{-2}$& $\del_1, \ \del_2, \ \del_3$\\
\hline $\fg_{-1}$& $\del_4, \ \del_5,\ \del_6,\ w_7=x_5\del_3+
 x_6\del_2+\del_7,$\\
&$w_8=x_4\del_3+x_6\del_1+\del_8, \ w_9=x_4\del_2+x_5\del_1+\del_9$\\
\hline

$\fg_{0}\cong$& $X_1^+=x_2\del_1+x_4\del_5+
 x_7\del_8,\ X_1^-=x_1\del_2+x_5\del_4+
 x_8\del_7, \ X_3^\pm=\big[X_1^\pm, X_2^\pm\big],$\\
$\fsl(3)\oplus \fgl(2) $&$ X_2^+= x_3\del_2+
 x_5\del_6+x_8\del_9, \ X_2^-=x_2\del_3+
 x_6\del_5+x_9\del_8, \ H_i=\big[X_i^+,X_i^-\big]$ \text{~for $i=1,2$};\\
&$\widetilde X_1^+=x_7\ x_8\del_3+ x_7\
 x_9\del_2+ x_8\
 x_9\del_1+x_7\del_4+x_8\del_5+
x_9\del_6,$\\
&$\widetilde X_1^-=x_4x_5\del_3 + x_4
 x_6\del_2+ x_5
 x_6\del_1+x_4\del_7+x_5\del_8+x_6\del_9,$\\
&$d=x_1\del_1+x_2\del_2+x_3\del_3+ x_4\del_4+
x_5\del_5+x_6\del_6,\ \widetilde H_1=\big[\widetilde X_1^+,
\widetilde X_1^-\big]$\\
\hline
\end{tabular}
\end{equation*}
The $\fg_0$-module $\fg_1$ has the following lowest- weight vectors:
\begin{gather*}
v_1=x_1 x_4\del_3+ x_1 x_6\del_1+x_2 x_5\del_3+ x_2 x_6\del_2+
 x_1\del_8+x_2\del_7+ x_4 x_6\del_4+ x_4 x_9\del_7\\
\hphantom{v_1=}{} +x_5 x_6\del_5+ x_5 x_9\del_8+ x_6 x_7\del_7+ x_6 x_8\del_8,\\
v_2= x_1 x_5\del_3+ x_1 x_6\del_2+x_1\del_7+ x_5 x_6\del_4+ x_5
x_9\del_7+ x_6 x_8\del_7.
\end{gather*}

\subsection{No simple partial prolongs} For $\un$ unconstrained,
$\dim\fg_1=18$. The module $V_1$ generated by $v_1$ is
6-dimensional, and the module $V_2$ generated by $v_2$ is
8-dimensional; $V_1\subset V_2$.

\textit{Critical coordinates of the unconstrained shearing vector}: $\widetilde\un_4,\dots, \widetilde\un_9$.

\section[A~description of $\widetilde{\fv\fle}\big(9;\widetilde{\un}\big):=
\textbf{F}(\fv\fle(5;\un\vert 4))$ for $p=2$]{A~description of $\boldsymbol{\widetilde{\fv\fle}\big(9;\widetilde{\un}\big):=
\textbf{F}(\fv\fle(5;\un\vert 4))}$ for $\boldsymbol{p=2}$}\label{SwidetildeVLE9}

Whenever possible in this section, we do not indicate the shearing vectors. This Lie superalgebra is the complete prolong of its
negative part, see Section~\ref{ssExDef}.

\underline{For $p=0$}, the $\fg_0$-action on $\fg_{-1}$ is that on the tensor
product of a~2-dimensional space on the space of semi--densities in 2
odd indeterminates. So it is not possible to just reduce modulo $2$ the formulas derived
for the characteristic~$0$. We have to understand how $\fg_0$ acts on
$\fg_{-1}$ when $p=2$. For this, we use the explicit form of elements
of $\fvle(4|3)$ for $p=2$, see equation~\eqref{D_{f,g}p=2}.

Note that the mapping $\Le_f\longmapsto D_{(f,0)}$ determines a
Lie superalgebra isomorphism between~$\fle(3)$ and its image in
$\fvle$. However, first, the mapping $D_{0,-}: g\longmapsto D_{(0, \ g)}$ has
the kernel:
\begin{equation}\label{kerg}
\Ker(D_{0,-})=\{g\in\cO(x;\underline M|\xi)\,|\,\deg_{\xi}g<2
\text{~and~} \Delta g=0\},
\end{equation}
and, second, certain coincidences $D_{(f,0)}=D_{(0, g)}$ might
occur. Formula \eqref{D_{f,g}p=2} makes it clear that such a~coincidence
takes place if and only if (recall formula \eqref{BandDelta} for $B_g$)
\begin{equation*}
\label{feqg1} B_f=0 ,\quad \Delta f=\Delta g=0,\quad \text{and}\quad
\Le_f=B_g.
\end{equation*}

Taking equation~\eqref{BandDelta} into account, these conditions are
equivalent to following conditions:
\begin{equation}
\label{feqg2} f=f(x), \quad g=\mathop{\sum}\limits_{(i,j,k)\in A_3}
\pderf{f}{x_i}\xi_j\xi_k.
\end{equation}

The grading of the Lie superalgebra we are interested in is induced
by the following grading of the space of generating functions:
\begin{equation}\label{gr} \deg \xi_1=0, \quad \deg x_1=2, \quad \deg \xi_2=\deg
\xi_3=\deg x_2=\deg x_3=1, \quad \deg y=0.
\end{equation}

Clearly, $\deg D_{(f, g)}=\deg f-2=\deg g-2$. Therefore (here we
introduce the 9 indeterminates $z_1$, $z_2$, $z_3$, $z_8$, $z_9$ (even) and
$z_4$, $z_5$, $z_6$, $z_7$ (odd) of the ambient Lie superalgebra
$\fvect(5;\un|4)$ containing our $\fg$) and $\del_i:=\del_{z_i}$
\begin{equation*}\label{indetsForVle9}
\footnotesize
\renewcommand{\arraystretch}{1.4}
\begin{tabular}{|l|l|l|}
\hline $\fg_{i}$&the basis elements&in terms of $\fvect(5;\un|4)$\\
\hline
&&\\[-1.5em]
\hline
$\fg_{-2}$&$D_{(f, 0)}, \text{~where~} f=\xi_1$&$\del_1$\\
\hline $\fg_{-1}$&$D_{(f, 0)}, \text{~where~} f=x_i, \
\xi_1\xi_i\,|\,\xi_i,\ \xi_1x_i \text{~for~}
i=2,3$&$x_i\lra\del_{2+i},\ \xi_1\xi_i\lra\del_{4+i}+z_{2+i}\del_1$ \\
&&$\xi_i\lra\del_{i},\ \xi_1x_i\lra \del_{6+i}+z_{i}\del_1$\\
\hline
\end{tabular}
\end{equation*}
because, for the nonzero vector fields of the form $D_{(0, g)}$
lying in $\fg_-$, we have, thanks to conditions~\eqref{feqg2}, the
following identifications:
\[
D_{(0, \xi_1\xi_2)}=D_{(x_3, 0)}, \quad D_{(0,\xi_1\xi_3)}=D_{(x_2, 0)}.
\]

Because the tautological representation of $\fsl(2)$ is isomorphic to its dual, we identify
\[
\fg_{-1}\simeq W=V\otimes \Lambda(\xi,\eta),\quad \text{where} \quad V=\Span(v_1,v_2)
\]
using the rules listed in Table \eqref{TabRules}. The table also
contains the explicit form of the vector fields $D_{(f,0)}\in\fg_{-1}$ needed to calculate the action of the fields of the
form $D_{(0,g)}\in \fg_0$ on $\fg_{-1}$ (the action of the fields
of the form $D_{(f,0)}\in \fg_0$ can be computed in terms of generating functions and the bracket in $\fle$).
\begin{equation}\label{TabRules}
\footnotesize
\renewcommand{\arraystretch}{1.4}
\begin{tabular}{|c|c|c||c|c|c|}
\hline
$f$&$D_{(f, 0)}$& the image in $W$ & $f$&$D_{(f,0)}$&
the image in $W$ \\
\hline
&&&&&\\[-1.5em]
\hline
$\xi_2$ & $\del_{x_2}$ & $v_1$& $\xi_1\xi_2$ & $\xi_2\del_{x_1}+
\xi_1\del_{x_2}+y\del_{\xi_3}$ & $ v_1\otimes \xi$\\
$\xi_3$ & $\del_{x_3}$ & $v_2$ & $\xi_1\xi_3$ & $\xi_3\del_{x_1}+
\xi_1\del_{x_3}+y\del_{\xi_2}$ & $ v_2\otimes \xi$\\
\hline
$x_2$ & $\del_{\xi_2}$ & $v_2\otimes \eta $& $\xi_1 x_2$ &
$\xi_1\del_{\xi_2}+x_2\del_{x_1}$ & $v_2\otimes \xi\eta $\\
$x_3$ & $\del_{\xi_3}$ & $ v_1\otimes \eta$& $\xi_1 x_3$ &
$\xi_1\del_{\xi_3}+x_3\del_{x_1}$ & $ v_1\otimes \xi\eta$\\
\hline
\end{tabular}
\end{equation}

\begin{Proposition}[the component $\fg_0$ of $\fg=\fv\fle(5;\un\vert 4):=\fvle(4;\un\vert 3;1)$ and its action on $\fg_{-1}$]\label{vle_0}
The component $\fg_0$ of $\fg=\fv\fle(5;\un\vert 4):=\fvle(4;\un\vert 3;1)$ consists of vector fields $D_{(f, \ g)}$,
where $\deg f=\deg g=2$ in the grading~\eqref{gr}. Table~\eqref{TabSummary} shows the correspondence between pair of generating functions $(f, g)$ and operators in $\End(W)$
\begin{equation}\label{TabSummary}
\tiny
\renewcommand{\arraystretch}{1.4}
\begin{tabular}{|c|c||c|c|}
\hline
$(f,g)$ & its image in $\End(W)$ & $(f,g)$ & its image in $\End(W)$\\
\hline
&&&\\[-1.5em]
\hline
$f=\mathop{\sum}\limits_{i,j=2,3}a_{ij}x_i\xi_j$, & $\begin{pmatrix}a_{22}&a_{23}\\a_{32}&a_{33}\end{pmatrix}\otimes \mathbbmss{1}$ &
$f=\xi_1\mathop{\sum}\limits_{i,j=2,3}a_{ij}x_i\xi_j$, & $\begin{pmatrix}a_{22}&a_{23}\\a_{32}&a_{33}\end{pmatrix}\otimes \xi$ \\
$\Delta(f)=0$,\ \ $g=0$ &$a_{22}+a_{33}=0$ & $\Delta(f)=0$,\ \ $g=0$ &$a_{22}+a_{33}=0$ \\
\hline
$\big(x_2^{(2)}, 0\big)$&$\begin{pmatrix}0&1\\0&0\end{pmatrix}\otimes \eta$&
$\big(\xi_1x_2^{(2)}, 0\big)$&$\begin{pmatrix}0&1\\
0&0\end{pmatrix}\otimes \xi\eta$\\
$\big(x_3^{(2)}, 0\big)$&$\begin{pmatrix}0&0\\1&0\end{pmatrix}\otimes \eta$&
$\big(\xi_1x_3^{(2)}, 0\big)$&$\begin{pmatrix}0&0\\
1&0\end{pmatrix}\otimes \xi\eta$\\
$(x_2x_3, 0)$&$\begin{pmatrix}1&0\\
0&1\end{pmatrix}\otimes \eta$&
$(\xi_1x_2x_3, 0)$&$\begin{pmatrix}1&0\\
0&1\end{pmatrix}\otimes \xi\eta$\\
\hline
$(x_1, 0)$ & $\del_\xi$ & $(\xi_1x_1, 0)$ & $\xi\del_\xi$\\
$(\xi_2\xi_3, 0)$ & $\del_\eta$ & $(\xi_1\xi_2\xi_3, 0)$ & $\xi\del_\eta$ \\
\hline
$(x_2\xi_2, 0)$ & $\begin{pmatrix}1&0\\0&0\end{pmatrix}\otimes \mathbbmss{1} +\eta\del_\eta$&
$(\xi_1x_2\xi_2, 0)$ & $\begin{pmatrix}1&0\\
0&0\end{pmatrix}\otimes \xi +\xi\eta\del_\eta$\\
\hline
$(0, \xi_1x_1)$ & $\eta\del_\xi $& $(0, \xi_1\xi_2\xi_3)$ & $\xi\del_\xi+\eta\del_\eta$\\
\hline

$(0, \xi_1\xi_2x_2)$ & $\begin{pmatrix}1&0\\
0&0\end{pmatrix}\otimes \eta +\xi\eta\del_\xi$ &&\\
\hline
\end{tabular}
\end{equation}
Thus, $\fg_0\simeq\fd(\fsl(2)\otimes \Lambda(2) \ltimes
\fvect(0|2))$,\index{$dmathfrak(\fg)$@$\fd(\fg)$} where the operator of outer derivation added to the
ideal $\fsl(2)\otimes \Lambda(2) \ltimes \fvect(0|2)$ is $\cD=D\otimes \One$, where $D=\left(\begin{smallmatrix}
1 & 0\\
0 & 0\\
\end{smallmatrix}\right)$, and the subalgebra $\fvect(0|2)$ acts as $X\mapsto
\One\otimes X + D\otimes \Div(X)$ for any $X\in \fvect(0|2)$. Hence, $\fg_{-1}\simeq\Vol(0|2)\oplus \Lambda(2)$, as $\fvect(0|2)$-module.
\end{Proposition}

\begin{proof} Let us begin with fields of the form $D_{(f,0)}$.

As we have already noted above, $[D_{(f_1, 0)}, D_{(f_2, 0)}]=D_{(\{f_1, f_2\},0)}$, and hence the action of such fields
can be described in terms of the generating functions and the
Buttin bracket.

If $f=\mathop{\sum}\limits_{i,j=2,3}a_{ij}x_i\xi_j$ and
$\Delta(f)=0$, then $f$ acts on $\fg_{-1}$ as
$\mathop{\sum}\limits_{i,j=2,3}a_{ij}(\xi_j\del_{\xi_i}+x_i\del_{x_j})$
which, thanks to our identification, corresponds to the action of
the operator
$\left(\begin{smallmatrix}a_{22}&a_{23}\\a_{32}&a_{33}\end{smallmatrix}\right)\otimes
\mathbbmss{1}$; i.e., the elements of this form span the subspace
$\fsl(V)\otimes \mathbbmss{1} \in \End(W)$.

Analogously, the functions of the form
$f=\xi_1\mathop{\sum}\limits_{i,j=2,3}a_{ij}x_i\xi_j$ such that
$\Delta(f)=0$ act on $\fg_{-1}$ as
$\xi_1\mathop{\sum}\limits_{i,j=2,3}a_{ij}(\xi_j\del_{\xi_i}+x_i\del_{x_j})$
which corresponds to the action of the operator
$\left(\begin{smallmatrix}a_{22}&a_{23}\\a_{32}&a_{33}\end{smallmatrix}\right)\otimes
\xi$; i.e., the elements of this form span the subspace
$\fsl(V)\otimes\xi\in \End(W)$.

The correspondence between the other operators $D_{(f,g)}$ and operators in $\End(W)$ is not so evident. To find it, in tables \eqref{t1}--\eqref{t3} and \eqref{tt1}--\eqref{tt2} we indicate nonzero actions of an operator $D_{(f,g)}$ in the first row of these tables and the same actions in the basis of the space $W$ in the second row. This allows us to reestablish the operator in $\End(W)$.

For example, Table \eqref{t1} describes the action of the generating function $f=x_2^{(2)}$. In terms of generating functions, $f$ acts as $\Le_f$, i.e., as the vector field (operator) $X=x_2\partial_{\xi_2}$. Looking at Table~\eqref{TabRules} we deduce that this $X$ acts as non-zero only on $\xi_2$ and $\xi_1\xi_2$. Explicitly, $X(\xi_2)=x_2$, and $X(\xi_1\xi_2)=\xi_1x_2$. This is precisely what is written in the first row of Table~\eqref{t1}.

In the 2nd row there is the same action in terms of the basis of $W$. Indeed, looking again at Table~\eqref{TabRules} we see that in
$W$ the incarnation of $\xi_2$ is denoted by $v_1$, whereas $x_2$ by $v_2\otimes \eta$; i.e., the equality $X(\xi_2)=x_2$ in the new notation takes the form
$X(v_1)= v_2\otimes \eta$:
\be\label{t1}
\begin{tabular}{|l|l|}
\hline
$\xi_2\longmapsto x_2$&$\xi_1\xi_2\longmapsto \xi_1x_2$\\
\hline

$v_1\longmapsto v_2\otimes \eta$&$v_1\otimes\xi\longmapsto
v_2\otimes \xi\eta$ \\ \hline \end{tabular}
\ee
Thus, $f=x_2^{(2)}$ acts as $x_2\del_{\xi_2}$, i.e., as $\left(\begin{smallmatrix}0&1\\
0&0\end{smallmatrix}\right)\otimes \eta\in\End(W)$.

Similarly, $f=x_3^{(2)}$ acts as $\left(\begin{smallmatrix}0&0\\ 1&0\end{smallmatrix}\right)\otimes \eta\in\End
W$, whereas $f=x_2x_3$ acts as the operator
$\left(\begin{smallmatrix}1&0\\ 0&1\end{smallmatrix}\right)\otimes \eta\in\End(W)$.

Thus, the functions of the form $f=f(x_2,x_3)$ span $\fsl(V)\otimes\eta\in \End(W)$.

Analogously, the functions of the form $f=\xi_1f(x_2,x_3)$ span
$\fsl(V)\otimes \xi\eta\in \End(W)$.

Clearly, $f=x_1$ acts on $\fg_{-1}$ as $\del_{\xi_1}$ which
corresponds to the operator $\del_\xi\in \End(W)$, and $f=\xi_1 x_1$
acts on $\fg_{-1}$ as $\xi_1 \del_{\xi_1}$ which corresponds to the
operator $\xi\del_\xi\in \End(W)$. The element $f=\xi_2\xi_3$ acts on
$\fg_{-1}$ as $\xi_2\del_{x_3}+\xi_3\del_{x_2}$, i.e., as $\del_\eta\in \End(W)$:
\begin{gather}\label{t2}
\begin{tabular}{|l|l|l|l|}
\hline $x_3\longmapsto \xi_2$&$x_2\longmapsto
\xi_3$&$\xi_1x_3\longmapsto \xi_1\xi_2$&
$\xi_1x_2\longmapsto \xi_1\xi_3$\\
\hline
$v_1\otimes \eta\longmapsto v_1$&$v_2\otimes \eta\longmapsto v_2$&
$v_1\otimes \xi\eta\longmapsto v_1\otimes \xi$& $v_2\otimes
\xi\eta\longmapsto v_2\otimes \xi$\\
\hline
\end{tabular}
\end{gather}
Analogously, $f=\xi_1\xi_2\xi_3$ acts as the
action of $\xi\del_\eta\in \End(W)$.

Finally, $f=x_2\xi_2$ acts as $\left(\begin{smallmatrix}
1&0\\0&0 \end{smallmatrix}\right)\otimes \mathbbmss{1} +\eta\del_\eta$:
\be\label{t3}
\begin{tabular}{|l|l|l|l|}
\hline $\xi_2\longmapsto \xi_2$&$\xi_1\xi_2\longmapsto \xi_1\xi_2$&$x_2\longmapsto
x_2$&
$\xi_1x_2\longmapsto \xi_1x_2$\\
\hline

$v_1\longmapsto v_1$&$v_1\otimes \xi\longmapsto v_1\otimes \xi$&
$v_2\otimes \eta\longmapsto v_2\otimes \eta,$& $v_2\otimes
\xi\eta\longmapsto v_2\otimes \xi\eta$\\
\hline
\end{tabular}
\ee
Analogously, $f=\xi_1x_2\xi_2$ acts as the
operator $\left(\begin{smallmatrix} 1&0\\0&0 \end{smallmatrix}\right)\otimes \xi
+\xi\eta\del_\eta$.

Now, let us describe the action of operators $D_{(0, g)}$. First,
taking equation~\eqref{feqg2} into account, we have
\begin{gather*}
D_{(0,\xi_2\xi_3)}=D_{(x_1, 0)}, \quad D_{(0,\xi_1\xi_2x_3)}=
D_{(x_3^{(2)},0)}, \\
D_{(0,\xi_1\xi_3x_2)}=D_{(x_2^{(2)},0)}, \quad D_{(0,\xi_1\xi_2x_2+\xi_1\xi_3x_3)}=D_{(x_2 x_3,0)}.
\end{gather*}
Now, taking the kernel \eqref{kerg} into account, we only have to
establish the three operators corresponding to the functions
$g=\xi_1x_1$, $\xi_1\xi_2\xi_3$, and $\xi_1\xi_2 x_2$.

For $g=\xi_1x_1$, the operator $D_{(0,g)}=\partial_y$ corresponds to $\eta\del_\xi\in\End(W)$:
\be\label{tt1}
\begin{tabular}{|l|l|}
\hline $\xi_1\xi_2\longmapsto x_3$&$\xi_1\xi_3\longmapsto
x_2$\\
\hline
$v_1\otimes\xi\longmapsto v_1\otimes\eta$&
$v_2\otimes\xi\longmapsto v_2\otimes\eta$\\
\hline
\end{tabular}
\ee

If $g=\xi_1\xi_2\xi_3$, then $D_{(0, g)}=
\mathop{\sum}\limits_{1\leq i\leq 3}\xi_i\del_i$,
which corresponds to 
$\xi\del_\xi+\eta\del_\eta\in \End(W)$.

Finally, let $D_{(0, \xi_1\xi_2
x_2)}=x_2\del_{\xi_3}+\xi_1\del_y$. This operator acts as follows:
\be\label{tt2}
\begin{tabular}{|l|l|}
\hline $\del_{x_2}=D_{(\xi_2, 0)}\longmapsto \del_{\xi_3}=
D_{(x_3, 0)}$
&$\xi_3\del_{x_1}+\xi_1\del_{x_3}+y\del_{\xi_2}=D_{(\xi_1\xi_3,0)}$\\
&$\longmapsto\xi_1\del_{\xi_2}+x_2\del_{x_1}=D_{(\xi_1x_2, 0)}$\\
\hline
$v_1\longmapsto v_1\otimes \eta$&
$v_2\otimes \xi\longmapsto
v_2\otimes \xi\eta$\\
\hline
\end{tabular}
\ee
which corresponds to the action of the operator $\left(\begin{smallmatrix}
1&0\\0&0\end{smallmatrix}\right)\otimes \eta+\xi\eta\del_\xi\in \End(W)$.

$D_{(x_1\xi_1+x_2\xi_2, \xi_1\xi_2\xi_3)}$
corresponds to the operator
$\cD:=\left(\begin{smallmatrix}1&0\\0&0\end{smallmatrix}\right)\otimes
\mathbbmss{1}\in\End(W)$, and $D_{(x_2\xi_2,\xi_1\xi_2\xi_3)}$
corresponds to the operator
$\left(\begin{smallmatrix}1&0\\0&0\end{smallmatrix}\right)\otimes \mathbbmss{1} +\xi\del_\xi\in\End(W)$.
\end{proof}

\textit{Critical coordinates of $\widetilde\un^u$}:
$\widetilde\un_1$, $\widetilde\un_2$, $\widetilde\un_3$, and
$\sdim\fg_1=20|20$. There are five lowest-weight vectors in $\fg_1$:
\begin{gather*}
v_1= z_1z_3\del_1 +z_1\del_9 +z_3 z_6\del_6 +z_3 z_7\del_7
+z_3z_8\del_8 +z_3z_9\del_9,\\
v_2= z_3^{(2)}\del_3 +z_3 z_4\del_4 +z_3 z_6\del_6 +z_3 z_9\del_9
+z_4 z_6\del_9 +z_1 z_3\del_1 +z_1\del_9 +z_2z_3\del_2 \\
\hphantom{v_2=}{}
  +z_2 z_4\del_5 +z_2 z_7\del_6 +z_2 z_8\del_9 +
z_3z_6\del_6 +z_3 z_7\del_7 +z_4z_7\del_8,\\
v_3=z_3^{(2)}\del_2 +z_3 z_4\del_5 +z_3 z_7\del_6 +z_3 z_8\del_9
+z_4z_7\del_9,\\
v_4=z_1z_2\del_1 +z_1\del_8 +z_2 z_6\del_6 +z_2 z_7\del_7
+z_2z_8\del_8 +z_2 z_9\del_9 +z_2^{(2)}\del_2\\
\hphantom{v_4=}{}
  +z_2 z_3\del_3 +z_2 z_7\del_7 +z_3 z_5\del_4 +
z_3z_6\del_7 +z_3 z_9\del_8 +z_5 z_6\del_9 +z_1z_2\del_1 +z_1\del_8\\
\hphantom{v_4=}{} +z_2 z_4\del_4 +z_2 z_7\del_7 +z_2 z_9\del_9
+z_4z_6\del_8 +z_2z_5\del_5 +z_2 z_6\del_6 +z_5z_7\del_8,\\
v_5= z_3^{(2)}\del_3 +z_3 z_4\del_4 +z_3 z_6\del_6 +z_3 z_9\del_9
+z_4 z_6\del_9 +z_3 z_5\del_5 +z_3 z_7\del_7
+z_3z_9\del_9 +z_5z_7\del_9.
\end{gather*}
\textbf{No simple partial prolongs.} The module generated by either one
of $v_1$, $v_2$ is $\fg_1$. The modules~$V_i$ generated by either of
$v_3$, $v_4$ are of dimension $4|4$ and $\sdim([\fg_{-1},
V_i])=4|4$; $\sdim V_5=8|8$ and $\sdim([\fg_{-1}, V_i])=8|10$, so there are no
new simple partial prolongs, see~\ref{noprol}.

\subsection{Desuperization} For $\un$ unconstrained, the critical
coordinates are those that correspond to the formerly odd indeterminates.

\section[$\fk\fle\big(15;\widetilde{\un}\big):= \textbf{F}(\fk\fle(5;\un\vert 10))$]{$\boldsymbol{\fk\fle\big(15;\widetilde{\un}\big):= \textbf{F}(\fk\fle(5;\un\vert 10))}$}\label{SKleGradK}

Whenever possible in this section, we do not indicate the shearing vectors.\index{$kmathfrak\fle\big(15;\widetilde{\un}\big)$@$\fk\fle\big(15;\widetilde{\un}\big)$} The Lie superalgebra $\fk\fle(5;\un\vert 10)$ is the complete prolong of its negative part, see Section~\ref{ssExDef}.

Recall that for $\fg=\fk\fle(5|10)$, we have $\fg_\ev=\fsvect(5|0)\simeq d\Omega^{3}(5|0)$ and
$\fg_\od=\Pi\big(d\Omega^{1}(5|0)\big)$ with the natural $\fg_\ev$-action on
$\fg_\od$ and the bracket of any two odd elements being their
product, where we identify
\begin{gather*}
dx_{i}\wedge dx_{j}\wedge dx_{k}\wedge dx_{l}\otimes
\vvol^{-1}=\sign(ijklm)\del_{x_{m}}\text{ for any permutation
$(ijklm)$ of $(12345)$}.
\end{gather*}

\looseness=1 Let $x_{i}$, where $1\leq i \leq 5$, be the even indeterminates,
$\del_i:=\del_{x_i}$. Let $\theta_{ab}$, where $1\leq a,b\leq 5$, be
an odd indeterminate such that $\theta_{ab}=-\theta_{ba}$; in
particular, $\theta_{aa}=0$ and we may assume that $a<b$. Let
$\delta_{ab}:=\del_{\theta_{ab}}$. Let $\fg_0=\fsl(5)=\fsl(V)$ act
on $\fg_{-2}$ as on its tautological 5-dimensional module $V$. Let
$E^2(V)$ be the 2nd exterior power of $V$. For a~basis of the
nonpositive components of $\textbf{F}(\mathfrak{kle}(5;\un|10))$ we
take the following elements (only Chevalley generators are given for
$\fg_0$):
\begin{equation*}\label{vleTilde9}
\tiny
\renewcommand{\arraystretch}{1.4}
\begin{tabular}{|c|l|}
\hline $\fg_{i}$&the basis elements \\ \hline
&\\[-1.5em]
\hline
$\fg_{-2}=V$& $\del_1, \del_2, \del_3, \del_4, \del_5$\\
\hline
$\fg_{-1}=E^2(V)$& $w_{12} =\delta_{12}+
\theta_{34}\partial_5+\theta_{45}\partial_3-
\theta_{35}\partial_4, \ w_{13} =\delta_{13}+
\theta_{25}\partial_4-\theta_{24}\partial_5-\theta_{45}\partial_2$, \\
&$ w_{14} =\delta_{14}+\theta_{23}\partial_5+
\theta_{35}\partial_2-\theta_{25}\partial_3,\ w_{15} =
\delta_{15}+\theta_{24}\partial_3-\theta_{23}\partial_4-
\theta_{34}\partial_2, $\\
&$ w_{23} =\delta_{23}+\theta_{45}\partial_1,\ w_{24} =
\delta_{24}+\theta_{35}\partial_1,\ w_{25} =\delta_{25}+
\theta_{34}\partial_1, $\\
&$\delta_{34},\ \ \delta_{35}, \ \ \delta_{45}$
\\
\hline
$\fg_{0}=\fsl(V)$& $ Z_1=x_1\partial_2+
 \theta_{23}\theta
_{24}\partial_5+
 \theta_{23}\theta
_{25}\partial_4+
 \theta_{24}\theta
_{25}\partial_3+
 \theta_{23}\delta
_{13} +\theta_{24}\delta_{14}+
 \theta_{25}\delta
_{15}, $\\
 & $ Z_2=x_2\partial_3+
 \theta_{34} \theta
_{35}\partial_1+
 \theta_{13}\delta
_{12}+ \theta_{34}\delta_{24}+
 \theta_{35}\delta
_{25},$\\
& $Z_3=x_3\partial_4+
 \theta_{14}\delta
_{13} + \theta_{24}\delta_{23}+
 \theta_{45}\delta
_{35},\ Z_{4}=x_4\partial_5+
 \theta_{15}\delta
_{14} + \theta_{25}\delta_{24}+
 \theta_{35}\delta
_{34}$\\
& $H_{1}=[Z_1,Y_1], \ H_{2}=[Z_2,Y_2], \ H_{3}=[Z_3,Y_3],\
 H_{4}=[Z_4,Y_4]$\\
& $Y_1=x_2\partial_1+
 \theta_{13}\theta
_{14}\partial_5+
 \theta_{13}\theta
_{15}\partial_4+
 \theta_{14}\theta
_{15}\partial_3+
 \theta_{13}\delta
_{23} + \theta_{14}\delta_{24}+
 \theta_{15}\delta
_{25}, $\\
 & $Y_{2}=x_3\partial_2+
 \theta_{24}\theta
_{25}\partial_1+
 \theta_{12}\delta
_{13} + \theta_{24}\delta_{34} +
 \theta_{25}\delta
_{35}$\\
&$Y_{3}=x_4\partial_3+
 \theta_{13}\delta
_{14} + \theta_{23}\delta_{24}+
 \theta_{35}\delta
_{45},\ Y_{4}=x_5\partial_4+
 \theta_{14}\delta
_{15} + \theta_{24}\delta_{25}+
 \theta_{34}\delta
_{35}, $\\
\hline
\end{tabular}
\end{equation*}

The $\fg_0$-module $\fg_1$ is irreducible of dimension 40. The
lowest-weight vector is
\begin{gather*}
v_1  = \theta_{12}\theta_{13}\delta_{15}+
\theta_{12}\theta_{23}\delta_{25}+\theta_{13}\theta_{23}\delta_{35}+
\theta_{14}\theta_{23}\delta_{45}+\theta_{12}\theta_{13}
\theta_{23}\partial_4+\theta_{12}\theta_{14}\theta_{23}\partial_3\\
\hphantom{v_1  =}{}
+\theta_{13}\theta_{14}\theta_{23}\partial_2+
\theta_{13}\theta_{23}\theta_{24}\partial_1+x_5\delta_{45}.
 \end{gather*}

\textit{No simple partial prolongs. Critical coordinates of the
shearing vector for
 $\fk\fle\big(15;\widetilde\un\big)$ are those
corresponding to the formerly odd indeterminates}.

\section[$\widetilde{\fk\fle}\big(15;\widetilde{\un}\big):=\textbf{F}(\fk\fle(9;\un\vert 6))$]{$\boldsymbol{\widetilde{\fk\fle}\big(15;\widetilde{\un}\big):=\textbf{F}(\fk\fle(9;\un\vert 6))}$}\label{SkleTilde}

The construction of $\fk\fle(9;\un|6)$, and its desuperization $\widetilde{\fk\fle}\big(15;\widetilde{\un}\big)$,\index{$kmathfrak\fle\big(15;\widetilde{\un}\big)$@$\widetilde{\fk\fle}\big(15;\widetilde{\un}\big)$}
resemble that of $\fkas$, see Section~\ref{Skas}.

However, $\fkas$ is a~partial prolong of $\fk(1|6)_{\leq 0}\oplus \fkas_1$, where $ \fkas_1$ is a~``half" of
 $\fk(1|6)_1$. Indeed, $\fkas_1$ is one of the two irreducible modules
whose direct sum is $\fk(1|6)_1$, and the Lie superalgebra $\fk\fle(9;\un|6)$ is the
prolong of $\fk(9;\un|6)_-\oplus\fk\fle(9;\un|6)_0$, the latter summand constituting a~half of $\fk(9;\un|6)_0$; this half
corresponds to either of the two possible embeddings
$\fsvect(0|4)\tto\fosp(6|8)$ corresponding to the representations of
$\fsvect(0|4)$ in the space $T_0^0(0|4):=\Vol_0(0|4)/\Kee\cdot \vvol$,
see~\eqref{Vol_0}, and in its dual.

This is why $\fk\fle(9;\un|6)$ is NOT the complete prolong of its negative part, see Section~\ref{ssExDef}.

To determine the component $\fg_0$ of $\fg$, we have to consider a~linear combination of two elements: the central element $Z$ commuting with the image of $\fsvect(0|4)$ in $\fosp(6|8)$ and an outer derivation, say $D=\xi_1\del_{\xi_1}\in \fvect(0|4)$.

\looseness=1 Let $\fG$ be the prolong of the nonpositive part where $\fG_0:=(\fsvect(0|4)\ltimes\Kee D)\bigoplus\Kee Z$ and $\fG_-:=\fg_-$.
Having computed $[\fG_1, \fg_{-1}]$ we determine the coefficients in the linear combination $aZ+bD$ that should belong to
$\fg_0:=\fsvect(0|4)\ltimes\Kee (aD+bZ)$ from the condition \mbox{$[\fG_1,\fg_{-1}]=\fg_0$}.

\looseness=1 To realize the Lie superalgebra $\fg$ by vector fields, we use the
representation of the even part of $\fg$ as $\fsvect(5;\uM)$ and its
odd part as $\Pi\big(d\Omega^1(5;\uM)\big)$: whatever the $\Zee$-grading of
$\fg$, the components $\fg_\ev$ and $\fg_\od$ have the needed
nonpositive part. For convenience, we use $\fgl(5)$-weights of the
elements of $\fg$, having added the outer derivation~-- the grading
operator~-- to $\fsvect(5;\uM)$.

Let $u_1, \dots, u_5$ be a~basis of the space $U$ we used to define
$\fsvect(U)$ and $d\Omega^1(U)$. In our grading, $\deg(u_5)=2$ and
$\deg(u_i)=1$ for $i<5$. Let $x_1, \dots, x_{15}$ be the desuperized indeterminates. Then,
\begin{gather}
\del_{x_1} + \cdots  \lra \del_{u_1},\quad \dots, \quad\del_{x_4} + \cdots \lra \del_{u_4},\nonumber\\
\del_{x_5} + \cdots  \lra \Pi(du_1 du_2), \quad \dots ,\quad \del_{x_{10}} \lra \Pi(du_3 du_4),\nonumber\\
\del_{x_{11}} \lra u_1 \del_{u_5},\quad \dots, \quad\del_{x_{14}}\lra
u_4 \del_{u_5},\quad \del_{x_{15}} \lra \del_{u_5}.\label{U}
\end{gather}
The functor $\Pi$ is interpreted as multiplication (tensoring) by the 1-dimensional
module whose generator $\Pi$ has the following weight $w$ to make
the weight and degree compatible:
\begin{gather*}
w(\Pi) = \left(-\nfrac12, \dots, -\nfrac12\right), \\
\deg(\Pi)=-\nfrac52.
\end{gather*}
To get rid of fractions, we multiply all weights by 2; assuming that $\deg du_i=\deg u_i$ we have
\begin{gather*}
w(x_5) = w(\Pi) + w(du_1) + w(du_2)  \\
\hphantom{w(x_5)}{}= (-1,-1,-1,-1,-1) + (2,0,0,0,0) + (0,2,0,0,0) \\
\hphantom{w(x_5)}{} = (1,1,-1,-1,-1).
\end{gather*}

\looseness=1 Now, the weights are symmetric in the sense that if there is an
element of weight $(2,0,0,0,0)$, there should be elements whose
weight have all coordinates but one equal to 0, one coordinate being
equal to~2. This symmetry helps to find correct expressions of the
vector fields in each component. Thus, the weights of the
indeterminates in the new grading are as follows:
\begin{alignat*}{4}
& x_1\to \{-2,0,0,0,0\}, \quad && x_6\to \{1,-1,1,-1,-1\}, \quad &&x_{11}\to \{2,0,0,0,-2\}, &\\
& x_2\to \{0,-2,0,0,0\}, \quad && x_7\to \{1,-1,-1,1,-1\}, \quad && x_{12}\to \{0,2,0,0,-2\}, & \\
& x_3\to \{0,0,-2,0,0\}, \quad && x_8\to \{-1,1,1,-1,-1\}, \quad && x_{13}\to \{0,0,2,0,-2\}, & \\
& x_4\to \{0,0,0,-2,0\}, \quad && x_9\to \{-1,1,-1,1,-1\}, \quad &&x_{14}\to \{0,0,0,2,-2\}, & \\
& x_5\to \{1,1,-1,-1,-1\}, \quad && x_{10}\to \{-1,-1,1,1,-1\}, \quad && x_{15}\to \{0,0,0,0,-2\}. &
\end{alignat*}
The \looseness=1 degree is equal to one half of (the sum of the first 4 coordinates plus the doubled fifth one).

In equation~\eqref{kle9/6} we give the basis of the negative part and
generators of the 0th component. It is possible to generate the
semi-simple part of $\fg_0$ by just 1 positive and 4 negative
generators, or 4 positive and 1 negative ones, but for symmetry we
give 4 and 4 of them. These 8 generators do not generate the element
$D+Z\in[\fg_1,\fg_{-1}]$ of weight $(0,0,0,0,0)$, so we give it
separately.
\begin{equation}\label{kle9/6} \footnotesize
\renewcommand{\arraystretch}{1.4}
\begin{tabular}{|l|l|} \hline
$\fg_{i}$&the generators \\ \hline
&\\[-1.5em]
\hline
$\fg_{-2}$&$\partial_{15}
$\\
\hline
$\fg_{-1}\simeq \textbf{F}\big(T_0^0(0|4)\big)$,&$\del_1+x_{11}\del_{15}, \
\del_2+x_{12}\del_{15}, \ \del_3+x_{13}\del_{15}, \
\del_4+x_{14}\del_{15}, \ \del_5+x_{10}\del_{15}$,\\
see \eqref{Vol_0}& $\del_6+x_9\del_{15}$, $\del_7+x_8\del_{15}$, $\del_8, \dots,
\del_{14}$\\
\hline
$\fg_{0}\simeq \Kee(D+Z)$& $\{-1,-1,-1,1,1\}\to
x_2\del_6+x_9\del_{12}+x_3\del_5+x_{10}\del_{13}+x_1\del_8$\\
$\rtimes\fsvect(4;\One)$& $\hphantom{\{-1,-1,-1,1,1\}\to}{}+
x_7\del_{11}+x_1x_7\del_{15}$\\
& $\{-1,-1,1,-1,1\}\to
x_2\del_7+x_8\del_{12}+x_4\del_5+x_{10}\del_{14}+
x_1\del_9$\\
& $\hphantom{\{-1,-1,1,-1,1\}\to}{}+x_6\del_{11}+x_1x_6\del_{15}$\\
&$\{-1,1,-1,-1,1\}\to x_3\del_7+x_8\del_{13}+x_4\del_6+x_9\del_{14}+
x_1\del_{10}$\\
& $\hphantom{\{-1,1,-1,-1,1\}\to}{} +x_5\del_{11}+x_1x_5\del_{15}$\\
&$\{1,-1,-1,-1,1\}\to x_2\del_{10}+x_5\del_{12}+
x_2x_5\del_{15}+x_3\del_9+x_6\del_{13}$\\
&$\hphantom{\{1,-1,-1,-1,1\}\to} +x_3x_6\del_{15}+x_4\del_8+x_7\del_{14}+
x_4x_7\del_{15}$\\
&$\{0,0,0,0,0\}\to
x_1\del_1+x_8\del_8+x_9\del_9+x_{10}\del_{10}+x_{12}\del_{12}+
x_{13}\del_{13}$\\
& $\hphantom{\{0,0,0,0,0\}\to}{}+x_{14}\del_{14}+x_{15}\del_{15}$\\
&$\{0,0,0,4,-2\}\to x_{14}\del_4+x_{14}^{(2)}\del_{15},\ \
\{0,0,4,0,-2\}\to x_{13}\del_3+x_{13}^{(2)}\del_{15}$\\
&$\{0,4,0,0,-2\}\to x_{12}\del_2+x_{12}^{(2)}\del_{15},\ \
\{4,0,0,0,-2\}\to x_{11}\del_1+x_{11}^{(2)}\del_{15}$\\
\hline
\end{tabular}
\end{equation}

\textit{No simple partial prolongs. Critical coordinates for
$\widetilde{\fk\fle}\big(15;\widetilde\un\big)$ are those corresponding to
formerly odd indeterminates}.

\section[$\fk\fle_3\big(20;\widetilde{\un}\big):= \textbf{F}(\fk\fle(9;\un\vert 11))$]{$\boldsymbol{\fk\fle_3\big(20;\widetilde{\un}\big):= \textbf{F}(\fk\fle(9;\un\vert 11))}$}\label{Skle_3}

Whenever possible in this section, we do not indicate the shearing vectors.\index{$kmathfrak\fle_3\big(20;\widetilde{\un}\big)$@$\fk\fle_3\big(20;\widetilde{\un}\big)$} The Lie superalgebra $\fk\fle(9;\un\vert 11)$ is the complete prolong of its
negative part, see Section~\ref{ssExDef}.\index{$kmathfrak\fle_3\big(20;\widetilde{\un}\big)$@$\fk\fle_3\big(20;\widetilde{\un}\big)$}

\subsection[Description of $\fk\fle(9;\un|11)_-$]{Description of $\boldsymbol{\fk\fle(9;\un|11)_-}$} We consider the
realization of $\fg=\fk\fle$ as the direct sum of
$\fg_\ev=\fsvect(U)$ and $\fg_\od=\Pi\big(d\Omega^1(U)\big)$, where
$U=\Span(u_1,\dots,u_5)$. Let $i,j=1, 2$, while $a,b,c= 3,4,5$. Let
$\{ijabc\}=\{12345\}$ as sets, $\del_\alpha:=\del_{u_\alpha}$ for
any index $\alpha$. Set, cf.~\eqref{regrad}:
\begin{equation}\label{SQ}
\deg u=(3,3,2,2,2), \quad \deg du=(0,0,-1,-1,-1), \quad \text{where $u=(u_1,\dots,u_5)$}.
\end{equation}
Then,
\begin{gather*}
\fg_{-3}=\Span(\del_1, \del_2),\\
\fg_{-2}=\Span(\del_a, du_a\wedge du_b \text{~~for any
$a,b=3,4,5)$},\\
\fg_{-1}=\Span(u_a\del_i, du_i\wedge du_a \text{~for any $i=1,2$,
$a=3,4,5$}).
\end{gather*}
The brackets are as in grading $K$, see \eqref{regrad}:
\begin{gather*}
[\fg_{-1},\ \fg_{-1}]\colon \ [du_1\wedge dux_a, du_2\wedge
du_b]=\del_c\text{~for $\{a,b,c\}=\{3,4,5\}$},\\
\hphantom{[\fg_{-1},\ \fg_{-1}]\colon}{} \ [du_i\wedge du_a,
u_b\del_i]=du_a\wedge du_b,\\
[\fg_{-1},\ \fg_{-2}]\colon \  [du_i\wedge du_a,\ du_b\wedge du_c]=\del_j,\
[\del_a,\ u_a\del_i]=\del_i.
\end{gather*}

\subsubsection[Description of $\fk\fle(9;\un|11)_-$ in terms of vector fields]{Description of $\boldsymbol{\fk\fle(9;\un|11)_-}$ in terms of vector fields}
We use the realization of Section~\ref{SkleTilde} with the same weights and degrees~\eqref{SQ}.

Recall that $\fg:=\fk\fle(9;\un|11)$ is the prolong of its
\textit{negative} part, and ${\fg_\ev=\fsvect(5; \underline{M})}$. For a~basis of the negative part we take the
following elements, where we denote the 20 indeterminates by~$x$,
set $\delta_i:=\del_{x_i}$:
\begin{equation*}\label{kle20_3} \footnotesize
\renewcommand{\arraystretch}{1.4}
\begin{tabular}{|l|l|} \hline
$\fg_{i}$&the generators \\ \hline
&\\[-1.5em]
\hline
$\fg_{-3}$&$ \{2,0,0,0,0\}\to \delta_{19}, \ \
\{0,2,0,0,0\}\to \delta_{20}$\\
\hline
$\fg_{-2}$&$\{0,0,2,0,0\}\to \delta_{13}+
x_7\delta_{19}+x_8\delta_{20}, \
\{0,0,0,2,0\}\to \delta_{14}+x_9\delta_{19}+x_{10}\delta_{20}$,\\
&$
\{0,0,0,0,2\}\to \delta_{15}+x_{11}\delta_{19}+x_{12}\delta_{20},$\\
&$\{1,1,-1,-1,1\}\to \delta_{16}, \ \
\{1,1,-1,1,-1\}\to \delta_{17}, \ \
\{1,1,1,-1,-1\}\to \delta_{18} $\\
\hline
$\fg_{-1}\simeq $&$\{-1,1,-1,1,1\}\to \delta_1+x_5\delta_{15}+
x_6\delta_{14}+x_9\delta_{16}+x_{11}\delta_{17}+x_{18}\delta_{20}$\\
&$\{-1,1,1,-1,1\}\to \delta_2+x_4\delta_{15}+x_6\delta_{13}+
x_7\delta_{16}+x_{11}\delta_{18}+x_{17}\delta_{20}$\\
&$\{-1,1,1,1,-1\}\to \delta_3+x_4\delta_{14}+x_5\delta_{13}+
x_7\delta_{17}+x_9\delta_{18}+x_{16}\delta_{20}$\\
&$\{1,-1,-1,1,1\}\to \delta_4+x_{10}\delta_{16}+x_{12}\delta_{17}+
x_{18}\delta_{19}$\\
&$\{1,-1,1,-1,1\}\to \delta_5+x_8\delta_{16}+x_{12}\delta_{18}+
x_{17}\delta_{19}$\\
&$\{1,-1,1,1,-1\}\to \delta_6+x_8\delta_{17}+x_{10}\delta_{18}+
x_{16}\delta_{19}$\\
&$\{2,0,-2,0,0\}\to \delta_7, \ \
\{0,2,-2,0,0\}\to \delta_8, \ \
\{2,0,0,-2,0\}\to \delta_9, \ \
\{0,2,0,-2,0\}\to \delta_{10}$\\
&$\{2,0,0,0,-2\}\to \delta_{11}, \ \
\{0,2,0,0,-2\}\to \delta_{12}$\\
\hline
\hline
\end{tabular}
\end{equation*}

\textit{No simple partial prolongs. Critical coordinates of $\widetilde\un^u$ for $\fk\fle_3\big(20;\widetilde\un\big)$ are those
corresponding to formerly odd indeterminates}.

\section[$\fk\fle_2\big(20;\widetilde{\un}\big):= \textbf{F}(\fk\fle(11;\un\vert 9))$]{$\boldsymbol{\fk\fle_2\big(20;\widetilde{\un}\big):= \textbf{F}(\fk\fle(11;\un\vert 9))}$}\label{Skle_2}

Whenever possible in this section, we do not indicate the shearing vectors.\index{$kmathfrak\fle_2\big(20;\widetilde{\un}\big)$@$\fk\fle_2\big(20;\widetilde{\un}\big)$} The Lie superalgebra $\fk\fle(11;\un\vert 9)$\index{$kmathfrak\fle_2\big(20;\widetilde{\un}\big)$@$\fk\fle_2\big(20;\widetilde{\un}\big)$} is the complete prolong of its negative part, see Section~\ref{ssExDef}.

\subsection[Description of $\fk\fle(11;\un|9)_-$]{Description of $\boldsymbol{\fk\fle(11;\un|9)_-}$} In \cite{ShP}, the Lie
superalgebra $\fk\fle$ was constructed from a~central extension of
$\fsle^{(1)}(4)$ with central element further denoted by $c$. The algebra
$\fsle^{(1)}(4)$ was considered in the grading where the degrees of the odd indeterminates
are all $0$. The regradings of this realization are
listed in equation~\eqref{regrad}. Let us give details.

Let the degrees of the generating functions of $\fsle^{(1)}(4)$ be determined as follows:
\[
\deg \xi_3=\deg \xi_4=0, \quad \deg q_3=\deg q_4=2, \quad \deg q_i=\deg
\xi_i=1 \quad \text{for $i=1,2$}.
\]
Then, (recall that the parities of the function are opposite to the
``natural'' ones, and $c$ is even)
\[
 \fg_{-2}=\Span(c,\xi_3,\xi_4,\xi_3\xi_4),\qquad
\fg_{-1}=\Span(\xi_1,\xi_2,q_1,q_2)\otimes\Lambda(\xi_3,\xi_4), \\
\]
with the nonzero brackets of the generating functions $f$ and $g$ in $\xi_3$ and $\xi_4$ being as follows:
\begin{gather*}
{}[f\xi_1,\ g\xi_2]=c \int_\xi(fg\xi_1\xi_2), \quad \text{where $\int_\xi F={}$ coeff. of $\xi_1\xi_2\xi_3\xi_4$ in the expansion of $F$},\\
{}[f\xi_i,\ gq_i]=\begin{cases}0&\text{if $f,g\in\Kee$},\\
fg&\text{otherwise} \end{cases} \quad \text{for $i=1,2$}.
\end{gather*}

\subsubsection[Description of $\fk\fle(11;\un|9)_-$ in terms of vector fields]{Description of $\boldsymbol{\fk\fle(11;\un|9)_-}$ in terms of vector fields}

The above was a~description easy to understand for humans. To compute with the aid of \textit{SuperLie}, we use the realization
of Section~\ref{Skle_3} with the same weights and the degrees given by (compare with~\eqref{SQ})
\begin{equation*}
\deg u=(2,2,2,1,1), \quad \deg du=(0,0,0,-1,-1), \quad \text{where $u=(u_1,\dots,u_5)$}.
\end{equation*}
Let us express the basis of $\fg_{-1}$ in terms of the $u_i$
introduced in~\eqref{U}:
\begin{gather}
\del_{x_1} + \cdots \lra \del_{u_4}, \quad
\del_{x_2} + \cdots \lra \del_{u_5},\nonumber\\
\del_{x_3} + \cdots \lra \Pi(du_1 du_4), \quad \dots ,\quad
\quad \ \del_{x_8} +\cdots \lra \Pi(du_3 du_5),\nonumber\\
\del_{x_9} \lra u_4 \del_{u_1}, \quad \dots ,\quad
\del_{x_{14}} \lra u_5 \del_{u_3},\nonumber\\
\del_{x_{15}} \lra \Pi(u_4 du_4 du_5), \quad
\del_{x_{16}} \lra \Pi(u_5 du_4 du_5),\nonumber\\
\del_{x_{17}} \lra \del_{u_1}, \quad \dots ,\quad \del_{x_{19}} \lra
\del_{u_3},\quad
\del_{x_{20}} \lra \Pi(du_4 du_5).\label{klePi}
\end{gather}
The Lie superalgebra $\fk\fle(11;\un|9)$ is the prolong of the
\textit{negative} part. For a~basis of the negative part we take the
following elements, see \eqref{kle20}. For their weights we take
\[
w(u_i) = w(du_i) = (0,\dots,2,\dots,0), \quad w(\Pi)=(-1,\dots,-1).
\]
We select the degree of $\Pi$ so as to ensure the correct degrees of
the $\del_{x_i}$, see \eqref{kle20}, where by abuse of notation
$\del_{i}:=\del_{x_i}$. Looking at the expression of
$\del_{x_{20}}$, see~\eqref{klePi}, we set $\deg(\Pi) = -4$.
Likewise, the weights of $\del_{15}$ and $\del_{16}$, see~\eqref{kle20}, are deduced from their expressions in terms of the
$u_i$, see \eqref{klePi}:
\begin{equation}\label{kle20} \footnotesize
\renewcommand{\arraystretch}{1.4}
\begin{tabular}{|l|l|} \hline
$\fg_{i}$&the generators \\ \hline
&\\[-1.5em]
\hline
$\fg_{-2}$&$ \{-2,0,0,0,0\}\to \del_{17}, \
\{0,-2,0,0,0\}\to \del_{18}, \
\{0,0,-2,0,0\}\to \del_{19}$, \\
& $\{-1,-1,-1,1,1\}\to \del_{20}$\\
\hline
$\fg_{-1}\simeq $&$ \{0,0,0,-2,0\}\to \del_1+
x_9\del_{17}+x_{10}\del_{18}+x_{11}\del_{19}+
x_{15}\del_{20}$,\\
&$\{0,0,0,0,-2\}\to \del_2+x_{12}\del_{17}+
x_{13}\del_{18}+x_{14}\del_{19}+x_{16}\del_{20}$,\\
&$\{1,-1,-1,1,-1\}\to \del_3+x_8\del_{18}+
x_6\del_{19}+x_{12}\del_{20}$,\\
& $\{1,-1,-1,-1,1\}\to \del_4+x_7\del_{18}+
x_5\del_{19}+x_9\del_{20}$,\\
&$\{-1,1,-1,1,-1\}\to \del_5+x_8\del_{17}+x_{13}\del_{20}$,\\
& $\{-1,1,-1,-1,1\}\to \del_6+x_7\del_{17}+x_{10}\del_{20}$,\\
&$\{-1,-1,1,1,-1\}\to \del_7+x_{14}\del_{20},\
\{-1,-1,1,-1,1\}\to \del_8+x_{11}\del_{20}$,\\
&$\{-2,0,0,2,0\}\to \del_9, \
\{0,-2,0,2,0\}\to \del_{10}, \
\{0,0,-2,2,0\}\to \del_{11}$\\
&$\{-2,0,0,0,2\}\to \del_{12}, \
\{0,-2,0,0,2\}\to \del_{13}, \
\{0,0,-2,0,2\}\to \del_{14},$\\
&$\{-1,-1,-1,3,1\}\to \del_{15}, \
\{-1,-1,-1,1,3\}\to \del_{16}$\\
\hline
\end{tabular}
\end{equation}

\textit{No simple partial prolongs. The critical coordinates of the shearing vector for
$\fk\fle_2\big(20;\widetilde\un\big)$ are those corresponding to the formerly
odd indeterminates}. Explicitly: noncritical coordinates of the
shearing vector correspond to $x_1$, $x_2$, $x_{17}$, $x_{18}$,
$x_{19}$.

\section[The Lie superalgebra $\fm\fb(4\vert 5)$ over $\Cee$]{The Lie superalgebra $\boldsymbol{\fm\fb(4\vert 5)}$ over $\boldsymbol{\Cee}$}\label{Smb}

In this section, we illustrate the algorithm presented in detail in~\cite{Shch}, verify and rectify one formula from \cite{CCK}.
This
algorithm allows one to describe vectorial Lie superalgebras by
means of differential equations. In \cite{Sh5, Sh14} the algorithm
was used to describe the exceptional simple vectorial Lie
superalgebras over $\Cee$.

The Lie superalgebra $\fm\fb(4|5)$ has three realizations as a
transitive\index{(Lie (super)algebra, transitive@Lie (super)algebra, transitive} and primitive\index{(Lie (super)algebra, primitive@Lie (super)algebra, primitive} (i.e., not preserving invariant foliations
on the space where it is realized by means of vector fields) vectorial Lie
superalgebra. Speaking algebraically, the requirement that it should be
transitive and primitive vectorial Lie superalgebra is the same as
to have a~$W$-filtration, so $\fm\fb(4|5)$ has three $W$-filtrations.

Two of these $W$-filtrations are of depth 2, and one is of depth 3.
In each realization this Lie superalgebra is the complete prolong of
its negative part, see Section~\ref{ssExDef}. In this section we
consider the case of depth 3 (the grading $K$); i.e., we explicitly
solve the differential equations singling out our Lie superalgebra.
We thus explicitly obtain the expressions for the elements of~$\fm\fb(4|5; K)$.

In this realization, the Lie superalgebra
$\fg=\fm\fb(3|8)=\fm\fb(4|5; K)$ is the complete prolong of its
negative part $\fg_-=\fg_{-3}\oplus \fg_{-2}\oplus \fg_{-1}$, where
\[
\sdim\fg_{-3}=0|2, \quad \sdim\fg_{-2}=3|0, \quad \sdim\fg_{-1}=0|6.
\]
We would like to embed $\fg_-$ into the Lie superalgebra
\begin{equation*}\label{u}
\fv:=\fvect(3|8)=\fder\Cee[u_1,u_2,u_3;
\eta_1,\eta_2,\eta_3,\zeta_1, \zeta_2,\zeta_3, \chi_1, \chi_2]
\end{equation*}
considered with the grading
\[
\deg \eta_i=\deg \xi_i=1, \quad \deg u_i=2, \quad \deg \chi_j=3
\quad \text{for any $i,j$}.
\]
According to the algorithm described in~\cite{Shch}, we find in
$\fv_-$ two mutually commuting families of elements: $X$-vectors
(the basis of $\fg_-$) and $Y$-vectors. The table of
correspondences, where $i=1,2,3$ and $j=1,2$:
\begin{equation*}
\label{111}\footnotesize
\begin{tabular}{|c|c|c|c|}
\hline $k$ & basis in $\fm\fb_{-k}$ & $X$&$Y$\\
\hline $-1$ & $q_1$, $q_2$, $q_3$ & $X_{\eta_i}$&$Y_{\eta_i}$\\
 & $\xi_2\xi_3$, $\xi_3\xi_1$, $\xi_1\xi_2$ & $X_{\zeta_i}$&$Y_{\zeta_i}$\\
\hline $-2$ & $\xi_1$, $\xi_2$, $\xi_3$ & $X_{u_i}$&$Y_{u_i}$\\
\hline $-3$ & $1$, $\widehat 1$ & $X_{\chi_j}$&$Y_{\chi_j}$\\
\hline
\end{tabular}
\end{equation*}
The nonzero commutation relations for the $X$-vectors are of the form ($(i,j,k)\in A_3$):
 \begin{gather*}
 {}[X_{\eta_i},  X_{u_i}]=-X_{\chi_1}, \quad
[X_{\zeta_i},  X_{u_i}]=-X_{\chi_2},\quad
{}[X_{\eta_i}, X_{\zeta_k}]=-X_{u_j}, \quad [X_{\eta_i},\ X_{\zeta_j}]=X_{u_k}.
\end{gather*}
The nonzero commutation relations for the $Y$-vectors correspond to the negative of the above structure constants:
 \begin{gather*}
 {}[Y_{\eta_i},  Y_{u_i}]=Y_{\chi_1}, \quad
[Y_{\zeta_i},  Y_{u_i}]=Y_{\chi_2},\quad
{}[Y_{\eta_i},  Y_{\zeta_k}]=Y_{u_j}, \quad [Y_{\eta_i},Y_{\zeta_j}]=-Y_{u_k}.
\end{gather*}
Let us represent an arbitrary vector field $D\in\fvect(3|8)$ in the form
\begin{equation} \label{fldD}
X=F_1Y_{\chi_1}+ F_2Y_{\chi_2}+\mathop{\sum}\limits_{1\leq i\leq
3} (f_{\zeta_i}Y_{\zeta_i}+f_{\eta_i}Y_{\eta_i}+f_{u_i}Y_{u_i} ).
\end{equation}
As it was shown in \cite{Shch}, any $X\in\fm\fb(3|8)$ is completely
determined by a~pair of functions $F_1$, $F_2$ by means of equations,
where $i=1,2,3$ and $(i,j,k)\in A_3$:
\begin{gather}
 \label{coor1}
Y_{\zeta_i}(F_1)=0,\quad Y_{\eta_i}(F_1)=-(-1)^{p(f_{u_i})}f_{u_i}=
Y_{\zeta_i}(F_2),\quad Y_{\eta_i}(F_2)=0,
\\
 \label{coor2}
Y_{\zeta_i}(f_{u_i})= Y_{\eta_i}(f_{u_i})=0,
\\
 \label{coor3}
Y_{\zeta_i}(f_{u_j})= (-1)^{p(f_{\eta_k})}f_{\eta_k}, \quad
Y_{\eta_i}(f_{u_j})=-(-1)^{p(f_{\zeta_k})}f_{\zeta_k},
\\
\label{coor4}
Y_{\zeta_i}(f_{u_k})= -(-1)^{p(f_{\eta_j})}f_{\eta_j}, \quad
Y_{\eta_i}(f_{u_k})=(-1)^{p(f_{\zeta_j})}f_{\zeta_j}.
\end{gather}
(Comment: since $(i,j,k)\in A_3$, i.e., is an even permutation, the
formulas \eqref{coor3} and \eqref{coor4} are different; of course
one can express the system by one formula inserting the sign of
permutation.) Therefore the functions $F_1$, $F_2$ must satisfy the
following three groups of equations:
\begin{equation}
 \label{eqF}
 Y_{\zeta_i}(F_1)=0,\quad Y_{\eta_i}(F_1)=
Y_{\zeta_i}(F_2),\quad Y_{\eta_i}(F_2)=0 \quad \text{for $i=1,2,3$}.
\end{equation}
The relations \eqref{coor1}, \eqref{coor3} determine the remaining
coordinates while the relations \eqref{coor2}, \eqref{coor4} follow
from \eqref{coor1}, \eqref{coor3} and the commutation relations that the $Y$-vectors
obey. Indeed, since $p(f_{u_j})=p(f_{u_i})=p(X)$, we have
\[
Y_{\zeta_i}(f_{u_j})=\begin{cases}-
(-1)^{p(f_{u_i})}Y_{\zeta_i}Y_{\zeta_i}(F_2)=0&\text{for $i=j$},\\
-(-1)^{p(f_{u_j})}Y_{\zeta_i}Y_{\zeta_j}(F_2)=
(-1)^{p(f_{u_i})}Y_{\zeta_j}Y_{\zeta_i}(F_2)=
-Y_{\zeta_j}(f_{u_i})&\text{for $i\neq j$}.
\end{cases}
\]
Besides,
\[
f_{\zeta_k}=-(-1)^{p(f_{\zeta_k})}Y_{\eta_i}(f_{u_j})=
-Y_{\eta_i}Y_{\zeta_j}(F_2)=Y_{\zeta_j}Y_{\eta_i}(F_2)+Y_{u_k}(F_2)=Y_{u_k}(F_2).
\]
We similarly get the expressions for the remaining coordinates:
\[f_{\eta_k}=Y_{u_k}(F_1).\]

Therefore, an arbitrary element $X\in\fm\fb(3|8)$ is of the form
\begin{gather}\label{XF}
X=X^{F}=F_{1}Y_{\chi_1}+F_{2}Y_{\chi_2}+ \sum \limits_{1\leq
i\leq 3}\big(Y_{u_i}(F_2)Y_{\zeta_i}
+Y_{u_i}(F_1)Y_{\eta_i}-(-1)^{p(X)}Y_{\zeta_i}(F_2)Y_{u_i} \big),\!\!\!
\end{gather}
where the pair of functions $F=\{F_1, F_2\}$ satisfies the system of
equations \eqref{eqF}.

We select the $Y$-vectors so that the equations \eqref{eqF} the
functions $F_1$, $F_2$ should satisfy were as simple as possible. For
example, take the following $Y$-vectors, where $i=1,2,3$,\ $s=1,2$,\
$(i,j,k)\in A_3$:
\begin{alignat*}{3}
& Y_{\eta_i}=\partial_{\eta_i}+\zeta_k\partial_{u_j}-
\zeta_j\partial_{u_k}+(\zeta_k\eta_j-\zeta_j\eta_k)
\partial_{\chi_1}-\zeta_j\zeta_k\partial_{\chi_2},\quad &&
Y_{\zeta_i}=\partial_{\zeta_i}, & \\
& Y_{u_i}=\partial_{u_i}+\eta_i\partial_{\chi_1}+\zeta_i\partial_{\chi_2},
\quad && Y_{\chi_s}=\partial_{\chi_s}. &
\end{alignat*}
Then, the corresponding $X$-vectors are of the form
\begin{alignat*}{3}
& X_{\eta_i}=\partial_{\eta_i}+u_i\partial_{\chi_1},\quad &&X_{\zeta_i}=
\partial_{\zeta_i}-\eta_j\partial_{u_k}+
\eta_k\partial_{u_j}-\eta_j\eta_k\partial_{\chi_1}+
u_i\partial_{\chi_2},& \\
& X_{u_i}=\partial_{u_i},\quad && X_{\chi_s}=\partial_{\chi_s}.&
\end{alignat*}

The Lie superalgebra $\fm\fb(3|8)$ consists of the vector fields
preserving the distribution determined by the following equations
for the vector field $D$ of the form~\eqref{fldD}:
\begin{equation}
\label{distr1} f_{u_1}=f_{u_2}=f_{u_3}=F_1=F_2=0.
\end{equation}

Let us express the coordinates $f$ of the field $D$ in the $Y$-basis
in terms of the standard coordinates in the basis of partial derivatives:
\[
D=g_{\chi_1}\partial_{\chi_1} +g_{\chi_2}\partial_{\chi_1}+
\mathop{\sum}\limits_{1\leq i\leq 3}(g_{\zeta_i}\partial_{\zeta_i}+
g_{\eta_i}\partial_{\eta_i} +g_{u_i}\partial_{u_i}).
\]
We get
\begin{gather*}
f_{u_i}=g_{u_i}+g_{\eta_j}\zeta_k-g_{\eta_k}\zeta_j \quad \text{for $1\leq i\leq 3$ and $(i,j,k)\in A_3$}, \nonumber\\
F_1=g_{\chi_1}-\mathop{\sum} g_{u_i}\eta_i, \quad
F_2=g_{\chi_2}-\mathop{\sum}
g_{u_i}\zeta_i-\mathop{\sum}\limits_{1\leq i\leq 3,\ (i,j,k)\in A_3}
g_{\eta_i}\zeta_j\zeta_k.
\end{gather*}
Therefore, in the standard coordinates, the distribution singled out by conditions~\eqref{distr1} is given by the equations:
\begin{gather} g_{u_i}+g_{\eta_j}\zeta_k-g_{\eta_k}\zeta_j=0 \quad \text{for $i=1,2,3$}, \nonumber\\
g_{\chi_1}-\mathop{\sum} g_{u_i}\eta_i=0, \quad
g_{\chi_2}-\mathop{\sum}
g_{u_i}\zeta_i-\mathop{\sum}\limits_{1\leq i\leq 3,\ (i,j,k)\in A_3}
g_{\eta_i}\zeta_j\zeta_k=0.  \label{distr2}
\end{gather}
The three equations
determined by the first line of \eqref{distr2} allow one to express $g_{u_i}$
and substitute into the third line to get
\[
g_{\chi_2}+\mathop{\sum}\limits_{1\leq i\leq 3,\ (i,j,k)\in A_3}
g_{\eta_i}\zeta_j\zeta_k=0.
\]

Assuming that the pairing of the space of vector fields with that of 1-forms is given by the formula
\[\langle
f\partial_\xi,\ gd\xi\rangle =(-1)^{p(g)}fg  \quad \text{for any $f,g\in\cF$},
\]
we see that the
distribution is singled out by Pfaff equations given by the
following 1-forms\footnote{We do not use the formulas thus obtained
in THIS text. However, they describe the algebra in meaningful
terms, ``as preserving a~distribution" and explicitly define this
distribution. So we provide these formulas, and keep them for future use.}:
\begin{gather*}
d u_i+\zeta_jd\eta_k-\zeta_kd\eta_j, \text{~~where $(i,j,k)\in A_3$}, \\
d\chi_1-\mathop{\sum} \eta_id u_i,\quad
d\chi_2+\mathop{\sum}\limits_{i\text{ such that }(i,j,k)\in A_3}
\zeta_j\zeta_kd\eta_i.\label{distr3}
\end{gather*}

Let us now solve the system \eqref{eqF}.

Since $Y_{\zeta_i}=\partial_{\zeta_i}$, the condition
$Y_{\zeta_i}(F_1)=0$ implies that $F_1=F_1(u, \eta,\chi)$. The
condition $Y_{\zeta_i}(F_2)=Y_{\eta_i}(F_1)$ takes the form:
\[
\pderf{F_2}{\zeta_i}=\pderf{F_1}{\eta_i}+\left(\zeta_k\pderf{F_1}{u_j}-
\zeta_j\pderf{F_1}{u_k}\right) +\left(\zeta_k\eta_j-
\zeta_j\eta_k\right)\pderf{F_1}{\chi_1}-\zeta_j\zeta_k\pderf{F_1}{\chi_2},
\]
wherefrom (since $F_1$ does not depend on $\xi$) we see that
\begin{equation} \label{F2}
F_2=\mathop{\sum}\limits_{1\leq i\leq 3}\zeta_i\pderf{F_1}{\eta_i} -
\mathop{\sum}\limits_{i=1,2,3,~~(i,j,k)\in
A_3}\zeta_i\zeta_j\left(\pderf{F_1}{u_k}+\eta_k\pderf{F_1}{\chi_1}\right)
- \zeta_1\zeta_2\zeta_3\pderf{F_1}{\chi_2} +\alpha_2,
\end{equation}
where $\alpha_2=\alpha_2(u,\eta,\chi)$, i.e., does not depend on
$\zeta$.

Let us consider the last group of equations \eqref{eqF}:
\begin{equation}\label{newEQ}
Y_{\eta_i}(F_2)=0\quad \text{for $i=1,2,3$}.
\end{equation}
To solve this system, take the expression \eqref{F2} for $F_2$ and
apply the operator~$Y_{\eta_i}$. As a~result, we get a~function
depending on various indeterminates, in particular, on $\zeta_j$. By
virtue of~\eqref{newEQ}, the coefficients of all monomials in
$\zeta$ should vanish. Observe that the coefficient of
$\zeta_1\zeta_2\zeta_3$ vanishes automatically. The terms of degree
$0$ in $\zeta$ are of the form:
\[
\pderf{\alpha_2}{\eta_i}=0 \Longrightarrow \alpha_2=\alpha_2(u,\chi).
\]
Now, let us look at the degree 1 terms in $\zeta$. To get them we
should either take the term independent of $\zeta$ in expression~\eqref{F2} for $F_2$ (and this is $\alpha_2$), and apply to it the
degree 1 terms in $\zeta$ of $Y_{\eta_i}$, i.e.,
\[
\zeta_k\del_{u_j}-\zeta_j\del_{u_k}+(\zeta_k\eta_j-\zeta_j\eta_k)\del_{\chi_1},
\]
 or, the other way round, take the degree 1 terms in $\zeta$ in
\eqref{F2}, i.e., $\sum_s\zeta_s\frac{\del F_1}{\del\eta_s}$, and apply to it the degree 0 in $\zeta$ term of the operator
$Y_{\eta_i}$, i.e., $\del_{\eta_i}$.

Therefore, the terms of degree $1$ in $\zeta$ are of the form:
\[
\zeta_k\pderf{\alpha_2}{u_j}- \zeta_j\pderf{\alpha_2}{u_k}
+\left(\zeta_k\eta_j- \zeta_j\eta_k\right)\pderf{\alpha_2}{\chi_1}=
\zeta_j\frac{\partial^2F_1}{\partial \eta_i\partial\eta_j}+
\zeta_k\frac{\partial^2F_1}{\partial \eta_i\partial\eta_k},
\]
implying that
\begin{equation*}
F_1=\mathop{\sum}\limits_{i=1,2,3, (i,j,k)\in
A_3}\eta_i\eta_j\pderf{\alpha_2}{u_k} +
\eta_1\eta_2\eta_3\pderf{\alpha_2}{\chi_1} +
\alpha_1(u,\chi)+\mathop{\sum}\limits_{1\leq i\leq 3}
f_i(u,\chi)\eta_i.
\end{equation*}
So, the functions $F_1$, $F_2$ are completely determined by the 5
functions $\alpha_1$, $\alpha_2$, $f_1$, $f_2$, $f_3$ that depend
only on $u$ and $\chi$.

The terms of degree $2$ in $\zeta$ follow from the same expression
\eqref{F2} and the same explanation as in the above paragraph leads to the
equation (the coefficient of $\zeta_j\zeta_k$):
\begin{equation} \label{dgr2xi}
\mathop{\sum}\limits_{1\leq s\leq 3}\left(
\frac{\partial^2F_1}{\partial u_s\partial\eta_s}-
\eta_s\pder{\eta_s}\pderf{F_1}{\chi_1}\right)+ \pderf{F_1}{\chi_1} +
\pderf{\alpha_2}{\chi_2} =0.
\end{equation}

Let us expand this equation in parts corresponding to degrees in $\eta$. In degree $0$ we have:
\begin{equation} \label{eqfF}
\mathop{\sum}\limits_{1\leq i\leq 3} (-1)^{p(f_i)}\pderf{f_i}{u_i}+
\pderf{\alpha_1}{\chi_1}+\pderf{\alpha_2}{\chi_2}=0.
\end{equation}
In degrees $1$, $2$, $3$ in $\eta$ the equation \eqref{dgr2xi} is automatically satisfied.

Let us express equation~\eqref{eqfF} in the following more lucid way. We
designate
\[
f_i:=f_i^0+f_i^1\chi_1+f_i^2\chi_2+f_i^{12}\chi_1\chi_2,
\quad
\alpha_s:=\alpha_s^0+\alpha_s^1\chi_1+\alpha_s^2\chi_2+
\alpha_s^{12}\chi_1\chi_2.
\]
The equation \eqref{eqfF} is equivalent to the following system of
four equations:
\begin{gather*}
\sum\limits_{1\leq i\leq 3} \pderf{f_i^{12}}{u_i}=0,\quad
\alpha_1^{12}- \!\sum \limits_{1\leq i\leq 3}
\pderf{f_i^2}{u_i}=0,\quad
\alpha_2^{12}+\!\sum\limits_{1\leq i\leq 3}
\pderf{f_i^1}{u_i}=0,\quad
\alpha_1^1+\alpha_2^2+\!\sum\limits_{1\leq i\leq 3}
\pderf{f_i^0}{u_i}=0.
\end{gather*}

Let us describe the commutation relations in $\fm\fb(3|8)$ more
explicitly. Let us represent the vector field \eqref{XF} as
\begin{gather} \label{XF1}
X^{F}=x^F+ \mathop{\sum}\limits_{1\leq i\leq
3}\big(f_{\zeta_i}Y_{\zeta_i}
+f_{\eta_i}Y_{\eta_i}-(-1)^{p(X)}f_{u_i}Y_{u_i} \big),
\quad \text{where $x^F= F_{1}\partial_{\chi_1}+F_{2}\partial_{\chi_2}$,}\!\!\!
\end{gather}
and observe that, taking relation \eqref{coor1} and \eqref{coor3} into account, we have
\begin{gather*}
\big[X^F, X^G\big]=X^H,\quad \text{where}\\
H_1=\big[x^F,x^G\big]_1+ \sum \limits_{1\leq i\leq
3}\big(f_{u_i}g_{\eta_i}-(-1)^{p(X^G)}f_{\eta_i}g_{u_i}\big),\\
H_2=\big[x^F,x^G\big]_2+\sum\limits_{1\leq i\leq 3}\big(f_{u_i}g_{\zeta_i}- (-1)^{p(X^G)}f_{\zeta_i}g_{u_i}
\big).\end{gather*}
Observe that it suffices to compute only the \textit{defining}
components of $F$, $G$, and $H$:
\begin{equation*}
\begin{tabular}{|l|l|}
\hline the pair&determined by the
set\\
\hline
$F$&$\{\alpha_s, f_i\,|\,s=1,2, \ i=1,2,3\}$\\
$G$& $\{\beta_s, g_i\,|\,s=1,2,\ i=1,2,3\}$\\
$H$&$\{\gamma_s, h_i \,|\,s=1,2,\
i=1,2,3\}$\\
\hline \end{tabular}
\end{equation*}
Then, we get
\begin{gather}
 \gamma_1 =
 \mathop{\sum}\limits_{1\leq i\leq 3}\left(-f_i\pderf{\beta_1}{u_i}+
(-1)^{p(X^G)}\pderf{\alpha_1}{u_i}g_i \right)\nonumber\\
\hphantom{\gamma_1 =}{}  +\left(\alpha_1\pderf{\beta_1}{\chi_1}+
\alpha_2\pderf{\beta_1}{\chi_2}\right)  -
 (-1)^{p(X^F)p(X^G)}\left(\beta_1\pderf{\alpha_1}{\chi_1}+
 \beta_2\pderf{\alpha_1}{\chi_2}\right),\nonumber\\
 \gamma_2  =
\sum\limits_{1\leq i\leq 3}\left(-f_i\pderf{\beta_2}{u_i}+
(-1)^{p(X^G)}\pderf{\alpha_2}{u_i}g_i
 \right)\nonumber\\
\hphantom{\gamma_2 =}{} +\left(\alpha_1\pderf{\beta_2}{\chi_1}+
\alpha_2\pderf{\beta_2}{\chi_2}\right) -
 (-1)^{p(X^F)p(X^G)}\left(\beta_1\pderf{\alpha_2}{\chi_1}+
 \beta_2\pderf{\alpha_2}{\chi_2}\right),\nonumber\\
 h_i  =-\mathop{\sum}\limits_{1\leq r\leq 3}
 f_r\pderf{g_i}{u_r}+
\mathop{\sum}\limits_{1\leq r\leq 3}
 \pderf{f_i}{u_r}g_r \nonumber\\
\hphantom{h_i  =}{}  -(-1)^{p(X^G)}
 \left(\pderf{\alpha_2}{u_j}\pderf{\beta_1}{u_k}-\pderf{\alpha_2}{u_k}
 \pderf{\beta_1}{u_j}
 -\pderf{\alpha_1}{u_j}\pderf{\beta_2}{u_k} +
 \pderf{\alpha_1}{u_k}\pderf{\beta_2}{u_j}\right) \nonumber\\
\hphantom{h_i  =}{}
 + \sum \limits_{s=1, 2} \alpha_s \pderf{g_i}{\chi_s} -
 (-1)^{p(X^F)p(X^G)} \mathop{\sum}\limits_{s=1, 2}
 \beta_s \pderf{f_i}{\chi_s}\quad \text{for
 $i=1,2,3$, $(i,j,k)\in A_3$}. \label{gamma1}
\end{gather}

In what follows we identify the vector field $X^F$ with the collection
\begin{equation}
 \label{coll}
\{\alpha_s,\ f_i\,|\,s=1,2,\  i=1,2,3\}.
\end{equation}
The bracket of vector fields corresponds to the bracket of such collections given by equations~\eqref{gamma1}.

Consider now the even part $\fm\fb(3|8)_\ev$ of our algebra. Since $p(F_1)=p(F_2)=\od$, it follows that $p(\alpha_s)=\od$ and
$p(f_i)=\ev$ for all $s$ and $i$. The component $\fm\fb(3|8)_\ev$ has the three subspaces:
\[
\fm\fb(3|8)_\ev=V_1\oplus V_2\oplus V_3.
\]

The subspace $V_1$ is determined by the collection \eqref{coll} such
that
\[
\{\alpha_1=\alpha_2=0,\; f_i=f_i(u)\chi_1\chi_2\,|\, \mathop{\sum}\limits_{i=1,2,3} \pderf{f_i}{u_i}=0\}.
\]
Equations~\eqref{gamma1} imply that the vector fields generated by such
functions form a~commutative ideal in $\fm\fb(3|8)_\ev$; we will
identify this ideal with $d\Omega^1(3)$:
\[
\{0,0,f_i\,|\,i=1,2,3\} \longmapsto -
\mathop{\sum}\limits_{i\text{~such that~}(ijk)\in A_3} f_idu_j\wedge
du_k.
\]

The subspace $V_2$ is determined by the collection \eqref{coll} such
that $f_i=0$ for $i=1,2,3$. We will identify this space with
$\Omega^0(3)\otimes \fsl(2)$, by setting
\[
\{\alpha(u)(a\chi_1+b\chi_2),\ \alpha(u)(c\chi_1-a\chi_2),\
0,0,0\}\longmapsto \alpha(u)\otimes \begin{pmatrix} a& c\\b&-a
\end{pmatrix},
\]
 where $\alpha\in\Omega^0(3)$, $a,b,c\in\Cee$. Equations~\eqref{gamma1} imply that the subspaces $V_1$ and $V_2$ commute
with each other whereas the brackets of two collections from $V_2$ is in our notation of the form
\[
[f\otimes A, g\otimes B]=fg\otimes [A,B]+df\wedge dg \cdot \tr AB.
\]

Concerning $V_3$, we have the following three natural ways to describe it: in all three cases we take $f_i=f_i(u)$ for all $i$, whereas
for the $\alpha_s$, we select one of the following:
\begin{alignat}{3}
& (a)\ \ \alpha_1=-\sum\pderf{f_i}{u_i}\chi_1,\quad &&\alpha_2=0, &\nonumber\\
& (b)\ \  \alpha_1=0,\quad &&
\alpha_2=-\sum\pderf{f_i}{u_i}\chi_2, &\nonumber \\
& (c)\ \ \alpha_1=-\frac 12 \sum
\pderf{f_i}{u_i}\chi_1,\quad &&  \alpha_2=-\frac 12
\sum\pderf{f_i}{u_i}\chi_2.& \label{3possib}
\end{alignat}

\underline{For $p\neq 2$}, the case (c) is more convenient to simplify
the brackets. Thus, we identify $V_3$ with $\fvect(3)$, by means of
the mapping
\begin{equation*}
\left\{-\frac 12 \sum\pderf{f_i}{u_i}\chi_1,\; -\frac 12
\sum\pderf{f_i}{u_i}\chi_2, \; f_1(u),\ f_2(u),\; f_3(u)\right\}
\longmapsto D_f=-\mathop{\sum} f_i(u)\partial_{u_i}.
\end{equation*}
The actions of $D_f$ on the subspace $V_1$ (as on the space of
$2$-forms) and $V_2$ (as on the space $\cF\otimes\fsl(2)$ of
$\fsl(2)$-valued functions) are natural. The bracket of two elements if the form
$D_f$ is, however, quite different from the usual bracket thanks
to an extra term:
\[
[D_f,  D_g]=D_fD_g-D_gD_f-\frac 12 d(\Div D_f)\wedge d(\Div D_g).
\]

Consider now the odd part: $\fm\fb(3|8)_\od$. We have
$p(F_1)=p(F_2)=\ev$, and hence
\[
p(\alpha_s)=\ev,\quad p(f_i)=\od.
\]

Let $V_4$ consist of collections \eqref{coll} with $f_i=0$. We
identify $V_4$ with $\Omega^0(3)\vvol^{-1/2}\otimes \Cee^2$, by
setting
\begin{equation*}
\{(\alpha(u)w_1,\; \alpha(u)w_2,\; 0,0,0)\}\longmapsto
\alpha(u)\vvol^{-1/2}\otimes \begin{pmatrix} w_2 \\
-w_1\end{pmatrix}.
\end{equation*}
Let $V_5$ consist of the collections \eqref{coll}, where
\begin{equation} \label{V5}
f_i=f_i(u)(v_1\chi_1+v_2\chi_2),\quad
\alpha_1=v_2\sum\pderf{f_i}{u_i}\chi_1\chi_2,\quad
\alpha_2=-v_1\sum\pderf{f_i}{u_i}\chi_1\chi_2.
\end{equation}
We identify $V_5$ with $\Omega^2(3)\vvol^{-1/2}\otimes \Cee^2$, by
assigning to the collection \eqref{V5} the element
\begin{equation*}
\omega\vvol^{-1/2}\otimes \begin{pmatrix} v_1\\v_2 \end{pmatrix},
\quad \text{where~} \omega= -\mathop{\sum}\limits_{i\text{~such
that~}(ijk)\in A_3} f_idu_j\wedge du_k.
\end{equation*}

Let us sum up a~description of the spaces $V_i$ and their elements, see Table~\eqref{tablsoo}.
\begin{table}[ht]\centering
\begin{equation}
\renewcommand{\arraystretch}{1.4}
\label{tablsoo} \footnotesize
\begin{tabular}{|@{\,}c@{\,}|@{\,}c@{\,}|@{\,}c@{\,}|@{\,}c@{\,}|@{\,}c@{\,}|}
 \hline
The space & $\alpha_1$ &$\alpha_2$ & $f_i$ & the element of $V_i$\\
 \hline
 $V_1\cong d\Omega^1(3)$ & $0$ & $0$ &
 $f_i(u)\chi_1\chi_2,\; $ & $\omega=\sum
 f_idu_j\wedge du_k, \; $\\
 & & & $\sum\pderf{f_i}{u_i}=0$ &$d\omega=0$\\
 \hline
 $\begin{matrix}V_2\cong \\
 \Omega^0(3)\otimes \fsl(2)\end{matrix}$ &
$ \alpha(u)(a\chi_1+b\chi_2)$ & $ \alpha(u)(c\chi_1-a\chi_2)$ & 0
& $\alpha(u)\otimes \begin{pmatrix} a&c\\b&-a \end{pmatrix}$\\
 \hline
 $V_3\cong $ & $-\frac 12 f(u)\chi_1$ & $-\frac 12 f(u)\chi_2$
 & $f_i(u)$ & $D=-\sum f_i(u)\partial_{u_i}$ \\
$ \fvect(3)$ & & & $f(u)=\sum\pderf{f_i}{u_i}$ & $\Div
D=-f(u)$\\
 \hline
&&&&\\[-1.5em]
 \hline
 $\begin{matrix}V_4\cong\\
 \Omega^0\vvol^{-1/2}\otimes \Cee^2\end{matrix}$& $ \alpha(u)w_1$ &
 $\alpha(u)w_2$ & $0$ &
 $\frac{\alpha(u)}{\vvol^{1/2}}\otimes \begin{pmatrix} w_2 \\ -w_1
 \end{pmatrix}$\\
 \hline
 $V_5\cong$ &$ v_2 f(u)\chi_1\chi_2$ & $ -v_1 f(u)\chi_1\chi_2$
 &$f_i(u)(v_1\chi_1+v_2\chi_2)$ & $\frac{\omega}{\vvol^{1/2}}\otimes
\begin{pmatrix} v_1 \\
 v_2 \end{pmatrix}$\\
$ \Omega^2\vvol^{-1/2}\otimes \Cee^2$ & &
&$f(u)=\sum\pderf{f_i}{u_i}$ & $\omega=\sum
 f_idu_j\wedge du_k$\\
 \hline
\hline
\end{tabular}
 \end{equation}\vspace{-4mm}
\end{table}

Having explicitly computed the brackets using expressions~\eqref{gamma1} and presenting the result by means of correspondences~\eqref{tablsoo}, we obtain the formulas almost identical to those offered in \cite{CCK}. The difference, however, is vital: the Jacobi identity either holds or not.

We have already given the brackets of the even elements. The brackets of elements of $\fm\fb_\ev$ and $\fm\fb_\od$ are of the form:
\begin{gather*}
 [V_1,V_4]\colon \ \big[\omega, \alpha\vvol^{-1/2}\otimes v\big]=\alpha\cdot
\omega\vvol^{-1/2}\otimes v \in V_5,\\
  [V_2,V_4]\colon \ \big[f\otimes A,\alpha\vvol^{-1/2}\otimes
v\big]=f\alpha\vvol^{-1/2}\otimes Av-df\wedge
d\alpha\vvol^{-1/2}\otimes Av \in V_4\oplus V_5,\\
 [V_3,V_4]\colon \ \big[D,\alpha\vvol^{-1/2}\otimes v\big]=(D(\alpha)-\frac 12
\Div D\cdot \alpha)\vvol^{-1/2}\otimes v\\
 \phantom{[V_3,V_4]\colon} \ {} +\frac 12 d(\Div D)\wedge d\alpha\cdot
\vvol^{-1/2}\otimes
v \in V_4\oplus V_5,\\
{}[V_1,V_5]  =0, \\
{}[V_2,V_5]\colon \  \big[f\otimes A,\omega \vvol^{-1/2}\otimes v\big]= f\omega
\vvol^{-1/2}\otimes Av\in V_5,\\
{}[V_3,V_5]\colon \ \big[D,\omega \vvol^{-1/2}\otimes v\big]= (L_D\omega-\frac 12
\Div D\cdot \omega) \vvol^{-1/2}\otimes v\in V_5.
\end{gather*}
To describe in these terms the bracket of two odd elements, perform
the following natural identifications:
\begin{gather*}
\frac{\Omega^2(3)}{\vvol}\cong \fvect(3)\colon \ \frac{\omega}{\vvol}
\lra D_\omega,\\
i_{D_\omega}(\vvol)=\omega,  \ \text{i.e.,} \
\mathop{\sum}\limits_{\{i,j,k\}=\{1,2,3\}\text{~such that~}(ijk)\in
A_3} f_idx_j\wedge dx_k
\lra \sum f_i\partial_{i},\\
\Lambda^2\Cee^2\cong \Cee\colon \  v\wedge w\lra \det
\begin{pmatrix}v_1 & w_1\\
v_2&w_2 \end{pmatrix}, \\
S^2(\Cee^2)\cong \fsl(2)\colon \  v\cdot w \lra \begin{pmatrix}
-v_1w_2-v_2w_1 & 2v_1w_1\\
-2v_2w_2 & v_1w_2+v_2w_1 \end{pmatrix}.
\end{gather*}
The bracket of two odd elements is of the form:
\begin{gather}
 [V_4,V_4]\colon \ \left[\frac{f}{\vvol^{1/2}}\otimes v, \frac
{g}{\vvol^{1/2}}\otimes w\right]=\frac{df\wedge dg \otimes v\wedge
w}{\vvol}\in\fvect(3)=V_3,\nonumber\\
 [V_5,V_5]\colon \ \left[\frac{\omega_1}{\vvol^{1/2}}\otimes v,
\frac{\omega_2}{\vvol^{1/2}}\otimes w\right]= (D_{\omega_1}(\omega_2)-(\Div
D_{\omega_2})\cdot \omega_1)v\wedge w\in V_1,\nonumber\\
 [V_4,V_5]\colon \ \left[\frac{f}{\vvol^{1/2}}\otimes
v,\frac{\omega}{\vvol^{1/2}}\otimes w\right]\nonumber\\
\hphantom{[V_4,V_5]\colon}{} \  =\frac{f\omega}{\vvol}\otimes v\wedge w - \frac
12(fd\omega-\omega\wedge df)\otimes v\cdot w +df\wedge d(\Div
D_\omega)\otimes v\wedge w.\label{lazha}
\end{gather}
In the last line above, the first summand lies in $V_3$, the second
one in $V_2$, and the third one in~$V_1$. \textit{The difference as
compared with} \cite{CCK}: the coefficient of the third summand in
the last line on \eqref{lazha} should be $1$ whereas in \cite{CCK}
it is equal to $\frac12$.

To verify, compute the Jacobi identity (it holds for 1 and not for
$\frac12$) for the triple
\[
u_3du_2\wedge du_3\in V_1,\quad \frac{u_1}{\vvol^{1/2}}\otimes e_1,
\quad \text{and} \quad \frac{u_2}{\vvol^{1/2}}\otimes e_2\in V_4,
\quad \text{where $e_1$, $e_2$ span $\Cee^2$.}
\]

\underline{For $p=2$}, when case (c) in \eqref{3possib} is not defined, we can select
any one of the cases (a) or (b), we take case (a) for definiteness. In these
cases~(a) and~(b), we get two embeddings $\fvect(3)\subset
\fm\fb(3|8)_\ev$.

\section[The Lie algebra $\textbf{F}(\mathfrak{mb}(3;\un\vert 8))$ is a~true deform of $\fsvect\big(5;\widetilde{\un}\big)$]{The Lie algebra $\boldsymbol{\textbf{F}(\mathfrak{mb}(3;\un\vert 8))}$ is a~true deform of $\boldsymbol{\fsvect\big(5;\widetilde{\un}\big)}\!$}\label{S_F(mb)}

In this section, we describe the analog of the complex Lie
superalgebra $\mathfrak{mb}(3|8)$ for $p=2$ and consider its
desuperization. For brevity, whenever possible we do not indicate the shearing vectors.\index{$F(\mathfrak{mb}(3;\un\vert 8))$@$\textbf{F}(\mathfrak{mb}(3;\un\vert 8))$}

In Section~\ref{Smb}, we showed that an arbitrary vector field
$X^F\in\fg$, where $\fg=\mathfrak{mb}(3|8)$, is of the form~\eqref{XF1} and is determined by 5 functions
$(\alpha_1,\alpha_2,f_1,f_2,f_3)$ in indeterminates
$\chi_1$, $\chi_2$, $u_1$, $u_2$,~$u_3$. Now, in discussing
$\textbf{F}(\mathfrak{mb}(3|8))$, we assume that all these
indeterminates are even.

For consistency we replace $\chi_i$ with $u_{3+i}$, and $\alpha_i$
with $f_{3+i}$. Accordingly we denote $X^F$ by~$X^f$, where
$f=(f_1,f_2,f_3,f_4,f_5)$. The equation \eqref{eqfF} takes the form
\begin{equation}
\label{sv} \mathop{\sum}\limits_{1\leq i\leq 5}\pderf{f_i}{u_i}=0.
\end{equation}

The equation \eqref{sv} is the only condition imposed on the
functions $f_i$, and hence there are no restrictions on the values
of coordinates of the shearing vector corresponding to the
indeterminates $u_i$, including $u_4$ and $u_5$.

Consider the mapping\index{$X^f$}\index{$D^f$}
\begin{equation}\label{imbSvect}
\vf\colon \ \fg\tto \fvect(5), \quad X^f \longmapsto D^f:=\sum
f_i\del_{u_i}.
\end{equation}
Clearly, this is a~linear injective mapping. Formula~(\ref{sv})
implies that $\vf(\fg)=\fsvect(5)$. The mapping $\vf$ is not,
however, an isomorphism of Lie algebras $\fg$ and $\fsvect(5)$.
Indeed, equations~(\ref{gamma1}) rewritten in new notation imply the
following equality (since $p=2$, we skip the signs):
\begin{equation}\label{bracket} \vf\big(\big[X^f, X^g\big]\big)=\big[D^f,
D^g\big]+\mathop{\sum}\limits_{(i,j,k)\in S_3}
\left(\pderf{f_4}{u_i}\pderf{g_5}{u_j} +
\pderf{f_5}{u_i}\pderf{g_4}{u_j}\right)\pder{u_k}.
\end{equation}

Realization of $\fk\fle$ convenient in what follows: for
$\fg = \fk\fle(5|10)$, we have $\fg_\ev= \fsvect(5|0)\simeq
d\Omega^3$ and $\fg_\od=\Pi\big(d\Omega^1\big)$ with the natural
$\fg_\ev$-action on $\fg_\od$, while the bracket of any two odd
elements is their product naturally identified with a
divergence-free vector field.

For any $D=\mathop{\sum}\limits_{1\leq i\leq 5} f_i\del_{u_i}\in
\fsvect(5)$, we define
\[
Z_i(D):=du_i\wedge df_i\in d\Omega^1(5)
\]
and construct the embedding (as a~vector space)
\begin{equation}
\label{kle} \psi\colon \ \fsvect(5|0)\tto \text{F}(\fk\fle), \quad D\longmapsto
D+Z_4(D)+Z_5(D).
\end{equation}

Let us compute the bracket of two fields of the form \eqref{kle}:
\[
\big[D^f+Z_4\big(D^f\big)+Z_5\big(D^f\big),  D^g+Z_4\big(D^g\big)+Z_5\big(D^g\big)\big].
\]
In order not to write too lengthy expressions, let us compute,
separately, the brackets of individual summands. First, let $i=1,2,3$, and $k=4,5$:
\begin{gather}
[f_i\partial_{u_i}, g_k\partial_{u_k}+Z_k(g_k\partial_{u_k})]=
[f_i\partial_{u_i}, g_k\partial_{u_k}+
du_k\wedge dg_k]\nonumber\\
\hphantom{[f_i\partial_{u_i}, g_k\partial_{u_k}+Z_k(g_k\partial_{u_k})]}{}
=\left(f_i\pderf{g_k}{u_i}\right)\pder{u_k}+ du_k\wedge
d\left(f_i\pderf{g_k}{u_i} \right)+
g_k\pderf{f_i}{u_k}\pder{u_i}\nonumber\\
\hphantom{[f_i\partial_{u_i}, g_k\partial_{u_k}+Z_k(g_k\partial_{u_k})]}{}
=[f_i\partial_i,
g_k\partial_k]+Z_k([f_i\partial_i, g_k\partial_k]).\label{br1}
\end{gather}
Here we applied the Leibniz formula for the action of a~vector field
on a~2-form, and the expressions for the Lie derivative along the
vector field $X$:
\[
L_{f_i\partial_i}(du_k)=0\quad \text{and} \quad L_X\circ d=d\circ L_X.
\]

Now, let $k=4$ or 5:
\begin{gather}
[f_k\partial_{u_k}+Z_k(f_k\partial_{u_k}),
g_k\partial_{u_k}+Z_k(g_k\partial_{u_k})]=[f_k\partial_{u_k}+du_k\wedge df_k, g_k\partial_{u_k}+du_k\wedge
dg_k]\nonumber\\
=[f_k\partial_{u_k},  g_k\partial_{u_k}]+
L_{f_k\partial_{u_k}}(du_k\wedge dg_k)+
L_{g_k\partial_{u_k}}(du_k\wedge df_k)+
[du_k\wedge df_k,  du_k\wedge dg_k]\nonumber\\
=\left(f_k\pderf{g_k}{u_k}+ g_k\pderf{f_k}{u_k}
\right)\partial_{u_k}+ df_k\wedge dg_k+du_k\wedge
d\left(f_k\pderf{g_k}{u_k} \right)+
dg_k\wedge df_k+du_k\wedge d\left(g_k\pderf{f_k}{u_k} \right)\nonumber\\
 =[f_k\partial_{u_k}, g_k\partial_{u_k}]+Z_k([f_k\partial_{u_k},g_k\partial_{u_k}]).\label{br2}
\end{gather}

Finally, let $i=4$ and $k=5$:
\begin{gather}
 [f_4\partial_{u_4}+Z_4(f_4\partial_{u_4}),g_5\partial_{u_5}+Z_5(g_5\partial_{u_5})]=
[f_4\partial_{u_4}+du_4\wedge df_4,g_5\partial_{u_5}+du_5\wedge dg_5]\nonumber\\
\qquad{}=\left(f_4\pderf{g_5}{u_4}\right)\pder{u_5}+ du_5\wedge
d\left(f_4\pderf{g_5}{u_4} \right)+ g_5\pderf{f_4}{u_5}\pder{u_4}+
du_4\wedge d\left(g_5\pderf{f_4}{u_5} \right)\nonumber\\
\qquad\quad{}  +\frac{du_4\wedge df_4\wedge du_5\wedge dg_5}{\vvol}\nonumber\\
\qquad{} = [f_4\partial_4, g_5\partial_5]+Z_4([f_4\partial_4,g_5\partial_5])+
Z_5([f_4\partial_4, g_5\partial_5])+ \mathop{\sum}\limits_{(i,j,k)\in
S_3}\pderf{f_4}{u_i} \pderf{g_5}{u_j}\pder{u_k}.\label{br3}
\end{gather}

The expressions \eqref{br1}, \eqref{br2}, and \eqref{br3} show that the through mapping $\psi\circ \vf$ determines an embedding $\fg\tto\fk\fle$, and hence the Lie algebra $\fg$ is isomorphic to the thus-constructed Lie subalgebra of $\fk\fle$.

\begin{Remark} Note that, thanks to formulas \eqref{br1} and
\eqref{br2}, the image of Lie algebra $\fsvect(5)$ under the
embedding
\[
\fsvect(5)\tto \fk\fle, \quad D\longmapsto D+Z_k(D) \quad \text{for any $k$}
\]
is isomorphic to $\fsvect(5)$. The image of the embedding with three
additional terms
\begin{equation}\label{IrasMap}
\fsvect(5)\tto \fk\fle, \quad D\longmapsto D+Z_3(D)+Z_4(D)+Z_5(D)
\end{equation}
is not a~proper subalgebra of $\fk\fle$: it generates the whole~$\fk\fle$. Indeed: take the bracket of the images of two fields of
the form~$f\partial_4, g\partial_5\in\fsvect(5)$; we see, thanks to equation~\eqref{br3}, that the image of $\fsvect(5)$ under the mapping~\eqref{IrasMap} must contain 2-forms such as $du_3\wedge dh$ for certain~$h$, and hence this image is not a~subalgebra. Since the $\fsvect(5)$-module $d\Omega^1(5)$ is irreducible, the image of~\eqref{IrasMap} generates the whole $\fk\fle$.
\end{Remark}

\subsection[The Lie algebra $\textbf{F}(\mathfrak{mb}(3;\un|8))$ is a true deform of $\fsvect\big(5;\widetilde \un\big)$]{The Lie algebra $\boldsymbol{\textbf{F}(\mathfrak{mb}(3;\un|8))}$ is a true deform of $\boldsymbol{\fsvect\big(5;\widetilde \un\big)}$}\label{ssNonSemiSvect}

Indeed, for the shearing vectors of the form $\un_\infty$, all
W-gradings of $\fm\fb$ are the same as over~$\Cee$. None of them has a~maximal subalgebra of codimension 5, whereas $\fsvect(5)$ has such a~subalgebra; cf.\ deforms described in~\cite{Tyu, W} as well.

We consider $\fg$ as a
deform of $\fsvect(5)$ with the grading
\[
\deg u_a=2, \quad \deg u_i=3, \quad \text{where~~} a=1,2,3,\ i=4,5
\]
and the new bracket \eqref{bracket} designated $[[- ,-]]$:\index{$([[- ,-]]$@$[[- ,-]]$}
\begin{equation}\label{brack1}
\big[\big[D^f,D^g\big]\big]=\big[D^f,D^g\big]+c\big(D^f,D^g\big),
\end{equation}
where $D^F=\sum f_i\del_i\in \fsvect(5)$,
$[- ,-]$ is the usual bracket of vector fields, and the
cocycle that determines the deform is
\[
c\big(D^f,D^g\big) = \mathop{\sum}\limits_{(i,j,k)\in S_3} \left(\pderf{f_4}{u_i}\pderf{g_5}{u_j} +
\pderf{f_5}{u_i}\pderf{g_4}{u_j}\right)\pder{u_k}.
\]

All calculations in this realization are rather simple. We have
(observe that thanks to formulas~\eqref{imbSvect} and~\eqref{brack1}
brackets between the elements of $\fg_{-1}$ are nontrivial, and
$\fg_{-1}$ generates the negative part)
\[
\renewcommand{\arraystretch}{1.4}
\begin{tabular}{|l|l|}
\hline $\fg_{i}$&its basis \\ \hline
&\\[-1.5em]
\hline
$\fg_{-3}$&$\del_4,\ \del_5$\\
\hline
$\fg_{-2}$&$\del_1,\ \del_2,\ \del_3$\\
\hline
$\fg_{-1}$&$u_a\del_i, \text{ where } a=1,2,3,\; i=4,5$\\
\hline
\end{tabular}
\]
We also have
\begin{gather*}
\fg_0=\fsl(3)\oplus \fsl(2)\oplus \Kee (u_1\del_1+u_5\del_5), \quad \text{where}\\
\fsl(3)=\left\{\mathop{\sum}\limits_{a,b=1,2,3}
\alpha_{ab}u_a\del_b\,|\,\mathop{\sum}\limits_{1\leq a\leq 3}
\alpha_{aa}=0\right \}, \quad \fsl(2)=\Span(u_4\del_5, u_5\del_4,
u_4\del_4+u_5\del_5).
\end{gather*}

\subsubsection[The deforms of $\fsvect(n;\un)$ for $p>3$]{The deforms of $\boldsymbol{\fsvect(n;\un)}$ for $\boldsymbol{p>3}$}\label{SvectP>3}

These deforms are described in \cite{W}.

\subsection{Partial prolongs} The Lie algebra
$\fg=\textbf{F}(\mathfrak{mb}(3|8))$ constructed above is the
\textit{complete} prolong of its negative part, see
Section~\ref{ssExDef}; let us investigate if there is
a~\textit{partial} prolong inside $\fg$. The component $\fg_1=V_1\oplus V_2$
is the direct sum of the
following $\fg_0$-invariant subspaces:
\[
V_1=\Span(u_i\del_a\,|\,i=4,5,\ a=1,2,3),\quad
V_2=\Span(u_au_b\del_i\,|\,a,b=1,2,3,\ i=4,5).
\]
The $\fg_0$-module $V_1$ is irreducible.

The $\fg_0$-module $V_2$ contains an irreducible $\fg_0$-submodule
$V_2^0= \Span(x_ax_b\del_i\,|\,a\ne b)\subset V_2$ and~$\fg_0$ acts
in the quotient space as follows: $\fsl(3)$ acts in $V_2/V_2^0$ by
zero and $\fsl(2)$ acts as $\id_{\fsl(2)}$ with multiplicity 3, so
$\dim V_2/V_2^0=8$.

Using (\ref{brack1}) it is easy to see that
\[
[[V_1,\fg_{-1}]]=\fg_0,\quad [[V_2,\fg_{-1}]]\subset \fsl(3).
\]

This means that only partial prolongs with
$\widetilde\fg_1\subset\fg_1$ containing $V_1$ can be simple.

For $\widetilde\fg_1=V_1\oplus V_2^0$, the partial prolong with the
unconstrained shearing vector which is of the form
$\un^u=(1,1,1,\infty,\infty)$ is a~deform of $\fsvect\big(5;\un^u\big)$.

For $\widetilde\fg_1=V_1\oplus V_2^0\oplus \Span\big(u_1^{(2)}\del_i\,|\,
i=4,5\big)$, the partial prolong is a~deform of $\fsvect\big(5;\un^u\big)$ with
$\un^u=(\infty,1,1,\infty,\infty)$.

For $\widetilde\fg_1=V_1\oplus V_2^0\oplus \Span\big(u_a^{(2)}\del_i\,|\,
a=1,2,\   i=4,5\big)$, the partial prolong is a~deform of
$\fsvect\big(5;\un^u\big)$ with $\un^u=(\infty,\infty,1,\infty,\infty)$.

The subspace $V_1$ is commutative and the partial prolong $\fh$ with
$V_1$ as the first component is trivial, i.e., $\fh=\fg_{-3}\oplus
\fg_{-2}\oplus \fg_{-1}\oplus \fg_{0}\oplus
\big(\widetilde\fg_{1}=V_1\big)$. Since $[[V_1,\fg_{-2}]]=0$, it follows
that $\fg_{-3}\oplus \fg_{-2}$ is an ideal in $\fh$. The simple
$24$-dimensional quotient obtained is isomorphic to $\fsl(5)$ with
the degrees of Chevalley generators being $(0,\pm1, 0, 0)$.

\textit{Conclusion.} There are no new algebras as partial prolongs.

\section[$\fm\fb\big(9;\widetilde{\uM}\big)$ and analogs of semi-densities for $p=2$]{$\fm\fb\boldsymbol{\big(9;\widetilde{\uM}\big)}$ and analogs of semi-densities for $\boldsymbol{p=2}$}\label{Smb4/5}

The Lie algebra $\fm\fb\big(9;\widetilde\uM\big)$\index{$mmathfrak\fb\big(9;\widetilde\uM\big)$@$\fm\fb\big(9;\widetilde\uM\big)$} is the desuperization of
$\fm\fb(4;\uM|5)$, the $p=2$ analog of $\fm\fb(4|5)$ over~$\Cee$. It
can be obtained from the Lie algebra $\textbf{F}(\fm\fb(3;\un|8))$,
a deform of $\fsvect\big(5;\widetilde\un\big)$, by regrading of the latter:
\[
\deg u_5=2, \quad \deg u_i=1 \quad\text{for}  \ i=1,2,3,4.
\]

Let us recall a~description of $\fm\fb(4|5)$ as the Lie superalgebra
that preserves something.

\underline{Over $\Cee$}, the Lie superalgebra $\fm\fb(4|5)$ was
initially constructed as follows. Consider the Lie superalgebra
$\fm(3;3)$: this is the regrading of $\fm(3)$ which preserves the distribution given by the Pfaff equation
with the form $d\tau+\sum q_id\xi_i$; this regrading is a~$\Zee$-grading of depth
1, see \eqref{nonstandgr}. We have (assuming $p\big(\vvol^{1/2}\big)=\od$)
\[
\fm=(\fm_{-1}, \fm_0)_*,
\]
where
\[
\fm_0=\fvect(\xi)\rtimes
\Lambda(\xi)\tau \quad \text{and} \quad  \fm_{-1}=\Lambda(\xi)\otimes\vvol^{1/2}
\stackrel{\text{as spaces}}{\simeq} \Pi(\Lambda(\xi)) .
\]
Here $\fvect(\xi)=\Span\big(\sum f_i(\xi)q_i\big)$. Denote $\fn:=\Lambda(\xi)\tau$.

Considering $\fm_{-1}$ as a~$\fvect(\xi)$-module, we preserve the
multiplication of the Grassmann algebra~$\Lambda(\xi)$; i.e., the
$\fvect(\xi)$-action satisfies the Leibniz rule, whereas the ideal
$\fn$ of $\fm_0$ does not preserve this multiplication. However,
there is an isomorphism of $\fvect(\xi)$-modules $\sigma\colon \fn\tto
\Pi(\fm_{-1})$ and the action of $\fn$ on $\fm_{-1}$ is accomplished
with the help of this isomorphism\footnote{Speaking informally,
although $\fn$ does not preserve the multiplication in $\fm_{-1}$
considered as the Grassmann algebra, $\fn$ ``remembers'' this
multiplication. And since $\fm$ is the Cartan prolong, it also
somehow ``remembers'' this structure.

The bilinear form $\omega$ with which we construct the central
extension $\fm_-=\fm_{-2}\oplus \fm_{-1}$ is the Berezin integral\index{(Berezin integral@Berezin integral}
(the coefficient of the highest term) of the product of the two
functions:
\[
\omega(g_1,g_2)=\int g_1g_2\vvol\quad \text{for any} \
g_1,g_2\in\fm_{-1};
\]
i.e., it also ``remembers'' the multiplication in $\fm_{-1}$.}:
\[
[f,g]=\sigma(f)\cdot g \quad \text{for any}  \ f\in \fn, \  g\in\fm_{-1}.
\]
The ``right'' question therefore is not ``which elements of $\fm_0$
preserve $\omega$?'', but rather
\begin{equation}\label{Q}
\text{``which elements of $\fm_0$ preserve $\omega$
\textit{conformally}, up to multiplication by a~scalar?''}
\end{equation}
It is precisely these elements which are derivations of the Lie
superalgebra $\fm_-$, and since $\fm$ is the maximal algebra that
``remembers'' the multiplication, it follows that the whole of
$\fder(\fm_-)$ lies inside~$\fm_0$.

Let us give an interpretation of the analog of the space of
semi-densities for $p=2$.

\underline{For $p\ne 2$}, we know the answer to the question
\eqref{Q}: these are elements of the two types:
\begin{enumerate}\itemsep=0pt
\item[(a)] the elements of $(\fb_{1/2}(3))_0$, the space of linear vector
fields preserving the form $\omega$, i.e., elements of the form
\[
D+\nfrac12\Div D, \quad \text{where $D\in\fvect(\xi)$;}
\]
\item[(b)] the elements of the form $c\cdot \tau\in\fn$ which multiply $\omega$ by $2c\in\Kee$.
\end{enumerate}

\underline{For $p=2$}, the elements of the form $c\cdot 1\in\fn$,
where $c\in\Kee$, multiply $\omega$ by $2c=0\in\Kee$. More\-over, any
function $f\in\fn$ preserves $\omega$ as well:
\[
f\omega(g_1,g_2)=\int ((fg_1)g_2+g_1(fg_2))\vvol=0 \quad \text{for any} \
f\in\fn, \ g_1,g_2\in\fm_{-1}.
\]
Thus, the form $\omega$ is preserved by $\fsvect(\xi)\rtimes
\Lambda(\xi)$ which is isomorphic to the subalgebra of linear
(degree 0) vector fields in $\fb_\infty$. This should have been
expected: since $2=0$, then $\frac12=\infty$.

The elements conformally preserving $\omega$ are precisely
$\xi_i\partial_i\longleftrightarrow q_i\xi_i$, so we have to add
their sum to the 0th part and calculate the Cartan prolong.

Now we are able to obtain the basis of the nonpositive components of
$\fm\fb(9;\uM)$. A~realization of the weight basis of the negative
components and generators of the 0th component by vector fields is
as follows, see Section~\ref{natGener} ($X_i^\pm$ are the Chevalley
generators of $\fsl(3)=\fsvect(0|3)_0$):
\begin{equation*}\label{w_24} \footnotesize
\renewcommand{\arraystretch}{1.4}
\begin{tabular}{|l|l|} \hline
$\fg_{i}$&the generators \\ \hline
&\\[-1.5em]
\hline
$\fg_{-2}$&$\partial_1 $\\
\hline
$\fg_{-1}\simeq $&$\xi_1\lra \partial_2,\ \xi_2\lra
\partial_3,\ \xi_3\lra \partial_4,\ \xi_1\xi_2 \xi_3\lra
\partial_9+x_8\partial_1,$
\\
$\cO(3;\One)$& $ 1\lra \partial_8,\ \xi_1\xi_2\lra \partial_5+
x_4\partial_1, \ \xi_1\xi_3\lra
\partial_6+x_3\partial_1,\ \xi_2\xi_3\lra
\partial_7+x_2\partial_1 $\\
\hline
$\fg_{0}\simeq \Kee D\rtimes $& $\del_{\xi_3}\lra x_6x_7\partial_1+
x_4\partial_8+ x_6\partial_2+ x_7\partial_3+
x_9\partial_5,$ \\
$(\fsvect(3;\One)\rtimes\cO(3;\One)$&
$X_1^-\lra x_3\partial_2+ x_7\partial_6,\ X_1^+\lra x_2\partial_3+ x_6\partial_7,$\\
&$X_2^-\lra x_4\partial_3+x_6\partial_5, \ X_2^+\lra x_3\partial_4+
x_5\partial_6,$\\
&$\xi_1\xi_2\del_{3}\lra x_4^{(2)}\partial_1+x_4\partial_5,\
\xi_1\xi_3\del_{2}\lra x_3^{(2)}\partial_1+x_3\partial_6$,\\
& $\xi_2\xi_3\del_{1}\lra x_2^{(2)}\partial_1+x_2\partial_7$, $\xi_1\xi_2 \xi_3\lra x_8^{(2)}\partial_1+x_8\partial_9,$\\
& $D= x_1\del_1+ x_2\del_2\ + x_5\del_5 + x_6\del_6+x_9\del_9$
 \\
\hline
\end{tabular}
\end{equation*}

For $\uM$ unconstrained, $\dim\fg_1=64$. The lowest-weight vectors
in $\fg_1$ are
\begin{gather*}
v_1 =  x_2x_5 x_6\partial_1+
 x_2^{(2)}\partial_8+ x_2 x_5\partial_3+x_2 x_6\partial_4 +
 x_2 x_9\partial_7 + x_5 x_6\partial_7,\\
v_2 =  x_2x_3 x_4\partial_1 + x_2 x_3\partial_5 + x_2 x_4\partial_6
+ x_3 x_4\partial_7 + x_2 x_7\partial_9 + x_2 x_8\partial_2 + x_4
x_5\partial_9+ x_4 x_8\partial_4 \\
\hphantom{v_2 =}{} + x_5 x_8\partial_5 +
 x_7 x_8\partial_7+
 x_1 x_8\partial_1+ x_1\partial_9 + x_2 x_7\partial_9 +
 x_2 x_8\partial_2 + x_6 x_8\partial_6 \\
\hphantom{v_2 =}{}+ x_8 x_9\partial_9 + x_3 x_6\partial_9 +
x_3 x_8\partial_3 + x_8^{(2)}\partial_8,\\
v_3= x_2x_5 x_7\partial_1+ x_3x_5 x_6\partial_1 + x_2 x_3\partial_8
+ x_2 x_5\partial_2 + x_2 x_7\partial_4 + x_2 x_9\partial_6 + x_3
x_5\partial_3 + x_3 x_6\partial_4 \\
\hphantom{v_3=}{}
 + x_3 x_9\partial_7 + x_5 x_6\partial_6 + x_5 x_7\partial_7,\\
v_4=x_1\del_8 + x_5x_6\del_2+ x_5 x_7\del_3+
 x_5x_9\del_5 + x_6 x_7\del_4 + x_6x_9\del_6 + x_7 x_9\del_7.
 \end{gather*}

\textit{Critical coordinates:} $\uM_5=\uM_6=\uM_7=\uM_9=1$. This Lie
algebra is a~regrading of $\fm\fb_3(11;\un)$.

\subsection[No simple partial prolongs with the whole $\fg_0$]{No simple partial prolongs with the whole $\boldsymbol{\fg_0}$} There are
remarkable elements in~$\fg_1$:
\begin{gather*}
v_1=x_2^{(2)}\del_7+x_2^{(3)}\del_1, \quad\!
v_2=x_3^{(2)}\del_6+x_3^{(3)}\del_1, \quad\!
v_3=x_4^{(2)}\del_5+x_4^{(3)}\del_1, \quad\!
v_4=x_8^{(2)}\del_9+x_8^{(3)}\del_1.
\end{gather*}
Each of the first three vectors generates a~submodule of $\dim=32$;
any two of the first three generate a~submodule of $\dim=40$; all three
together generate a~submodule of $\dim=48$. The last one generates a
submodule of $\dim=8$. All 4 together generate a~submodule of
$\widetilde\fg_1$ of $\dim=56$. The quotient $\fg_1/\widetilde\fg_1$ is an
irreducible $\fg_0$-module. We have $\dim ([\fg_{-1}, \widetilde\fg_1]) = 25$ while $\dim \fg_0 = 26$; absent is the vector of
weight~0:{\samepage
\[
x_1\del_1 + x_2\del_2 + x_5\del_5 + x_6\del_6 + x_9\del_9.
\]
Note that $[\fg_{-1},\fg_1] = \fg_0$.}

For the 24-dimensional intersection $\fg_{1i}$ of the 32-dimensional
submodules, we see that $\fg_{-1}$ is irreducible over $[\fg_{-1},\fg_{1i}]$, and hence over $[\fg_{-1},\widetilde\fg_1]$; we have
$\dim ([\fg_{-1},\fg_{1i}]) = 21$.

The elements absent in $[\fg_{-1}, \fg_{1i}]$ as compared with $\fg_0$:
\begin{gather*}
x_8^2\del_1 + x_8\del_9, \quad x_2^2\del_1 + x_2\del_7, \quad
x_3^2\del_1 + x_3\del_6, \quad
x_4^2\del_1 + x_4\del_5.
\end{gather*}

\section[A~description of $\fm\fb_2\big(11;\widetilde{\un}\big):=\textbf{F}(\fm\fb(5;\un\vert 6))$]{A~description of $\boldsymbol{\fm\fb_2\big(11;\widetilde{\un}\big):=\textbf{F}(\fm\fb(5;\un\vert 6))}$}\label{Smb11_2}

Whenever possible in this section, we do not indicate the shearing  vectors.\index{$mmathfrak\fb_2\big(11;\widetilde{\un}\big)$@$\fm\fb_2\big(11;\widetilde{\un}\big)$} The Lie superalgebra $\fm\fb(5;\un\vert 6)$ is the complete prolong of its negative part, see Section~\ref{ssExDef}.

Let us consider $\fg:=\textbf{F}(\fm\fb(5;\un|6))$ as a~deform of
$\fsvect(5)$ with the grading
\[
 \deg u_1 =\deg u_2 =1,\quad \deg u_3 =\deg u_4 =\deg u_5 =2.
\]

Let us describe the complete prolong of this negative part of this
Lie superalgebra, see Section~\ref{ssExDef}. We deduce the form of
the vector fields forming a~basis of the negative part of
$\fm\fb(5;\un|6)$ from nonzero commutation relations between
$\del_k$ and $x_i\del_a$, where $k=1,\dots,5$, $a=3,4,5$, and
$i=1,2$, cf.~\eqref{bracket}, \eqref{brack1}, considered as elements
of $\textbf{F}(\fm\fb(5;\un|6))$:
\begin{equation}\label{weightMB}
[[\del_i,  u_i\del_a]]=\del_a, \quad [[u_1\del_4,u_2\del_5]]=[[u_1\del_5, u_2\del_4]]=\del_3.
\end{equation}
For a~basis we take realization in vector fields in 5 indeterminates
$z_k$, where $k=1,\dots,5$, and 6 indeterminates $z_{ia}$, where
$a=3,4,5$ and $i=1,2$, of which $z_1$, $z_2$, $z_3$, $z_{13}$,
$z_{23}$ are even while $z_4$, $z_5$, $z_{14}$, $z_{15}$, $z_{24}$,
$z_{25}$ are odd and $\delta_i:=\del_{z_i}$:
\begin{equation*}\label{mb_2(11)} \footnotesize
\renewcommand{\arraystretch}{1.4}
\begin{tabular}{|l|l|} \hline
$\fg_{i}$&the generators (even$\,|\,$odd)\\ \hline
&\\[-1.5em]
\hline

$\fg_{-2}$&$\delta_3 \,|\,\delta_4,\ \delta_5$ \\
\hline

$\fg_{-1} $&$\delta_1, \ \delta_2 ,\ \delta_{13}+z_1\delta_3,\
\delta_{23}+z_2\delta_3\,|\,\
\delta_{14}+z_1\delta_4+z_{25}\delta_3, \ \delta_{24}+z_2\delta_4, \
\delta_{15}+z_1\delta_5+z_{24}\delta_3, \ \delta_{25}+z_2\delta_5$\\
\hline
\end{tabular}
\end{equation*}
Because the bracket \eqref{weightMB} is a~deformation that does not
preserve the grading given by the torus in $\fgl(5)$, we consider the
part of the weights that is salvaged, namely, we just exclude the
3rd coordinate of the weight; whereas the weight of $x_3$ is defined
to be equal to $(-1,-1,1,1)$.

The dimension of $\fg_0$ is the same for all $p$; it is the
expressions of the elements that differ. The raising operators in $\fg_0$ are those
of weight $(1,-1,0,0)$ or $(0,0,1,-1)$, and those with a~positive sum
of coordinates of the weight, $\dim(\fg_0^+) = 13$; we skip their
explicit description (it is commented with \% marks in the \TeX\ file available in arXiv).

The lowering operators in $\fg_0$ are those of weight $(-1,1,0,0)$ or
$(0,0,-1,1)$, and those with a~negative sum of coordinates of the
weight; $\dim(\fg_0^{-}) = 4$:
\begin{gather}
\fbox{\{-1,-1,0,1\}}\to \{z_{13}\delta_{14}+z_{23}\delta_{24}+z_{15}\delta_2+ z_{25}\delta_1+z_3\delta_4+
z_{15}z_{23}\delta_3+z_{15}z_{24}\delta_4+z_{15}z_{25}\delta_5 \},\nonumber \\[1mm]
\fbox{\{-1,-1,1,0\}}\to  \{z_{13}\delta_{15}+z_{23}\delta_{25}+z_{14}\delta_2+
z_{24}\delta_1+z_3\delta_5+z_{14}z_{23}\delta_3+ z_{14}z_{24}\delta_4+z_{14}z_{25}\delta_5 \},\nonumber \\[1mm]
\fbox{\{-1,1,0,0\}}\to  \{z_2\delta_1+z_{13}\delta_{23}+z_{14}\delta_{24}+z_{15}\delta_{25}+ z_{14} z_{15}\delta_3 \},\nonumber \\[1mm]
\fbox{\{0,0,1,-1\}}\to  \{z_{14}\delta_{15}+z_{24}\delta_{25}+z_4\delta_5 \}.\label{mbLast}
\end{gather}

\textit{Noncritical coordinates}: $N_1$, $N_2$, $N_3$.

For the unconstrained shearing vector, $\sdim\fg_1=20|20$ with the
lowest-weight vector
\begin{gather*}
v_1=x_3\delta_1 +x_6x_7\delta_7 +x_6x_8\delta_8
+x_6x_{10}\delta_{10} +x_7x_8\delta_9 +x_7x_{10}\delta_{11}\\
\hphantom{v_1=}{}
+x_8x_{10}\delta_2 +x_8x_9x_{10}\delta_4 +x_8x_{10}x_{11}\delta_5
\end{gather*}
generating the whole $\fg_0$-module $\fg_1$. All other highest- and
lowest-weight vectors together ge\-nerate a~submodule $V$ of $\fg_1$
of superdimension $16|16$ such that $[\fg_{-1}, V]$ is a
20-dimensional subalgebra of $\fg_0$. The $\fg_0$-module $\fg_1/V$
is irreducible.

\subsection{Desuperization}
Its 0th component is the same as in \eqref{mbLast} with parities forgotten.

\textit{Noncritical coordinates}: $N_1, \dots, N_5$.

\section[On analogs of $\fk\fas$ for $p=2$]{On analogs of $\boldsymbol{\fk\fas}$ for $\boldsymbol{p=2}$}\label{Skas}

In this section, whenever
possible, we do not indicate the shearing vectors. All computations in this section are performed for $p=2$;
however, for comparison, we also recall expressions obtained earlier
over $\Cee$ in \cite{Sh5, Sh14, ShP}. These expressions do not differ, usually,
from those for $p>2$.

The Lie superalgebra $\fk\fas$ over $\Cee$ was the last example
needed to complete the list of simple W-graded vectorial Lie
superalgebras, see \cite{Sh5, Sh14}. Its nonpositive part is the
same as that of $\fg:=\fk(1|6)$ (generated by the functions in the
even $t$ and 6 odd indeterminates) in its standard $\Zee$-grading
while the component $\fg_1$ is exceptional, as a~$\fg_0$-module, among
various $\fk(1|n)$: only for $n=6$ does $\fg_1$ split into 3 irreducible
components: one depends on $t$, the other two are dual to each
other. For any $p\neq 2$, we define two copies of $\fk\fas$; each of them is the partial prolong
generated by the nonpositive part, and the two submodules of
$\fg_1$: the one that depends on $t$, and any one of the other two submodules.

To distinguish between these two isomorphic copies of $\fkas$, we
denote by $\fkas^\xi$ the one whose space of generating functions
contains the product $\xi_1\xi_2\xi_3$; let $\fkas^\eta$ be the
one whose space of generating functions contains the product
$\eta_1\eta_2\eta_3$. We always consider only $\fkas^\xi$, see
\eqref{kas}, so we skip the superscript.

\underline{For $p=2$}, the structure of $\fg_1$ as a~$\fg_0$-module is rather
complicated, and it is not clear what should one take for an analog
of $\fk\fas$. Let us investigate.

\subsection[The component $\fg_1$ as $\fg_0$-module for $\fg:=\fk(1|6)$ if $p=2$]{The component $\boldsymbol{\fg_1}$ as $\boldsymbol{\fg_0}$-module for $\boldsymbol{\fg:=\fk(1|6)}$ if $\boldsymbol{p=2}$}

 Let $\fg=\fk(1|6)$ be described in terms of generating
functions of $\xi$, $\eta$, $t$, where $\xi=(\xi_1,\xi_2, \xi_3)$ and
$\eta=(\eta_1,\eta_2,\eta_3)$, with the bracket{\samepage
\begin{equation}\label{K1_6}
\{f,g\}_{\rm k.b.}=\pderf{f}{t}(1+E')(g)+(1+E')(f)\pderf{g}{t}+
\mathop{\sum}\limits_{1\leq i\leq 3}
\left(\pderf{f}{\xi_i}\pderf{g}{\eta_{i}}+
\pderf{f}{\eta_{i}}\pderf{g}{\xi_{i}}\right),
\end{equation}
where $E'=\sum \xi_i\partial_{\xi_i}$, and the standard grading
$\deg t=2$, $\deg\xi_i=\deg\eta_i=1$.}

We have $ \fg_0\simeq \fd\big(\fo_\Pi^{(1)}(4)\big)$.\index{$dmathfrak(\fg)$@$\fd(\fg)$} However, since
$\fG_0\simeq \fd(\fo_\Pi(4))$ for $\fG:=\fk(7;\un)$, we have to
investigate how can we enlarge $\fd\big(\fo_\Pi^{(1)}(4)\big)$ to get a~``correct'' version of $\textbf{F}(\fk\fas)_0$.

Consider the subalgebra
$\fh:=\fo_\Pi^{(1)}(6)\simeq\Lambda^2(\xi,\eta)=\fh_{-2}\oplus\fh_0\oplus\fh_2$
of $\fg_0=\fd(\fo_\Pi^{(1)}(6))$, where
\begin{equation}\label{h_i}
\fh_{-2}=\Lambda^2(\xi),\quad \fh_{0}=\Span(\xi_i\eta_j\,|\,i,j=1,
2,3)\simeq\fgl(3),\quad \fh_{2}=\Lambda^2(\eta).
\end{equation}

Set $\Phi:=\sum\xi_i\eta_i$.\index{((ZZPhi@$\Phi$} For $p\neq 2$ and the contact bracket
\eqref{KB}, we have $\ad_\Phi|_{\fh_i}=i\id$, hence the grading in~\eqref{h_i}. For $p=2$, the elements $\Phi$ and $t$ interchange
their roles: $\Phi$ commutes with $\fg_0=\fc\fo_\Pi^{(1)}(6)$, while
$t$ is a~grading operator on $\fg_0=\fh_0\oplus\fh_1$, where
$\fh_1=\fh_{-2}\oplus\fh_2$.

\subsection{Desuperization} Under desuperization $\fk(1;\un|6)$ turns into
$\fG:=\fk\big(7;\widetilde\un\big)$, whereas the Lie algebra $\fh$, see \eqref{h_i},
turns into $\fH:=\fo_\Pi(6)=S^2(\xi,\eta)\subset\fG_0\simeq
\fd(\fo_\Pi(4))$, where
\begin{equation*}
\fH_{-2}=S^2(\xi),\quad \fH_{0}=\Span(\xi_i\eta_j\,|\,i,j=1,
2,3)\simeq\fgl(3),\quad \fH_{2}=S^2(\eta).
\end{equation*}

The highest-weight vectors of the $\fG_0$-module $S^3(\xi,\eta)$ are
as follows (in parentheses are the dimensions of the respective
$\fG_0$-modules these vectors generate)
\[\xi_1^{(3)} (26), \quad \xi_1\xi_2^{(2)} (26), \quad \xi_1\xi_2\xi_3 (14),\quad
\xi_1(\xi_2\eta_2+\xi_3\eta_3) (6).
\]
The lowest-weight vectors and the dimensions of the $\fG_0$-modules
these vectors generate are the same with the replacement $\xi \lra
\eta$. However, since the modules generated by lowest or
highest-weight vectors do not span the whole of $\fG_1$ if $p=2$, it is more
natural to describe this component differently, as follows.

Bracketing $\xi_i^{(3)}$ (resp.~$\eta_j^{(3)}$) with $\fg_{-1}$
yields $\xi_i^{(2)}$ (resp.~$\eta_j^{(2)}$), and \textit{since each
of the $26$-dimen\-sio\-nal modules generated by any cube contains only
one cube, to have all squares in $\textbf{F}(\fk\fas)_0$, we have to
take for $\textbf{F}(\fk\fas)_1$ the module generated by all cubes.
But we cannot do this: the prolong of the module containing all cubes is equal to~$\fk(7)$.} Let us establish which cubes should be absent in the correct version of $\textbf{F}(\fk\fas)_1$ and how many versions are there.

At this stage we do not yet know what shall we eliminate in $\fG_0$ to get a~correct version of~$\textbf{F}(\fk\fas)_0$, so we consider
modules over $\fg_0$.

\underline{The $\fg_0$-submodules of $E^3(\xi, \eta)$ and $\fg_1$}. The submodule $V=\Span(\xi_i\Phi, \eta_i\Phi)_{i=1}^3\subset E^3(\xi, \eta)$ is the smallest; observe that in $E^3(\xi, \eta)$ all squares vanish. By adding any of the following 8~one-dimensional modules
\begin{equation}\label{8dim}
\text{spanned by expressions $x_1x_2x_3$, where $x_i$ is any of $\xi_i$ or $\eta_i$},
\end{equation}
we can enlarge $V$ and still have a~$\fg_0$-submodule. Together these modules span a~14-dimensional submodule $W$. The quotients $E^3(\xi, \eta)/W\simeq V^*$ and $W/V$ are irreducible $\fg_0$-modules.

Now, let us involve $t$. Set $P:=\Span(t\xi_j,
\eta_i(t+\xi_j\eta_j)\,|\,i\ne j)$. As is easy to see, $\dim V\cap
P=3$. The $\fg_0$-module generated by $t\xi_i$ is of dimension 16,
as space it is the direct sum
\[
(V +  P)\oplus  \text{a 4-dimensional subspace of the 8-dimensional space \eqref{8dim}.}
\]
The $\fg_0$-module generated by $t\xi_i$ and $t\eta_j$ is of dimension 26; as a~vector space it is the direct sum
\[
(V +  P)\oplus \text{the 8-dimensional space \eqref{8dim}.}
\]

\underline{The $\fg_0$-submodules of $S^3(\xi, \eta)$ and $\fG_1$}.
There are 6 modules of dimension 26 each, each of them is generated
by one cube $\big(\xi_i^{(3)}$ or $\eta_j^{(3)}\big)$. The intersection of
$\geq 2$ of these modules is a~20-dimensional module
$E^3(\xi,\eta)$. Unions of several of these modules form $32$-,
$38$-, $44$-, and $50$-dimensional submodules, or the whole
$S^3(\xi,\eta)$.

The module generated by $\xi_1^{(3)}$ is of dimension 16, and contains 5 elements with $\xi_1^{(2)}$, $V$ and 4 elements of the form
$\xi_1x_2x_3$, see~\eqref{8dim}. The intersection of all 6 such
modules generated by cubes is equal to $W=V\oplus$the 8-dimensional
space~\eqref{8dim}.

The dimension of the union of the modules generated by $\xi_i^{(3)}$ and
$\xi_j^{(3)}$ for $i\neq j$ is equal to 24.

The dimension of the union of the modules generated by $\xi_i^{(3)}$ and
$\eta_j^{(3)}$ for any $i$ and $j$ is equal to 26.

The dimension of the union of the modules generated by $\xi_i^{(3)}$ and
$\eta_j^{(3)}$ for all $i$ and $j$ is equal to 50, it is all $\fG_1$
except for $V^*$. \textit{Verdict:}\index{$F(\fk\fas(1;\un\vert 6))$@$\textbf{F}(\fk\fas(1;\un\vert 6))$}
\begin{gather}
\text{$\textbf{F}(\fk\fas(1;\un|6))$ is the prolong of $\fk(7)_-$
and $Y\oplus (V+ P)$, }\nonumber\\
\text{where $Y$ is the 8-dimensional module \eqref{8dim}},\nonumber\\
V:=\Span(\eta_i\Phi,  \xi_j\Phi)_{i,j=1}^3,\quad P:=\Span(t\xi_j, \eta_i(t+\xi_j\eta_j)\,|\,i\ne j).\label{kas}
\end{gather}

\subsection{Remark: useful formulas for manual computations} The
lowest-weight vectors of the $\fg_0$-module $\fg_1$ are as follows:
\begin{gather*}
\eta_1 \eta_2 \xi_3, \quad \eta_1 \eta_2 \eta_3, \quad t\eta_1
\quad \text{for $p=0$},\\
\eta_1 \eta_2 \xi_3, \quad \eta_1 \eta_2 \eta_3, \quad
\eta_1\eta_2\xi_2+\eta_1\eta_3 \xi_3=\eta_1\Phi \quad \text{for $p=2$}.
\end{gather*}
Clearly, $\Lambda^3(\xi,\eta)$ is a~$\fg_0$-submodule. Let us describe it.

Let $X_0:=\eta_1\eta_2\eta_3$. The subalgebra $\fh_2\subset \fg_0$
commutes with $X_0$, and $\fh_0$ acts on $X_0$ by scalar operators,
so $U(\fg_0)X_0=U(\fh_{-2})X_0$. Denote $V_1:=\Span(\eta_1\Phi,\eta_2\Phi,\eta_3\Phi)$.

We have
\[
\{\xi_1\xi_2,X_0\}=\xi_2\eta_2\eta_3+\xi_1\eta_1\eta_3=\eta_3(\xi_1\eta_1+\xi_2\eta_2)=\eta_3\Phi.
\]
Similar computations show that $[\fh_{-2},X_0]=V_1$.

Let us now describe $[\fh_{-2},V_1]$. Clearly, see \eqref{K1_6},
\[
\{\xi_i\xi_j, \eta_\alpha\Phi\}=\{\xi_i\xi_j, \eta_\alpha\}\cdot\Phi+
\eta_\alpha\{\xi_i\xi_j, \Phi\}.
\]
We have
\[
\{\xi_i\xi_j, \Phi\}=\begin{cases}2\xi_i\xi_j&\text{for $p\neq 2$,}\\
0&\text{for $p=2$,}\end{cases}
\]
and respectively we have
\begin{gather*}
\{\xi_1\xi_2,\eta_1\Phi\}=
\begin{cases}\xi_2(-\xi_1\eta_1+\xi_3\eta_3), \\
\xi_2\Phi, \end{cases}\quad \{\xi_1\xi_3,\eta_1\Phi\}=
\begin{cases}\xi_3(-\xi_1\eta_1+\xi_2\eta_2), \\
\xi_3\Phi, \end{cases}\\
\{\xi_2\xi_3,\eta_1\Phi\}=\begin{cases}2\eta_1\xi_2\xi_3, \\
0. \end{cases}
\end{gather*}

We set $V_2:=[\fh_{-2},V_1]=\Span(\xi_1\Phi,\xi_2\Phi,\xi_3\Phi)$ and see that $[\fh_{-2},V_2]=0$.

Therefore, $\Kee X_0\oplus V_1\oplus V_2$ is a~$\fg_0$-submodule.

We have $\{\eta_i\eta_j, \eta_\alpha\Phi\}=0$; i.e., $V=V_1\oplus
V_2$ is a~submodule to which $\Kee X_0$ us glued ``from above''.

Absolutely analogously, if $Y_0:=\xi_1\xi_2\xi_3$, then $Y_0$
generates the submodule $\Kee Y_0\oplus V_1\oplus V_2$, and $\Kee
Y_0$ is now a~submodule glued to $V$ ``from below".

Now, let $W_1:=\Span(\xi_2\xi_3\eta_1, \xi_1\xi_3\eta_2,\xi_1\xi_2\eta_3)$ and $W_2:=\Span(\xi_1\eta_2\eta_3,
\xi_2\eta_1\eta_3, \xi_3\eta_1\eta_2)$. Then,
\begin{gather*}
\xi_i\xi_j\colon \ W_1\tto 0,\quad \xi_i\eta_j\colon \ W_1\tto
\begin{cases}V_2&\text{for $i\neq j$},\\
W_1&\text{for $i=j$},
\end{cases}\quad
\eta_i\eta_j\colon \ W_1\tto V_1,
\\
\xi_i\xi_j\colon \  W_2\tto V_2,\quad \xi_i\eta_j\colon \ W_2\tto
\begin{cases}V_1&\text{for $i\neq j$},\\
W_2&\text{for $i=j$},
\end{cases} \quad
\eta_i\eta_j\colon \ W_2\tto 0.
\end{gather*}

We see that $U:=V\oplus W_1\oplus W_2$ is a~submodule and $\fh$ annihilates the quotient $U/V$, and hence it is possible to glue any element of $W_1\oplus W_2$ to $V$.

Finally, the tautological representation of $\fo_\Pi^{(1)}(6)$ is realized in the 6-dimensional quotient of~$\Lambda^3(\xi,\eta)$; it
cannot be, however, singled out as a~SUBmodule, moreover, it is glued to the whole submodule $U$, including the elements~$X_0$ and~$Y_0$.

Now, look at the elements of $\fg_1$ whose expressions contain~$t$.

The Lie algebra $\fh_0=\fgl(3)$ acts in the same way as for $p=0$ (as on
the direct sum of the tautological $\fgl(3)$-module and its dual). Whereas
\[
\{\eta_i\eta_j, t\eta_k\}=
\begin{cases}\eta_i\eta_j\eta_k=X_0&\text{if the indices $i$, $j$, $k$ are distinct},\\
0&\text{otherwise}.\end{cases}
\]
Further,
\[
\{\xi_1\xi_2, t\eta_1\}=\xi_2(t+\xi_1\eta_1), \quad
\{\xi_1\xi_3, t\eta_1\}=\xi_3(t+\xi_1\eta_1), \quad
\{\xi_2\xi_3, t\eta_1\}=\xi_2\xi_3\eta_1,
\]
and $\{\xi_i\xi_j,\xi_2(t+\xi_1\eta_1)\}=0$ for all $i$, $j$.

Thus, under
the action of $\fh$ the space $Q_1=\Span(t\eta_1, t\eta_2, t\eta_3)$ generates the space
\[
Q_2=\Span(\xi_1(t+\xi_2\eta_2),\xi_2(t+\xi_3\eta_3), \xi_3(t+\xi_1\eta_1)),
\]
as well as $V_2$, $W_1$, and $\Kee X_0$.

We can try to twist the elements $t\eta_i$ by adding something to them to enable the subalgebra~$\fh_2$ annihilate them. Such twisted elements span the space $P_1:=\Span(\eta_i(t+\xi_j\eta_j)\,|\,i\ne j)$. Under the action of
$\fg_{-2}$ we obtain from $P_1$ the spaces
$P_2:=\Span(t\xi_1,t\xi_2,t\xi_3)$, $W_2$, and $\Kee Y_0$.

\subsection[The simple ideal $\fk\fas^{(1)}(1;\un|6)$ in $\fk\fas$]{The simple ideal $\boldsymbol{\fk\fas^{(1)}(1;\un|6)}$ in $\fk\fas$}\label{ssSimpleIdeal}

Let $\fg:=\fk\fas(1;\un|6)$ considered with the standard $\Zee$-grading. In $\fg_0$, the subalgebra $\fh$, see \eqref{h_i}, is not
simple: it contains the ideal $\fo_\Pi^{(2)}(6)$ of codimension 1,
consisting of matrices $\left(\begin{smallmatrix} A&B\\C&A^{\rm t}\end{smallmatrix}\right)$
with zero-diagonal symmetric matrices $B$ and~$C$ and $A\in \fsl(6)$
whereas $\fh$ consists of the same type matrices with $A\in
\fgl(6)$. Therefore, $\fg=\fk\fas(1;\un|6)$ contains a~simple ideal $\fk\fas^{(1)}(1;\un|6)$\index{$kmathfrak\fas^{(1)}(1;\un\vert 6)$@$\fk\fas^{(1)}(1;\un\vert 6)$}
of codimension 1, its outer derivation being the outer derivation of
$\fg_0$. This derivation is present in the versions of $\fk\fas$
considered in the next three sections.

\subsection[Desuperizations of $\fk\fas(1;\un|6)$]{Desuperizations of $\boldsymbol{\fk\fas(1;\un|6)}$}\label{ssDesupKAS}{}~{}

\textit{For one of the W-gradings of $\textbf{F}(\fk\fas)$, we do not require presence of all squares in $\textbf{F}(\fk\fas)_0$, but rather require their absence}; this affects the number of parameters the shearing vector depends on.

\textit{Critical coordinates.} The shearing vector $\widetilde\un$ of the desuperization $\fk:=\fk\big(7;\widetilde\un\big)$, the ambient of the desuperized $\fk\fas(1;\un|6)$, has no critical coordinates.

\section[$\widetilde{\fkas}(7;\uM):= \textbf{F}\big(\fkas\big(4;\widetilde{\un}\vert 3\big)\big)$, where $\fkas(4\vert 3):=\fkas(1\vert 6;3\eta)$]{$\boldsymbol{\widetilde{\fkas}(7;\uM):= \textbf{F}\big(\fkas\big(4;\widetilde{\un}\vert 3\big)\big)}$, where $\boldsymbol{\fkas(4\vert 3):=\fkas(1\vert 6;3\eta)}$}\label{Skas7}

\begin{equation}\label{nonyet0}
\renewcommand{\arraystretch}{1.4}
\footnotesize
\begin{tabular}{|c|c|c|}
\hline
element of &its action & the corresponding \\
$\fg_0$ &on $\fg_{-1}\simeq\Vol_0$ &vector field $\in\fvect(4;\un|3)$\\
\hline
$\xi_1$ & $\partial_{\eta_1}$ & $z_1\del_0+z_{12}\del_2+z_{13}\del_3$\\
$\xi_2$ & $\partial_{\eta_2}$ & $z_2\del_0+z_{12}\del_1+z_{23}\del_3$\\
$\xi_3$ & \fbox{$\partial_{\eta_3}$} & $z_3\del_0+z_{13}\del_1+z_{23}\del_2$\\
\hline
\hline
\multicolumn{3}{|c|}{$t\longleftrightarrow 1, \ \xi_i\eta_i
\longleftrightarrow \eta_i\del_{\eta_i}\Longrightarrow t+
\xi_i\eta_i\longleftrightarrow \eta_i\del_{\eta_i}+
\Div(\eta_i\del_{\eta_i})$}\\
\hline
$t$ & $1$ & $\sum z_s\del_s$ for any index $s$\\
$\xi_1\eta_1$& $\eta_1\del_{\eta_1}$&$z_1\del_1+z_{12}\del_{12}+
z_{13}\del_{13}$\\
$\xi_2\eta_2$ &$\eta_2\del_{\eta_2}$&$z_2\del_2+z_{12}\del_{12}+
z_{23}\del_{23}$\\
$\xi_3\eta_3$& $\eta_3\del_{\eta_3}$&$z_3\del_3+z_{13}\del_{13}+
z_{23}\del_{23}$\\
$\xi_1\eta_2$ & \fbox{$\eta_2\del_{\eta_1}$} & $z_1\del_2+z_{13}\del_{23}$\\
$\xi_1\eta_3$ & $\eta_3\del_{\eta_1}$ & $z_1\del_3+z_{12}\del_{23}$\\
$\xi_2\eta_1$ & $\eta_1\del_{\eta_2}$ & $z_2\del_1+z_{23}\del_{13}$\\
$\xi_2\eta_3$ & \fbox{$\eta_3\del_{\eta_2}$} & $z_2\del_3+z_{12}\del_{13}$\\
$\xi_3\eta_1$ & $\eta_1\del_{\eta_3}$ & $z_3\del_1+z_{23}\del_{12}$\\
$\xi_3\eta_2$ & $\eta_2\del_{\eta_3}$ & $z_3\del_2+z_{13}\del_{12}$\\
\hline
\hline
\multicolumn{3}{|c|}{$\xi_i\eta_j\eta_k\longleftrightarrow
\eta_j\eta_k\del_{\eta_i}\in\fsvect(\eta),\ \ i\ne j\ne k $}\\
\hline
$\xi_1\eta_2\eta_3$ & $\eta_2\eta_3\del_{\eta_1}$ & $z_1\del_{23}$\\
$\xi_2\eta_1\eta_3$ & $\eta_1\eta_3\del_{\eta_2}$ & $z_2\del_{13}$\\
$\xi_3\eta_1\eta_2$ & $\eta_1\eta_2\del_{\eta_3}$ & $z_3\del_{12}$\\
\hline
\hline
\multicolumn{3}{|c|}{$\eta_i\Phi\longleftrightarrow
\eta_i(\eta_j\del_{\eta_j}+\eta_k\del_{\eta_k})
\in\fsvect(\eta),\ \ i\ne j\ne k $}\\
\hline
$\eta_1\Phi$ & $\eta_1(\eta_2\del_{\eta_2}+\eta_3\del_{\eta_3})$ &
$z_2\del_{12}+z_3\del_{13}$\\
$\eta_2\Phi$ & $\eta_2(\eta_1\del_{\eta_1}+\eta_3\del_{\eta_3})$ &
$z_1\del_{12}+z_3\del_{23}$\\
$\eta_3\Phi$ & $\eta_3(\eta_2\del_{\eta_2}+\eta_1\del_{\eta_1})$ &
$z_1\del_{13}+z_2\del_{23}$\\
\hline
\hline
\multicolumn{3}{|c|}{$\eta_i(t+\xi_j\eta_j)\longleftrightarrow
\eta_i\eta_j\del_{\eta_j}+\eta_i= \eta_i\eta_j\del_{\eta_j}+
\Div(\eta_i\eta_j\del_{\eta_j}) $}\\
\hline
$\eta_1(t+\xi_2\eta_2)$ & $\eta_1\eta_2\del_{\eta_2}+\eta_1$ &
$z_0\del_1+z_3\del_{13}$\\
$\eta_2(t+\xi_3\eta_3)$ & $\eta_2\eta_3\del_{\eta_3}+\eta_2$ &
$z_0\del_2+z_1\del_{12}$\\
$\eta_3(t+\xi_1\eta_1)$ & $\eta_1\eta_3\del_{\eta_1}+\eta_3$ &
$z_0\del_3+z_2\del_{23}$\\
\hline
\hline
\multicolumn{3}{|c|}{$\eta_i\eta_i(t+\xi_k\eta_k)\longleftrightarrow
\eta_1\eta_2\eta_3\del_{\eta_k}+\eta_i\eta_j=
\eta_1\eta_2\eta_3\del_{\eta_k}+
\Div(\eta_1\eta_2\eta_3\del_{\eta_k}) $}\\
\hline
$\eta_1\eta_2(t+\xi_3\eta_3)$ & $\eta_1\eta_2\eta_3\del_{\eta_3}+
\eta_1\eta_2$ & $z_0\del_{12}$\\
$\eta_1\eta_3(t+\xi_2\eta_2)$ & $\eta_1\eta_2\eta_3\del_{\eta_2}+
\eta_1\eta_3$ & $z_0\del_{13}$\\
$\eta_2\eta_3(t+\xi_1\eta_1)$ & $\eta_1\eta_2\eta_3\del_{\eta_1}+
\eta_2\eta_3$ & $z_0\del_{23}$\\
\hline
\end{tabular}
\end{equation}

For bases in $\fg_{-1}$ and $\fg_0$ we take the following elements:\index{$kasmathfrak(7;\uM)$@$\widetilde{\fkas}(7;\uM)$}
\begin{equation}\label{nonyet}
\fg_{-1}\colon \ \begin{array}{l} \del_0\lra 1,\quad \del_{12}\lra
\eta_1\eta_2,\quad  \del_{13}\lra \eta_1\eta_3, \quad  \del_{23}\lra
\eta_2\eta_3,\\ \del_1\lra \eta_1,\quad  \del_2\lra \eta_2,\quad  \del_3\lra
\eta_3.\end{array}
\end{equation}

Let $\fk(1;\un|6)$ be considered as preserving the distribution
given by the form $dt+\sum \xi_id\eta_i$ with the contact bracket~\eqref{K1_6} and the grading of the generating functions given by,
see Table~\eqref{regrad}:
\[
\deg t=\deg\xi_i=1, \quad \deg\eta_i=0 \quad \text{for $i=1,2,3$}, \quad \text{hence
$\deg_{\rm Lie}(f)=\deg(f)-1$}.
\]

For $\fg:=\fkas\big(4;\widetilde\un|3\big)$, we have $\fg_{-1}\simeq
\Vol_0=\Span(f(\eta)\,|\,\int f=0)$, i.e, all polynomials of $\eta$
without the product of the three of them. In \eqref{nonyet},
\eqref{nonyet0} we express these fields in terms of the 7
indeterminates $z$; we set $\del_i:=\del_{z_i}$,
$\del_{ij}:=\del_{z_{ij}}$. We have (recall the definition of
$\Vol_0$, see~\eqref{Vol_0})
\[
\fg_{0}\simeq \fc(\fvect(0|3)),
\]
with $\fvect(0|3)$ acting on $\fg_{-1}$ as on the space of volume
forms, i.e., $D\longmapsto D+\Div(D)$, and the element $t$ generating the
center of $\fg_0$ acts on $\fg$ as the grading operator. To simplify notation, we redenote the indeterminates as
follows:
\begin{gather*}
z_0\lra 1\lra x_1,\quad z_{12}\lra \eta_1\eta_2\lra x_2,\quad z_{13}\lra
\eta_1\eta_3\lra x_3,\quad z_{23}\lra \eta_2\eta_3\lra x_4,\\
z_1\lra \eta_1\lra x_5,\quad z_2\lra \eta_2\lra x_6,\quad z_3\lra \eta_3\lra
x_7. \end{gather*}

\subsection{Partial prolongs} The unconstrained shearing vector only depends on $N_1$,
we have $\sdim\fg_1=16|18$.

Let $V_i$ be the $\fg_0$-submodule in
$\fg_1$ generated by $v_i$. For the unconstrained shearing vector,
we have $\sdim V_1=8|6$, $\sdim V_2=12|9$ with $[\fg_{-1},
V_i]=\fg_0$, and $V_1\subset V_2$.

The prolong in the direction of $V_1$, see \eqref{convent}, is
trivial, namely $\fg_2^{(V_1)}=0$.

The prolong in the direction of $V_2$, see \eqref{convent}, gives
$\sdim \fg_2^{(V_2)} =4|4$, and $\sdim \fg_2^{(V_2)}=0|1$.
$\sdim([V_2, V_2]) = 3|4$, while $[V_2, \fg_2^{(V_2)}] = 0$.

There are also 3 highest-weight vectors that generate nested modules $W_1 \subset W_2 \subset W_3$; we have
\[
W_1 = V_1, \quad \sdim(W_2)=8|7  \quad \text{and} \quad \sdim(W_3)=12|12.
\]

The prolong in the direction of $W_2$ is trivial, as is the prolong in the direction of $V_1$.

The prolong in the direction of the $12|10$-dimensional module $V_2+W_2$ is equal to the prolong in the direction of $V_2$.

The prolong in the direction of $W_3$: $\sdim \fg_i^{(W_3)} =12|12$ for every $i>1$; and $\un$ depends on one parameter: $N_1$.

The prolong in the direction of the $16|15$-dimensional module $V_2+W_3$ is the same as for the whole of $\fg_1$, and $\sdim \fg_i^{(V_2+W_3)} =16|16$ for every $i>1$; and hence $\un$ depends on one parameter:~$N_1$.

The lowest-weight vectors in $\fg_1$ are
\begin{gather*}
v_1=x_1x_2\del_6 +x_1x_3\del_7 +x_1x_5\del_1 +x_2x_5\del_2 + x_3x_5\del_3 +x_4x_5\del_4 \\
\phantom{v_1=}{}+x_2x_7\del_4 +x_3x_6\del_4 +x_5x_6\del_6 +x_5x_7\del_7,\\
v_2=x_2x_3\del_3 +x_2x_4\del_4 +x_2x_5\del_5 +x_2x_6\del_6 +x_3x_6\del_7 +x_4x_5\del_7 +x_5x_6\del_1.
\end{gather*}

The highest-weight vectors in $\fg_1$ are
\begin{gather*}
w_1=x_1^{(2)}\del_5 +x_1x_6\del_2 +x_1^{(2)}\del_5 +x_1x_7\del_3,\quad
w_2=x_1x_7\del_2,\quad
w_3=x_1^{(2)}\del_2.
\end{gather*}

\subsection{Desuperization} For the unconstrained shearing vector, we have $\dim\fg_1=55$ with three lowest-weight vectors. The first two are as above, and the third one is
\[
v_3=x_2 x_3\del_4 +x_2x_5\del_6 +x_3x_5\del_7 +x_5^{(2)}\del_1.
\]
The unconstrained shearing vector is of the form $\uM=(m,1,1,1,n,s,t)$.

\subsubsection{Partial prolongs} For the unconstrained shearing vector, we have
\[
\dim V_1=14, \quad \dim V_2=21, \quad \dim V_3=31 \quad \text{and} \quad \dim(V_3+ W_3)=40
\]
 with $[\fg_{-1}, V_i]=\fg_0$ for every $i$, and $V_1\subset V_2\subset V_3$, $W_2 \subset V_3$.

The unconstrained shearing vector for the prolong in the direction of $V_3$ depends on 1 parameter $N_5$, and $\dim \fg_i^{(V_3)} =32$ for every $i>1$.

The unconstrained shearing vector for the prolong in the direction of $V_3+ W_3$ depends on 2 parameters $N_1$, $N_5$, and
\[
\dim \fg_2^{(V_3+ W_3)} =56, \quad \dim \fg_3^{(V_3+ W_3)} =72, \quad \dim \fg_4^{(V_3+ W_3)} =88.
\]

\section[$\fkas(8;\uM):=\textbf{F}(\fkas(4;\un\vert 4))$]{$\boldsymbol{\fkas(8;\uM):=\textbf{F}(\fkas(4;\un\vert 4))}$}\label{Skas8}
Let $\fk(1;\un|6)$ \index{$kasmathfrak(8;\uM)$@$\fkas(8;\uM)$} be considered as preserving the distribution
given by the form $dt+\sum \xi_id\eta_i$ with the contact bracket
\eqref{K1_6} and the grading of the generating functions given by,
see Table \eqref{regrad}:
\[
\deg t=\deg\eta_i=1, \quad \deg\xi_i=0 \quad \text{for $i=1,2,3$}, \quad \text{hence
$\deg_{\rm Lie}(f)=\deg(f)-1$}.
\]

For the subalgebra $\fg=\fkas\big(4;\widetilde\un|4\big)$ of $\fk(1;\un|6;3\xi)$, see \eqref{kas}, we have
\[
\fg_0\simeq\begin{cases}
\fsl(1|3)\rtimes \cO(0|3), \quad \text{where $\fsl(1|3)\subset\fvect(0|3)$},
&\text{for $p\neq 2$},\\
\fd(\fsvect^{(1)}(0|3))\rtimes \cO(0|3), \quad \text{see \eqref{nonpos22}},
&\text{for $p=2$},
\end{cases}
\]
with the natural action of $\fsl(1|3)$ if $p\neq 2$, or $\fd(\fsvect^{(1)}(0|3))$ if $p=2$, on $\cO(0|3)$.

Indeed, the element $(t+\Phi)f(\xi)\in \fg_0$ acts on $\fg_{-1}$ as the operator of multiplication by $f(\xi)$. Additionally $\fg_0$
contains the following operators:
\begin{equation}\label{kas3xi}
\eta_i\lra
\del_{\xi_i}, \quad \xi_i\eta_j\lra \xi_i \del_{\xi_j},\quad \xi_i\Phi\lra \xi_i (\xi_j\del_{\xi_j}+\xi_k\del_{\xi_k}).
\end{equation}

\underline{For $p\neq 2$}, the last 3 elements in \eqref{kas3xi} do not belong
to $\fsvect(0|3)$ and the elements \eqref{kas3xi} generate $\fsl(1|3)$.

\underline{For $p= 2$}, the last 3 elements in \eqref{kas3xi} do belong to
$\fsvect(0|3)$, whereas the elements $\xi_i\eta_i$ do not; it is the
sum of any two of them that belongs to $\fsvect(0|3)$. So the
elements \eqref{kas3xi} generate a~subalgebra
$\fd(\fsvect^{(1)}(0|3))$ in the Lie algebra $\fder(\fsvect^{(1)}(0|3))$
of all derivations of $\fsvect^{(1)}(0|3)$.

In \eqref{nonpos22} we introduce 8 indeterminates $x$ needed to
express $\fkas(4;\un|4)$ and its desuperization in terms of vector
fields as the prolong; $\del_i:=\del_{x_i}$. For a~basis of the
nonpositive part (the $X_i^\pm$ are the Chevalley generators of what
is $\fsl(1|3)$ for $p\neq 2$ and turns into $\fsvect^{(1)}(0|3)$ for
$p=2$) we take:
\begin{equation}\label{nonpos22}
\footnotesize
\renewcommand{\arraystretch}{1.4}
\begin{tabular}{|l|l|}
\hline
$\fg_{-1}$ &even: $1\lra \del_1,\ \xi_1\xi_2\lra \del_2,\
\xi_1\xi_3\lra\del_3,\ \xi_2\xi_3\lra \del_4$\\

&odd: $\xi_1\lra \del_5, \ \xi_2\lra \del_6,\
\xi_3\lra \del_7, \ \xi_1 \xi_2 \xi_3\lra \del_8$\\
\hline

& $t+\Phi\lra \sum x_i\del_i,\
(t+\Phi)\xi_1\lra x_1\del_5 +x_4\del_8+x_6\del_2 +x_7\del_3,$\\

& $(t+\Phi)\xi_2\lra x_1\del_6 +x_3\del_8+x_5\del_2 +x_7\del_4,$\\

& $(t+\Phi)\xi_3\lra x_1\del_7 +x_2\del_8+x_5\del_3 +x_6\del_4,$\\

& $(t+\Phi)\xi_1\xi_2\lra x_1\del_2 +x_7\del_8, \
(t+\Phi)\xi_1\xi_3\lra x_1\del_3 +x_6\del_8,$\\

$\fg_{0}\simeq\cO(0|3)$ & $(t+\Phi)\xi_2\xi_3\lra x_1\del_4 +
 x_5\del_8,\ \ \ t\xi_1\xi_2\xi_3\lra x_1\del_8$\\

\cline{2-2} $\ltimes\fd(\fsvect^{(1)}(0|3))$&$\eta_1\lra
x_2\del_6 +x_3\del_7 +x_5\del_1 +x_8\del_4,\ \eta_2\lra x_2\del_5
+x_4\del_7 +x_6\del_1+x_8\del_3,$\\

& $X_1^-=\eta_3\lra x_3\del_5+x_4\del_6 +x_7\del_1 +x_8\del_2,$\\

&$X_1^+=\xi_1\Phi\lra x_6\del_2 +x_7\del_3, \ \xi_2\Phi\lra
x_5\del_2 +
x_7\del_4,\ \xi_3\Phi\lra x_5\del_3 + x_6\del_4,$\\

& $X_2^-=\eta_1\xi_2\lra x_3\del_4 +x_5\del_6, \ \eta_1\xi_3\lra
x_2\del_4 +x_5\del_7$,\\
& $ X_2^+=\eta_2\xi_1\lra x_4\del_3 +x_6\del_5$\\

&$X_3^-=\eta_2\xi_3\lra x_2\del_3 +x_6\del_7, \ \eta_3\xi_1\lra
x_4\del_2
+x_7\del_5$, \\
& $X_3^+=\eta_3\xi_2\lra x_3\del_2 +x_7\del_6$\\

&$\eta_1\xi_1\lra x_2\del_2+x_3\del_3 +x_5\del_5 +x_8\del_8, \
\eta_2\xi_2\lra x_2\del_2 +x_4\del_4 +x_6\del_6 +x_8\del_8$\\

&$\eta_3\xi_3\lra x_3\del_3 +x_4\del_4 +x_7\del_7 +x_8\del_8$\\

\hline
\end{tabular}
\end{equation}

For the unconstrained shearing vector, $\sdim \fg_1=16|16$. The
$\fg_0$-module $\fg_1$ splits into the 2~irreducible submodules:
$\fg_1= V_1\oplus V_2$, where $\sdim(V_1)=12|12$, and $\sdim(V_2)=4|4$.
There are the 2 highest-weight vectors in $\fg_1$:
\[
\begin{array}{l}
h_1=x_1 x_6\del_2+x_1 x_7\del_3+x_6x_7\del_8,\quad
h_2=x_1^{(2)}\del_8;
\end{array}
\]
and the 2 lowest-weight vectors:
\begin{gather*}
v_1=x_2x_3\del_3 +x_2x_4\del_4 +x_2x_5\del_5 +x_2x_6\del_6 +x_3x_6\del_7 +
x_4x_5\del_7 +x_5x_6\del_1\\
\hphantom{v_1=}{}
+x_5x_8\del_3 +x_6x_8\del_4,\\
v_2=x_1^{(2)}\del_1 +x_1x_2\del_2 +x_1x_7\del_7 +x_1x_8\del_8 +x_2x_7\del_8
+x_1^{(2)}\del_1 +x_1x_3\del_3 +x_1x_6\del_6\\
\phantom{v_2=}{} +x_1x_8\del_8 +x_3x_6\del_8 + x_1^{(2)}\del_1 +x_1x_4\del_4 +x_1x_5\del_5 +x_1x_8\del_8 +
x_4x_5\del_8 +x_5x_6\del_2 \\
\phantom{v_2=}{}+x_5x_7\del_3 +x_6x_7\del_4.
\end{gather*}

\subsection{Partial prolongs}\label{ss59}
We have $\sdim ([\fg_{-1},\fg_1]) = 12|10$, as it should be (having in mind
the outer derivation of $\fkas$);
for its representative, we can take $x_1\partial_1+
x_2\partial_2+ x_3\partial_3+ x_4\partial_4$.

In the direction of $V_1$, we have $\sdim ([\fg_{-1},V_1]) = 12|9$
(apart from the outer derivative, $x_1\partial_8$ is absent); the
$[\fg_{-1},V_1]$-module $\fg_{-1}$ is irreducible,
$\sdim(\fg_2^{(V_1)})=4|7$, and $\fg_3^{(V_1)}=0$. More precisely,
$\sdim ([V_1,V_1]) = 3|3$ and the $[\fg_{-1},V_1]$-modules $V_1$
and $[V_1,V_1]$ are irreducible. Thus, the superdimension of this
simple prolong is $31|28$.

In the direction of $V_2$ we have $\sdim ([\fg_{-1},V_2]) = 4|4$,
the $[\fg_{-1}, V_1]$-module $\fg_{-1}$ is not irreducible, so no
new simple partial prolongs exist in this direction.

\textit{Critical coordinate}: only $N_1$.

\subsection{Desuperization} The same as above, with dimension $a+b$ instead of superdimension~$a|b$.

\section[The Lie superalgebra $\fkas\big(5;\widetilde{\un}\vert 5\big) \subset\fk(1;\un\vert 6;1\xi)$]{The Lie superalgebra $\boldsymbol{\fkas\big(5;\widetilde{\un}\vert 5\big) \subset\fk(1;\un\vert 6;1\xi)}$}\label{Skas55}

Whenever possible we do not indicate the shearing vector.

Let $\fk(1;\un|6)$ be the Lie superalgebra which preserves the distribution given by the form
$dt+\sum \xi_id\eta_i$. Then, $\fk(1;\un|6)$ is endowed with the contact bracket \eqref{K1_6}; set $\deg K_f=\deg(f)-2$, where the grading of the generating functions is given by
\begin{equation*}
\deg t=\deg\eta_1=2, \quad \deg\xi_1=0, \quad
\deg\eta_i=\deg\xi_i=1\quad \text{for $i=2,3$.}
\end{equation*}

We identify $\fg_{-1}$ with $V(\Lambda)\otimes W$, where
$\Lambda:=\Lambda[\xi_1]$ and $V(\Lambda):=V\otimes \Lambda$; let
$V=\Span(v_1, v_2)$, $W=\Span(w_1, w_2)$. For a~basis of the
nonpositive part of $\fg$, we take the elements listed in~\eqref{nonposKas55}.

The component
\begin{equation}\label{FunnyG0}
\fg_{0}\cong \begin{cases}\fd\big(\big(\widetilde{\fsl}(W)\oplus (\fgl(V;\Lambda)\ltimes
\fvect(0|1))/\Kee Z\big)\big), & \\
\quad \text{where $\fd=\Kee D$, see
\eqref{nonposKas55}}&\text{if $p= 2$},\\
\widetilde{\fsl}(W)\oplus (\fgl(V;\Lambda)\ltimes
\fvect(0|1))&\text{if $p\neq 2$}
\end{cases}
\end{equation}
of the subalgebra $\fg:=\fkas\big(5;\widetilde\un|5\big)\subset
\fk(1;\un|6;1\xi)$, see \eqref{kas}, is rather complicated for $p=2$. To
describe this component, we compare it with the complete prolong of the negative
part, see Section~\ref{ssExDef}. The 0th component of this prolong
is equal to the 0th component of $\fk(1;\un|6; 1\xi)$. Its 3
elements that do not belong to $\fg_0$ are easy to find from the
description of $\fkas$ given in Section~\ref{Skas} (they are boxed):
\begin{equation}\label{nonposKas55}
\footnotesize
\renewcommand{\arraystretch}{1.4}
\begin{tabular}{|c|l|}
\hline $\fg_{i}$&the basis elements \\ \hline
&\\[-1.5em]
\hline
 $\fg_{-2}\simeq\Lambda$& $1\lra\del_2\,|\,\xi_1\lra \del_1$\\
\hline
 $\fg_{-1}\simeq\id_{\fsl(W)}$&$\xi_1\xi_2\lra \xi_1v_1\otimes w_1\lra \del_3,\
 \xi_1\xi_3\lra \xi_1v_1\otimes w_2\lra \del_4,$ \\
$\otimes\id_{\fgl(V;\Lambda)}$& $\xi_1\eta_2\lra
\xi_1v_2\otimes w_2\lra \del_5, \
\xi_1\eta_3\lra \xi_1v_2\otimes w_1\lra \del_6\,|$\\
&$\xi_2\lra v_1\otimes w_1\lra
x_5\del_1+\del_7, \ \xi_3\lra v_1\otimes w_2\lra
x_6\del_1+\del_8,$\\
& $\eta_2\lra v_2\otimes w_2\lra x_3\del_1+x_7\del_2+\del_9,\ \eta_3\lra v_2\otimes w_1\lra
x_4\del_1+x_8\del_2+
\del_{10}$\\
\hline
& $\xi_1E\otimes\One\lra \xi_1\Phi\lra
x_7\del_3+x_9\del_5+x_8\del_4
 +
x_{10}\del_6+x_7x_9\del_1 +x_8 x_{10}\del_1 $\\
&$\left(\begin{smallmatrix}0&0\\
0&\xi_1\end{smallmatrix}\right)\otimes\One\lra t\xi_1\lra x_7\del_3+x_8\del_4+x_{2}\del_1$\\

&$\xi_1X^-\otimes\One\lra \xi_1\eta_2\eta_3\lra x_7\del_6
+x_8\del_5+x_7x_8\del_1$\\

&$\xi_1X^+\otimes\One\lra
\xi_1\xi_2\xi_3\lra x_9\del_4 +x_{10}\del_3
+x_9x_{10}\del_1$\\

&$X^+\lra \xi_2\xi_3\lra
x_5\del_4 +x_6\del_3 +x_9\del_8
+x_{10}\del_7 +x_9x_{10}\del_2$\\

&$\widetilde X^+\lra \xi_3\eta_2\lra x_3\del_4 +x_6\del_5
+x_7\del_8
+x_{10}\del_9$\\

$\fg_{0}$&\fbox{$X^-$}$\lra \eta_2\eta_3\lra x_3\del_6
+x_4\del_5
+x_7\del_{10} +x_8\del_9 +x_7x_8\del_2$\\

&\fbox{$\widetilde X^-$}$\lra \xi_2\eta_3\lra x_4\del_3
+x_5\del_6 +x_8\del_7
+x_9\del_{10}$\\

&$D:=\left(\begin{smallmatrix}1&0\\
0&0\end{smallmatrix}\right)\otimes\One+\One\otimes\left(\begin{smallmatrix}0&0\\
0&1\end{smallmatrix}\right)\lra\xi_2\eta_2\lra x_3\del_3+x_5\del_5+x_7\del_7+x_9\del_9$\\

&$E\lra \xi_2\eta_2+\xi_3\eta_3\lra x_3\del_3+x_4\del_4+ x_5\del_5+x_6\del_6+x_7\del_7$\\
&$\quad{} +x_8\del_8+x_9\del_9+x_{10}\del_{10}$\\

&$\xi_1\del_{\xi_1}\lra \xi_1\eta_1\lra x_4\del_4+x_6\del_6+x_8\del_8 +x_{10}\del_{10}$\\

&$\left(\begin{smallmatrix}0&0\\
0&1\end{smallmatrix}\right)\otimes\One\lra t+\xi_1\eta_1\lra x_2\del_2 +x_3\del_3 +x_4\del_4
+x_9\del_9 +x_{10}\del_{10}$\\

&\fbox{$\del_{\xi_1}$}$\lra \eta_1\lra x_3 x_5\del_1 +x_4x_6\del_1 +x_5
x_7\del_2
+x_6x_8\del_2$\\
&$\quad{}+x_1\del_2 +x_3\del_7 +x_4\del_8
+x_5\del_9 +x_6\del_{10}$\\
\hline
\end{tabular}
\end{equation}

The component $\fg_0$ contains two copies of $\fsl(2)$; to
distinguish them, we endow one of them with a~tilde:
$\widetilde{\fsl}(2)=\fsl(W)$ generated by $\widetilde X^+$ and
$\widetilde X^-$, the other copy being $\fsl(V)$ generated by $X^+$
and $X^-$. These two copies of $\fsl(2)$ are ``glued"; their glued
sum has a~common center spanned by $E$; i.e., their direct sum is
factorized by a~1-dimensional subalgebra $\Kee Z$ in their
2-dimensional center, the explicit form of $Z$ is
inessential for us at the moment. Observe that $D,
\xi_1\del_{\xi_1}\notin [\fg_{-1}, \fg_1]$; only their sum $D+
\xi_1\del_{\xi_1}\lra \xi_1\eta_1+\xi_2\eta_2$ belongs to the
commutant.

In \eqref{nonposKas55}, we expressed the nonpositive part of $\fg$
by means of vector fields in 10 indeterminates~$x$ setting
$\del_i:=\del_{x_i}$.

The reader wishing to verify our computations will, of course, use
the contact bracket and generating functions to compute inside
$\fg_0$. The realization by vector fields is only needed to compute
$\fg_i$ for $i>0$ (with computer's aid to speed up the process).

\textit{The only noncritical coordinate of the shearing vector
$\un$} is $N_2$; it corresponds to what used to be $t$.

For the unconstrained shearing vector, we have $\sdim\fg_1=8|8$. The
only lowest-weight vector (w.r.t.\ the boxed operators) of $\fg_1$
that generates $\fg_1$ as a~$\fg_0$-module is
\begin{gather*}
u_1=x_1\del_9 +x_3x_4\del_4 +x_3x_6\del_6 +x_3x_7\del_7 +x_3x_9\del_9
+x_4x_6\del_5 +x_4x_7\del_8\\
\hphantom{u_1=}{}+x_4x_{10}\del_9 +x_6x_7\del_{10}
+x_6x_8\del_9 +x_1x_3\del_1 +x_1x_7\del_2 +x_6x_7x_8\del_2.
\end{gather*}
The other lowest-weight vector and the only highest-weight vector (together
and separately) generate a~submodule $V$ which, together with $\fg_{-1}$, generate
an 8-dimensional part of $\fg_0$. The quotient $\fg_1/V$ is an irreducible
$\fg_0$-module.

\subsection[Desuperization of $\fkas\big(5;\widetilde\un|5\big)$]{Desuperization of $\boldsymbol{\fkas\big(5;\widetilde\un|5\big)}$} The only critical coordinates are $N_1$ and $N_2$. (For the unconstrained shearing vector, $\dim \fg_1 =16$, $\dim(\fg_2)=20$, $\dim(\fg_3)=24$, $\dim(\fg_4)=28$.)

\section[A~description of $\widetilde{\fsb}\big(2^{n}-1;\widetilde{\un}\big)$ for $p=2$]{A~description of $\boldsymbol{\widetilde{\fsb}\big(2^{n}-1;\widetilde{\un}\big)}$ for $\boldsymbol{p=2}$}\label{StildeSB}
\index{$sbmathfrak\big(2^{n}-1;\widetilde{\un}\big)$@$\widetilde{\fsb}\big(2^{n}-1;\widetilde{\un}\big)$}

\subsection[Recapitulation: $p=0$, $n$ even]{Recapitulation: $\boldsymbol{p=0}$, $\boldsymbol{n}$ even} Let $q=(q_1,\dots, q_n)$ and
$\xi=(\xi_1,\dots, \xi_n)$. We consider the subsuperspace of functions $\Cee[q,\xi]$ of the form\index{((ZZXi:=\xi_1\cdots\xi_n$@$\Xi:=\xi_1\cdots\xi_n$}
\[
\{(1+\Xi)f(q,\xi)\,|\,\Delta(f)=0\quad \text{and} \quad \int_\xi
f\vvol_\xi=0\}, \quad \text{where $\Xi:=\xi_1 \cdots \xi_n$},
\]
with the Buttin bracket. In this section we use only this bracket and omit the index ``B.b''.

Let us compute the bracket in $(1+\Xi)\fsb^{(1)}(n;n)$ realized by
elements of $\fsb^{(1)}(n)$. We have
\begin{gather*}
\{(1+\Xi)f,(1+\Xi)g\}\\
\quad{}=
\begin{cases} \{f, g\}=(1+\Xi)\{f, g\}&\text{if
$\deg_\xi (f),
\deg_\xi (g)>0$};\\
\{(1+\Xi)f, g\}=(1+\Xi)\{f, g\}&\text{if $\deg_\xi (f)=0$,
$\deg_\xi (g)>1$};\\
\{(1+\Xi)f, g\}=&\\
(1+\Xi)\{f, g\}+\sum\partial_{\xi_i}\Xi
f\partial_{q_i}(g_i)\xi_i&\text{if $\deg_\xi (f)=0$, $\deg_\xi
(g)=1$,}\\
\stackrel{\text{(since $\sum\partial_{\xi_i}\Xi
f\partial_{q_i}(g_i)\xi_i=\Xi f\sum\partial_{q_i}(g_i)=0$)}}{=}
&\\
(1+\Xi)\{f, g\}&\text{if $g=\sum g_i(q)\xi_i$ and
$\Delta(g)=0$};\\
\sum\partial_{\xi_i}\Xi (\partial_{q_i}(f)g-\partial_{q_i}(g)f)&
\text{if $\deg_\xi (f)=\deg_\xi (g)=0$}.
\end{cases}
\end{gather*}

In the $\Zee$-grading of $\fg=\widetilde{\fsb}(n;n)$ by degrees
of the $q$ shifted by $-1$, we have:
\begin{itemize}\itemsep=0pt
\item $\fg_{-1}$ is spanned by monomials in $\xi$ of degrees 1 through
$n-1$, and by $1+\Xi$;
\item $\fg_{0}$ is spanned by functions of the form $g=(1+\Xi)\sum
g_i(\xi)q_i$, where $\sum \partial_{\xi_i}g_i=0$.
\end{itemize}

The $\fg_0$-action on $\fg_{-1}$ is as follows. If
$\deg_\xi(g_i)>0$, then we can ignore $\Xi$ in the factor $1+\Xi$
since $\Xi$ annihilates $g_i$, and hence $\ad_g$ acts on $\fg_{-1}$
as the vector field $\sum g_i\partial_{\xi_i}$ acts on the space of functions in $\xi$.

If $g=(1+\Xi)q_i$, then the $\ad_g$-acts on $\fg_{-1}$
precisely as an element of $\widetilde{\fsvect}(0|n)$ acts on the space $\Vol_\xi$:
\begin{gather*}
\{g, \xi_j\}=(1+\Xi)\delta_{ij},\quad
\{g, 1+\Xi\}=\partial_{\xi_i}(\Xi),\quad
\ad_g(f)=\partial_{\xi_i}(f) \quad \text{for monomials} \ f=f(\xi).
\end{gather*}
Since $(1-\Xi)\vvol$ is the invariant subspace in $\Vol_\xi$, it follows that, in the quotient space, we can take for a~basis elements of the form $f(\xi)\vvol_\xi$, where monomials $f$ differ from $1$ and $\Xi$, and either $1$ or $\Xi$. For reasons unknown, \textit{SuperLie} selected $\Xi$, not~1.

\subsubsection[Recapitulation: $p=0$, $n$ odd]{Recapitulation: $\boldsymbol{p=0}$, $\boldsymbol{n}$ odd} Everything is as above but
with $\tau\Xi$, where $\tau$ an odd parameter, instead of $\Xi$.

\subsection[$\textbf{F}\big(\widetilde{\fsb}\big(2^{n-1};\un|2^{n-1}-1\big)\big)$ for $n$ odd, $p=2$]{$\boldsymbol{\textbf{F}\big(\widetilde{\fsb}\big(2^{n-1};\un|2^{n-1}-1\big)\big)}$ for $\boldsymbol{n}$ odd, $\boldsymbol{p=2}$}\label{SSnOdd}

For $p=2$, it is possible to desuperize deforms with odd parameters and consider them in the category of
superspaces, see \cite{BLLS2}. We assume that $p(\vvol_\xi)\equiv n\mod 2$.

\subsubsection[Example: $n=3$]{Example: $\boldsymbol{n=3}$} For a~basis in $\fg_{-1}$, where $\del_i:=\del_{x_i}$, we take:
\begin{gather*}
\partial_1=\xi_1\vvol_\xi,\quad \partial_2=\xi_2\vvol_\xi,\quad \partial_3=\xi_3\vvol_\xi,\quad \partial_4=\xi_1\xi_2\vvol_\xi,\quad
\partial_5=\xi_1\xi_3\vvol_\xi,\nonumber\\
\partial_6=\xi_2\xi_3\vvol_\xi,\quad \partial_7=\tau\xi_1\xi_2\xi_3\vvol_\xi.
\end{gather*}

For a~basis of $\fg_0$, where $\delta_i:=\del_{\xi_i}$ we take the following elements, where
the $\fg_0$-action on $\fg_{-1}$ is given by realizations on the right of the $\lra$:
\begin{equation*}
\begin{split}
&(1+\tau\xi_1\xi_2\xi_3)\delta_1\lra x_1\del_7+x_4\del_2+x_5\del_3+x_7\del_6, \\
&(1+\tau\xi_1\xi_2\xi_3)\delta_2\lra x_2\del_7+x_4\del_1+x_6\del_3+x_7\del_5, \\
&(1+\tau\xi_1\xi_2\xi_3)\delta_3\lra x_3\del_7+x_5\del_1+x_6\del_2+x_7\del_4, \\
&\xi_1\delta_2\lra x_2\del_1+x_6\del_5, \\
&\xi_1\delta_3\lra x_3\del_1+x_6\del_4, \\
&\xi_2\delta_1\lra x_1\del_2+x_5\del_6, \\
&\xi_1\delta_1+\xi_2\delta_2\lra x_1\del_1+x_2\del_2+x_5\del_5+x_6\del_6, \\
&\xi_1\delta_1+\xi_3\delta_3\lra x_1\del_1+x_3\del_3+x_4\del_4+x_6\del_6,
\end{split} \quad
\begin{split}
&\xi_2\delta_3\lra x_3\del_2+x_5\del_4, \\
&\xi_3\delta_1\lra x_1\del_3+x_4\del_6, \\
&\xi_3\delta_2\lra x_2\del_3+x_4\del_5, \\
&\xi_1\xi_2\delta_3\lra x_3\del_4, \\
&\xi_1\xi_3\delta_2\lra x_2\del_5, \\
&\xi_1\xi_2\delta_2+\xi_1\xi_3\delta_3\lra x_2\del_4+x_3\del_5, \\
&\xi_2\xi_3\delta_1\lra x_1\del_6, \\
&\xi_1\xi_3\delta_1+\xi_2\xi_3\delta_2\lra x_1\del_5+x_2\del_6, \\
&\xi_1\xi_2\delta_1+\xi_2\xi_3\delta_3\lra x_1\del_4+x_3\del_6.
\end{split}
\end{equation*}

The weights are considered with respect to
$\fsl(3)\subset\textbf{F}\big(\widetilde\fsvect(0|3)\big)$, i.e.,
\[
w(\xi_1)=(1,0), \quad w(\xi_2)=(-1,1), \quad  w(\xi_3)=(0,-1).
\] The raising
elements are those for which either $w_1+w_2>0$ or $w_1=-w_2>0$;
the lowering elements are those for which either $w_1+w_2<0$ or
$w_1=-w_2<0$. (To find lowering and raising operators, we could have
considered a~$\Zee$-grading of $\widetilde\fsvect(0|n)$ by setting $\deg\xi_n=-n+1$
 and $\deg\xi_1=\dots=\deg\xi_{n-1}=1$ with ensuing
natural division into ``positive'' and ``negative'' parts.)

The highest-weight vectors of the $\fg_0$-module $\fg_1$ are
\begin{gather*}
w_1=x_2x_3\del_7 +x_2x_5\del_1 +x_2x_6\del_2 +x_2x_7\del_4 +x_3x_4\del_1 +x_3 x_6\del_3\nonumber\\
\hphantom{w_1=}{} +x_3x_7\del_5 +x_4x_6\del_4 +x_5x_6\del_5,\nonumber\\
w_2=x_2x_3\del_1 +x_2x_6\del_4 +x_3x_6\del_5,\nonumber\\
w_3=x_3^{(2)} \del_1 + x_3x_6\del_4.
\end{gather*}
The lowest-weight vectors of the $\fg_0$-module $\fg_1$ are
\begin{gather*}
v_1=x_1^{(2)}\del_1 +x_1x_2\del_2 +x_1 x_5\del_5 +x_1x_6\del_6
+ x_2x_5\del_6 +x_1^{(2)}\del_1 +x_1x_3\del_3\nonumber\\
\hphantom{v_1=}{} +x_1x_4\del_4 +x_1x_6\del_6 +x_3x_4\del_6,\nonumber\\
v_2=x_1^{(2)}\del_1 +x_1x_2\del_2 +x_1 x_5\del_5 +x_1x_6\del_6 +
x_2x_5\del_6 +x_1^{(2)}\del_1 +x_1x_3\del_3\nonumber\\
\hphantom{v_2=}{} +x_1x_4\del_4 +x_1x_6\del_6 +x_3x_4\del_6,\nonumber\\
v_3=x_1^{(2)} \del_6,
\end{gather*}

\textit{Partial prolongs}: The elements of $\fg_0$ absent in $\widetilde\fg_0:=[\fg_1,\fg_{-1}]$ are $\xi_1\xi_2\delta_3$, $\xi_1\xi_3\delta_2$, $\xi_2\xi_3\delta_1$. The $\widetilde\fg_0$-module $\fg_{-1}$ is irreducible.

Let $V_i$ and $W_i$ denote the $\widetilde\fg_0$-modules generated by $v_i$ and $w_i$, respectively.
We have
\begin{gather*}
\dim \fg_1 = 31,\quad \dim\fg_2 = 49, \quad \dim\fg_3 = 71,\nonumber\\
\dim V_1 = \dim W_1 =7,\quad \dim V_2 = \dim W_2 =8,\quad \dim V_3 = \dim W_3 =16,\nonumber\\
V_1 = W_1,\quad V_1\subset V_2 \subset V_3,\quad W_1\subset W_2 \subset W_3,\nonumber\\
\dim(V_2+W_2) = 9,\quad \dim(V_2+W_3) = \dim(V_3+W_2) = 17,\quad \dim(V_3+W_3) = 24.
\end{gather*}

The brackets with $\fg_{-1}$:
\begin{gather*}
\dim ([\fg_{-1}, V_1]) = \dim [\fg_{-1}, V_2+W_2] = 14,\nonumber\\
\dim ([\fg_{-1}, V_3]) = 15 \quad \text{(absent are $\xi_1\xi_2\delta_3$, $\xi_1\xi_3\delta_2$)},\nonumber\\
\dim ([\fg_{-1}, W_3]) = 15 \quad \text{(absent are $\xi_1\xi_3\delta_2$, $\xi_2\xi_3\delta_1$)},\nonumber\\
\dim ([\fg_{-1}, V_3+ W_3]) = 16 \quad \text{(absent is $\xi_1\xi_3\delta_2$)}.
\end{gather*}

Therefore (recall the convention \eqref{convent})
\begin{equation}\label{dimVin=31}
\footnotesize
\begin{tabular}{|l|l|}
\hline
Partial prolongs in the direction of&dimensions\\
\hline
$V_1 \text{~~or~~}V_2 \text{~~or~~} V_2+W_2$&$\dim \fg_2 = 1,\ \fg_3 = 0$\\
$V_3 \text{~~or~~}V_3+W_2$&$\dim\fg_2 = \dim \fg_3 = 16$\\
$V_3+W_3$&$\dim\fg_2 = 32, \ \ \dim \fg_3 = 40$\\
$[\fg_{-1}, V_1] \text{~~or~~} [\fg_{-1}, V_2+W_2]$&$\dim \fg_1 = 10
\text{~~absent are $v_3$ and $w_3$}, \ \dim \fg_2 = 1, \ \fg_3 = 0$.\\
\hline
\end{tabular}
\end{equation}

\textit{Critical coordinates of $\widetilde{\fsb}\big(7;\widetilde\un\big)$} are $N_4$, $N_5$, $N_6$, and $N_7$, as follows from \eqref{dimVin=31}.

\subsection[$\textbf{F}\big(\widetilde{\fsb}\big(2^{n-1}-1;\un|2^{n-1}\big)\big)$ for $n$ even, $p=2$]{$\boldsymbol{\textbf{F}(\widetilde{\fsb}(2^{n-1}-1;\un|2^{n-1}))}$ for $\boldsymbol{n}$ even, $\boldsymbol{p=2}$}
For the unconstrained shearing
vector~$\un^u$, the dimensions of homogeneous
components of $\fg=\widetilde{\fsb}(2^n-1;\un^u)$ are the same as
those of $\fsb^{(1)}(n)$ in the nonstandard grading
$\fsb^{(1)}(n;n)$ for $p=0$. The main idea:
$\fsb^{(1)}(n)=\IM\Delta|_{\fb(n)}$, where $\Delta=\sum\frac{\del^2
}{\del q_i \del\xi_i}$. The dimensions of homogeneous components for $n$ even are:
\begin{equation*}\label{dimsB}
\footnotesize
\begin{tabular}{|l|l|l|}
\hline $i$&$\sdim\fg_i$&$\sdim\fsb^{(1)}(n;n)_i$ \tsep{1pt} \\\hline
$-1$&$2^{n-1}|2^{n-1}$&$2^{n-1}|2^{n-1}-1$\tsep{2pt}\\
$0$&$n\big(2^{n-1}|2^{n-1}\big)$&$(n-1)\big(2^{n-1}|2^{n-1}\big)+1|0$\\
$1$&$\frac12n(n+1)\big(2^{n-1}|2^{n-1}\big)$&
$\frac12\big(n^2-n+2\big)\big(2^{n-1}|2^{n-1}\big)-0|1$\bsep{2pt}\\
\hline
\end{tabular}
\end{equation*}

Let the weights of $\xi_i$ be $w(\xi_i)=(0,\dots,0,1,0,\dots,0)$ with a~1 on the $i$th place for $i<n$ and $w(\xi_n)=(-1,\dots,-1)$.

\subsubsection[Example: $n=4$]{Example: $\boldsymbol{n=4}$} 
For a~basis (even $\mid$ odd) of the $\fg_{0}$-module $\fg_{-1}\simeq
\frac{\Pi(\Vol(0|n))}{\Kee(1+\xi_1\cdots\xi_4)\vvol_\xi}$, where
$\fg_{0}\simeq \widetilde{\fsvect}(0|4))$, we take:
\begin{equation*}
\begin{array}{@{}ll|ll}\del_1:=\xi_1\vvol_\xi,&\del_2:=\xi_2\vvol_\xi,&
\del_9:=\xi_1\xi_2\vvol_\xi,&\del_{10}:=\xi_1\xi_3\vvol_\xi,\\
\del_3=\xi_3\vvol_\xi,&\del_4=\xi_4\vvol_\xi, &
\del_{11}:=\xi_1\xi_4\vvol_\xi,&\del_{12}:=\xi_2\xi_3\vvol_\xi,
\\
\del_5=\xi_1\xi_2\xi_3\vvol_\xi, &\del_6= \xi_1\xi_2\xi_4\vvol_\xi, &
\del_{13}:=\xi_2\xi_4\vvol_\xi,&
\del_{14}:=\xi_3\xi_4\vvol_\xi, \\
\del_7=\xi_1\xi_3\xi_4\vvol_\xi,& \del_8=\xi_2 \xi_3
\xi_4\vvol_\xi,& \del_{15}:=\xi_1 \xi_2\xi_3\xi_4\vvol_\xi.&
\end{array}
\end{equation*}

\textit{Critical coordinates of $\widetilde{\fsb}\big(15;\widetilde\un\big)=
\textbf{F}\big(\widetilde{\fsb}(8;\un|7)\big)$} are the same as those of
$\widetilde{\fsb}(8;\un|7)$: $N_5=N_6=N_7=N_8=1$, and also all those corresponding to
the formerly odd indeterminates.

\subsubsection{Partial prolongs} We have
$\sdim \fg_0 = 25|24$, and $\fg_0$ contains a~simple ideal of $\sdim = 21|24$,
the quotient is commutative; $\fg_{-1}$ is irreducible over this ideal.
We have $\sdim \fg_1 = 56|55$, there are 3 highest-weight vectors and 2 lowest-weight vectors in $\fg_1$;
\begin{gather*}
V_1=W_1, \quad V_1\subset V_2, \quad W_1\subset W_2 \subset W_3, \nonumber\\
\sdim V_1 = 24|21, \quad \sdim V_2 = \sdim W_3 = 32|31, \quad \sdim W_2 = 24|22,\nonumber\\
 \sdim (V_2+W_2) = 32|32,\quad \sdim (V_2+W_3) = 40|40,\quad
\sdim \fg_2 = 105|104.
\end{gather*}

The highest-weight vector of $W_3$ is $w_3=x_4^{(2)}\del_5$. This answer seems strange:
the algebra is symmetric with respect to the permutation of the $\xi_i$ while the list of highest-weight vectors is not. Performing
all possible permutations we obtain similar vectors $x_1^{(2)}\del_8$, $x_2^{(2)}\del_7$, $x_3^{(2)}\del_6$ (which are
not highest/lowest with respect to the division into positive/negative weights we have
selected first), but generate similar submodules $Y_1$, $Y_2$, $Y_3$ (and $Y_4=W_3$).

We have $Y_1+Y_2+Y_3+Y_4 = \fg_1$ and $\sdim(Y_1+Y_2+Y_3) = 48|48$.

Other highest-weight vectors:
\begin{gather}
w_1=x_3x_4\del_{10}+x_3x_{13}\del_5+x_4^{(2)}\del_{11}+x_4x_{12}\del_5+x_4x_{13}\del_6\nonumber\\
\hphantom{w_1=}{} +
x_2x_4\del_9+x_2x_{14}\del_5+x_4^{(2)}\del_{11}+x_4x_{12}\del_5+x_4x_{14}\del_7,\nonumber\\
w_2=x_3x_4\del_9+x_3x_{14}\del_5 +x_4x_{14}\del_6.\label{HWV87}
\end{gather}

The lowest-weight vectors:
\begin{gather*}
v_1=x_1x_2\del_4+x_1x_5\del_7+x_1x_9\del_{11}+x_1x_{12}\del_{14}+x_2x_5\del_8\nonumber\\
\hphantom{v_1=}{} + x_2x_9\del_{13} +x_2x_{10}\del_{14}+x_9x_{10}\del_7+x_9 x_{12}\del_8,\nonumber\\
v_2=x_1^{(2)}\del_4+x_1x_5\del_8+x_1x_9\del_{13}+x_1x_{10}\del_{14}+x_9x_{10}\del_8.
\end{gather*}
We have
\begin{gather}
  \sdim ([\fg_{-1}, V_1]) = 21|24,\quad   \sdim ([\fg_{-1}, Y_i]) = 22|24,\nonumber\\
 \sdim ([\fg_{-1}, Y_i+Y_j]) = 23|24 \quad \text{for $i\neq j$},\nonumber\\
 \sdim ([\fg_{-1},Y_i+Y_j+Y_k]) = 24|24 \quad \text{for $i\neq j\neq k\neq i$}. \label{sdims87}
\end{gather}

\underline{Partial prolongs of $\fg_0$ and the following parts of $\fg_1$}:
\begin{itemize}\itemsep=0pt
\item from $V_1=W_1$ and $W_2$: $\sdim \fg_2 = 11|8$, $\sdim \fg_3 =0|1$, no parameters;
\item from $V_2$: $\sdim \fg_2 = 33|32$, 1 parameter: $N_1$ (same for
$W_3$, parameter $N_4$);
\item from $V_2+W_3$: $\sdim \fg_2 = 56|56$, 2 parameters: $N_1$ and
$N_4$, similar for $Y_i+Y_j$;
\item from $Y_1+Y_2+Y_3$: $\sdim \fg_2 = 80|80$, 3 parameters: $N_1$,
$N_2$ and $N_3$.
\end{itemize}

\underline{Partial prolongs of the following parts of $\fg_0$, see equation~\eqref{sdims87}}:
\begin{itemize}\itemsep=0pt
\item from $(21|24)$: $\sdim \fg_1=24|27$, $\sdim \fg_2 = 11|8$, no parameters;
\item from $(22|24)$: $\sdim \fg_1=32|34$, $\sdim \fg_2 = 33|32$, 1 parameter;
\item from $(23|24)$: $\sdim \fg_1=40|41$, $\sdim \fg_2 = 56|56$, 2 parameters;
\item from $(24|24)$: $\sdim \fg_1=48|48$, $\sdim \fg_2 = 80|80$, 3 parameters.
\end{itemize}

\section[$\fv\fas(4; \un\vert 4)$]{$\boldsymbol{\fv\fas(4; \un\vert 4)}$}\label{Svas}
In this section, we can omit $\un$ when the arguments do not depend on it.

\subsection[For $p\neq 2$]{For $\boldsymbol{p\neq 2}$}\index{$vmathfrak\fas(4\vert 4)$@$\fv\fas(4\vert 4)$}For $\fg=\fv\fas(4|4)$ described in Table~\ref{ssExcAsCTS} as the Cartan prolong of the pair $(\id_{\fa\fs}, \fa\fs)$, we have another description:
$\fg_\ev=\fvect(4|0)$ and $\fg_\od=\Omega^{1}(4|0)\otimes_{\Omega^0(4|0)} \Vol^{-1/2}(4|0)$ with the
natural $\fg_\ev$-action on $\fg_\od$, and the bracket of odd elements given by
\[
\left[\frac{\omega_{1}}{\sqrt{\vvol}},
\frac{\omega_{2}}{\sqrt{\vvol}}\right]=
\frac{d\omega_{1}\wedge\omega_{2}+ \omega_{1}\wedge d\omega_{2}}{
\vvol},
\]
where we identify
\begin{gather}\label{ident}
\frac{dx_{i}\wedge dx_{j}\wedge dx_{k}}{
\vvol}=\sign(ijkl)\partial_{x_{l}}\quad \text{for any permutation $(ijkl)$ of $(1234)$}.
\end{gather}

\subsection[For $p=2$]{For $\boldsymbol{p=2}$} The first impression is that the characteristic-2 version of the Lie superalgebra $\fv\fas$ does not exist: the cocycle that determines the central extension $\fa\fs$ of $\fspe(4)$ is trivial, see~\cite{BGLL3}.
The following problem is most natural.

\begin{Problem}[on analogs of $\fas$ for $p=2$]\label{OPas} \index{(Problem, open@Problem, open} For $p = 2$, there are $8$ analogs of $\fpe(n)$ and 8 analogs of $\fspe(n)$, and lots of their nontrivial central extensions, see~{\rm \cite{BGLL3}}. Is there a~nontrivial central extension $\fe$ of one of these Lie (super)algebras, and an irreducible $\fe$-module $M$ such that $(M, \fe)_{1}\neq 0$? \end{Problem}

The above-mentioned ``first impression'' was, however, too hasty. Define \textit{a~character-$2$ analog of $\fvas$}
in the form $\fg=\fg_\ev\oplus \fg_\od$, where $\fg_\ev
= \fsvect(4|0)\rtimes \cF(4|0)$ and $\fg_\od=\Omega^1(4|0)$ with the natural $\fg_\ev$-action on $\fg_\od$. Define the bracket of odd elements by the formula
\[
[\omega_1,\omega_2]=\frac{d(\omega_1\wedge \omega_2)}{\vvol}+\Div\left(\frac{d\omega_1\wedge \omega_2}{\vvol}\right)
\]
subject to identification \eqref{ident}. Define the square of every
$\omega = \sum f_i dx_i\in\Omega^1(4|0)$ as follows, where $(i,j,k)$ is a~permutation of indices $(1,2,3)$: \index{$vasmathfrak$@$\fvas$}
\[
\omega^2 := \frac{d\omega \wedge \omega}{\vvol} + \sum\limits_{(i,j,k)\in S_3} \frac{\del f_i}{\del x_j}\left( \frac{\del f_k}{\del x_4} + \frac{\del f_4}{\del x_k} \right).
\]

Let $\fsvect(4|0)=\fsvect(y)$, where $y=(y_1,y_2,y_3,y_4)$. Consider the $\Zee$-grading of $\fg$ of depth~$1$, by setting
\[
\deg y_i=1,\quad \deg(dy_i)=-1.
\]
We get an embedding $\fg\tto\fvect(4|4)$. Let us describe the non-positive components of the embedded
algebra. Let the coordinates of the ambient be $x$ and $\xi$, and let us identify the basis elements of $\fg_{-1}$ with the following vector fields in $\fvect(4|4)=\fvect(x|\xi)$
\[
\del_{y_i}\longleftrightarrow \del_{x_i}, \quad dy_i \longleftrightarrow \del_{\xi_i}.
\]

Then, $\fg_{-1}=\Span\{\del_{y_i}, dy_i\} $ for $ i=1,\dots,4$, and $(\fg_0)_\ev$ consists of the pairs $(D, c)$, where $D=\sum_{i,j} a_{ij}y_i\del_{y_j}$ is any vector field such that $\Div D=0$, and $c\in\Kee$, whereas $(\fg_0)_\od$ consists of 1-forms $y_idy_j$:
\begin{equation}\label{IT}\footnotesize
\renewcommand{\arraystretch}{1.4}
\begin{tabular}{|c|c|c|}
\hline
Element of $\fg_0$ & its non-zero action & the corresponding vector field\\
\hline
&&\\[-1.5em]
\hline
$y_i\del_{y_j}, i\ne j$ & $\del_{y_i}\mapsto \del_{y_j},\; dy_j\mapsto dy_i$ & $x_i\del_{x_j}+\xi_j\del_{\xi_i}, i\ne j$\\
\hline
$y_i\del_{y_i}+y_j\del_{y_j}$, & $\del_{y_i}\mapsto \del_{y_i},\;\del_{y_j}\mapsto \del_{y_j},$ & $x_i\del_{x_i}+x_j\del_{x_j}+\xi_i\del_{\xi_i}+\xi_j\del_{\xi_j}$,\\
$i<j$& $\; dy_i\mapsto dy_i, \; dy_j\mapsto dy_j$ &$i< j$\\
\hline
$1$&$\begin{cases}\id&\text{on $(\fg_{-1})_\od$}\\
0&\text{on $(\fg_{-1})_\ev$}\\
\end{cases}$&$\sum \xi_i\partial_{\xi_i}$\\
\hline
&&\\[-1.5em]
\hline
$y_idy_i $ & $\del_{y_i}\mapsto dy_i$ & $x_i\del_{\xi_i}$\\
\hline
$y_idy_j,$ & $\del_{y_i}\mapsto dy_j,\; dy_k\mapsto \del_{y_l}, \;dy_l\mapsto \del_{y_k}, $ & $x_i\del_{\xi_j}+\xi_k\del_{x_l}+\xi_l\del_{x_k}$\\
$ i\ne j$& $(i,j,k,l)\in S_4$ & \\
\hline
\end{tabular}
\end{equation}

This $\fg_0$ is a~characteristic-2 analog of $\fas$.\index{$asmathfrak$@$\fas$}
In the basis $\del_{x_1}, \dots, \del_{x_4}, \del_{\xi_1},\dots, \del_{\xi_4}$ the $\fg_0$-action in $\fg_{-1}$ is given by the following (super)matrix whose correspondence to vector fields we give explicitly only for $(\fg_0)_\ev$ since the correspondence with $(\fg_0)_\od$ is too cumbrous to describe; moreover, is not worth the trouble thanks to the explicit table~\eqref{IT}:
\begin{gather*}
\begin{pmatrix} a& c\\
b+\tilde c & a^{\rm t} \end{pmatrix}+d\diag(0_4, 1_4), \quad \text{where} \ b^{\rm t}=b,\ c\in ZD, \ a\in\fsl(4), \nonumber\\
\tilde c_{ij}=E_{kl} \quad \text{for $k<l$ and $d\in\Kee$},\quad
\text{corresponding to $\sum a_{ij}y_i\del_{y_j} +d$}.
\end{gather*}

\begin{Claim}[description of the simple part of $\fvas$]\quad
\begin{enumerate}\itemsep=0pt
\item[$1)$] The Lie superalgebra $\fvas^{(1)}(4;\un|4)$ is simple; its even part is $\fsvect^{(1)}(4;\un|0)\rtimes\Vol_0(4;\un|0)$, see Section~{\rm \ref{ssVasNotSimple}}.
\item[$2)$] The critical coordinates of the shearing vector for the simple Lie algebra $\fvas^{(1)}\big(8;\widetilde\un\big)$~-- the desuperization of $\fvas^{(1)}(4;\un|4)$~-- are the ones that correspond to formerly odd indeterminates.
\end{enumerate}
\end{Claim}

\subsubsection{Partial prolongs} In order to investigate possible partial prolongs, we have to consider the $\fg_0$-submodules $V_i$ of $\fg_1$ such that $[\fvas_{-1}, V_i]=\fg_0$. Since it is not clear what is a~lowest/highest-weight vector with respect to $\fg_0$, we consider the lowest-weight vectors with respect to $(\fg_0)_\ev$, and build the $\fg_0$-submodules from them.

\begin{Claim}[lowest-weight vectors] \label{sssParPr}
There are $7$ lowest-weight vectors LWVs:
\begin{equation}\label{v}\footnotesize
\renewcommand{\arraystretch}{1.4}
\begin{tabular}{|c|c|c|}
\hline
 LWV & its image in $(x,\xi)$-model & its image in $y$-model\\
\hline
\hline
$v_1$ &$x_2x_3\del_{\xi_4}+x_2x_4\del_{\xi_3}+x_3x_4\del_{\xi_2}$ &$y_2y_3dy_4+y_2y_4dy_3+y_3y_4dy_2$ \\
\hline
$v_2$ &$x_3x_4\del_{x_1}+x_3\xi_1\del_{\xi_4}+x_4\xi_1\del_{\xi_3}$ &$y_3y_4\del_{y_1}$ \\
\hline
$v_3$ &$x_3^{(2)}\del_{\xi_4}+x_3x_4\del_{\xi_3}$ &$y_3^{(2)}dy_4+y_3y_4dy_3$ \\
\hline
$v_4$ &$x_4^{(2)}\del_{\xi_4}$ &$y_4^{(2)}dy_4$ \\
\hline
$v_5$ &$x_4^{(2)}\del_{x_1}+x_4\xi_1\del_{\xi_4}$ &$y_4^{(2)}\del_{y_1}$ \\
\hline
$v_6$ &$x_3x_4\del_{\xi_4}+x_4\xi_1\del_{x_2}+x_4\xi_2\del_{x_1}+\xi_1\xi_2\del_{\xi_4} $&$y_3y_4dy_4$ \\
\hline
$v_7$ &$x_4\sum_{i=1}^4\xi_i\del_{\xi_i}+\xi_1\xi_2\del_{x_3}+\xi_1\xi_3\del_{x_2}+\xi_2\xi_3\del_{x_1}$ &$y_4$ \\
\hline
\hline
\end{tabular}
\end{equation}
\end{Claim}

\begin{Claim}[no partial prolongs] \label{sssParPr1} Let $\fg:=\fvas(4;\un|4)$.\quad
\begin{enumerate}\itemsep=0pt
\item[$1)$] Let $V_i$ be the $\fg_0$-submodule of $\fg_1$ generated by $v_i$, see \eqref{v}. Then, $[V_i, \fg_{-1}]=\fg_0$ for all~$i$.
For $\un$ unconstrained, we have $\sdim(\fg_1)=40|40$, and
\begin{gather*}
V_7=V_6=V_2=V_1= V_3\cap V_4, \quad \text{this $\fg_0$-module is irreducible}, \quad V_4=V_5,\\
\sdim({V_1})= 24|24,\quad \sdim({V_3})=\sdim({V_4})= 28|28.
\end{gather*}
\item[$2)$]  Let $\fg^{V_i}:=(\fg_-, \fg_0, V_i)_*$ be the prolong in the direction of $V_i$. Then, $\fg^{V_1}= \fvas(4;\One|4)$.
We have $\sdim \fvas^{(1)}(4;\One|4)=60|64$.
\item[$3)$] In the quotient $\fg_1/V_1$, to each $i\in \{1,2,3,4\}$ there corresponds a~$4|4$-dimensional submodule $M_i$ spanned by the images of $y_i^{(2)}\partial_{y_j}$ and $y_i^{(2)}dy_j$ for $j = 1,2,3,4$. Each $M_i$ is irreducible, and the images of $M_i$ and $M_j$ in $\fg_1/V_1$ do not intersect for $i\neq j$. Thus, $\fg_1$ contains a~submodule~$V_1$ corresponding to $\un=\One$, and up to four modules $M_i$ glued to $V_1$ if $\un\neq\One$. The partial prolongation in the direction of $(\oplus_{i\in I\subset\{1,2,3,4\}}M_i)\rtimes V_1$ is $\fvas(4;\un|4)$, where
    \[ N_i = \begin{cases}\infty&\text{if $i\in I$},\\
1&\text{if $i\not\in I$.}\end{cases}
\]
\end{enumerate}
\end{Claim}

\begin{proof}[Idea of the proof] Since there is no complete reducibility, to prove item 3) we have to consider also highest-weight vectors (HWV) with respect to $(\fg_0)_\ev$. Then, we are able to find the two quotients modules $M_i$ invisible in table~\eqref{v} since their LWVs go to $V_1$ under $(\fg_0)_\ev$. We have already encountered similar phenomenon in previous sections considering LWVs and HWVs with respect to the whole $\fg_0$ for respective $\fg$. We skip the table of HWVs analogous to~\eqref{v}.
\end{proof}

\section[Cartan prolongs of the Shen algebra; Melikyan algebras for $p=2$]{Cartan prolongs of the Shen algebra; Melikyan algebras\\ for $\boldsymbol{p=2}$}\label{SShen}

\subsection{Brown's version of the Melikyan algebra in characteristic~2}\label{ssBrMe}
Brown~\cite{Bro} described characteristic-2 analog of the Melikyan algebra\index{(Melikyan algebra@Melikyan algebra}
as follows.
As spaces, and $\Zee/3$-graded Lie algebras, let
\begin{equation*}
L(\un):=\fg_\ev \oplus\fg_{\bar 1}\oplus\fg_{\bar 2} \simeq
\fvect(2; \un)\oplus \Vol(2; \un)\oplus \cO(2; \un) .
\end{equation*}
The $\fg_\ev$-action on the $\fg_{\bar i}$ is natural (adjoint, on
volume forms, and functions, respectively); $\cO(2; \un)=\Kee[u_1,
u_2; \un]$ is the space of functions; $\Vol(2;\un)$ is the
space of volume forms with volume element $\vvol:= \vvol(u)$, where $u=(u_1,
u_2)$. Let the
multiplication in $L(\un)$ be given, for any $f,g\in\cO(2; \un)$, by
the following formulas:
\begin{equation*}
{}[f\vvol, g\vvol]=0,\quad [f\vvol, g]=fH_g,\quad  [f, g]:=
H_f(g)\vvol,
\end{equation*}
where
\begin{equation*}
 H_f=\pderf{f}{u_1}\del_{u_2}+\pderf{f}{u_2}\del_{u_1} .
\end{equation*}

Define a~$\Zee$-grading of $L(\un)$ by setting
\begin{equation*}
\deg u^{(\underline{r})}\del_{u_i} =3|\underline{r}|-3,\quad \deg
u^{(\underline{r})}\vvol =3|\underline{r}|-2,\quad
\deg u^{(\underline{r})} =3|\underline{r}|-4.\\
\end{equation*}
Now, set $\fme(5;\un):=L(\un)/L(\un)_{-4}$, where $ L(\un)_{-4}$ is
the center (the space of constants). The algebra $\fme(5;\un)$ is
not simple, because $\Vol(2; \un)$ has a~submodule of codimension 1;
but $\fme^{(1)}(5;\un)$ is simple; in~\cite{Ei}, Eick denoted what
we denote $\fme^{(1)}(5;\One)$\index{$memathfrak^{(1)}(5;\One)$@$\fme^{(1)}(5;\One)$} by $\text{Bro}_2(1,1)$. This algebra
was discovered by Shen Guangyu, see~\cite{Shen}, and should be
denoted somehow to commemorate his wonderful discovery, we suggest
to designate this Shen's analog of $\fg(2)$ by $\fg\fs(2)$.\index{$gmathfrak\fs(2)$, Shen's analog of $\fg(2)$@$\fg\fs(2)$, Shen's analog of $\fg(2)$}
\index{(Shen algebra@Shen algebra}

There are two $\Zee$-gradings of $\fg(2)$ with \textbf{one}
pair of Chevalley generators of degree $\pm 1$ (the other generators
being of degree 0): one $\Zee$-grading of depth 2 and the other
one of depth 3. As is easy to see, for the grading of depth 3, the
nonpositive parts of $\fg(2)$ over fields $\Kee$ of characteristic~
$p\neq 3$ and those of $\fme(5;\un)$ are isomorphic.
Remarkably, this description holds for any $p\neq 3$, see
\cite{Shch}. For $p=3$, the positive parts of the prolongation have
the same dimensions as those of $\fg(2)$ for $p\neq 2, 5$, but
$[\fg_1, \fg_{-1}]=\Kee 1_2$, the center of $\fgl(2)$. (By the way,
the realization of the nonpositive components of $\fg(2)$, see
equation~\eqref{Zgrad102bis}, that works for $p\neq 3$, should be
modified for $p=3$, but we skip this since neither the complete prolong
nor any partial prolong is simple.)

Let $U[k]$ be
the $\fgl(V)$-module which is $U$ as $\fsl(V)$-module, and let the
central element $z\in \fgl(V)$ represented by the unit matrix, which acts on
$U[k]$ as $k\id$, where $k$ should be understood modulo $p$. Then,
the grading of depth 3 is of the form
\begin{equation*}
\footnotesize
\renewcommand{\arraystretch}{1.4}
\begin{tabular}{|c|c|c|c|}
\hline $\fg_0$&$\fg_{-1}$&$\fg_{-2}$&$\fg_{-3}$\cr \hline
$\fgl(2)\simeq\fgl(V)$&$V=V[-1]$&$E^2(V)[-2]$&$V[-3]$\cr
\hline\end{tabular}\quad\text{~~for $\Char\Kee\neq 3$}.
\end{equation*}

Set $\del_i:=\del_{x_i}$
to distinguish it from $\del_{u_i}$; we use both representations in
terms of $x$ and~$u$, whichever is more convenient. Here is the (borrowed from \cite{Shch}) description of nonpositive components of $\fme(5;\un)$,
which are the same as those of $\fg\fs(2)$ and $\fg(2)$, by means of
vector fields:
\begin{equation}\label{Zgrad102bis}
\footnotesize \renewcommand{\arraystretch}{1.4}
\begin{tabular}{|l|l|}
\hline
$\fg_{i}$&the basis elements\\ \hline
\hline
$\fg_{-3}$ &
$\del_{u_1}\longleftrightarrow\del_{1}, \
\del_{u_2}\longleftrightarrow\del_{2}$\\
\hline

$\fg_{-2}$ & $\vvol\longleftrightarrow\del_{3}$\\
\hline

$\fg_{-1}$ & $X_2^-:=u_1\longleftrightarrow
(x_3+x_4x_5)\del_{2}+\del_4,\quad u_2\longleftrightarrow
x_3\del_{1}+x_4\del_3+ \del_5$
 \\
\hline

& $u_1\del_{u_1}\longleftrightarrow
x_1\del_{1}+x_3\del_3+x_4\del_4$,\\

$\fg_{0}\simeq$& $X_1^+:=u_1\del_{u_2}\longleftrightarrow
x_5^{(3)}\del_{1}+(x_1+
x_4x_5^{(2)})\del_{2}+x_5^{(2)}\del_3+x_5\del_4$\\
$\fgl(2)$& $X_1^-:=u_2\del_{u_1}\longleftrightarrow(x_2+
x_4^{(2)}x_5)\del_{1}+x_4^{(3)}\del_{2}+ x_4^{(2)}\del_3+
x_4\del_5$
 \\
& $u_2\del_{u_2}\longleftrightarrow
x_2\del_{2}+x_3\del_3+
x_5\del_5$\\
\hline
\end{tabular}
\end{equation}
The highest-weight vector in $\fg_{-1}$ is $X_2^-:=u_1$. Consider
the positive part of $\fg=\fg\fs(2)$. The lowest-weight vector in
$\fg_1$ is given by the vector field
\[
X_2^+:=x_4^{(3)}x_5\del_{2}+\big(x_2+ x_4^{(2)}x_5\big)\del_3+x_4 x_5\del_5\quad (=u_2\vvol).
\]
So far, the generators and the dimensions of the components look like
their namesakes in $\fg(2)$ for $p>3$; however, the relations are
different: To facilitate comparison with presentations in terms of
Chevalley generators, set $H_i:=[X_i^+, X_i^-]$; i.e.,
\begin{gather*}
H_1=x_1\del_{1}+x_2\del_{2}+x_4\del_4+x_5\del_5\quad
(=u_1\del_{u_1}+u_2\del_{u_2}),\nonumber\\
H_2=x_2\del_{2}+x_3\del_3+
x_5\del_5\quad  (= u_2\del_{u_2}). 
\end{gather*}
Clearly, $H_1$ is the central element of $\fg_0$; for its grading element we take $u_1\del_{u_1}$, see~\cite{BGL1}.

\begin{Lemma}[the multiplication tables in
$\fg\fs(2)$ and $\fg(2)$]\label{LemmaShen} The multiplication tables in
$\fg\fs(2)$ and $\fg(2)$ are as follows \textup{(for $\fg(2)$, we get
$[H_i,X_j^\pm]=\pm A_{ij} X_j^\pm$, not $[H_i,X_j^\pm]=\pm
A_{ji} X_j^\pm$; let $X_3^\pm:=[X_1^\pm, X_2^\pm]$)}
\begin{equation*}
\footnotesize
\renewcommand{\arraystretch}{1.4}
\begin{tabular}{|@{\,}l@{\,}|@{\,}c@{\,}||@{\,}l@{\,}|@{\,}c@{\,}||@{\,}l@{\,}|@{\,}c@{\,}|}
\hline
in $\fg\fs(2)$&in $\fg(2)$& in $\fg\fs(2)$&in $\fg(2)$
& in $\fg\fs(2)$&in $\fg(2)$\\
\hline
$[H_1,X_1^+]=0$&$2 X_1^+$&$[H_2,X_1^+]= X_1^+$&
$-3 X_1^+$&$[H_1,H_2]=0$&$0$\\
$[H_1,X_2^+]= X_2^+$&$-X_2^+$&$[H_2,X_2^+]=0$&
$2 X_2^+$&$[X_1^-, X_2^-]=x_3\del_{1}+x_4\del_3+ \del_5$&$X_3^-$\\

$[H_1,X_1^-]=0$&$-2 X_1^-$&$[H_2,X_1^-]=X^-_1$&
$3X^-_1$&$[X_1^+, X_2^+]= u_1\vvol $&$X_3^+$\\

$[H_1,X_2^-]=X^-_2$&$ X^-_2$&$[H_2,X_2^-]=0$&
$-2 X_2^-$&$[X_1^\pm,X_2^\mp]=0$&$0$\\
\hline
\end{tabular}
\end{equation*}
\end{Lemma}

\textit{Critical coordinates of $\fme(5;\un)$}: $\un_3=1$.

The $\fg_0$-module $\fg_1$ is generated by the lowest-weight vector
$X_2^+$; we have $\dim\fg_1=2$. Since~$X_1^\pm$ and $X_2^+$ contain
$x_4$ and $x_5$ in degrees~2 and~3, see equation~\eqref{Zgrad102bis}, the
corresponding coordinates of the shearing vector in the generic case are
$\geq 2$; for the shearing vector with the smallest coordinates still
ensuring simplicity; i.e., for
$\un=(1,1,1,2,2)$, the prolong $\fg$ is of dimension
17; it has ideals of dimension 14, 15, 16. The ideal of dimension
14 is simple, see \cite{Bro,Ei, Shen}.

\section{Miscellaneous remarks}\label{Ssvect}

\subsection[Desuperizations that are nonsimple if $N_i<\infty$ for all $i$]{Desuperizations that are nonsimple if $\boldsymbol{N_i<\infty}$ for all $\boldsymbol{i}$}\label{ssSvect}

In Section~\ref{ssSimpleIdeal}, the simple derived algebras of various
W-graded versions of $\fk\fas$ are described; this is new. The results of
\textit{this} section are not new (although they were usually considered for $p>2$); see, e.g., Lemma~2.4 in~\cite{Kfil};
we present them for completeness, see also equation~\eqref{specValOfLambForP2} and Section~\ref{ssBrMe} on $\fm\fe^{(1)}$.

\subsubsection[$\fg=\fsvect(n;\un)$]{$\boldsymbol{\fg=\fsvect(n;\un)}$}\label{ssSvectNotSimple} Let us
prove that the elements of the form
\[
D_k=\left(\prod\limits_{i\in\{1,\dots,n\},\ i\neq k}
x_i^{(2^{N_i}-1)}\right)\del_k
\]
do not lie in $\fg^{(1)}$. In what follows we assume that $k=n$, for
definiteness. As $\fg$ is a~sum of its $\Zee^n$-weighted components,
it suffices to show that $D_n$ cannot be obtained as the bracket of
two elements homogenous with respect to the grading by the weight. As the $x_n$-weight
(i.e., weight with respect to $x_n\del_n$) of $D_n$ is equal to $-1$, which
is also the minimal possible $x_n$-weight in $\fg$, it follows that,
in order to obtain $D_n$ as a~bracket, one of the factors (we say
``factor" speaking about the Lie bracket, just as we do it for an
associative multiplication) has to have weight~${-1}$ as well. Then,
if this factor is homogenous w.r.t. the $\Zee^n$-weight, it must be a
monomial of the form $a=\Big(\mathop{\prod}\limits_{1\leq i\leq
n-1}x_i^{(r_i)}\Big)\del_n$ up to a~scalar multiplier, where
${0\leq r_i< 2^{N_i}}$. Then, from the weight considerations, the
other factor must be of the form
\begin{gather*}
b = \mathop{\sum}\limits_{1\leq i<n,\
r_i>0}c_i\left(\prod\limits_{1\leq j<n, \  j\neq i}
x_j^{(2^{N_j}-1-r_j)}\right)x_i^{(2^{N_i}-r_i)}\del_i +
c_n\left(\prod\limits_{1\leq j\leq
n-1}x_j^{(2^{N_j}-1-r_j)}\right)x_n\del_n.
\end{gather*}
Clearly,
\begin{gather*}
{}[a,b] = \left( \mathop{\sum}\limits_{1\leq i<n\text{~such that~}
r_i>0,\ i=n}c_i\right) D_n,\\
\Div b = \left( \mathop{\sum}\limits_{1\leq i<n\text{~such that~}
r_i>0, \ i=n}c_i\right) \left(\prod\limits_{1\leq j\leq
n-1}x_j^{(2^{N_j}-1-r_j)}\right).
\end{gather*}
So $b\in\fg$ if and only if $[a,b]=0$, hence $\fg^{(1)}$ contains no
elements of the same weight as $D_n$.

\subsubsection[$\fvas$ for $p=2$]{$\boldsymbol{\fvas}$ for $\boldsymbol{p=2}$}\label{ssVasNotSimple} In this case, the even part of the Lie superalgebra $\fvas(4;\un|4)$, and of its $\Zee/2$-graded desuperization, should be diminished to get a~simple Lie algebra, namely
\[
\fvas^{(1)}(4;\un|4)_\ev=\fsvect^{(1)}(4;\un|0)\rtimes \Vol_0(4;\un|0).
\]

\subsubsection{The Lie (super)algebra of contact vector fields}\label{nonSimp} \underline{Let $p\neq 2$}. As follows from
equation~\eqref{div2}, if $2n+2-m\equiv 0\mod p$, then the Lie superalgebra $\fk(2n+1;\un|m)$ is
divergence-free, its derived algebra is simple.

If $2n+2-m\equiv -2\mod p$, then $\fk(2n+1;\un|m)\simeq \Vol$, and hence not simple; it contains a~codimension 1 ideal, $\fk^{(1)}(2n+1;\un|m)$.

 \underline{Let $p=2$}. If $(n, m)\neq(0,0)$, then the Lie (super)algebra $\fk(2n+1;\un|2m)$ is
divergence-free if $n+m+1\equiv 0\mod 2$, see equation~\eqref{K_f2}.

The \textit{Zassenhaus algebra} $\fvect(1;\un)$ for $p=2$ is not simple; observe that
$\fvect(1;\un)\simeq\fk(1;\un)$.

\subsection[On deforms of $\fsvect$ and $\fh$. Quantizations]{On deforms of $\boldsymbol{\fsvect}$ and $\boldsymbol{\fh}$. Quantizations} \label{ssTyu}

\begin{itemize}\itemsep=0pt
\item In \cite{Tyu}, Tyurin described non-isomorphic
filtered deforms of
the Lie algebras of series $\fsvect$ for $p>3$ considered in the
\textit{standard} $\Zee$-grading. There are three
statements in~\cite{Tyu} that should be corrected.

First, in the introduction to~\cite{Tyu}, Tyurin wrote that in
\cite{Kfil} Kac proved that
all deforms of $\fsvect$ for $p>3$ are filtered. Kac did not claim
this in~\cite{Kfil}. Moreover, Kac did not claim he described \textit{all}
filtered deformations, either; Kac writes only about filtered deformations
associated with the \textit{standard} $\Zee$-gradings.

Today, when the simple modular Lie algebras are classified for $p>3$, the list of all
their deforms is not needed for \textit{classification},
but is a~useful part of \textit{interpretation} of the algebras found, see, e.g., \cite{Sk0,Sk1}; this is
of independent interest, like knowledge about ``occasional isomorphisms"
$\fo(3)\simeq\fsl(2)$ or $\fo(6)\simeq\fsl(4)$, or $\fvect(1|1)\simeq \fm(1)\simeq\fk(1|2)$, as abstract Lie superalgebras.

Second, for any $p$, a~particular deformation~-- called \textit{quantization} in physical literature~-- of the Poisson Lie algebra on~2 indeterminates, induces a~deform
of $\fsvect(2;\un)\simeq \fh(2;\un)$, at least for $\un$ of the form
$(a,a)$ for any $a\geq1$, cf.~\cite{BLLS1}. Therefore, in~\cite{Tyu}, the claims describing all
deformations of $\fsvect(m;\un)$ should have been confined to $m>2$
and, moreover, Tyurin's main theorem should only claim a~complete description of
non-isomorphic \textit{filtered} deforms related to the \textit{standard} $\Zee$-grading; for examples of filtered deforms of
$\fsvect^{(1)}(3;\One)\simeq \fh^{(1)}(4;\One)$ corresponding to distinct $\Zee$-gradings, see~\cite{ChKu}.

Although other deforms of $\fh(2n;\un)$ do not provide us with new Lie algebras, they do provide us with new deforms, non-isomorphic to the filtered deforms.

Third, Wilson~\cite{W} corrected the main result of Tyurin who found all normal shapes of volume forms for $p>2$, but missed an isomorphism. Wilson wrote only about normal shapes of volume forms, thus avoiding any discussion of deforms of $\fsvect$.

\item The deform of $\fsvect(5;\un)$ we describe here is a~completely new, exceptional, simple vectorial Lie algebra. It exists only for $p=2$, the case neither Tyurin nor Wilson considered.

\item The characteristic-2 analogs of exceptional deformations of $\fh$ and $\fb$ described in \cite{LSh} can have both even and odd parameters. The complete description of the deformations is unknown.
\end{itemize}

\section{Tables}\label{STables}

\subsection[Series of vectorial Lie superalgebras over $\Cee$; conditions for their simplicity]{Series of vectorial Lie superalgebras over $\boldsymbol{\Cee}$;\\ conditions for their simplicity}
In Table \eqref{1.4}, FD indicates finite dimension.
\begin{equation}\label{1.4}\footnotesize
\renewcommand{\arraystretch}{1.4}
\begin{tabular}{|l|l|}
\hline $N$&the family and conditions for its simplicity\cr \hline

1&$\fvect(m|n; r)$ for $m\geq 1$ and $0\leq r\leq n$ \cr \hline

2 &
$\fvect(0|n; r)$ for $n>1$ and $0\leq r\leq n$ (FD)\cr \hline

3&$\fsvect(m|n; r)$ for $m>1$, $0\leq r\leq n$ \cr \hline

4 &
$\fsvect(0|n; r)$ for $n>2$ and $0\leq r\leq n$ (FD)\cr \hline

5&$\fsvect^{(1)}(1|n; r)$ for $n>1$, $0\leq r\leq n$\cr \hline

6 &$\widetilde{\fsvect}(0|n)$ for $n>2$ (FD)\cr \hline
&\\[-1.5em]
\hline 
7&$\fk(2m+1|n; r)$ for $0\leq r\leq [\frac{n}{2}]$ unless
$(m|n)=(0|2k)$\cr &$\fk(1|2k; r)$ for $0\leq r\leq k$ except
$r=k-1$\cr \hline

8&$\fh(2m|n; r)$ for $m>0$ and $0\leq r\leq [\frac{n}{2}]$\cr
\hline
9&$\fh_{\lambda}(2|2; r)$ for $\lambda\neq -2, -\frac32, -1, \frac12,
0, 1, \infty$, and \cr &$r=0$, $1$ and $\Reg_{\fh}$ (see Sect.~$1.3.1$ in \cite{LSh})\cr \hline

10&$\fh^{(1)}(0|n)$ for $n>3$ (FD)\cr \hline

\hline
11&$\fm(n|n+1; r)$ for $0\leq r\leq n$ except $r= n-1$\cr
 \hline

12&$\fsm(n|n+1; r)$ for $n>1$, but $n\neq 3$ and $0\leq r\leq n$
except $r= n-1$\cr \hline

13&$\fb_{\lambda}(n|n+1; r)$ for $n>2$, where $\lambda\neq 0, 1,
\infty$ and $0\leq r\leq n$ except $r= n-1$\cr
\hline

14&$\fb_{1}^{(1)}(n|n+1; r)$ for $n>2$ and $0\leq r\leq n$ except $r=
n-1$\cr \hline

15&$\fb_{\infty}^{(1)}(n|n+1; r)$ for $n>2$ and $0\leq r\leq n$ except
$r= n-1$\cr \hline

16&$\fle(n|n; r)$ for $n>1$ and $0\leq r\leq n$ except $r= n-1$\cr
\hline

17&$\fsle^{(1)}(n|n; r)$ for $n>2$ and $0\leq r\leq n$ except $r= n-1$\cr
\hline

18&$\widetilde{\fsb}_{\mu}(2^{n-1}-1|2^{n-1})$ for $\mu\neq 0$ and
$n>2$\cr \hline
\end{tabular}
 \end{equation}

\begin{landscape}

\subsection[Lie algebras $\textbf{F}(\fg)$ over $\Kee$ ($p=2$) analogous to serial
vectorial Lie\\ superalgebras~$\fg$ over $\Cee$ and names of
both]{Lie algebras $\boldsymbol{\textbf{F}(\fg)}$ over $\boldsymbol{\Kee}$ ($\boldsymbol{p=2}$) analogous to serial
vectorial Lie superalgebras $\boldsymbol{\fg}$ over $\boldsymbol{\Cee}$ and names of
both}\label{forgetseries}
\begin{equation}\label{cts-prol}
\footnotesize
\renewcommand{\arraystretch}{1.3}
\begin{tabular}{|c|c|c|c|c|c|c|}
\hline
$N$&$\fg$&$\fg_{-2}$&$\fg_{-1}$&$\fg_0$&$\textbf{F}(\fg_0)$&
$\textbf{F}(\fg)$\cr
\hline
&&&&&&\\[-1.4em] \hline

$1$&$\begin{array}{l}\fvect(n|m)\\
\text{for $mn\neq 0$, $n > 1$ or $m = 0$, $n >
2$}\end{array}$&$-$&$\id\simeq V
$&$\fgl(n|m)\simeq\fgl(V)$&$\fgl(n+m)$&$\fvect(n+m;\widetilde\un)$\cr
\hline
&&&&&&\\[-1.4em]
 \hline $2$&$\fsvect(n|m)$ for $m,n\neq 1$&$-$&$\id\simeq V
$&$\fsl(n|m)\simeq\fsl(V)$&$\fsl(n+m)$&$\fsvect(n+m;\widetilde\un)$\cr
\hline
&&&&&&\\[-1.4em]
\hline $3$&$\fh_\cB(2n|m)$, where $mn\neq 0$, $n >
1$,&$-$&$\id$&$\fosp_\cB(m|2n)$&$\fo_{F(\cB)}(m+2n)$&
$\fh_{F(B)}(2n+m;\widetilde\un)$\cr \hline
&&&&&&\\[-1.4em]
\hline $4$&$\fk(2n+1|m)$ for $mn\neq 0$ and $m$
even&$\Fee$&$\id\simeq
V$&$\fc\fosp_\cB(m|2n)\simeq\fc\fosp(V)$&$\fc\fo_{F(\cB)}(m+2n)$&
$\fk(2n+m+1;\widetilde\un)$ \cr \hline
&&&&&&\\[-1.4em]
 \hline

$5$&$\fm(n):=\fm(n|n+1)$ for $n>1$&$\Pi(\Fee)$&$\id\simeq
V$&$\fc\fpe(n)\simeq\fc\fpe(V)$&$\fc\fpe(n)$&$\fk(2n+1;\widetilde\un)$\cr
\hline

$6$&$\fb_{\lambda}(n;n)$ for $n>1$&$-$&$\Pi(\Vol^{\lambda}(0|n))$&
$\fvect(0|n)$&$\fvect(n;\One)$&
 $\fpo_{\lambda}(2n+1;\widetilde\un)$\cr

$6_1$&$\fb_{1}^{(1)}(n;n)$ for $n>1$&$-$&$\Pi(\Vol_0(0|n))$&
$\fvect(0|n)$&$\fvect(n;\One)$&
 $\fpo_{1}^{(1)}(2n+1;\widetilde\un)$\cr

 $6_\infty$&$\fb_{\infty}^{(1)}(n;n)\simeq\fsm(n;n)$ for $n>1$&$-$&$\Pi(\Vol_0(0|n))$&
$\fsvect(0|n)\rtimes \Vol_0(0|n)$&$\fsvect(n;\One)\rtimes \Vol_0(n;\One)$&
 $\fpo_{\infty}^{(1)}(2n+1;\widetilde\un)$\cr

$7$&$\fb_{a, b}(n)$ for $n>1$&$\Pi(\Fee)$&$\id$& $\fspe(n)_{a, b}$&
$\fspe(n)_{a, b}$& $\fpo_{a, b}(2n;\widetilde\un)$\cr \hline
&&&&&&\\[-1.4em]
 \hline

$8$&$\fle(n):=\fle(n|n)$ for $n>1$&$-$&$\id\simeq V$&$\fpe(n)\simeq
\fpe(V)$&$\fpe(n)$&$\fh_\Pi(2n;\widetilde\un)$\cr \hline

$9$&$\fsle(n):=\fsle(n|n)$ for $n>1$&$-$&$\id\simeq
V$&$\fspe(n)\simeq
\fspe(V)$&$\fspe(n)$&$\fs\fh_\Pi(2n;\widetilde\un)$\cr \hline
&&&&&&\\[-1.4em]
\hline

$10$&$\widetilde{\fsb}_{\mu}(2^{n-1}-1|2^{n-1})$ or $\widetilde{\fsb}_{\mu}(2^{n-1}|2^{n-1}-1)$&$-$
&$\frac{\Pi(\Vol(0|n))}{\Fee(1-\mu\xi_1\cdots\xi_n)\vvol_\xi}$&
$\widetilde{\fsvect}_\mu (0|n)$&$\widetilde{\fsvect}_\mu
(n;\One)$&$\widetilde{\fsb}_{\mu}(2^{n}-1;\widetilde\un)$\cr
\hline\end{tabular}
\end{equation}

\subsubsection{Remarks}\label{TabRem} In all lines $\Par\widetilde\un=\dim\widetilde\un$, except for the bottom line, see Section~\ref{StildeSB}.
To save space, we skip most of the conditions for simplicity in
Table~\eqref{cts-prol}. In columns $\fg_i$ for $i<0$, obviously, $\Fee$ is $\Cee$ or
$\Kee$. In lines $N=6, 7$, we have $\lambda =\frac{2a}{n(a-b)}$ for
$p\neq2$ and $\lambda =\frac{a}{b}$ for $p=2$. In line $6_\infty$, we identify $\Vol_0$ with a~subspace of the space of functions $\cF$. In line 10, the Lie superalgebra
\[
\widetilde{\fsvect}_\mu (0|n):=(1+\mu \xi_1\cdots \xi_n)\fsvect(0|n)
\quad \text{preserves the volume element~~$(1-\mu \xi_1\cdots
\xi_n)\vvol_\xi$,~~where $p(\mu)\equiv n\mod 2$.}
\]
For $n$ even, we can (and do) set $\mu=1$, whereas $\mu$ odd should
be considered as an additional indeterminate on which the
coefficients depend. The Lie superalgebras $\widetilde{\fsvect}_\mu
(0|n)$ are isomorphic for nonzero $\mu$'s; and therefore so are the
algebras
\begin{gather*}
\widetilde{\fsb}_{\mu}\big(2^{n-1}-1|2^{n-1}\big):=(1+\mu \xi_1\cdots
\xi_n)\fsb(n;n) \text{for $n$ even},\quad
\widetilde{\fsb}_{\mu}\big(2^{n-1}|2^{n-1}-1\big):=(1+\mu \xi_1\cdots
\xi_n)\fsb(n;n) \text{for $n$ odd},
\end{gather*}
Recall the definition of $\fspe(n)_{a, b}$ in Section~\ref{AandB}.

\textit{To be specified}: Some of the Lie superalgebras in Table
(\ref{cts-prol}) are not simple; it is their quotients modulo their centers
or ideals of codimension~1 which are simple (such are $\fsvect(1|m)$,
$\fh(0|m)$, $\fb_{\lambda}(n)$ for certain values of $\lambda$, and
$\fsle(n)$); some small values of superdimension should be excluded
(like $(1|1)$ and $(0|m)$, where $m\leq 2$, for the $\fsvect$ series;
$(0|m)$, where $m\leq 3$, for the $\fh$ series; etc.)
\end{landscape}

\subsection[Exceptional vectorial Lie superalgebras over $\Cee$]{Exceptional vectorial Lie superalgebras over $\boldsymbol{\Cee}$}\label{excepC}
\[\footnotesize
\renewcommand{\arraystretch}{1.3}
\begin{tabular}{|c|c|c|c|c|}
\hline $\fg$&$\fg_{-2}$&$\fg_{-1}$&$\fg_0$&$\fg(\sdim\fg_{-})$\cr
\hline
&&&&\\[-1.4em]
 \hline

$\fv\fle(4|3)$&$-$&$\Pi(\Lambda(3)/\Cee
1)$&$\fc(\fvect(0|3))$&$\fv\fle(4|3)$\cr \hline

$\fv\fle(4|3; 1)$&$\Cee[-2]$& $\id_{\fsl(2;\Lambda
(2))}\boxtimes \vvol^{1/2}$& $\fc(\fsl(2;\Lambda(2))\ltimes
T^{1/2}(\fvect(0|2))$&$\fv\fle(5|4)$\cr \hline

$\fv\fle(4|3; K)$&$\id_{\fsl(3)}\boxtimes
\Cee[-2]$&$\id^*_{\fsl(3)}\boxtimes \id_{\fsl(2)}\boxtimes
\Cee[-1]$&$\fsl(3)\oplus\fsl(2)\oplus\Cee z$&$\fv\fle(3|6)$\cr
\hline
&&&&\\[-1.4em]
 \hline

$\fv\fas(4|4)$&$-$&$\spin$&$\fas$&$\fv\fas(4|4)$\cr \hline
&&&&\\[-1.4em]
 \hline

$\fk\fas$&$\Cee[-2]$&$\Pi(\id)$& $\fc\fo(6)$&$\fk\fas(1|6)$\cr
\hline

$\fk\fas(; 1\xi)$
&$\Lambda(1)$&$\id_{\fsl(2)}\boxtimes \id_{\fgl(2;\Lambda(1))}$&
$\fsl(2)\oplus[\fgl(2;\Lambda(1))\ltimes
\fvect(0|1)]$&$\fk\fas(5|5)$\cr \hline

$\fk\fas(; 3\xi)$&$-$&
$\Lambda(3)$&$\Lambda(3)\oplus\fsl(1|3)$&$\fk\fas(4|4)$\cr \hline

$\fk\fas(; 3\eta)$&$-$&$\Vol_{0}(0|3)$&
$\fc(\fvect(0|3))$&$\fk\fas(4|3)$\cr

\hline
&&&&\\[-1.4em]
 \hline

$\fm\fb(4|5)$&$\Pi(\Cee[-2])$&$\Vol^{1/2}
(0|3)$&$\fc(\fvect(0|3))$&$\fm\fb(4|5)$\cr \hline

$\fm\fb(4|5; 1)$&$\Lambda(2)/\Cee 1$
&$\id_{\fsl(2;\Lambda(2))}\boxtimes \vvol^{1/2}$
&$\fc(\fsl(2;\Lambda(2))\ltimes
T^{1/2}(\fvect(0|2))$&$\fm\fb(5|6)$\cr \hline

$\fm\fb(4|5; K)$&$\id_{\fsl(3)}\boxtimes
\Cee[-2]$&$\Pi(\id^*_{\fsl(3)}\boxtimes \id_{\fsl(2)}\boxtimes
\Cee[-1])$&$\fsl(3)\oplus\fsl(2)\oplus\Cee z$&$\fm\fb(3|8)$\cr

\hline
&&&&\\[-1.4em]
 \hline

$\fk\fle(9|6)$&$\Cee[-2]$&$\Pi(T^0_{0}(\vec 0))$&$\fsvect_{3,
4}(0|4)$&$\fk\fle(9|6)$\cr \hline

$\fk\fle(9|6; 2)$&$\Pi(\id_{\fsl(1|3)})$&$\id_{\fsl(2;\Lambda(3))}$&
$\fsl(2;\Lambda(3)) \ltimes \fsl(1|3)$&$\fk\fle(11|9)$\cr \hline

$\fk\fle(9|6;
K)$&$\id$&$\Pi(\Lambda^2(\id^*))$&$\fsl(5)$&$\fk\fle(5|10)$\cr
\hline

$\fk\fle(9|6;
CK)$&$\id_{\fsl(3;\Lambda(1))}^*$&$\id_{\fsl(2)}\boxtimes
\id_{\fsl(3;\Lambda(1))}$&$\fsl(2)\oplus\fsl(3;\Lambda(1)) \ltimes
\fvect(0|1)$&$\fk\fle(9|11)$\cr \hline
\end{tabular}
\]
\textit{Depth $3$}: None of the simple $W$-graded vectorial Lie
superalgebras over $\Cee$ is of depth $>3$ and only two
superalgebras are of depth 3:
\begin{equation}\label{depth3}
\fm\fb(3|8)_{-3}=\Pi(\Cee\boxtimes \id_{\fsl(2)}\boxtimes \Cee[-3]),\quad
\fk\fle(9|11)_{-3}= \Pi(\id_{\fsl(2)}\boxtimes \Cee[-3]).
 \end{equation}
For the definition of the module $\Vol_0$, see \eqref{Vol_0}. Here,
$T^{1/2}$ is the representation of $\fvect$ in the module of
$\frac12$-densities, and $\fas$ is the nontrivial central extension
of $\fspe(4)$, cf.~\cite{BGLL3}. For the definition of $\fsvect_{3,
4}(0|4)$, see Section~\ref{AandB}. In the 0th term $\left
(\fsl(2)\boxtimes \Lambda(3)\right) \ltimes \fsl(1|3)$ of
$\fg=\fk\fle(11;\un|9)$, we consider $\fsl(1|3)$ naturally embedded
into $\fvect(0|3)$ with its tautological action on the space
$\Lambda(3)$ of ``functions''.

For the notation $\Cee[i]$,\index{$Cmathbb[i]$@$\Cee[i]$} see Section~\ref{Conv4.2.1}.

\subsection[The exceptional simple vectorial Lie superalgebras over $\Cee$
as Cartan prolongs]{The exceptional simple vectorial Lie superalgebras over $\boldsymbol{\Cee}$\\
as Cartan prolongs}\label{ssExcAsCTS} For depth $2$, we sometimes write $(\fg_{-2},
\fg_{-1}, \fg_{0})_{*}$ for clarity.

In Table~\eqref{regrad}, there are given indeterminates and their respective
degrees in the re\-gra\-ding~$R(r)$.

\begin{equation*}\footnotesize 
\renewcommand{\arraystretch}{1.4}
\begin{tabular}{|ll|l|}
\hline $\fv\fle(4|3; r)=(\Pi(\Lambda(3))/\Cee\cdot 1,
\fc\fvect(0|3))_{*}$& $\subset\fvect(4|3; R(r))$& $r= 0, 1, K$\cr
\hline

$\fv\fas(4|4)=(\spin, \fas)_{*}$& $\subset \fvect(4|4)$&\cr \hline

$\fk\fas(1|6; r)$& $\subset \fk(1|6; r)$&$r=0, 1\xi, 3\xi$\cr
$\fk\fas(1|6; 3\eta)=(\Vol_{0}(0|3), \fc(\fvect(0|3)))_{*}$&
$\subset \fsvect(4|3)$&$r=3\eta$\cr \hline

$\fm\fb(4|5; r)=(\fba(4), \fc\fvect(0|3))_{*}$& $\subset \fm(4|5;
R(r))$& $r=0, 1, K$\cr \hline

$\fk\fle(9|6; r)=(\fhei(8|6), \fsvect_{3, 4}(0|4))_{*}$& $\subset
\fk(9|6; r)$& $r=0, 2$, CK\cr

$\fk\fle(9|6; K)=(\id_{\fsl(5)}, \Lambda^2(\id^*_{\fsl(5)}),
\fsl(5))_{*}$& $\subset \fsvect(5|10; R(K)), r=K$&\cr

\hline
\end{tabular}
\end{equation*}

\begin{equation}\label{regrad}\footnotesize
\setcounter{MaxMatrixCols}{16}
 \tabcolsep=3pt
\begin{tabular}{|l|l|}
\hlx{hv}
$\fv\fle(4|3)$& $ R(0)=\begin{pmatrix}
x_1&x_2&x_3&y&|&\xi_1&\xi_2&\xi_3\\ 1&1&1&1&|&1&1&1\\ \end{pmatrix}$\cr \hlx{vv}

$\fv\fle(5|4)$&$R(1)=\begin{pmatrix}
x_1&x_2&x_3&y&|&\xi_1&\xi_2&\xi_3\\ 2&1&1&0&|&0&1&1\\ \end{pmatrix}$\cr \hlx{vv}

$\fv\fle(3|6)$&$R(K)=\begin{pmatrix}
x_1&x_2&x_3&y&|&\xi_1&\xi_2&\xi_3\\
2&2&2&0&|&1&1&1\\\end{pmatrix}$\cr \hlx{vhv}

$\fm\fb(4|5)$&$R(0)=\begin{pmatrix}
x_0&x_1&x_2&x_3&|&\xi_0&\xi_1&\xi_2&\xi_3&\tau \\
1&1&1&1&|&1&1&1&1;&2\\\end{pmatrix}$\cr
\hlx{vv}
$\fm\fb(5|6)$& $R(1)=\begin{pmatrix}
x_0&x_1&x_2&x_3&|&\xi_0&\xi_1&\xi_2&\xi_3&\tau \\
0&2&1&1&|&2&0&1&1;&2\\\end{pmatrix}$\cr
\hlx{vv}
$\fm\fb(3|8)$&$R(K)=\begin{pmatrix}
x_0&x_1&x_2&x_3&|&\xi_0&\xi_1&\xi_2&\xi_3&\tau \\
 0&2&2&2&|&3&1&1&1;&3\\\end{pmatrix}$\cr

\hlx{vhv}
$\fk\fas(1|6)$&$R(0)=\begin{pmatrix}
t&|&\xi_1&\xi_2&\xi_3&\eta_1&\eta_2&\eta_3 \\
 2&|&1&1&1&1&1&1\\\end{pmatrix}$\cr
\hlx{vv}
$\fk\fas(5|5)$&$R(1\xi)=\begin{pmatrix}
t&|&\xi_1&\xi_2&\xi_3&\eta_1&\eta_2&\eta_3 \\
 2&|&0&1&1&2&1&1\\\end{pmatrix}$\cr
\hlx{vv}
$\fk\fas(4|4)$&$R(3\xi)=\begin{pmatrix}
t&|&\xi_1&\xi_2&\xi_3&\eta_1&\eta_2&\eta_3 \\
1&|&0&0&0&1&1&1\\\end{pmatrix}$\cr
\hlx{vv}
$\fk\fas(4|3)$&$R(3\eta)=\begin{pmatrix}
t&|&\xi_1&\xi_2&\xi_3&\eta_1&\eta_2&\eta_3 \\
 1&|&1&1&1&0&0&0\\\end{pmatrix}$\cr
\hlx{vhv}
$\fk\fle(9|6)$&$R(0)=\begin{pmatrix}
q_1&q_2&q_3&q_4&p_1&p_2&p_3&p_4&t&|&\xi_1&\xi_2&\xi_3&
\eta_1&\eta_2&\eta_3\\
1&1&1&1&1&1&1&1;&2&|&1&1&1&1&1&1\\
\end{pmatrix}$\cr
\hlx{vv}
$\fk\fle(11|9)$&$R(2)=\begin{pmatrix}
q_1&q_2&q_3&q_4&p_1&p_2&p_3&p_4&t&|&\xi_1&\xi_2&\xi_3&
\eta_1&\eta_2&\eta_3\\
1&1&2&2&1&1&0&0;&2&|&0&1&1&2&1&1\\
\end{pmatrix}$\cr
\hlx{vv}
$\fk\fle(5|10)$&$R(K)=\begin{pmatrix}
q_1&q_2&q_3&q_4&p_1&p_2&p_3&p_4&t&|&\xi_1&\xi_2&\xi_3&
\eta_1&\eta_2&\eta_3\\
2&2&2&2&1&1&1&1;&2&|&1&1&1&1&1&1\\
\end{pmatrix}$\cr
\hlx{vv}
$\fk\fle(9|11)$&$R(CK)=\begin{pmatrix}
q_1&q_2&q_3&q_4&p_1&p_2&p_3&p_4&t&|&\xi_1&\xi_2&\xi_3&
\eta_1&\eta_2&\eta_3\\
3&2&2&2&0&1&1&1;&3&|&2&2&2&1&1&1\\
\end{pmatrix}$\cr
\hlx{vh}
\end{tabular}
\end{equation}

\begin{landscape}
\subsection{Exceptional vectorial Lie superalgebras over $\Kee$ and their
desuperizations}\label{forgetexcept1}
\[\footnotesize
\renewcommand{\arraystretch}{1.3}
\begin{tabular}{|c|c|c|c|c|c|}
\hline $\fg$&$\textbf{F}(\fg_{-2})$&$\textbf{F}(\fg_{-1})$&
$\textbf{F}(\fg_0)$&$\textbf{F}(\fg)$&$\Par\widetilde\un$\cr \hline
&&&&&\\[-1.4em]
 \hline

$\fv\fle(4;\un|3)$&$-$&$\cO(3;\One)/\Kee
1$&$\fc(\fvect(3;\One))$&$\fv\fle\big(7;\widetilde\un\big)$&$3$\cr \hline

$\fv\fle(3;\un|6)$&$\id_{\fsl(3)}\boxtimes
\Kee[*]$&$\id^*_{\fsl(3)}\boxtimes \id_{\fsl(2)}\boxtimes
\Kee[*]$&$\fsl(3)\oplus\fsl(2)\oplus\Kee
z$&$\fv\fle\big(9;\widetilde\un\big)$&$3$\cr \hline
&&&&&\\[-1.4em]
 \hline

$\fk\fas(1;\un|6)$&$\Kee[*]$&$\id$&
$\fc\fo_\Pi^{(1)}(6)$&$\fk\fas\big(7;\widetilde\un\big)$&$7$\cr \hline

$\fk\fas(5;\un|5)$
&$\cO(1;\One)$&$\id_{\fsl(2)}\boxtimes \id_{\fgl(2)}\boxtimes \cO(1;\One)$&
$\fd((\widetilde{\fsl}(W)\oplus (\fgl(V;\cal O)(1;\underline{1}))\ltimes \fvect(1;\One))/\Kee Z)$, see \eqref{FunnyG0}&$\fk\fas\big(10;\widetilde\un\big)$&$7$\cr \hline

$\fk\fas(4;\un|4)$&$-$& $\cO(3;\One)$&$\cO(3;\One)\ltimes
\fd(\fsvect^{(1)}(3;\One))$, see \eqref{nonpos22}&
$\fk\fas\big(8;\widetilde\un\big)$&$7$\cr \hline

$\fk\fas(4;\un|3)$&$-$&$\Vol_{0}(3;\One)$& $\fc(\fvect(3;\One))$&
$\widetilde{\fkas}\big(7;\widetilde\un\big)$&$3$\cr

\hline
&&&&&\\[-1.4em]
 \hline

$\fm\fb(4;\un|5)$&$\Kee[*]$&$\cO(3;\One)$& $\fsvect(3;\One)\rtimes
\cO(3;\One)$&$\fm\fb\big(9;\widetilde\un\big)$&$5$\cr \hline

$\fm\fb(3;\un|8)$&$\id_{\fsl(3)}\boxtimes
\Kee[*]$&$\id^*_{\fsl(3)}\boxtimes \id_{\fsl(2)}\boxtimes
\Kee[*]$&$\fsl(3)\oplus\fsl(2)\oplus\Kee
z$&$\fm\fb_3\big(11;\widetilde\un\big)$&$5$\cr
\hline

$\fm\fb(5;\un|6)$&$\cO(2;\One)/\Kee 1$
&$\id_{\fsl(2)}\boxtimes \cO(2;\One)$
&$\fc(\fsl(2)\boxtimes \cO(2;\One)\ltimes
T^{\infty}(\fsvect(2;\One)\rtimes \cO(2;\One))$&
$\fm\fb_2\big(11;\widetilde\un\big)$&$5$\cr \hline

\hline
&&&&&\\[-1.4em]
 \hline

$\fk\fle(5;\un|10)$&$\id$&$\Lambda^2(\id^*)$&$\fsl(5)$&
$\fk\fle\big(15;\widetilde\un\big)$&$5$\cr \hline

$\fk\fle(11;\un|9)$&$\id_{\fsl(4)}$&$\id_{\fsl(2)}
\boxtimes \cO(3;\One)$ & $\left (\fsl(2)\boxtimes \cO(3;\One)\right)
\ltimes \fpgl(4)$&$\fk\fle_2\big(20;\widetilde\un\big)$&$5$\cr \hline

$\fk\fle(9;\un|11)$&$\id_{\fsl(3)}^*\boxtimes \cO(1;\One)$&
$\id_{\fsl(2)}\boxtimes \left
(\id_{\fsl(3)}\boxtimes \cO(1;\One)\right)$&$\fsl(2)\oplus\left
(\fsl(3)\boxtimes \cO(1;\One) \ltimes
\fvect(1;\One)\right)$&$\fk\fle_3\big(20;\widetilde\un\big)$&$5$\cr
\hline
&&&&&\\[-1.4em]
\hline

$\fv\fle(5;\un|4)$&$\Kee[*]$&$\id\boxtimes \cO(2;\One)$&
$\fc(\fsl(2)\boxtimes \cO(2;\One)\ltimes
T^{\infty}(\fvect(2;\One))$&$\widetilde{\fv\fle}\big(9;\widetilde\un\big)$&$3$\cr
\hline

$\fk\fle(9;\un|6)$&$\Kee[*]$&$T^0_0$&$\fsvect(4;\One)\ltimes\Kee(D+Z)$,
see \eqref{kle9/6} &$\widetilde{\fk\fle}\big(15;\widetilde\un\big)$&$5$\cr
\hline
&&&&&\\[-1.4em]
\hline
$\fv\fas(4;\un|4)$&$-$&$\id_{\textbf{F}(\fas)}$&
$\textbf{F}(\fas)$&$\fvas\big(8;\widetilde\un\big)$&$4$\cr
\hline \hline

\end{tabular}
\]
Recall the definition of the module $\Vol_0$, see \eqref{Vol_0} and
\eqref{Kee*}; before desuperization we replace \eqref{depth3} with
\begin{equation*}
\fm\fb(3|8)_{-3}=\Pi(\Kee\boxtimes \id_{\fsl(2)}\boxtimes \Kee[*]),\quad
\fk\fle(9|11)_{-3}= \Pi(\id_{\fsl(2)}\boxtimes \Kee[*]).
\end{equation*}

To distinguish the two desuperizations of $\fk\fle$ realized by
vector fields on the spaces of the same dimension, we indicate by an
index the depths of these algebras, e.g., \index{$Kmathbb[*]$@$\Kee[*]$}
$\fk\fle_2\big(20;\widetilde\un\big)$; if both algebras are of the same
depth, we cover one of the desuperizations with a~tilde. Clearly,
under the desuperization we should ignore the change of parity in
the negative components of $\textbf{F}(\fg)$.
\end{landscape}

\subsection*{Acknowledgements}

We thank S.~Skryabin for providing us with \cite{Sk0} and elucidations. We heartily thank the referees, especially the one who wrote 26 pages of constructive comments, for their help; their suggestions considerably improved and clarified the exposition.
S.B.\ and D.L.\ were partly supported by the grant
AD 065 NYUAD.

\LastPageEnding


\addcontentsline{toc}{section}{References}
\begin{thebibliography}{99}
\footnotesize\itemsep=0pt

\bibitem{ALSh}
Alekseevsky D., Leites D., Shchepochkina I., Examples of simple {L}ie
 superalgebras of vector fields, \textit{C.~R. Acad. Bulg. Sci.} \textbf{34}
 (1080), 1187--1190.

\bibitem{BBH}
Benayadi S., Bouarroudj S., Hajli M., Double extensions of restricted {L}ie
 (super)algebras, \href{https://doi.org/10.1007/s40598-020-00149-5}{\textit{Arnold Math.~J.}} \textbf{6} (2020), 231--269,
 \href{https://arxiv.org/abs/1810.03086}{arXiv:1810.03086}.

\bibitem{BGP}
Benkart G., Gregory T., Premet A., The recognition theorem for graded {L}ie
 algebras in prime characteristic, \href{https://doi.org/10.1090/memo/0920}{\textit{Mem. Amer. Math. Soc.}} \textbf{197}
 (2009), xii+145~pages, \href{https://arxiv.org/abs/math.RA/0508373}{arXiv:math.RA/0508373}.

\bibitem{BW}
Block R.E., Wilson R.L., Classification of the restricted simple {L}ie
 algebras, \href{https://doi.org/10.1016/0021-8693(88)90216-5}{\textit{J.~Algebra}} \textbf{114} (1988), 115--259.

\bibitem{BGLL1}
Bouarroudj S., Grozman P., Lebedev A., Leites D., Divided power (co)homology.
 {P}resentations of simple finite dimensional modular {L}ie superalgebras with
 {C}artan matrix, \href{https://doi.org/10.4310/HHA.2010.v12.n1.a13}{\textit{Homology Homotopy Appl.}} \textbf{12} (2010),
 237--278, \href{https://arxiv.org/abs/0911.0243}{arXiv:0911.0243}.

 \bibitem{BGLL3}
Bouarroudj S., Grozman P., Lebedev A., Leites D., Derivations and central
 extensions of simple modular {L}ie algebras and superalgebras,
 \href{https://arxiv.org/abs/1307.1858}{arXiv:1307.1858}.

\bibitem{BGLLS}
Bouarroudj S., Grozman P., Lebedev A., Leites D., Shchepochkina I., New simple
 {L}ie algebras in characteristic~2, \href{https://doi.org/10.1093/imrn/rnv327}{\textit{Int. Math. Res. Not.}}
 \textbf{2016} (2016), 5695--5726, \href{https://arxiv.org/abs/1307.1551}{arXiv:1307.1551}.

\bibitem{BGL4}
Bouarroudj S., Grozman P., Leites D., Deforms of the symmetric modular {L}ie
 superalgebras, \href{https://arxiv.org/abs/0807.3054}{arXiv:0807.3054}.

\bibitem{BGL3}
Bouarroudj S., Grozman P., Leites D., New simple modular {L}ie superalgebras as
 generalized prolongations, \href{https://doi.org/10.1007/s10688-008-0025-3}{\textit{Funct. Anal. Appl.}} \textbf{42} (2008),
 161--168, \href{https://arxiv.org/abs/0704.0130}{arXiv:0704.0130}.

\bibitem{BGL1}
Bouarroudj S., Grozman P., Leites D., Classification of finite dimensional
 modular {L}ie superalgebras with indecomposable {C}artan matrix,
 \href{https://doi.org/10.3842/SIGMA.2009.060}{\textit{SIGMA}} \textbf{5} (2009), 060, 63~pages, \href{https://arxiv.org/abs/0710.5149}{arXiv:0710.5149}.

\bibitem{BKLLS}
Bouarroudj S., Krutov A., Lebedev A., Leites D., Shchepochkina I., Restricted
 {L}ie (super)algebras in characteristic~3, \href{https://doi.org/10.1007/s10688-018-0206-7}{\textit{Funct. Anal. Appl.}}
 \textbf{52} (2018), 49--52, \href{https://arxiv.org/abs/1809.08582}{arXiv:1809.08582}.

\bibitem{BLLS2}
Bouarroudj S., Lebedev A., Leites D., Shchepochkina I., Classifications of
 simple {L}ie superalgebras in characteristic~$2$, \href{https://arxiv.org/abs/1407.1695}{arXiv:1407.1695}.

\bibitem{BLLS1}
Bouarroudj S., Lebedev A., Leites D., Shchepochkina I., Lie algebra
 deformations in characteristic~2, \href{https://doi.org/10.4310/MRL.2015.v22.n2.a3}{\textit{Math. Res. Lett.}} \textbf{22}
 (2015), 353--402, \href{https://arxiv.org/abs/1301.2781}{arXiv:1301.2781}.


\bibitem{BLW}
Bouarroudj S., Lebedev A., Wagemann F., Deformations of the {L}ie algebra
 {${\mathfrak o}(5)$} in characteristics~3 and~2, \href{https://doi.org/10.1134/S0001434611050191}{\textit{Math. Notes}}
 \textbf{89} (2011), 777--791, \href{https://arxiv.org/abs/0909.3572}{arXiv:0909.3572}.

\bibitem{BL}
Bouarroudj S., Leites D., Simple {L}ie superalgebras and nonintegrable
 distributions in characteristic~{$p$}, \href{https://doi.org/10.1007/s10958-007-0046-0}{\textit{J.~Math. Sci.}} \textbf{141}
 (2007), 1390--1398, \href{https://arxiv.org/abs/math.RT/0606682}{arXiv:math.RT/0606682}.

\bibitem{Bro}
Brown G., Families of simple {L}ie algebras of characteristic two,
 \href{https://doi.org/10.1080/00927879508825259}{\textit{Comm. Algebra}} \textbf{23} (1995), 941--954.

\bibitem{CCK}
Cantarini N., Cheng S.-J., Kac V., Errata to: ``{S}tructure of some {$\mathbb
 Z$}-graded {L}ie superalgebras of vector fields'' [\href{https://doi.org/10.1007/BF01237358}{\textit{{T}ransform.
 {G}roups}} {\bf 4} (1999), 219--272] by {C}heng and {K}ac, \href{https://doi.org/10.1007/s00031-004-9005-8}{\textit{Transform.
 Groups}} \textbf{9} (2004), 399--400.

\bibitem{ChKu}
Chebochko N.G., Kuznetsov M.I., Integrable cocycles and global deformations of
 {L}ie algebra of type~{$G_2$} in characteristic~2, \href{https://doi.org/10.1080/00927872.2016.1233241}{\textit{Comm. Algebra}}
 \textbf{45} (2017), 2969--2977.

\bibitem{CK1}
Cheng S.-J., Kac V.G., Generalized {S}pencer cohomology and filtered
 deformations of {${\bf Z}$}-graded {L}ie superalgebras, \href{https://doi.org/10.4310/ATMP.1998.v2.n5.a7}{\textit{Adv. Theor.
 Math. Phys.}} \textbf{2} (1998), 1141--1182, \href{https://arxiv.org/abs/math.RT/9805039}{arXiv:math.RT/9805039}.

\bibitem{CK2}
Cheng S.-J., Kac V.G., Structure of some {$\mathbb Z$}-graded {L}ie
 superalgebras of vector fields, \href{https://doi.org/10.1007/BF01237358}{\textit{Transform. Groups}} \textbf{4} (1999), 219--272.

\bibitem{CK1a}
Cheng S.-J., Kac V.G., Addendum: ``{G}eneralized {S}pencer cohomology and
 filtered deformations of {${\mathbb Z}$}-graded {L}ie superalgebras''
 [\href{https://doi.org/10.4310/ATMP.1998.v2.n5.a7}{\textit{{A}dv. {T}heor. {M}ath. {P}hys.}} {\bf 2} (1998), 1141--1182],
 \href{https://doi.org/10.4310/ATMP.2004.v8.n4.a2}{\textit{Adv. Theor. Math. Phys.}} \textbf{8} (2004), 697--709.

\bibitem{Ei}
Eick B., Some new simple {L}ie algebras in characteristic 2,
 \href{https://doi.org/10.1016/j.jsc.2010.05.003}{\textit{J.~Symbolic Comput.}} \textbf{45} (2010), 943--951.

\bibitem{Et}
Etingof P., Koszul duality and the {PBW} theorem in symmetric tensor categories
 in positive characteristic, \href{https://doi.org/10.1016/j.aim.2017.06.014}{\textit{Adv. Math.}} \textbf{327} (2018),
 128--160, \href{https://arxiv.org/abs/1603.08133}{arXiv:1603.08133}.

\bibitem{GPS}
Gomis J., Par\'{\i}s J., Samuel S., Antibracket, antifields and gauge-theory
 quantization, \href{https://doi.org/10.1016/0370-1573(94)00112-G}{\textit{Phys. Rep.}} \textbf{259} (1995), 1--145,
 \href{https://arxiv.org/abs/hep-th/9412228}{arXiv:hep-th/9412228}.

\bibitem{GZ}
Grishkov A., Zusmanovich P., Deformations of current {L}ie algebras.
 {I}.~{S}mall algebras in characteristic~2, \href{https://doi.org/10.1016/j.jalgebra.2016.11.024}{\textit{J.~Algebra}} \textbf{473}
 (2017), 513--544, \href{https://arxiv.org/abs/1410.3645}{arXiv:1410.3645}.

\bibitem{Gr}
Grozman P., SuperLie (2013), \url{http://www.equaonline.com/math/SuperLie}.

\bibitem{GL1}
Grozman P., Leites D., Defining relations for {L}ie superalgebras with {C}artan
 matrix, \href{https://doi.org/10.1023/A:1026642004008}{\textit{Czechoslovak~J. Phys.}} \textbf{51} (2001), 1--21,
 \href{https://arxiv.org/abs/hep-th/9702073}{arXiv:hep-th/9702073}.

\bibitem{GL4}
Grozman P., Leites D., Structures of {$G(2)$} type and nonintegrable
 distributions in characteristic~{$p$}, \href{https://doi.org/10.1007/s11005-005-0026-6}{\textit{Lett. Math. Phys.}} \textbf{74}
 (2005), 229--262, \href{https://arxiv.org/abs/math.RT/0509400}{arXiv:math.RT/0509400}.

\bibitem{GLP}
Grozman P., Leites D., Poletaeva E., Defining relations for classical {L}ie
 superalgebras without {C}artan matrices, \href{https://doi.org/10.4310/hha.2002.v4.n2.a12}{\textit{Homology Homotopy Appl.}} \textbf{4} (2002), 259--275,
 \href{https://arxiv.org/abs/math.RT/0202152}{arXiv:math.RT/0202152}.

\bibitem{ILL}
Iyer U.N., Lebedev A., Leites D., Prolongs of (ortho-)orthogonal {L}ie
 (super)algebras in characteristic~2, \href{https://doi.org/10.1142/S1402925110000866}{\textit{J.~Nonlinear Math. Phys.}}
 \textbf{17} (2010), suppl.~1, 253--309.

\bibitem{ILMS}
Iyer U.N., Leites D., Messaoudene M., Shchepochkina I., Examples of simple
 vectorial {L}ie algebras in characteristic~2, \href{https://doi.org/10.1142/S1402925110000878}{\textit{J.~Nonlinear Math.
 Phys.}} \textbf{17} (2010), suppl.~1, 311--374.

\bibitem{Kfil}
Kac V.G., Description of filtered {L}ie algebras with which graded {L}ie
 algebras of {C}artan type are associated, \href{https://doi.org/10.1070/IM1974v008n04ABEH002128}{\textit{Math. USSR Izv.}} \textbf{8}
 (1974), 801--835.

\bibitem{K2}
Kac V.G., Lie superalgebras, \href{https://doi.org/10.1016/0001-8708(77)90017-2}{\textit{Adv. Math.}} \textbf{26} (1977), 8--96.

\bibitem{KWK}
Kac V.G., Corrections to: ``Exponentials in {L}ie algebras of
 characteristic~$p$''
 [\href{https://doi.org/10.1070/IM1971v005n04ABEH001116}{\textit{Math. USSR
 Izv.}} {\bf 5} (1971), 777--803], \href{https://doi.org/10.1070/IM1995v045n01ABEH001632}{\textit{Russian Acad. Sci. Izv. Math.}}
 \textbf{45} (1995), 229.

\bibitem{K}
Kac V.G., Classification of infinite-dimensional simple linearly compact {L}ie
 superalgebras, \href{https://doi.org/10.1006/aima.1998.1756}{\textit{Adv. Math.}} \textbf{139} (1998), 1--55.

\bibitem{K10}
Kac V.G., Classification of supersymmetries, in Proceedings of the
 {I}nternational {C}ongress of {M}athematicians, {V}ol.~{I} ({B}eijing, 2002),
 Higher Ed. Press, Beijing, 2002, 319--344, \href{https://arxiv.org/abs/math-ph/0302016}{arXiv:math-ph/0302016}.

\bibitem{Kapp}
Kaplansky I., Graded {L}ie algebras, {U}niversity Chicago, Chicago, 1975,
 available at \url{http://www1.osu.cz/~zusmanovich/links/files/kaplansky/}.

\bibitem{Kir}
Kirillov S.A., Sandwich algebras in simple finite-dimensional {L}ie algebras,
 Ph.D.~Thesis, {N}izhny Novgorod, 1992, for a~short summary, see Kirillov~S.A., The sandwich subalgebra in {L}ie algebras of {C}artan type, \textit{Russian Math. (Iz. VUZ)} \textbf{36} (1992), 16--23.

\bibitem{KL}
Kochetkov Yu., Leites D., Simple finite-dimensional {L}ie algebras in
 characteristic 2 related to superalgebras and on a notion of finite simple
 group, in Proceedings of the International Conference on Algebra, Part~1
 (Novosibirsk, August 1989), \textit{Contemp. Math.}, Vol.~131, Editors L.A.~Bokut', Yu.L.~Ershov, A.I.~Kostrikin, Amer. Math. Soc., Providence, RI, 1992, 59--67.

\bibitem{Kos1}
Kostrikin A.I., The beginnings of modular {L}ie algebra theory, in Group
 Theory, Algebra, and Number Theory ({S}aarbr\"{u}cken, 1993), de Gruyter,
 Berlin, 1996, 13--52.

\bibitem{KD}
Kostrikin A.I., Dzhumadil'daev A.S., Modular {L}ie algebras: new trends, in
 Algebra ({M}oscow, 1998), Editor Yu.~Bahturin, de Gruyter, Berlin, 2000,
 181--203.

\bibitem{Koch}
Kotchetkoff Yu.Yu., D\'eformations des superalg\`ebres de {B}uttin et
 quantification, \textit{C.~R.~Acad. Sci. Paris S\'er.~I Math.} \textbf{299}
 (1984), 643--645, for details, see id., Deformations of {L}ie superalgebras,
 VINITI depositions, 1984, 24~pages (in Russian), Letter to {L}eites, 1985.

\bibitem{KrLe}
Krutov A., Lebedev A., On gradings modulo 2 of simple {L}ie algebras in
 characteristic~2, \href{https://doi.org/10.3842/SIGMA.2018.130}{\textit{SIGMA}} \textbf{14} (2018), 130, 27~pages,
 \href{https://arxiv.org/abs/1711.00638}{arXiv:1711.00638}.

\bibitem{KLLS}
Krutov A., Lebedev A., Leites D., Shchepochkina I., Non-degenerate invariant
 symmetric bilinear forms on simple {L}ie superalgebras in characteristic~$2$,
 Oberwolfach Preprint OWP 2020-02, available at
 \url{https://publications.mfo.de/handle/mfo/3697}.

\bibitem{LaF}
Ladilova A.A., Filtered deformations of the {F}rank algebras, \href{https://doi.org/10.3103/S1066369X09080076}{\textit{Russian
 Math.}} \textbf{53} (2009), 43--45.

\bibitem{LaY}
Ladilova A.A., Filtered deformations of {L}ie algebras of series~{$Y$},
 \href{https://doi.org/10.1007/s10958-009-9739-x}{\textit{J.~Math. Sci.}} \textbf{164} (2010), 91--94.

\bibitem{LaZ}
Ladilova A.A., Filtered deformations of {L}ie algebras of series~$Z$, in
 Proceedings of the International Con\-fe\-rence Dedicated to 100-th Anniversary
 of V.V.~Morozov (Kazan, September 25--30, 2011) and Youth School-Conference
 ``Modern Problems of Algebra and Mathematical Logic'' (Kazan Federal
 University, 2011), Kazan Federal University, Russia, 2011, 123--124.

\bibitem{Leb}
Lebedev A., Analogs of the orthogonal, {H}amiltonian, {P}oisson, and contact
 {L}ie superalgebras in characteris\-tic~2, \href{https://doi.org/10.1142/S1402925110000854}{\textit{J.~Nonlinear Math. Phys.}}
 \textbf{17} (2010), suppl.~1, 217--251.

\bibitem{LL}
Lebedev A., Leites D., On realizations of the {S}teenrod algebras (with an
 appendix by P.~Deligne), \textit{J.~Prime Res. Math.} \textbf{2} (2006),
 101--112.

\bibitem{Le2}
Leites D., New {L}ie superalgebras and mechanics, \textit{Soviet Math. Doklady}
 \textbf{18} (1977), 1277--1280.

\bibitem{Lei}
Leites D., Towards classification of simple finite dimensional modular {L}ie
 superalgebras, \textit{J.~Prime Res. Math.} \textbf{3} (2007), 101--110,
 \href{https://arxiv.org/abs/0710.5638}{arXiv:0710.5638}.

\bibitem{Lsos}
Leites D. (Editor), Seminar on supersymmetries, Vol.~1, Algebra and calculus on
 supermanifolds, MCCME, Moscow, 2011, a~draft in English is available.

\bibitem{Lsos2}
Leites D. (Editor), Seminar on supersymmetries, Vol.~2, Algebra and calculus on
 supermanifolds. Additional chapters, MCCME, Moscow, {t}o appear.

\bibitem{LSerg}
Leites D., Serganova V., Metasymmetry and {V}olichenko algebras, \href{https://doi.org/10.1016/0370-2693(90)91086-Q}{\textit{Phys.
 Lett.~B}} \textbf{252} (1990), 91--96, for a~correction and additions
 see~\cite{Lsos2}.

\bibitem{LS}
Leites D., Shchepochkina I., The classification of simple {L}ie superalgebras
 of vector fields, Preprint MPIM-2003-28, available at
 \url{http://www.mpim-bonn.mpg.de/preblob/2178}, for a~short version, see
 Toward classification of simple vectorial {L}ie superalgebra, in Seminar on
 Supermanifolds (Reports of Stockholm University, 31/1988-14), 235--278, and
 in Differential Geometric Methods in Theoretical Physics (Davis, CA, 1988),
 Editors W.~Nahm, L.~Chau, \textit{NATO Adv. Sci. Inst. Ser.~B Phys.},
 Vol.~245, Plenum, New York, 1990, 633--651.

\bibitem{LS1}
Leites D., Shchepochkina I., On the classification of simple {L}ie
 superalgebras of vector fields with polynomial coefficients, an updated and
 detailed version of~\cite{LS}.

\bibitem{LS0}
Leites D., Shchepochkina I., Quivers and {L}ie superalgebras, in Commutative
 Algebra, Representation Theory and Combinatorics, Conference in Honor of
 D.~Buchsbaum (Northeastern Uversity, Boston, October 18--20, 1997), Editors
 D.~Eisenbud, A.~Martsinkovsky, J.~Weyman, Northeastern Uversity, Boston,
 1997, 67.

\bibitem{LSh}
Leites D., Shchepochkina I., How should an antibracket be quantized?,
 \href{https://doi.org/10.1023/A:1010312700129}{\textit{Theoret. and Math. Phys.}} \textbf{126} (2001), 281--306,
 \href{https://arxiv.org/abs/math-ph/0510048}{arXiv:math-ph/0510048}.

\bibitem{Ma}
Manin Yu.I., Introduction to the theory of schemes, \href{https://doi.org/10.1007/978-3-319-74316-5}{\textit{Moscow Lectures}},
 Vol.~1, Springer, Cham, 2018.

\bibitem{Mo}
Molotkov V., Infinite-dimensional and colored supermanifolds,
 \href{https://doi.org/10.1142/S140292511000088X}{\textit{J.~Nonlinear Math. Phys.}} \textbf{17} (2010), suppl.~1, 375--446.

\bibitem{Ri}
Richardson Jr. R.W., On the rigidity of semi-direct products of {L}ie algebras,
 \href{https://doi.org/10.2140/pjm.1967.22.339}{\textit{Pacific~J. Math.}} \textbf{22} (1967), 339--344.

\bibitem{Ru}
Rudakov A.N., Deformations of simple {L}ie algebras, \href{https://doi.org/10.1070/IM1971v005n05ABEH001204}{\textit{Math. USSR-Izv.}}
 \textbf{5} (1971), 1120--1126.

\bibitem{Sh3}
Shchepochkina I., Exceptional simple infinite-dimensional {L}ie superalgebras,
 \textit{C.~R. Acad. Bulg. Sci.} \textbf{36} (1983), 313--314.

\bibitem{Sh5}
Shchepochkina I., Five simple exceptional {L}ie superalgebras of vector fields,
 \href{https://doi.org/10.1007/BF02465205}{\textit{Funct. Anal. Appl.}} \textbf{33} (1999), 208--219,
 \href{https://arxiv.org/abs/hep-th/9702121}{arXiv:hep-th/9702121}.

\bibitem{Sh14}
Shchepochkina I., The five exceptional simple {L}ie superalgebras of vector
 fields and their fourteen regradings, \href{https://doi.org/10.1090/S1088-4165-99-00012-6}{\textit{Represent. Theory}} \textbf{3}
 (1999), 373--415, \href{https://arxiv.org/abs/hep-th/9702121}{arXiv:hep-th/9702121}.

\bibitem{Shch}
Shchepochkina I., How to realize a {L}ie algebra by vector fields,
 \href{https://doi.org/10.1007/s11232-006-0078-5}{\textit{Theoret. and Math. Phys.}} \textbf{147} (2006), 821--838,
 \href{https://arxiv.org/abs/math.RT/0509472}{arXiv:math.RT/0509472}.

\bibitem{ShP}
Shchepochkina I., Post G., Explicit bracket in an exceptional simple {L}ie
 superalgebra, \href{https://doi.org/10.1142/S0218196798000235}{\textit{Internat.~J. Algebra Comput.}} \textbf{8} (1998),
 479--495, \href{https://arxiv.org/abs/physics/9703022}{arXiv:physics/9703022}.

\bibitem{Shen}
Shen G.Y., Variations of the classical {L}ie algebra {$G_2$} in low
 characteristics, \textit{Nova~J. Algebra Geom.} \textbf{2} (1993), 217--243.

\bibitem{Sk0}
Skryabin S., The normal shapes of symplectic and contact forms over algebras of
 divided powers, VINITI deposition, 1986, 8504-B86 (in Russian),
 \href{https://arxiv.org/abs/1906.11496}{arXiv:1906.11496}.

\bibitem{Sk1}
Skryabin S., Classification of Hamiltonian forms over divided power algebras,
 \href{https://doi.org/10.1070/SM1991v069n01ABEH001232}{\textit{Math. USSR Sb.}} \textbf{69} (1991), 121--141.

\bibitem{Sk2}
Skryabin S., A contragredient {L}ie algebra of dimension 29 over a field of
 characteristic~3, \href{https://doi.org/10.1007/BF00971230}{\textit{Sib. Math.~J.}} \textbf{34} (1993), 548--554.

\bibitem{SkT1}
Skryabin S., Toral rank one simple {L}ie algebras of low characteristics,
 \href{https://doi.org/10.1006/jabr.1997.7231}{\textit{J.~Algebra}} \textbf{200} (1998), 650--700.

\bibitem{S}
Strade H., Simple {L}ie algebras over fields of positive characteristic.
 {I}.~{S}tructure theory, \textit{De Gruyter Expositions in Mathematics},
 Vol.~38, \href{https://doi.org/10.1515/9783110197945}{Walter de Gruyter \& Co.}, Berlin, 2004.

\bibitem{Tyu}
Tyurin S.A., The classification of deformations of a special {L}ie algebra of
 {C}artan type, \href{https://doi.org/10.1007/BF01140026}{\textit{Math. Notes}} \textbf{24} (1988), 948--954.

\bibitem{WK}
Weisfeiler B.Ju., Kac V.G., Exponentials in {L}ie algebras of
 characteristic~{$p$}, \href{https://doi.org/10.1070/IM1971v005n04ABEH001116}{\textit{Math. USSR Izv.}} \textbf{5} (1971), 777--803.

\bibitem{W}
Wilson R.L., Simple {L}ie algebras of type~{$S$}, \href{https://doi.org/10.1016/0021-8693(80)90182-9}{\textit{J.~Algebra}}
 \textbf{62} (1980), 292--298.

\bibitem{Y}
Yamaguchi K., Differential systems associated with simple graded {L}ie
 algebras, in Progress in Differential Geometry, \textit{Adv. Stud. Pure
 Math.}, Vol.~22, \href{https://doi.org/10.2969/aspm/02210413}{Math. Soc. Japan}, Tokyo, 1993, 413--494.

\end{thebibliography}
\end{document}